\documentclass[11pt]{article}
\usepackage{lscape}
 \usepackage[ansinew]{inputenc}
 \usepackage[dvips]{graphicx}
\usepackage[psamsfonts]{amssymb} 
\usepackage[all]{xy}   
 \usepackage[psamsfonts]{eucal}
\usepackage{amssymb, amsmath}
\usepackage{amsthm}
\usepackage{enumerate}
\newcommand{\be}{\begin{equation}}
\newcommand{\ee}{\end{equation}}

\newcommand{\am}{ \, \widehat{\amalg}\, }
\newcommand{\bd}{\begin{displaymath}}
\newcommand{\ed}{\end{displaymath}}
\newcommand{\ben}{\begin{enumerate}}
\newcommand{\een}{\end{enumerate}}
\newtheorem{Theorem}{Theorem} 

\newtheorem{Prop}{Proposition}[section]

\newtheorem{lem}{Lemma}[section]

\title{Enriched Differential Lie Algebras in Topology}
 \author{Yves F\'elix and Steve Halperin}
\begin{document}

\maketitle

The context of this article is the Sullivan rationalization
$$X\to X_{\mathbb Q}$$
of path connected topological spaces. This is constructed by Sullivan via his Sullivan algebras, which are commutative differential graded algebras (cdga's), $(\land V,d)$, in which: (i) $V= V^{\geq 1}$, (ii) $\land V$ is the free graded commutative algebra generated by $V$, and (iii) $d$ satisfies a certain "Sullivan condition". As Sullivan realized from the outset, the component, $d_1 : V\to V\land V$ of $d$, defines by duality a Lie bracket in 
$$L_V :=s^{-1} \mbox{Hom}(V;\mathbb Q);$$
($s^{-1}$ is inverse suspension). Moreover, \emph{$L_V$ is a graded Lie algebra, the homotopy Lie algebra of $(\land V,d)$}.

We recall now that each connected space admits a quasi-isomorphism
$$(\land V,d)\stackrel{\simeq}{\longrightarrow} A_{PL}(X)$$
from a Sullivan algebra, where $A_{PL} :$ Top $\longrightarrow$ Cdga is a functor defined by Sullivan that induces an isomorphism $H(A_{PL}(X))\cong H^*(X;\mathbb Q)$. By definition, this identifies $(\land V,d)$ as a Sullivan model for $X$. For these models (\cite{S})
 $$H^*(X;\mathbb Q) \cong H(\land V,d)\hspace{3mm}\mbox{and there is a bijection of sets }sL_V\stackrel{\cong}{\longrightarrow} \pi_*(X_{\mathbb Q}).$$
  Thus a  Sullivan model is an algebraic object linking homotopy invariants of the topological spaces $X$ and $X_{\mathbb Q}$. 
  
  Our focus throughout this article is to establish new
properties and examples of these invariants, many of which together with the details of the constructions, are provided in Part III.

Sullivan algebras,  $(\land V,d)$, have an additional fundamental (but nonetheless not generally appreciated) property :
$$V = \varinjlim_\alpha V_\alpha,$$
where $\{V_\alpha\}$ is the family of finite dimensional subspaces of $V$ for which $\land V_\alpha$ is preserved by $d$. Dualizing exhibits $L_V$ as the inverse limit
$$L_V= \varprojlim_\alpha L_\alpha, \hspace{1cm}  
 L_\alpha = L_{V_\alpha}.$$ In particular, each $L_\alpha$ is finite dimensional and nilpotent.
 
  This motivates our introduction of the category of enriched graded Lie algebras:

\vspace{2mm}\noindent {\bf Definition.} An \emph{enriched graded Lie algebra} is a graded Lie algebra, $L= L_{\geq 0}$, together with an inverse system of surjections $\rho_\alpha : L\to L_\alpha$ onto finite dimensional nilpotent Lie algebras. If $L\stackrel{\cong}{\longrightarrow} \varprojlim_\alpha L_\alpha$, then $L$ is complete.

\vspace{2mm} It follows from \cite{Bou} that the functor $\varprojlim$ is exact when applied to inverse systems of short exact sequences of finite dimensional vector spaces.  This makes the enriched structure of a homotopy Lie algebra an essential property on which most of this article   relies. 

In particular, just as the introduction of differentials in a number of algebraic categories permits the construction of a homotopy theory for those categories, we introduce here the category Edgl of enriched differential graded Lie algebras $(L, \partial)$: 

\vspace{2mm} \noindent {\bf Definition.} \emph{An enriched differential graded Lie algebra}  (edgl) is a differential graded Lie algebra with an enriched structure given by
$$(L, \partial) = \varprojlim_\alpha (L_\alpha, \partial_\alpha).$$

\vspace{2mm}This category is then equipped (Part IV) with a homotopy theory   analogous to that developed by Sullivan for cdga's. Specifically, analogous to Sullivan algebras are the \emph{profree} edgl's $(\overline{\mathbb L}_T, \partial)$ in which $\overline{\mathbb L}_T$ is the completion of the free graded Lie algebra $\mathbb L_T$ generated by $T$ (Part II). Just as a Sullivan algebra $(\land V,d)$ is minimal if $d: V\to \land^{\geq 2}V$ so $(\overline{\mathbb L}_T, \partial)$ is minimal if $\partial : T\to $ completion of $[\overline{\mathbb L}_T,\overline{\mathbb L}_T]$. Similarly with the Sullivan model of a cdga, any edgl $(L, \partial)$ admits a quasi-isomorphism
$$(\overline{\mathbb L}_T,\partial) \stackrel{\simeq}{\longrightarrow} (L,\partial)$$
from a minimal profree edgl, and this determines $(\overline{\mathbb L}_T,\partial)$ up to isomorphism.

Beyond the analogy, there is an explicit connection between edgl's and Sullivan algebras. An enriched edgl $(L, \partial) = \varprojlim_\alpha (L_\alpha, \partial_\alpha)$ determines the Sullivan algebra $(\land V, d_0+d_1)= \varinjlim_\alpha (\land V_\alpha, d_{0,\alpha}+d_{1,\alpha})$ in which $\partial_\alpha$ dualizes to $d_{0,\alpha}$ and $d_{1,\alpha}$ dualizes to the Lie bracket $L_\alpha \times L_\alpha \stackrel{[\,\,]}{\longrightarrow} L_\alpha$.

Now suppose $X$ is a connected space. An \emph{edgl model} of $X$ is a profree edgl $(\overline{\mathbb L}_T, \partial)$ for which the corresponding Sullivan algebra is a Sullivan model for $X$. Moreover, (Part V) $X$ has a minimal edgl model $(L,\partial)$ and this is unique up to isomorphism. It satisfies
$$H(L,\partial) = L_X,  \mbox{the homotopy Lie algebra of $X$, and } sT \hspace{2mm}\mbox{ is the dual of }  H^*(X;\mathbb Q) .$$

This introduction of edgl models fully places Sullivan's approach to rational homotopy theory in the long standing tradition, initiated by Quillen \cite{Q} in 1970, of modelling spaces by differential graded Lie algebras. That tradition remains very active in recent work, for example in \cite{LM} and \cite{4A}.
 
 Moreover, whereas Sullivan models provide simple descriptions of the fibre of a fibration, edgl models provide simple descriptions of the cofibre of a cofibration such as a space of the form $X\cup_YZ$. In both cases the appropriate choice of model can enable a computational approach as well as the construction of interesting examples. (Unsurprisingly Sullivan models and edgl models do not provide simple descriptions of cofibres, respectively fibres). 
 
For example, if $\gamma : (\overline{\mathbb L}_T,0) \to (L, \partial)$ is an edgl representative of a continuous map
 $$g : \vee_\alpha S^{n_\alpha}\to X$$ then an edgl model of $X\cup_g \vee_\alpha D^{n_\alpha +1}$ is given by the free Lie product
 $$(L\, \widehat{\amalg}\, \overline{\mathbb L}_{sT}, \partial + \partial_\gamma)$$
 where $\partial_\gamma (L) = 0$ and $\partial_\gamma (sx) = \gamma (x)$. This construction extends to provide an edgl model of a CW complex, $X$ of the form $(\overline{\mathbb L}_T, \partial)$ in which for each $k$
 \begin{enumerate}
 \item[$\bullet$] $(\overline{\mathbb L}_{T_{<k}}, \partial)$ is an edgl model of the $k$-skeleton, $X_k$, and
 \item[$\bullet$] The inclusion $(\overline{\mathbb L}_{T_{<k}}, \partial)\to (\overline{\mathbb L}_{T_{<k}}\, \widehat{\amalg}\, \overline{\mathbb L}_{T_k}, \partial)$ is obtained from the attaching map constructing $X_{k+1}= X_k\cup_g \vee_\alpha D^{k+1}_\alpha$.
 \end{enumerate}
 
  The interplay between Sullivan models and edgl models is a powerful tool for understanding the impact of simple topological constructions on rational homotopy groups. Because of the natural morphism
  $$\pi_*(X)\to \pi_*(X_{\mathbb Q})$$ this can also provide insight into the structure of $\pi_*(X)$ and, in particular, on $\pi_1(X_{\mathbb Q})$ and $\pi_1(X)$. 
 In Part VI we consider a particular question : when does a cell attachment 
  $$\iota_X : X\to X\cup_g \vee_\alpha D^{n_\alpha + 1}$$
  result in a surjection
  $$\pi_*(\iota_X)_{\mathbb Q} : \pi_*(X_{\mathbb Q})\to \pi_*\left( X\cup_g \vee_\alpha D^{n_\alpha +1}\right)_{\mathbb Q} ?$$
  In this case we say $g$ is \emph{rationally inert}. Even in some relatively simple examples it is unknown what conditions force $\pi_*(\iota_X)$ to be surjective. By contrast, we provide four specific characterizations of rational inertness, one specific to Sullivan models and another in terms of edgl models. For instance, extending results obtained in \cite{HL} for simply connected spaces, we show that 
  \begin{enumerate}
  \item[$\bullet$] The attachment of a single cell to a wedge of odd spheres is rationally inert if the attaching map is not rationally trivial, and
  \item[$\bullet$] If $Y^n$ is a finite CW complex, and its cohomology algebra $H^*(Y;\mathbb Q)$ satisfies Poincar\'e duality and has at least two generators, then $Y$ is obtained from the $(n-1)$-skeleton by a rationally inert attachment.
  \end{enumerate}
  
  In summary, this is an example of the application of edgl models in combination with cdga models to understand the morphism
  $$\pi_*(X)\to \pi_*(X_{\mathbb Q}) $$
  in particular situations. More generally, our focus throughout is best summarized by the title of Part VII: Applications in Topology, which contains further applications and examples.
  
Rational homotopy theory begins with three fundamental contributions: Quillen's \cite{Q} for simply connected spaces, that of Bousfield-Kan \cite{BK}, and Sullivan's \cite{S}. We note that for CW complexes of finite type the latter two are equivalent \cite{BG}. More recently,   there are other approaches to a rational homotopy for non-simply connected spaces based on dgl's. The two closest ones are due to   Lazarev and Markl (\cite{LM}) and to Buijs, F\'elix, Murillo and Tanr\'e (\cite{4A}). In the second one, the authors associate   to each finite simplicial set $X$ a pronilpotent dgl $M_X$.   For a more general space $X$ their chosen model is   the limit $\varinjlim_a  M_{(X_a)}$,  limit  taken over all the finite subcomplexes $X_a\subset X$.   For finite simplicial spaces our models are quasi-isomorphic to these.

There are of course many other approaches to rational homotopy for non-simply connected spaces. In (\cite{HGT}) the authors associate to a space a simplicial model of its universal cover together with a cofibrant replacement of the action of the fundamental group. In (\cite{Pridham}), the author associates to a simplicial space $X$ the $\mathbb Q$-pro-algebraic completion of its simplicial group $GX$, denoted $X^{alg} $. This model has also nice properties, and in particular an isomorphism between $\pi_0(X^{alg}) $ and the Malcev completion of $\pi_1(X)$. 
  
  \vspace{5mm}\noindent {\bf Conventions}
  \begin{enumerate}
  \item[$\bullet$] A \emph{topological space} is either a simplicial set or a CW complex, and the two homotopy categories are identified by the singular simplex and Steenrod realization functors \cite{May}. In either case a connected space is path connected.
  \item[$\bullet$] By a sphere we mean a \emph{connected sphere} $S^n$, $n\geq 1$.
  \item[$\bullet$] If $X$ is a topological space then for simplicity, $H(X) = H^*(X;\mathbb Q)$ denotes the rational singular cohomology algebra of $X$. We use the standard notation for homology nd for singular cohomology with other coefficients.
  \item[$\bullet$] Graded objects are families $Q = \{Q^n\}_{n\in \mathbb Z}$, where some of the $Q^n$ may be equal to $0$. We adopt the convention $\{Q^n\}_{n\in \mathbb Z}= \{Q_n\}_{n\in \mathbb Z}$ with $Q_n = Q^{-n}$. If $x\in Q$ then $(-1)^{deg\, x} = (-1)^n$ whether $x\in Q^n$ or $x\in Q_n$.
  \item[$\bullet$] If $Q= \{Q^n\}$ is a graded object then its suspension $sQ$ is the graded object defined by $(sQ)_n = Q_{n-1}$. 
  \item[$\bullet$] Because the notation $(\,\,)^*$ is used everywhere, and because dual vector spaces appear here so frequently, we simplify the notation Hom$(V,\mathbb Q)$ and write
  $$(V^\vee)_n = \mbox{Hom}(V^n, \mathbb Q)$$
  for any graded vector space $V$.
  \item[$\bullet$] A \emph{graded Lie algebra} is a graded vector space $L$ together with a linear map of degree $0$, $[\,,\,]: L\otimes L\to L$, satisfying
  $$[x,y]+ (-1)^{deg\,x\,deg\, y} [y,x] = 0 \hspace{3mm}\mbox{and } [x,[y,z]]= [[x,y],z] + (-1)^{deg\, x\,deg\, y}[y,[x,z]].$$

  \item[$\bullet$] By \emph{a commutative graded differential algebra (cdga)}, $(A,d)$ we mean a graded algebra $A= \{A^n\}_{n\geq 0}$ satisfying $ab= (-1)^{deg\, a\, deg\, b}\, ba$ and $H^0(A) = \mathbb Q$. 
  
  An \emph{augmentation} is a morphism $\varepsilon : (A,d)\to (\mathbb Q,0)$.

  \item[$\bullet$] A \emph{Sullivan algebra} is a cdga $(\land V,d)$ in which $\land V$ is a free graded commutative algebra. We denote by $\land^pV$ the space spanned by the products of length $p$ of elements in $V$. To be a Sullivan algebra the differential must satisfy $V = \cup_{n\geq 0} V_n$ where $V_0 = V\cap \mbox{ker}\, d$ and $V_{n+1} = V\cap d^{-1}(\land V_n)$.
  \item[$\bullet$] A \emph{quadratic Sullivan algebra} is a Sullivan algebra $(\land V,d)$ in which $d : V\to \land^2V$.
  \item[$\bullet$] A \emph{$\Lambda$-extension} is a sequence of cdga morphisms
  $$\xymatrix{A\ar[rr]^\eta&& A\otimes \land Z\ar[rr]^\rho && \land Z}$$
  in which: $\eta (a)= a\otimes 1$, $\rho= \varepsilon_A\otimes id$ with $\varepsilon_A$ an augmentation, and the quotient $\land Z$ is a Sullivan algebra.
  \end{enumerate}

\newpage
\tableofcontents

\newpage
\part{Enriched Lie algebras}

\section{Introduction}

By an \emph{enriched Lie algebra} we mean a graded Lie algebra, $L=L_{\geq 0}$, equipped with a family   ${\mathcal I}= \{\rho_\alpha : L\to L_\alpha\}$ of surjective morphisms satisfying the following conditions
\begin{enumerate}
\item[(i)] Each $L_\alpha  $ is a finite dimensional nilpotent Lie algebra. 
\item[(ii)] The index set, $\{\alpha\}$ is a directed set under the partial order given by
$$\alpha \geq \beta \Longleftrightarrow \mbox{ker}\, \rho_\alpha \subset \mbox{ker}\, \rho_\beta:$$
i.e. for each $\alpha, \beta$ there is a $\gamma$ such that $\mbox{ker}\, \rho_\gamma \subset \mbox{ker}\, \rho_\alpha \cap \mbox{ker}\, \rho_\beta$.
\item[(iii)] $\cap_\alpha \mbox{ker}\, \rho_\alpha = 0$.
\end{enumerate}
 
Moreover, we shall use ${\mathcal I}$ to denote both the family $\{\rho_\alpha\}$, and the index set $\{\alpha\}$, and denote ker$\, \rho_\alpha$ by $I_\alpha$.

\vspace{3mm}\noindent {\bf Example.} In any graded Lie algebra $L= L_{\geq 0}$ the set of all surjective morphisms onto   finite dimensional and nilpotent graded Lie algebras forms a directed system as above. If the intersection of the kernels is zero then this defines   an enriched structure.
However, in principle, $L$ may admit two distinct such structures.

Finally, as we shall show in \S 8, a minimal Sullivan algebra induces a natural enriched structure in its homotopy Lie algebra.

\vspace{3mm}
\noindent {\bf Definition.} A \emph{morphism} $\varphi : (E, {\mathcal I}_E) \to (F, {\mathcal I}_F)$ of enriched Lie algebras is a morphism of graded Lie algebras such that for each $\alpha \in {\mathcal I}_F$ there is a $\beta\in {\mathcal I}_E$ such that $\varphi : \mbox{ker}\, \rho_\beta \to \mbox{ker}\, \rho_\alpha$.

\subsection{Completions}

\vspace{3mm}\noindent {\bf Definition.} An enriched Lie algebra $(L, {\mathcal I}_L)$ is \emph{complete} if $L=\varprojlim_\alpha L_\alpha$.

\vspace{3mm}
If $(L, {\mathcal I}_L)$ is an enriched Lie algebra we define $(\overline{L}, {\mathcal I}_{\overline{L}})$ by setting $\overline{L} = \varprojlim_\alpha L_\alpha$ and ${\mathcal I}_{\overline{L}}= \{\overline{\rho_\alpha}: \overline{L}\to L_\alpha\}$ where the $\overline{\rho_\alpha}$ are the projections of the inverse limit.

It follows from the defining properties of enriched Lie algebras that $(\overline{L}, {\mathcal I}_{\overline{L}}$ is complete, that the map $\lambda_L : L\to \overline{L}$ induced by the $\rho_\alpha$ is injective and that
 $$\lambda_L : (L, {\mathcal I}_L)\to (\overline{L},  {\mathcal I}_{\overline{L}})$$
is a morphism of enriched Lie algebras.

\vspace{3mm}\noindent {\bf Definition.} The morphism $\lambda_L$ is the \emph{completion} of $(L, {\mathcal I}_L)$.

\vspace{3mm}In particular,  $(L, {\mathcal I}_L)$ is \emph{complete} if $\lambda_L$ is an isomorphism.  
Moreover, almost by definition, we have
\begin{eqnarray}
\label{i1} 
L\cap \mbox{ker}\, \overline{\rho_\alpha} = \mbox{ker}\, \rho_\alpha
\end{eqnarray}
and that $\overline{L}$ itself is complete.

Finally, observe that any morphism $\psi :  (E, {\mathcal I}_E) \to (F, {\mathcal I}_F)$ of enriched Lie algebras extends naturally to a morphism $\overline{\psi}$ between their completions. In fact, set ${\mathcal I}_{E, \alpha}= \{\beta\, \vert\, \mbox{ker}\, \rho_\beta \subset \psi^{-1}(\mbox{ker}\, \rho_\alpha)\}$. Then $\psi$ induces morphisms
$$\psi_\alpha : \varprojlim_{\beta\in {\mathcal I}_{E, \alpha}} E_\beta \to F_\alpha.$$
Moreover, $\overline{E} = \varprojlim_\alpha \varprojlim_{\beta \in {\mathcal I}_{E, \alpha} } E_\beta$ and $\psi$ extends to the morphism
$$\overline{\psi} = \varprojlim_\alpha \psi_\alpha.$$ It is straightforward to check that $\overline{\psi}$ is the unique morphism extending $\psi$.

\vspace{3mm}
Now observe that a sub Lie algebra $E\subset L$ of an enriched Lie algebra, $(L, \{\rho_\alpha\})$, inherits the enriched structure given by the surjections  $\rho_\alpha\vert_E : E\to \rho_\alpha(E)\subset \rho_\alpha (L)$. Thus the inclusion of $E$ in $L$ is a morphism of enriched Lie algebras. Moreover, if $E$ is an ideal then the surjections $L/E\to \rho_\alpha L/\rho_\alpha E$ define an enriched structure ${\mathcal I}_{L/E}$ in $L/E$.

\vspace{3mm}\noindent {\bf Definition.} $(E, {\mathcal I}_E)$ and $(L/E, {\mathcal I}_{L/E})$ are respectively the \emph{sub and quotient enriched Lie algebras} determines by the inclusion $E\to L$.   

\vspace{3mm}\noindent {\bf Definition.} A subset $S\subset L$ in a complete enriched Lie algebra is \emph{closed} if $S=\overline{S}= \varprojlim_\alpha \rho_\alpha (S)$. 

\vspace{3mm}\noindent {\bf Remark.} It is immediate that if $E$ is a sub Lie algebra of a complete enriched Lie algebra $(L, {\mathcal I}_L)$ then the completion $(\overline{E}, {\mathcal I}_{\overline{E}})$ is a sub enriched Lie algebra of $(\overline{L}, {\mathcal I}_{\overline{L}})$.

\vspace{3mm}
Next, if ${\mathcal I}$ and ${\mathcal I}'$ are families   satisfying the conditions above for $L$, we say that ${\mathcal I}$ and ${\mathcal I}'$ are \emph{equivalent} (${\mathcal I}\sim {\mathcal I}'$) if the identities, $(L, {\mathcal I})\to (L, {\mathcal I}')$ and $(L, {\mathcal I}')\to (L, {\mathcal I})$  are morphisms. This is an equivalence relation, and if $\varphi : (F, {\mathcal G})\to (E, {\mathcal I})$ is a morphism, then $\varphi : (F, {\mathcal G}') \to (E, {\mathcal I}') $ is also a morphism whenever ${\mathcal G}\sim {\mathcal G}'$ and ${\mathcal I}\sim {\mathcal I}'$. In particular, if ${\mathcal I}$ and ${\mathcal I}'$ are equivalent families, then the identity induces an isomorphism
$$\varprojlim_{\alpha \in \mathcal I} L_\alpha \stackrel{\cong}{\longrightarrow} \varprojlim_{\beta\in {\mathcal I}'} L_\beta.$$
Thus the completion, $\lambda_L: L\to \overline{L}$ is independent of the choice of ${\mathcal I}$ in its equivalence class.

\vspace{3mm}
Now we define the category, ${\mathcal C}$, of enriched Lie algebras as follows:
\begin{enumerate}
\item[$\bullet$] The objects of ${\mathcal C}$ are the pairs $(L, \{{\mathcal I}\})$ in which $L=L_{\geq 0}$ is a graded Lie algebra and $\{{\mathcal I}\}$ is an equivalence class  satisfying the conditions above.
\item[$\bullet$] The morphism of ${\mathcal C}$ are the morphisms $\varphi : (F,{\mathcal G})\to (E, {\mathcal I})$.
\end{enumerate}
In particular $L\leadsto \overline{L}$ is a functor from ${\mathcal C}$ to the sub category of complete enriched Lie algebras, the inclusion $L\to \overline{L}$ is a morphism, and any morphism $(E, {\mathcal J})\to (F, {\mathcal I})$ extends uniquely to a morphism $\overline{E}\to \overline{F}$.

\vspace{3mm}\noindent {\bf Examples.} (1). Let $(L,\{\rho_\alpha\})$ be an enriched Lie algebra and let $(\rho_\beta)$ be a subfamily such that for each $\alpha $ there is a $\beta$ with $\mbox{ker}\, \rho_\beta\subset \mbox{ker}\, \rho_\alpha$. Then $(\rho_\alpha)$ and $(\rho_\beta)$ are equivalent families.

(2). \emph{Full enriched Lie algebras}. An enriched Lie algebra $(L, {\mathcal I}) $ is \emph{full} if whenever an ideal $I$ contains some $\mbox{ker}\, \rho_\beta$ then $\overline{L}= \mbox{ker}\, \rho_\alpha$ for some $\alpha$. If $(L, {\mathcal J}) $ is any enriched Lie algebra then ${\mathcal J} $ extends to a unique full enriched family ${\mathcal I}$ consisting of all the ideals containing some $\mbox{ker}\, \rho_\beta$. Evidently ${\mathcal I}\sim {\mathcal J}$. Moreover if ${\mathcal I}$ and ${\mathcal J}$ are equivalent full enriched structure, then they coincide.

(3). Suppose that $L$ is the inverse limit of an inverse system of nilpotent finite dimensional Lie algebras $Q_\alpha$ : $L=\varprojlim_\alpha Q_\alpha$. Denote $q_\alpha : L\to Q_\alpha$ the corresponding projections. Then $(L, \{  \rho_\alpha\})$ is a complete enriched Lie algebra with $p_\alpha : L\to q_\alpha (L)$ coinciding with $q_\alpha$.

\vspace{3mm} A key aspect of an enriched Lie algebra $(L, {\mathcal I})$ is that ${\mathcal I}$ makes the Lie algebra accessible to finite dimensional arguments, even when $L$ is not a graded vector space of finite type. In particular, enriched Lie algebras behave well with respect to inverse limits, as follows from the following special case of a Theorem in \cite{Bou}, for which we give a short proof in the next section :
 
\begin{Prop} (Bourbaki).
\label{p1.1} Suppose $0\to A_\alpha\to B_\alpha\to C_\alpha \to 0$ is a directed system of short exact sequences of graded vector spaces of finite type. Then
\begin{eqnarray}
\label{i2}
0\to \varprojlim_\alpha A_\alpha \to \varprojlim_\alpha B_\alpha \to \varprojlim_\alpha C_\alpha \to 0
\end{eqnarray}
is also short exact. In particular, for chain complexes of finite type, homology commutes with inverse limits.
\end{Prop}

The first properties of complete enriched Lie algebras are given by the following Lemma

\begin{lem}
\label{l1.1}
Let $L$ be a complete enriched Lie algebra. Then,
\begin{enumerate}
\item[(i)] Let $I\subset L$ be a closed ideal in a complete enriched Lie algebra, then $L/I$ is a complete enriched Lie algebra.
\item[(ii)] A finite sum of closed subspaces of $L$ is closed.  
\item[(iii)] An arbitrary intersection of closed subspaces of $L$ is closed.
\item[(iv)] If $S\subset L$ has finite type, then $[S,L]$ is closed.
\end{enumerate}\end{lem}

\vspace{3mm}\noindent {\sl proof.} (i) By Proposition \ref{p1.1}, we have
$$\overline{L/I }= \varprojlim_\alpha L_\alpha / \rho_\alpha (I) = \varprojlim_\alpha L_\alpha \, /\, \varprojlim_\alpha \rho_\alpha (I) = L/I.$$

(ii) It is sufficient to show that if $S$ and $R$ are closed then $S+R$ is closed. By Proposition \ref{p1.1},
$$\renewcommand{\arraystretch}{1.7}
\begin{array}{ll}
\overline{S+R}/\overline{S} &= \displaystyle\varprojlim_\alpha \rho_\alpha (S+R) / \rho_\alpha (S) = \displaystyle\varprojlim_\alpha \displaystyle\frac{\rho_\alpha (S) + \rho_\alpha (R)}{\rho_\alpha (S)}
\\
& = \displaystyle\varprojlim_\alpha \displaystyle\frac{\rho_\alpha (R)}{\rho_\alpha (S) \cap \rho_\alpha (R)} 
  = \displaystyle \overline{R} / \displaystyle \varprojlim_\alpha \rho_\alpha (S) \cap \rho_\alpha (R)
\end{array}
\renewcommand{\arraystretch}{1}$$
In particular, it follows that the inclusions of $\overline{S}$ and $\overline{R}$ in $\overline{S+R}$ define a surjection $\overline{S}+ \overline{R}\to \overline{S+R}$. Since $S$ and $R$ are closed, $S+R\to \overline{S+R}$
is an isomorphism.

\vspace{3mm} (iii) Suppose $x\in \overline{\cap_\sigma S_\sigma}$, where each $S_\sigma$ is a closed subspace of $T$. Then $$\rho_\alpha x\in \rho_\alpha (\cap S_\sigma) \subset \cap_\sigma \rho_\alpha (S_\sigma).$$
Therefore $\rho_\alpha x \in \rho_\alpha (S_\sigma)$ for each $\sigma$. Thus $x = (\rho_\alpha x)\in \varprojlim_\alpha \rho_\alpha (S_\sigma) = \overline{S_\sigma} = S_\sigma$. Thus $x\in \cap_\sigma S_\sigma$ and $\cap_\sigma S_\sigma = \overline{\cap_\sigma S_\sigma}$.

\vspace{3mm} (iv) In view of (ii) it is sufficient to prove that for any $x\in L$, $[x,L]$ is closed. Let $C_\alpha = \{y_\alpha \in L_\alpha \,\vert\, [\rho_\alpha x, y_\alpha]= 0\}.$ It is immediate from the definition that $\{C_\alpha\}$ is an inverse system, and we set
$$C = \varprojlim_\alpha C_\alpha.$$

Now we show that $C = \{y\in L\, \vert\, [x,y]= 0\}.$  Indeed, if $y\in C$, then for all $\alpha$,
$$\rho_\alpha [x,y] = [\rho_\alpha x, \rho_\alpha y] \in [\rho_\alpha x, C_\alpha] = 0,$$
and hence $ [x,y]= 0$. On the other hand, if $[x,y]= 0$ then $[\rho_\alpha x, \rho_\alpha y]= \rho_\alpha [x,y]= 0$ and so each $\rho_\alpha y\in C_\alpha$. Thus $y\in C$. 

Finally, apply (i) to obtain the commutative diagram
$$\xymatrix{L/C = \varprojlim_\alpha L_\alpha \, / \, \varprojlim_\alpha C_\alpha\ar[d]_\cong^{[x,-]} \ar[rrr]^-= &&&\varprojlim_\alpha L_\alpha/C_\alpha\ar[d]_\cong^{[x,-]}\\
[x,L] \ar[rrr]^-= &&&\varprojlim_\alpha [\rho_\alpha x, L_\alpha]}$$
Since $\varprojlim_\alpha [\rho_\alpha x, L_\alpha]= \overline{[x,L]}$, (v) is established. \hfill$\square$

\vspace{3mm} Complete enriched Lie algebras arise naturally in Sullivan's approach to rational homotopy theory :   the homotopy Lie algebra of a Sullivan model is naturally a complete enriched Lie algebra and this is an essential aspect when the theory is extended to all connected CW complexes. Thus our objective here  is  to develop the properties of enriched Lie algebras, in particular with a view to future topological applications. In particular, while many of the basic properties of ${\mathcal C}$ are parallel those for the category of all graded Lie algebras, the enrichment gives the objects in ${\mathcal C}$ a distinct "homotopy flavour". In particular each enriched Lie algebra $(L, {\mathcal I})$ determines a quadratic Sullivan algebra $\land Z$ and a simplicial set $\langle \land Z\rangle$.

\subsection{Universal enveloping algebras and the group $G_L$}
 
\vspace{1mm} Recall next that the classical completion of the universal enveloping algebra $UE$ of a graded Lie algebra $E$ is the inverse limit
 $$\widehat{UE} = \varprojlim_n UE/J^n,$$
 $J^n$ denoting the $n^{th}$ power of the augmentation ideal. In particular
 $$\widehat{UE} = \mathbb Q \oplus \widehat{J},$$
 and $\widehat{J} := \varprojlim_n J/J^n$ is the augmentation ideal for $\widehat{UE}$.

 \vspace{3mm}\noindent {\bf Definition.} 
The \emph{completion} $\overline{UL}$ of an enriched Lie algebra $(L, \{\rho_\alpha\})$ is defined  by
 $$\overline{UL}:= \varprojlim_\alpha \widehat{UL_\alpha} = \varprojlim_{n,\alpha} UL_\alpha/J_\alpha^n.$$

 \vspace{2mm}Passing to inverse limits shows that the inclusions $L_\alpha \to L_\alpha/J_\alpha^n$ define an inclusion
 $$\overline{L}\to \overline{UL}.$$ Note also that   $$ \overline{UL} = \mathbb Q \oplus \overline{J},$$
 where $\overline{J} = \varprojlim_{\alpha, n} J_\alpha/J_\alpha^n$; $\overline{J}$ is the augmentation ideal for $\overline{UL}$. It is also immediate that 
  a morphism $\varphi : E\to L$ of enriched Lie algebras extends to a morphism $\overline{U\varphi}: \overline{UE}\to \overline{UL}$.

 Finally, in analogy with the filtration of $UL$ by the ideals $J^n$ we filter $\overline{UL}$ by the ideals
 $$\overline{J}^{(n)}:= \varprojlim_{k,\alpha} J_\alpha^n/J_\alpha^{k+n} \subset \overline{J}= \overline{J}^{(1)}.$$
 Then $\overline{UL}$ is complete with respect to this filtration:
 \begin{eqnarray}\label{i3}
 \overline{UL}= \varprojlim_n \overline{UL}/\overline{J}^{(n)}.
 \end{eqnarray}
 In fact, (\ref{i2}) yields
 $$
 \renewcommand{\arraystretch}{1.4}
 \begin{array}{ll}
 \overline{UL} = \varprojlim_{n,\alpha} UL_\alpha /J_\alpha^n & = \varprojlim_k \varprojlim_{n,\alpha} (UL_\alpha/J_\alpha^{n+k})/ (J_\alpha^n/J_\alpha^{n+k})\\
 & = \varprojlim_n \varprojlim_{k,\alpha} (UL_\alpha/J_\alpha^{n+k})/ \varprojlim_{k,\alpha} (J_\alpha^n/J_\alpha^{n+k})\\
 & = \varprojlim_n \overline{UL}/\overline{J}^{(n)}.
 \end{array}
 \renewcommand{\arraystretch}{1}
$$ 
 
\vspace{3mm}\noindent {\bf Example.} 
 When  $\, L/[L,L]$ has finite type, then $\overline{UL}= \widehat{UL}$.
 
\vspace{3mm}\noindent {\bf Remark.} In analogy with complete enriched Lie algebra, a \emph{complete enriched algebra} is the projective limit of an inductive systems of finite dimensional nilpotent graded algebras. The graded commutator makes a complete enriched algebra a complete enriched graded Lie algebra. The induced functor from complete enriched algebras to complete enriched Lie algebras has a left adjoint functor that associates to a Lie algebra $L$ its universal enveloping algebra $\overline{UL}$.

 \vspace{3mm} Associated with a complete enriched Lie algebra, $L$, is a group $G_L$ defined as follows: for each $\alpha$, the diagonal map $L_\alpha \to L_\alpha \times L_\alpha$ extends to a morphism
 $$\Delta_\alpha : \widehat{U}L_\alpha \to \widehat{U}L_\alpha \, \widehat{\otimes}\, \widehat{U}L_\alpha:= \widehat{U}(L_\alpha \times L_\alpha).$$
 An element $x\in \widehat{U}L_\alpha$ is called \emph{group-like} if $\Delta_\alpha x= x\otimes x$, and the group-like elements form a group $G_\alpha$ with multiplication given by the product in the algebra $\widehat{U}L_\alpha$.
 
 \vspace{3mm}\noindent {\bf Definition.} The \emph{fundamental group}, $G_L$, of a complete enriched Lie algebra, $L$, is the group
 $$G_L = \varprojlim_\alpha G_\alpha.$$
 
 \vspace{3mm} Note that, since each $L_\alpha$ is finite dimensional and nilpotent, the classical exponential map $L_\alpha \to \widehat{U}L_\alpha$ is well defined and it is straightforward to verify that it is a bijection onto $G_\alpha$. The inverse is then given by the standard power series for $\log$, which is trivially also convergent. Passing to inverse limits this yields bijections which we denote by 
 $$\xymatrix{L \ar@/^/[rr]^\exp && G_L.\ar@/^/[ll]^\log}$$

\section{Enriched vector spaces}

\vspace{3mm}\noindent {\bf Definitions.}   \begin{enumerate}
\item[(1)] An \emph{enriched vector space} is a graded vector space $T$, together with a family ${\mathcal I}_T$ of surjective linear maps $\rho_\alpha : T\to T_\alpha$ of degree zero, and satisfying
\begin{enumerate}
\item[(i)] Each dim$\, T_\alpha<\infty$.
\item[(ii)] The index set $\{\alpha\}$ is a directed set under the partial order given by
$$\alpha \geq \beta \, \Longleftrightarrow \mbox{ker}\, \rho_\alpha \subset \mbox{ker}\, \rho_\beta.$$
\item[(iii)] $\cap_\alpha \mbox{ker}\, \rho_\alpha = 0$. 
\end{enumerate}
\item[(2)] A linear map $\varphi : S\to T$ (of any degree) is \emph{coherent} with respect to enriched structures ${\mathcal I}_S = \{\rho_\beta\}$ and ${\mathcal I}_{T} = \{\rho'_\alpha\}$ if for each $\alpha$ the is a $\beta$ such that $\varphi(\mbox{ker}\,\rho_\beta)\subset \mbox{ker}\, \rho'_\alpha$. If $\varphi$ has degree $0$, it is a \emph{morphism} of enriched vector spaces.
\item[(3)] Two enriched structures $\{\rho_\alpha\}$ and $\{\rho'_\beta\}$ in a vector space $T$ are \emph{equivalent} if $id_T$ is coherent as a map $(T, \{\rho_\alpha\})\to (T, \{\rho'_\beta\})$ and as a map $(T, \{\rho'_\beta\})\to (T, \{\rho_\alpha\})$.
\end{enumerate}

As with enriched Lie algebras we shall use ${\mathcal I}_T$ to denote both the family $\{\rho_\alpha\}$ and the index set $\{\alpha\}$, and we note that the underlying vector space of an enriched Lie algebra inherits the 
structure of an enriched vector space.

\vspace{3mm}\noindent {\bf Important remark}. The same enriched vector space may have two inequivalent enriched Lie algebra structures.

\vspace{3mm} Now for enriched vector spaces   the following   definitions extend those in $\S 1$.

\vspace{3mm}\noindent {\bf Definitions.}
\begin{enumerate}
\item[(1)] The \emph{completion} of an enriched vector space $(T, {\mathcal I}_T)$ is the inclusion
$$\lambda_T : (T, {\mathcal I}_T)\to (\overline{T}, {\mathcal I}_{\overline{T}})$$
in which $\overline{T} = \varprojlim_\alpha T_\alpha$ and ${\mathcal I}_{\overline{T}}$ is the set of surjections $\varprojlim_\alpha T_\alpha \to T_\alpha$. In particular, $(T, {\mathcal I}_T)$ is complete if $\lambda_T : T\to \overline{T}$ is an isomorphism, and $(\overline{T}, {\mathcal I}_{\overline{T}})$ is always complete.
\item[(2)] \emph{Sub and quotient spaces.} Suppose $(T, {\mathcal I}_T)$ is an enriched vector space and $S\subset T$ is any subspace. The enriched structure in $T$ endows $S$ and $T/S$ with the enriched structures given respectively by
$${\mathcal I}_S = \{\rho_\alpha : S\to \rho_\alpha(S)\, \vert\, \alpha \in {\mathcal I}_T\} \hspace{3mm}\mbox{and } {\mathcal I}_{T/S}=\{\overline{\rho_\alpha}: T/S\to \rho_\alpha T/\rho_\alpha S\, \vert \alpha \in {\mathcal I}_T\}.$$
They are called, respectively, a \emph{sub enriched vector space} and a \emph{quotient enriched vector space}.
\item[(3)] A subspace $S\subset T$ in a complete enriched vector space is \emph{closed} if $S= \overline{S}$. 
\end{enumerate}

\vspace{3mm}\noindent {\bf Example.} The families of ideals $x^n$ and $(1-x)^n$ define two structures of enriched vector space on $\mathbb Q[x]$. The identity is not a coherent morphism of enriched vector space. However the morphism $\varphi : \mathbb Q[x]\to \mathbb Q[x]$ defined by $\varphi (x) = 1-x$ is a coherent morphism inducing an isomorphism of completions $\mathbb Q[[x]]\to \mathbb Q[[1-x]]$.

\vspace{3mm}\noindent {\bf Remarks}.
\begin{enumerate}
\item[(1)] If $(S, {\mathcal I}_S)\subset (T, {\mathcal I}_T)$ is an enriched subspace, then $\overline{S}\subset \overline{T}$ and the completion $(\overline{S}, {\mathcal I}_{\overline{S}})$ is an enriched subspace of $(\overline{T}, {\mathcal I}_{\overline{T}})$. 
\item[(2)]  A subspace $S\subset T$ is closed if and only if each subspace $S_k$ of degree $k$ is closed.
\end{enumerate}

\vspace{3mm} Suppose $(T, {\mathcal I}_T)$ is an enriched vector space. The advantage of the condition dim$\, T_\alpha<\infty$ is that for these spaces $T_\alpha \to (T_\alpha^\vee)^\vee$ is an isomorphism. Moreover the construction $(\,\, )^\vee$ converts inverse and inductive systems respectively to inductive and inverse systems. In particular we have the

\vspace{3mm}\noindent {\bf Definition.} The \emph{predual} of an enriched vector space $(T, {\mathcal I}_T)$ is the space, $$V:= \varinjlim_\alpha V_\alpha,$$
where $V_\alpha = T_\alpha^\vee$.

\vspace{3mm}\noindent {\bf Remarks.} 
\begin{enumerate}
\item[(1)] Let $V= \varinjlim_\alpha V_\alpha$ be the predual of an enriched vector space $(T, {\mathcal I}_T)$. Since $(\,)^\vee$ converts direct limits to inverse limits,
$$V^\vee = \varprojlim_\alpha V_\alpha^\vee = \varprojlim_\alpha T_\alpha = \overline{T}.$$
In particular, $T$ is complete if and only if it is the dual of its predual. 
\item[(2)]   This also provides the promised proof of Proposition \ref{p1.1}, since $(\, )^\vee$ preserves exact sequences. In fact a short exact sequence
$$0 \to A_\alpha \to B_\alpha \to C_\alpha\to 0$$
of graded vector spaces of finite type dualizes to a short exact sequence, and an inverse system of such short exact sequences is the dual of the inductive systems of their duals. Since direct limits preserve short exact sequences it follows that
$$0 \leftarrow \varinjlim_\alpha A_\alpha^\vee \leftarrow \varinjlim_\alpha B_\alpha^\vee \leftarrow \varinjlim_\alpha C_\alpha^\vee \leftarrow 0$$
is short exact. Dualizing again shows that
$$0\to \varprojlim_\alpha A_\alpha \to \varprojlim_\alpha B_\alpha \to \varprojlim_\alpha C_\alpha\to 0$$
is short exact.
\end{enumerate}

\begin{Prop}
\label{p2.1} Suppose a graded vector space $T$ has the form $V^\vee$.
\begin{enumerate}
\item[(i)] The inclusions of all the finite dimensional subspaces $V_\alpha$ of $V$ dualize to define a complete enriched structure in $T$ for which $V$ is the predual.
\item[(ii)] Any complete enriched structure in $T$ for which $V$ is the predual is equivalent to the structure in (i).
\item[(iii)] If $\varphi : S\to T$ is a coherent linear map between    complete enriched vector spaces, then $\varphi$ is the dual of a map $\psi : W\leftarrow V$ between their preduals, and any linear map $\psi : W\leftarrow V$ dualizes to a coherent linear map. 
\item[(iv)] A sequence $Z\leftarrow W\leftarrow V$ of linear maps is exact at $W$ if and only if the dual sequence is exact. In particular, if $\varphi : S\to T$ is a coherent linear map between complete enriched vector spaces then Im$\, \varphi$ and Ker$\, \varphi$ are complete.
\item[(v)]    Suppose $(N, \partial)$ is a complete enriched $\mathbb Z$-graded chain complex, equipped with a decreasing filtration $N(q)$ for which $\cap_q N(q)= 0$. The corresponding spectral sequence is then convergent.
\end{enumerate}
\end{Prop}

\vspace{3mm}\noindent {\sl proof.} (i) This is immediate from the definition.

(ii) Fix some enriched structure $\{\rho_\beta\}$ in $T$ with predual $V$ then $(\rho_\beta T)^\vee = V_\beta$ for some finite dimensional subspace $V_\beta \subset V$, and $V = \varinjlim_\beta V_\beta$. It follows from
$$V = \varinjlim_\beta \,(\varinjlim_{V_\alpha \subset V_\beta}\, V_\alpha)$$
that $\{\rho_\beta\}$ is equivalent to the $\{\rho_\alpha\}$ of (i). 

(iii) and (iv) Let $\{V_\alpha\}$ be the finite dimensional subspaces of $V$ and let $\{W_\gamma\}$ be the finite dimensional subspaces of $W$. Then set ${\mathcal I}_\gamma = \{\alpha \vert \psi (V_\alpha)\subset W_\gamma\}$.
From $V = \varinjlim_\gamma \varinjlim_{\alpha \in {\mathcal I}_\gamma} \, V_\alpha$ and the fact that $\psi : \varinjlim_{\alpha \in {\mathcal I}_\gamma} V_\alpha  \to W_\gamma $ it follows that $\psi^\vee : S\to T$ is coherent. On the other hand, if $\varphi $ is coherent, it follows that for each $\alpha$ there is a $\gamma$ such that $\varphi$ factors to give a map $S_\gamma \to T_\alpha$. Dualizing gives a map $V_\alpha \to W_\gamma$ and passing to direct limits defines $\psi : W\leftarrow V$ with $\psi^\vee = \varphi$.

Finally, elementary linear algebra gives Im$\, \varphi= (W/\mbox{ker}\, \psi)^\vee$ and Ker$\, \varphi= (V/ \psi (W))^\vee$, and so these subspaces are closed.

 (v) 
 In fact $N = M^\vee$ and $N/N(q) = M(q)^\vee$ where $N(q)\to N\to N/N(q)$ is the dual of $M(q)\to M\to M/M(q)$. The $M(q)$ provide an increasing filtration with $M=\cup_q M(q)$ so that the associated spectral sequence is convergent. The spectral sequence for $N$ is the dual spectral sequence.
 
 \hfill$\square$

\begin{lem}
\label{l2.1}
Let $(V(r), f_{rs} : V(r)\to V(s))$ be a   diagram of complete enriched vector spaces. Then $V= \varprojlim_r V(r)$ is a complete enriched vector space and is the limit of the diagram in the category of complete enriched vector spaces. 
\end{lem}

\vspace{3mm}\noindent {\sl proof.} Denote by $(U_r, g_{rs}: U_s\to U_r)$ the predual of the diagram and let $U = \varinjlim_r U_r$. Then $V = \varprojlim_r V_r = U^\vee$ is an enriched vector spaces and a straightforward computation shows that this is the limit in the category of complete enriched vector spaces.
\hfill$\square$

\begin{lem}
\label{l2.2} Suppose $S\subset T$ and $(T, \{\rho_\alpha\})$ is an enriched vector space. Then   
$$\overline{T}/\overline{S} = \varprojlim_\alpha \rho_\alpha T/\rho_\alpha S = \overline{(T/S)}.$$
\end{lem}

\vspace{3mm}\noindent {\sl proof:} Since $\rho_\alpha S\subset \rho_\alpha T$,
$$0\to \rho_\alpha S\to \rho_\alpha T\to \rho_\alpha T/\rho_\alpha S\to 0$$
is a short exact sequence of finite dimensional spaces. Thus the Lemma follows from Proposition \ref{p1.1}.

 \hfill$\square$

\begin{lem}
\label{l2.3}  Suppose $(T, \{\rho_\alpha\})$ is a complete enriched vector space.
\begin{enumerate}
\item[(i)] A finite sum of closed subspaces of $T$ is closed.  
\item[(ii)] An arbitrary intersection of closed subspaces of $T$ is closed.
\item[(iii)] For any $S\subset T$, $\overline{S}= \cap_\alpha (S+\mbox{ker}\, \rho_\alpha)$.
\item[(iv)] If $S\subset T$   is a graded space of finite type, then $S$ is closed.
 
\end{enumerate}\end{lem}

\vspace{3mm}\noindent {\sl proof.} (i) It is sufficient to show that if $S$ and $R$ are closed then $S+R$ is closed. By Lemma \ref{l2.2},
$$\renewcommand{\arraystretch}{1.7}
\begin{array}{ll}
\overline{S+R}/\overline{S} &= \displaystyle\varprojlim_\alpha \rho_\alpha (S+R) / \rho_\alpha (S) = \displaystyle\varprojlim_\alpha \displaystyle\frac{\rho_\alpha (S) + \rho_\alpha (R)}{\rho_\alpha (S)}
\\
& = \displaystyle\varprojlim_\alpha \displaystyle\frac{\rho_\alpha (R)}{\rho_\alpha (S) \cap \rho_\alpha (R)}\\
& = \displaystyle \overline{R} / \displaystyle \varprojlim_\alpha \rho_\alpha (S) \cap \rho_\alpha (R)
\end{array}
\renewcommand{\arraystretch}{1}$$
In particular, it follows that the inclusions of $\overline{S}$ and $\overline{R}$ in $\overline{S+R}$ define a surjection $\overline{S}+ \overline{R}\to \overline{S+R}$. Since $S$ and $R$ are closed, $S+R\to \overline{S+R}$
is an isomorphism.

\vspace{3mm} (ii) Suppose $x\in \overline{\cap_\sigma S_\sigma}$, where each $S_\sigma$ is a closed subspace of $T$. Then $$\rho_\alpha x\in \rho_\alpha (\cap S_\sigma) \subset \cap_\sigma \rho_\alpha (S_\sigma).$$
Therefore $\rho_\alpha x \in \rho_\alpha (S_\sigma)$ for each $\sigma$. Thus $x = (\rho_\alpha x)\in \varprojlim_\alpha \rho_\alpha (S_\sigma) = \overline{S_\sigma} = S_\sigma$. Thus $x\in \cap_\sigma S_\sigma$ and $\cap_\sigma S_\sigma = \overline{\cap_\sigma S_\sigma}$. 

\vspace{3mm} (iii) If $x\in \overline{S}$ then $\rho_\alpha x\in \rho_\alpha (S)$. Thus for some $x_\alpha \in S$, 
$$\rho_\alpha (x-x_\alpha) = 0$$
and so $x-x_\alpha \in \mbox{ker}\, \rho_\alpha$. Therefore $x\in S+\mbox{ker}\, \rho_\alpha$ for all $\alpha$.

In the reverse direction, suppose $x\in \cap_\alpha (S+\mbox{ker}\, \rho_\alpha)$. Then
$$\rho_\alpha x\in \rho_\alpha S.$$
Thus the coherent family $x = (\rho_\alpha x) \in \overline{S}$. 

\vspace{2mm} (iv) It follows from (i) that each $S_k$ is closed and therefore that $S$ is closed.

  \hfill$\square$

\vspace{3mm}\noindent {\bf Remark.} Suppose $(x_k)_{k\geq 1}$ is a sequence of elements in an enriched vector space, $(T, \{\rho_\alpha\})$. If for each $\alpha $ there is some $r(\alpha)$ such that
$$\sum_{k=r(\alpha) + 1}^m \, x_k \in \mbox{ker}\, \rho_\alpha, \hspace{1cm} m>r(\alpha),$$
then for $m>r(\alpha)$ the elements $\rho_\alpha (\sum_{k=1}^{m} x_k) \in \rho_\alpha (T)$ define a single element $y\in \overline{T}$,  
 and we write $y = \sum_k x_k.$

\vspace{3mm} \noindent {\bf Example.} A vector space can have enriched structures with non isomorphic completions. For instance let $V = \mathbb Q[[x]] = \varprojlim_n \mathbb Q[x]/x^n$. Then $V$ is complete for the ideals generated by the projections $\rho_n : V\to \mathbb Q[x]/x^n$. Denote by $\omega= \sum a_n x^n$ be a series with $a_n\neq 0$ for all $n$,   let $S\subset V$ be a direct summand of $\mathbb Q \omega$ in $V$, and let $\varepsilon : V\to \mathbb Q \omega$ be the projection with kernel $S$. 
 
 The morphisms $\sigma_n = \varepsilon \oplus \rho_n : V\to \mathbb Q \oplus \mathbb Q[x]/x^n$ define an enriched structure on $V$, that is not complete because the coherent sequence $(0, \rho_n(\omega))$ is not in $V$. 
 
\subsection{Weight decompositions}

\vspace{3mm}\noindent   A \emph{weight decomposition} in a graded Lie algebra is a decomposition
$$L= \oplus_{k\geq 0} L(k)$$
in which $[L(k), L(\ell)] \subset L(k+\ell)$. A \emph{weighted subspace} of $L$ is a subspace $S\subset L$ such that $S = \oplus_k S\cap L(k)$. We denote $S(k) = S\cap L(k)$.

A \emph{weighted enriched Lie algebra} is a weighted Lie algebra with a defining set of ideals of the form $I_\alpha = \oplus_k I_\alpha (k)$. 

\vspace{2mm} \emph{Henceforth in this example, $L= \oplus_{k\geq 1}L(k)$ denotes a fixed weighted   enriched Lie algebra with a defining set of weighted ideals $I_\alpha$.}

\vspace{2mm} There follows the

\begin{Prop} \label{p2.2}
With the hypotheses and notation above,  
\begin{enumerate}
\item[(i)] Each $L(k)$ is closed and 
$\overline{L} = \prod_{k}  {L(k)}.$
\item[(ii)] If $S\subset L$ is a weighted subspace then $\overline{S} = \prod_k \overline{S(k)}$. In particular, $\prod_k S(k)$ is closed if and only if each $S(k)$ is closed, in which case any subspace of the form $\prod_i S(k_i)$ is closed.
\item[(iii)] If each $L(k)$ is finite dimensional, then  $\overline{L} = \varprojlim_r L/\, \prod_{k>r} L(k).$
\end{enumerate}\end{Prop}

 \vspace{3mm}\noindent {\sl proof.} (i) Denote as usual by $\rho_\alpha : L\to   L_\alpha= L/I_\alpha$ the projection on the finite dimensional Lie algebra $L_\alpha$. Then
 $$L_\alpha = \oplus_k \rho_\alpha (L(k)) = \prod_k \rho_\alpha (L(k)),$$
and
 $$L=\varprojlim_\alpha L_\alpha = \prod_k \lim_\alpha \rho_\alpha (L(k)) = \prod_k \varprojlim_\alpha  L(k)/I_\alpha (k) = \prod_k \overline{L(k)}.$$
 
 (ii) Similar proof than the one for (i), and (iii) is a direct consequence. 
 
 \hfill$\square$

\section{Lower central series}
\emph{Throughout this entire section $(L, {\mathcal I})$ denotes a fixed complete enriched Lie algebra, with quotient maps $\rho_\alpha: L\to L/I_\alpha:=L_\alpha$.} 

\vspace{2mm} A key invariant for any graded Lie algebra $E$ is its lower central series
$$E= E^1\supset \dots \supset E^k\supset \dots$$
of ideals in which $E^k$ is the linear span of iterated commutators of length $k$ of elements in $E$. 
The \emph{classical completion} of $E$ is the inverse limit $$\widehat{E}=\varprojlim_n E/E^n,$$ and $E$ is \emph{pronilpotent} if $E\stackrel{\cong}{\to} \widehat{E}$.

\vspace{3mm}The analogue of the lower central series for   $L$ is the sequence of ideals
$$L=L^{(1)}\supset \dots \supset L^{(k)}\supset \dots$$
defined by $L^{(k)}= \overline{L^k}$. It is immediate that $L\leadsto L^{(k)}$ is a functor. Moreover, since $L^k\subset L^{(k)}$, \emph{each $L/L^{(k)}$ is a nilpotent Lie algebra : $\left(L/L^{(k)}\right)^k = 0$.}

\begin{lem}
\label{l3.1}  Let $L$ be a complete enriched Lie algebra. Then
$$L = \varprojlim_n L/L^{(n)}.$$
More precisely,
\begin{enumerate}
\item[(i)] $ L^{(k)} = \varprojlim_\ell L^{(k)}/L^{(k+\ell)}.$ 
\item[(ii)] $[L^{(k)}, L^{(\ell)}] \subset L^{(k+\ell)}.$
\item[(iii)] $L$ is a retract of $\widehat{L} =\varprojlim_n L/L^n$.
\end{enumerate}
\end{lem}

\vspace{1mm}\noindent {\sl proof.} (i) Since each $L_\alpha$ is nilpotent, for some $\ell = \ell (\alpha)$, $L_\alpha^{k+\ell}= 0$. Thus by Lemma \ref{l2.2},
$$\varprojlim_\ell \overline{L^k}/\overline{L^{k+\ell}} = \varprojlim_{\ell, \alpha} L_\alpha^k/L_\alpha^{k+\ell}= \varprojlim_\alpha L_\alpha^k = \overline{L^k}.$$

(ii) From Lemma \ref{l2.3} we obtain $\rho_\alpha (L^{(k)})= \rho_\alpha (\lim_{\beta \geq \alpha} \rho_\beta (L^k))= \rho_\alpha (L^k) = L_\alpha^k$. It follows that
$$\rho_\alpha [L^{(k)}, L^{(\ell)}] \subset [L_\alpha^k, L_\alpha^\ell] \subset L_\alpha^{k+\ell},$$
and so $[L^{(k)}, L^{(\ell)}] \subset \overline{[L^{(k)}, L^{(\ell)}]} \subset L^{(k+\ell)}.$
 
 (iii) For each $n$ and $\alpha$, the surjection $L\to L_\alpha/L_\alpha^n$ (induced from $\rho_\alpha$) factors as
 $$L\to L/L^n \to L_\alpha/L_\alpha^n.$$
 Thus we obtain
 $$L\to \varprojlim_n L/L^n \to \varprojlim_{n,\alpha} L_\alpha /L_\alpha^n = \varprojlim_\alpha L_\alpha = L,$$
 which decomposes id$_L$ as $L\stackrel{\varphi}{\to} \widehat{L}\stackrel{\psi}{\to} L.$
 
 \hfill$\square$
 
\vspace{3mm}\noindent {\bf Corollary.} A morphism $E\to F$ of complete enriched Lie algebras is an isomorphism if and only if it induces isomorphisms $E/E^{(n)} \stackrel{\cong}{\to} F/F^{(n)}$, $n\geq 2.$

\vspace{3mm}\noindent {\bf Remark.} While $L$ may not be pronilpotent, Lemma \ref{l3.1} identifies the sequence $L\supset \dots \supset L^{(k)}\supset \dots$ as an $N$-suite as defined by Lazard in \cite{La}.

\vspace{3mm}Finally, the classical lower central series for a group $G$ is the sequence $(G^n)$ of subgroups generated by iterated commutators $aba^{-1}b^{-1}$ of length $n$. In analogy with the series $L^{(n)}$ we define
$$G_L^{(n)} = \varprojlim_\alpha G_\alpha^n,$$
where $G_\alpha$ is the group in $\widehat{U}L_\alpha$ defined in \S 1.2. In particular it follows that $\exp$ and $\log$ restrict to bijections
$$\xymatrix{L^{(n)}\ar@/^/[rr]^\exp &&G_L^{(n)}\ar@/^/[ll]^\log}$$
since these are established for each $L_\alpha, G_\alpha$ in \cite[Theorem 2.2]{RHTII}. Then in view of Proposition \ref{p1.1} these factor to yield bijections
$$L^{(n)}/L^{(k)} \stackrel{\cong}{\longrightarrow} G_L^{(n)}/G_L^{(k)} $$
since the corresponding bijections for each $L_\alpha$ are established in \cite[Proposition 2.5]{RHTII}. In particular (\cite[Corollary 2.4]{RHTII}) for each $n$, the bijection
$$L^{(n)}/L^{(n+1)} \stackrel{\cong}{\longrightarrow} G_L^{(n)}/G_L^{(n+1)}$$
is an isomorphism of rational vector spaces.

 \vspace{3mm}\noindent {\bf Example.} Let $\{ L^{(n)}\}_{n\geq 0}$ denote the lower central series of a complete enriched Lie algebra $L$. Then (as follows from Lemma \ref{l3.1}) the direct sum
 $$g(L) := \oplus_{n\geq 0} L^{(n)}/L^{(n+1)}$$
 is an enriched weighted Lie algebra, with Lie bracket $L^{(n)}/L^{(n+1)} \times L^{(k)}/L^{(k+1)} \to L^{(n+k)}/L^{(n+k+1)}$ induced from the Lie bracket in $L$.

 \vspace{2mm}
\begin{lem}
\label{l3.2} Suppose $E\subset L$ is a sub Lie algebra.
\begin{enumerate}
\item[(i)] If $E+ L^{(2)} = L$, then $\overline{E}= L$.
\item[(ii)]   $\overline{E^n} = \overline{E}^{(n)}$, $n\geq 1$. \end{enumerate}\end{lem}

\vspace{3mm}\noindent {\sl proof.} (i) By hypothesis $\rho_\alpha L = \rho_\alpha E+ \rho_\alpha L^{(2)}= \rho_\alpha E + [\rho_\alpha L, \rho_\alpha L]$. Since $\rho_\alpha L$ is nilpotent this gives $\rho_\alpha E = \rho_\alpha L$ and $\overline{E}= L$.

(ii) Since $E\subset \overline{E}$, $\rho_\alpha (\overline{E})= \rho_\alpha (E)$. Therefore
$$\renewcommand{\arraystretch}{1.4}
\begin{array}{ll} \overline{E^n} & = \varprojlim_\alpha \rho_\alpha (E^n) = \varprojlim_\alpha (\rho_\alpha (E))^n\\ & = \varprojlim_\alpha (\rho_\alpha (\overline{E})^n = \varprojlim_\alpha \rho_\alpha (\overline{E}^n)= \overline{E}^{(n)}.\end{array}\renewcommand{\arraystretch}{1}$$ \hfill$\square$

\begin{Prop}
\label{p3.1} Suppose $L$ is a complete enriched graded Lie algebra. If each dim$\, (L/L^{(2)})_p<\infty$, then
\begin{enumerate} 
\item[(i)] For any $p$ and any integer $k$ there is an $\alpha$ such that the projection $(L/L^{(k)})_p\to (L_\alpha /L_\alpha^k)_p$ is an isomorphism.
\item[(ii)] $L^k= L^{(k)}, \hspace{3mm}  k\geq 1$, and so $L$ is pronilpotent.
\end{enumerate}
\end{Prop}

\vspace{3mm}\noindent {\sl proof.} (i).  By Lemma \ref{l2.2} there is an $\alpha_0$, such that for each $\beta\geq \alpha_0 $  the projection $p_\beta : L/L^{(2)}\to L_\beta/L_\beta^2$ is an isomorphism in degree $p$. This shows that for each $\alpha $ the dimension of $L_\alpha/L_\alpha^2$ in degree $p$ is bounded by the integer $N = \mbox{dim}\, (L_{\alpha_0}/L_{\alpha_0}^2)_p$. It follows that for any $r\geq 1$, and any $\alpha$ the dimension of $L_\alpha^r/L_\alpha^{r+1}$ in degree $p$ is bounded by $N^r$. Thus it follows from Lemma \ref{l2.2} that   $L^{(r)}/L^{(r+1)}= \varprojlim_\alpha L_\alpha^r/L_\alpha^{r+1}$ is finite dimensional in degree $p$. Therefore for any integer $k$ there is an index $\alpha_k$ such that the projection 
$$(L/L^{(k)})_p\to (L_{\alpha_k}/L_{\alpha_k}^k)_p$$
is an isomorphism.

(ii) 
It follows from Lemma \ref{l2.2} that, given $p$, for some $\alpha$ the Lie bracket in $L$  
\begin{eqnarray}\label{i4}
L/L^{(2)} \otimes L^{(r)}/L^{(r+1)} \to L^{(r+1)}/L^{(r+2)}, \hspace{1cm} r\geq 1,
\end{eqnarray}
in degree $p$  may be identified with the corresponding surjection
$$L_\alpha/L_\alpha^2\otimes L_\alpha^r/L_\alpha^{r+1} \to L_\alpha^{r+1} / L_\alpha^{r+2}.$$
This implies that the Lie bracket (\ref{i4}) is   surjective in each degree $p$.

Now let $x_1, \dots , x_m\in L$ represent a basis of  $(L/L^{(2)})_{\leq p}$. Then for any
$y\in L^{(r+1)}_p$ it follows   that there are elements $y_i\in L^{(r)}$ such that
$$y-\sum_{i=1}^r [x_i, y_i]\in L^{(r+2)}.$$
Given $x\in L^{(n)}$, this yields an inductive construction  of elements $y_i(\ell)\in L^{(\ell-1)}$, $\ell\geq k$, such that
$$x-\sum_{i=1}^m\sum_{\ell = k-1}^r [x_i, y_i(\ell)]\in L^{(r)}.$$
Set $y_i=\sum_\ell y_i(\ell)$. Then by construction
$$x= \sum_{i=1}^m [x_i, y_i].$$

This shows that $L^{(k)}\subset [L, L^{(k-1)}]$. Induction on $k$ gives
$$L^{(k)}\subset L^k,$$
and the reverse inclusion is obvious.

\hfill$\square$

\vspace{3mm}\noindent {\bf Corollary 1.}  Let $E$ be an enriched Lie algebra with dim$\, E/E^2<\infty$. Then dim$\,\overline{E}/\overline{E^k}<\infty$, $k\geq 1$, and 
the lower central series for $\overline{E}$ satisfies
$$\overline{E}^{(k)}= \overline{E}^k = \overline{E^k}.$$ 

\vspace{3mm}\noindent {\sl proof.} There is an 
$\alpha_0$ such that for each $\alpha \geq \alpha_0$, the projection $E_\alpha /E_\alpha^2\to E_{\alpha_0}/E_{\alpha_0}^2$ is an isomorphism. Thus by Lemma \ref{l2.2} $\overline{E}/\overline{E}^{(2)}$  is isomorphic to $  E_{\alpha_0}/E_{\alpha_0}^2$, and so is   finite dimensional. Proposition \ref{p3.1} gives then an isomorphism $\overline{E}/\overline{E}^{(2)}= \overline{E}/\overline{E}^2=E_{\alpha_0}/E_{\alpha_0}^2$ In particular, $E+\overline{E}^{(2)}= \overline{E}$. 
The result is then a direct consequence of Proposition \ref{p3.1}(ii) and  Lemma \ref{l3.2}(ii).  

\hfill$\square$

\vspace{3mm}\noindent {\bf Corollary 2.}  An complete enriched Lie algebra $L$ is the direct limit of closed   sub Lie algebras, $E$, satisfying dim$\, E/E^2<\infty$.

\vspace{3mm}\noindent {\sl proof.} Any Lie algebra is the direct limit of its finitely generated Lie algebras, $F$. Thus $L$ is the direct limit of the completions, $\overline{F}$. By Proposition \ref{p3.1} and Corollary 1, each $F$   satisfies dim$\, \overline{F}/\overline{F}^2<\infty$. 

\hfill $\square$

\vspace{3mm}Proposition \ref{p3.1} has the following analogue:

\begin{Prop}
\label{p3.2} If $E\subset L$ is a sub Lie algebra and $E/E^2$ is a graded vector space of finite type then for each $k$
\begin{enumerate}
\item[(i)]$E/E^k \stackrel{\cong}{\to} \overline{E}/\overline{E^k}$.
\item[(ii)] $\overline{E}^k = \overline{E^k} = \overline{E}^{(k)}$.
\item[(iii)] $\widehat{E}= \overline{E}= \widehat{\overline{E}}$.
\end{enumerate}
In particular, $\overline{E}$ is pronilpotent.
\end{Prop}

\vspace{3mm}\noindent {\sl proof:} (i) Since $E/E^2$ has finite type each $E/E^k$ has finite type. Thus (i) follows from Lemma \ref{l2.2}.

(ii) This follows exactly as in Corollary 1 to Proposition \ref{p3.1}.

(iii) This follows immediately from (i) and (ii). \hfill$\square$

\begin{Prop}
\label{p3.3} 
 Suppose $L$ is a complete enriched Lie algebra such that
 $L/L^2 $  has finite type.
Then $L$ has a unique (up to equivalence) structure as an enriched Lie algebra $(L, {\mathcal I})$, namely that given by
$${\mathcal I}= \{I_k\}_{k\geq 1} \hspace{5mm}\mbox{with } I_k = L^k+ L_{\geq k}.$$
\end{Prop}

\vspace{3mm}\noindent {\sl proof:} Suppose a second family ${\mathcal J}= \{J_\beta\}$ of ideals in $L$ also makes $L$ into an enriched Lie algebra. Since $L/J_\beta$ is finite dimensional and nilpotent, it is also immediate that $J_\beta \supset I_k$ for some $k$.

In the reverse direction, apply Proposition \ref{p3.1}(i) to conclude that for fixed $k$ and $n$ there is some $\beta = \beta (k,n)$ such that 
$$\rho_\beta : L_n /(L^k)_n \stackrel{\cong}{\longrightarrow} (L_\beta)_n /(L_\beta^k)_n.$$
It follows that 
$$(L^k)_n = \rho_\beta^{-1}(L_\beta^k)_n \supset (J_\beta)_n.$$
Choose $\gamma \in {\mathcal J}$ so that 
$$J_\gamma \subset \cap_{n=1}^k J_{\beta (k,n)}.$$
Then $$(L^k)_{\leq n} \supset (J_\gamma)_{\leq n},$$
and so $$I_k\supset J_\gamma.$$ \hfill$\square$

 \vspace{3mm}\noindent {\bf Corollary.} Let $L$ be a complete enriched Lie algebra and $E\subset L$ be a subalgebra such that $E/E^2$ has finite type. Then $E$ is closed if and only if $E$ is pronilpotent.
 
 \vspace{3mm}\noindent {\sl proof.} By Proposition \ref{p3.3}, the structure systems $\rho_\alpha E$ and $E/E^n$ are equivalent. This implies the result.\hfill$\square$

\vspace{3mm}\noindent {\bf Example.} Suppose a graded Lie algebra $E$ satisfies
$$\cap_k E^k= 0\,, \hspace{3mm}\mbox{dim}\, E/E^2<\infty, \hspace{3mm}\mbox{and dim}\, E= \infty.$$
Let $z\in \overline{E}$ satisfy $E\oplus \mathbb Q z \subset \overline{E}$ and let $F$ be the sub Lie algebra generated by $E$ and $\mathbb Q z$. Then $E\subset F\subset \overline{E}$, and it follows that $\overline{E}= \overline{F}$. However the dimension of $F/F^2$ depends on the choice of $z$. 

For instance, suppose $E$ is the completion of the free Lie algebra $\mathbb L(x,y)$ generated by $x,y$, and with the unique enriched structure described in Proposition \ref{p3.3}. Let $F$ be the sub Lie algebra generated by $x$, $y$ and $z= e^{ad\, x}(y)$. Denote by $t$ the operator $\mbox{ad}_x$. The series $t$ and $e^t$ are algebraically independent in $\mathbb Q[[t]]$, i.e., if $P(a,b)$ is a polynomial in 2 variables $a$ and $b$ and $P(t,e^t)= 0$ then $P$ is identically zero. It follows that if we have a polynomial relation of the form
 $$\sum_n \alpha_n \,\mbox{ad}_x^n (y) + \sum_n \beta_n \,\mbox{ad}^n_x\, e^{\mbox{ad}_x}(y) = 0$$
 then $(\sum \alpha_nt^n + \sum\beta_n t^n e^t) (y)= 0$, and $\alpha_n=\beta_n= 0$ for all $n$.   Therefore the classes of $x,y$ and $z$ are linearly independent in $F/F^2$, and $F/F^3$ has   dimension 3.  
 
 On the other hand, with the same $E$, if $z= \sum_{n\geq 0} \mbox{ad}^n_x(y)$, then $z-y= [x,z]$ and $F/F^2 = E/E^2$.

\section{The quadratic model of $(L, {\mathcal I})$}

\emph{In this section $(L, {\mathcal I})$ denotes a complete enriched Lie algebra. We denote
$${\mathcal I} = \{\rho_\alpha \, : \, L\to L_\alpha\},$$
so that $L= \varprojlim_\alpha L_\alpha$. }

\vspace{3mm}
On the other hand, a \emph{quadratic Sullivan algebra} is a Sullivan algebra $(\land V,d)$ in which $d : V\to \land^2V$. This endows $H(\land V)$ with the decomposition
$$H(\land V) = \oplus_k H^{[k]}(\land V)$$
in which $H^{[k]}(\land V)$ is the subspace represented by $\land^kV\cap \mbox{ker}\, d$. Moreover, $V = \varinjlim_n V_n$, where 
$$V_0= V\cap \mbox{ker}\, d\hspace{3mm}\mbox{and } V_{n+1}= V\cap d^{-1}(\land^2V_n).$$
Then a \emph{morphism} $\varphi : \land V\to \land W$   of quadratic Sullivan algebras is a cdga morphism which restricts to a linear map $V\to W$. Thus it restricts to linear maps 
$$\varphi_n : V_n \to W_n.$$
The filtration $(V_n)_{n\geq 0}$ will be called the \emph{standard filtration}.

On the other hand, recall that
the classical functor $E \leadsto C_*(E)$ from graded Lie algebras to cocommutative chain coalgebras is given by $C_*(E)= \land sE$, with differential determined by the condition
$$\partial  (sx\land sy) = (-1)^{deg\, x+1} s[x,y].$$
The dual, $C^*(E)= \mbox{Hom}(C_*(E), \mathbb Q)$ is, with the differential forms on a manifold, one of the earliest examples of a  cdga. In particular, if dim$\, E<\infty$ then
$C^*(E) = \land (sE)^\vee.$ 

Now set $V_\alpha = (sL_\alpha)^\vee$. Then $C_*(L_\alpha)= (\land V_\alpha, d_\alpha)$ and, since $L_\alpha$ is finite dimensional and nilpotent, $(\land V_\alpha, d_\alpha)$ is a quadratic Sullivan algebra in which 
\begin{eqnarray}\label{i5}
<d_\alpha v, sx,sy> = (-1)^{deg\, y +1} <v, s[x,y]>, \hspace{1cm} v\in V_\alpha, x,y\in L_\alpha.
\end{eqnarray}
Moreover, the $\{V_\alpha\}$ are the inductive system whose dual is the inverse system $\{sL_\alpha\}$ (Proposition \ref{p2.1}) and by that Proposition 
\begin{eqnarray}
\label{i6}\renewcommand{\arraystretch}{1.6}\left.\begin{array}{c}  sL =   \varprojlim_\alpha sL_\alpha =\varprojlim_\alpha V_\alpha^\vee = V^\vee, \hspace{2mm}\mbox{and }\\
<dv, sx,sy> = (-1)^{1+ deg\, y} <v, s[x,y]>, \hspace{3mm} v\in V, x,y\in L .\end{array}\right\}\renewcommand{\arraystretch}{1}\end{eqnarray}

\vspace{3mm}\noindent {\bf Definition.} $\land V$ is the quadratic model of $(L, {\mathcal I})$.

\vspace{3mm}\noindent {\bf Remarks.} 1. It follows from (\ref{i5}) that $L$ is the homotopy Lie algebra (\cite[\S2.1]{RHTII}) of $\land V$.

2. It follows (Proposition \ref{p2.1}) from (\ref{i6}) that $V$ is the predual of $L$ so that the enriched structure of $L$ is defined by the isomorphism $sL = V^\vee$. 

3.  The natural morphism $\land V= \varinjlim C^*(L_\alpha) \to C^*(L)$ will not, in general, be a quasi-isomorphism. For instance, if $V= V^1$ has zero differential and countably infinite dimension then $\land V$ is countable but $L$ and $C^*(L)$ are uncountable, and since $L$ is abelian, $C^*(L)= H(C^*(L))$. 

 \vspace{3mm}\noindent {\bf Example. Weighted Lie algebras.}  Let $E= \oplus_k  E(k)$ be a weighted enriched Lie algebra with defining ideals $I_\alpha = \oplus_k I_\alpha (k)$. The weighting then induces another gradation $V = \oplus_k V(k)$ in the generating space of the corresponding quadratic model:
 $$V(k) = \varinjlim_\alpha s\left( E(k)/I_\alpha(k)\right)^\vee.$$
 It is immediate that the differential preserves the induced gradation in $\land V$.

\vspace{2mm}
\begin{Prop}
\label{p4.1} (i) The correspondence $(L, {\mathcal I})\leadsto \land V$ is a contravariant isomorphism from the category of complete enriched Lie algebras to the category of quadratic Sullivan algebras.

(ii) Two morphisms $\varphi, \psi : \land V\to \land W$ of quadratic Sullivan algebras are based homotopic if and only if $\varphi = \psi$ or, equivalently, if they induce the same morphism between the homotopy Lie algebras.
\end{Prop}

\vspace{3mm} \noindent {\sl proof.} It follows from Proposition \ref{p2.1} that coherent linear maps $\varphi : E\to L$ between complete enriched Lie algebras are exactly the duals of linear maps $\psi : V_E\leftarrow V_L$ between the preduals. It follows from (\ref{i5}) that $\varphi$ is a morphism of Lie algebras if and only if $\psi$ extends to a morphism $\land V_E\leftarrow \land V_L$ of quadratic Sullivan algebras. This establishes (i), and (ii) is straightforward from the definitions. \hfill$\square$

 \vspace{3mm}\noindent {\bf Remark.} Proposition \ref{p4.1} couples the categories of complete enriched Lie algebras, $L$, and quadratic Sullivan algebras, $\land V$, as pairs $(L, \land V)$ where $\land V$ is the quadratic Sullivan model of $L$ and, equivalently, $L$ is the homotopy Lie algebra of $\land V$. 
 
\vspace{3mm}  The correspondence of Proposition \ref{p4.1} is reflected in the next Lemma.

\begin{lem}
\label{l4.1} (i) The degree 1 identification $L\stackrel{\cong}{\to} V^\vee$ induces isomorphisms
$$L/L^{(n+2)} \stackrel{\cong}{\longrightarrow} V_n^\vee, \hspace{1cm} n\geq 0,$$
and therefore identifies $L^{(n+2)} = \{x\in L\, \vert\, <V_p, sx>= 0\}.$

(ii) If $T \oplus L^{(2)}= L$ then $L$ is the closure of the sub Lie algebra, $E$, generated by $T$. \end{lem}

\vspace{3mm}\noindent {\sl proof.} (i)  This is straightforward when $L$ and $V$ are replaced by $L_\alpha$ and $V_\alpha$. In the general case by Lemma \ref{l2.2} we have the following sequence of isomorphisms
$$
\renewcommand{\arraystretch}{1.6}
\begin{array}{ll} L/L^{(n+2)} &= \overline{L}/\overline{L^{n+2}}= \displaystyle\varprojlim_\alpha \rho_\alpha (L)/ \rho_\alpha (L^{n+2})= \displaystyle\varprojlim_\alpha L_\alpha / L_\alpha^{n+2}\\
& = \displaystyle \varprojlim_\alpha  {(V_\alpha)_n}^\vee = \left( \displaystyle\varinjlim_\alpha  (V_\alpha)_n\right)^\vee = V_n^\vee.
\end{array}
\renewcommand{\arraystretch}{1}
$$

(ii) The inclusion $sE \to sL$ is the dual of a surjection $W \stackrel{\varphi}{\longleftarrow} V$ where $\land W$ is the quadratic model of $E$. By construction $E/E^{(2)} \stackrel{\cong}{\longrightarrow} L/L^{(2)}$, and so $W_0 \stackrel{\cong}{\longleftarrow} V_0$. If $\xymatrix{ W_n & V_n\ar[l]^\cong_\varphi}$ and $v\in V_{n+1}$ then $\varphi v= 0$ $\Rightarrow$ $dv=0$ $\Rightarrow$ $v\in V_0$ $\Rightarrow$ $v=0$. 
\hfill$\square$

\vspace{3mm} \noindent {\bf Remark.}  A morphism $\psi : E\to L$ of complete enriched Lie algebras induces morphisms $\psi(n) : E/E^{(n)}\to L/L^{(n)}$, $n\geq 1$. If $\varphi : \land V_L\to \land V_E$ is the corresponding morphism of quadratic models then $\psi(n+2)$ is dual to to the linear map $\varphi_n : (V_L)_n \to (V_E)_n$.

\begin{lem}
\label{l4.2} The closure of any subspace $S\subset L$ is the image of the inclusion $(V/K)^\vee \to V^\vee = L$, where 
$$K = \{ v\in V \, \vert\, <v,sS>= 0\}.$$
In particular, any closed subspace of $L$ has a closed direct summand.
\end{lem}

\vspace{3mm}\noindent {\sl proof.} For $x\in L$ and $v\in V_\alpha$  we have
$<v, sx> = <v, s\rho_\alpha x>.$
 Thus 
$$V_\alpha \cap K = \{v\in V_\alpha \, \vert \, <v,s\rho_\alpha S>= 0\}.$$
Set $P_\alpha = V_\alpha /V_\alpha \cap K$. Then, because $V_\alpha$ and $L_\alpha$ are finite dimensional, this gives
$$\rho_\alpha S = \{x\in L_\alpha \, \vert\, <V_\alpha \cap K, sx> = 0\} = (sP_\alpha)^\vee.$$
Thus denoting $V\cap K= \varinjlim_\alpha P_\alpha = P$ we obtain
 $$\overline{S}= \varprojlim_\alpha \rho_\alpha S = \varprojlim_\alpha (sP_\alpha)^\vee = (\varinjlim_\alpha sP_\alpha)^\vee = (sP)^\vee.$$ 

Finally the inclusion $S\to L$ is the dual of a surjection $V\to V/K$. Dualizing the inclusion $K\to V$ provides a surjection $L \to (sK)^\vee$ onto a closed direct summand of $S$. 
\hfill$\square$

\subsection{Closed subalgebras and ideals}
 
\vspace{1mm} Suppose $E\subset L$ is a closed subalgebra of a complete enriched Lie algebra. Then the inclusion corresponds via duality to a surjection $\xymatrix{\land V \ar[r]^\rho & \land Z}$ of the respective quadratic models, and by Lemma \ref{l4.2}, 
 $$E = \{x\in L\,\vert\, <\mbox{ker}\, \rho_{\vert V}, sx>= 0\}.$$
 Moreover (\cite[Lemma 10.4]{RHTII}), $d: \mbox{ker}\, \rho_{\vert V}\to \mbox{ker}\, \rho_{\vert V}\, {\scriptstyle \land}\, V$. Recall also that every such surjection induces an inclusion of a closed sub algebra of $L$.
 
 Now suppose $I\subset L$ is a closed ideal. Let $\rho : \land V\to \land Z$ be the corresponding surjection, and denote $W = \mbox{ker}\, \rho_{\vert V}$. In this case (\cite[Lemma 10.4]{RHTII}) the condition that $I$ be an ideal is equivalent to the condition $d: W\to \land^2W$. Thus $\land W$ is a sub quadratic Sullivan algebra and $I$ decomposes $\land V$ as the Sullivan extension
 $$\xymatrix{\land W\ar[r]^-\lambda & \land W\otimes \land Z = \land V \ar[r]^-\rho & \land V\otimes_{\land W}\mathbb Q = \land Z,}$$
 which dualizes to
 $$L/I \leftarrow L\leftarrow I.$$
 In particular this identifies $L/I$ as the homotopy Lie algebra of $\land W$.

  \begin{Prop}
\label{p4.2} Let $\varphi : E\to L$ be a morphism of complete enriched Lie algebras.
\begin{enumerate}
\item[(i)] Both Ker$\, \varphi$ and Im$\, \varphi$ are closed subspaces, respectively in $E$ and in $L$.
\item[(ii)] If $\varphi$ induces a surjection $E/E^{(2)}\to L/L^{(2)}$ then $\varphi$ is surjective.
\end{enumerate}
\end{Prop}

\vspace{3mm}\noindent {\sl proof.}  (i) is a special case of Proposition \ref{p2.1} (iii).

(ii) Denote by $\psi : \land V_E\to \land V_L$ the corresponding morphism given by Proposition \ref{p4.1}. By Lemma \ref{l4.1}, $(V_E\cap \mbox{ker}\, d)^\vee \cong L/L^{(2)}$ and $(V_L\cap \mbox{ker}\, d)^\vee \cong E/E^{(2)}$. It follows that $\psi : V_E\cap \mbox{ker}\, d\to V_L\cap \mbox{ker}\, d$ is injective. We suppose that for some integer $n$, $\psi : (V_E)_n\to (V_L)_n$ is injective. Then $\psi : \land (V_E)_n\to \land (V_L)_n$ is injective. Suppose that $v\in (V_E)_{n+1}$ is in the kernel of $\psi$, then $\psi (dv)= 0$, and so $dv= 0$. Therefore $v\in (V_E)_0$ which implies that $v=0$. It follows that the dual map $\varphi$ is surjective.

\hfill$\square$

\vspace{3mm}\noindent {\bf Corollary.} 
  Let $\varphi : E\to L$ be a surjective coherent morphism between complete enriched abelian Lie algebras. Then $\varphi$ admits a coherent section, $\sigma$.

\vspace{3mm}\noindent {\sl proof.}  By Proposition \ref{p4.2}, Ker$\, \varphi$ is a closed subspace. Now by Lemma \ref{l4.2}(iii), Ker$\, \varphi\subset L$ admits a closed direct summand $S$. It follows that $\varphi_{\vert S} : S\to $ Im$\, \varphi$ is an isomorphism, and we define $\sigma = (\varphi_{\vert S})^{-1}$. \hfill$\square$

 \begin{Prop}
\label{p4.3}
Suppose $\{L(\sigma), \varphi_{\sigma, \tau}: L(\tau)\to L(\sigma)\}$ is an inverse system of morphisms of complete enriched Lie algebras, and suppose that either the morphisms are surjective or else that the system is finite. Then
$$L:= \varprojlim_\sigma L(\sigma)$$
is naturally a complete enriched Lie algebra.
\end{Prop}

\vspace{3mm}\noindent {\sl proof.} Suppose first that the morphisms are surjective. 
Under the correspondence of Proposition \ref{p4.1}, the inverse system is the dual of a directed system $\{\land V(\sigma)\}$ of   injections between the quadratic Sullivan models of the $L(\sigma)$. But then $\land V := \varinjlim_\sigma \land V(\sigma)$ is a quadratic Sullivan algebra, and $L$ is the homotopy Lie algebra of $\land V$. It is immediate that $L$ is the inverse limit in the category of complete enriched Lie algebras.

When the system is finite, we replace the system $L(\sigma)$ by the sub system formed by the $G(\sigma)$ where $G(\tau)$ is the intersection of the images of morphisms $L(\sigma)\to L(\tau)$, coming from all $L(\sigma)$.  By Proposition \ref{p4.2} the images are closed sub algebras, and by Lemma \ref{l2.2} their intersection is also a closed sub Lie algebra. Now by construction the induced maps between the $G(\sigma)$ are surjective. Since $\varprojlim_\sigma L(\sigma) = \varprojlim_\sigma G(\sigma)$, the inverse limit is a complete enriched Lie algebra.

\hfill$\square$

\vspace{3mm}\noindent {\bf Remark.}   ZFC - set theory alone does not permit us to extend Proposition \ref{p4.1} to all graded Lie algebras, since   it is consistent with the ZFC-axioms that the same graded Lie algebra can support two inequivalent enriched structures.

In fact, let $L=L_0$ be a vector space and suppose (consistent with the ZFC axioms, \cite{Solo}) that there are isomorphisms 
$$V^\vee \cong sL\cong W^\vee$$
in which card$\,V\neq$ card$\, W$.
Let $\land V$ and $\land W$ be the quadratic Sullivan algebras with zero differential. Then $L$, regarded as an abelian Lie algebra is the homotopy Lie algebra of both $\land V$ and $\land W$, but since $\land V \not\cong \land W$ the corresponding enriched structures in $L$ are not equivalent.

Remark that it is also consistent with ZFC to admit the   Generalized Continuum Hypothesis (GCH) with the consequence  that an isomorphism $V^\vee \cong W^\vee$ implies an isomorphism $V\cong W$  (\cite{Jech})

\section{Representations}
\emph{Throughout this entire section $(L, {\mathcal I})$ denotes a fixed complete enriched Lie algebra, with quotient maps $\rho_\alpha: L\to L/I_\alpha:=L_\alpha$  and quadratic model $\land V$.}

\vspace{3mm} First recall that the left and right adjoint representations of $L$ in $L$ are given by
\begin{eqnarray}
\label{i7}
\rm{ad}_y(x)= [y,x] \hspace{5mm}\mbox{and } \mbox{r}\,\mbox{ad}_y(x)= (-1)^{deg\, y} [x,y].
\end{eqnarray}
As well, if $L\otimes M\to M$ is any representation of $L$ ($M$ is a \emph{left $L$-module}) then the dual, $M^\vee$, of $M$ inherits the right $L$-module structure given by
\begin{eqnarray}
\label{i8}
<m, a\cdot y> = (-1)^{deg\, y}<y\cdot m,a>, \hspace{1cm} m\in M, a\in M^\vee, y\in L.
\end{eqnarray}

\vspace{3mm}\noindent {\bf Definition.} (i) An \emph{elementary $L$-module} is a finite dimensional left (resp. right) nilpotent $L$-module, $Q$, for which some $I_\alpha \cdot Q = 0$ (resp. $Q\cdot I_\alpha = 0$).

(ii) An \emph{enriched $L$-module} is an $L$-module, $N$, together with a decomposition,
$$N = \varprojlim_{\tau} N_\tau,$$
of $N$ as an inverse limit of elementary $L$-modules.

(iii) A \emph{locally elementary $L$-module} is an $L$-module, $M$, together with a decomposition,
$M=\varinjlim_\sigma M_\sigma,$
of $M$ as a direct limit of elementary $L$-modules.

(iv) A \emph{morphism} of enriched (resp. locally elementary) $L$-modules is a morphism of $L$-modules which is the inverse limit (resp. direct limit) of morphisms of elementary $L$-modules.

 \begin{Prop}
 \label{p5.1}
 \begin{enumerate}
 \item[(i)] The correspondence $M\leadsto M^\vee$ is a contravariant isomorphism from the category of locally elementary left $L$-modules to the category of enriched right $L$-modules.
 
 \item[(ii)] Both the left and right adjoint representations of $L$ in $L$ convert $L$ into an enriched $L$-module.
 \item[(iii)] If $M$ is either a locally elementary or an enriched $L$-module, then that representation of $L$ extends uniquely to a representation of $\overline{UL}$. 
 \item[(iv)] The isomorphism $sL= V^\vee$ identifies the right adjoint representation of $L$ with the dual of the locally elementary representation in $V$ given by
 $$<y\cdot v, sx> = -< dv, sx,sy>.$$
 \end{enumerate}
 \end{Prop}

 \vspace{3mm}\noindent {\sl proof.}. (i), (ii), and (iii) follows immediately from the definitions and the enriched structure $L=\varprojlim_\alpha L_\alpha$, in which $L_\alpha$ is finite dimensional and nilpotent. Finally, (iv) follows from this and (\ref{i5}). \hfill$\square$

\vspace{3mm}\noindent {\bf Example.}
Let $L$ be the free Lie algebra on one generator $x$ in degree $0$. The space $M$ of finite sequences of rational numbers $(a_1, \dots , a_n)$ equipped with the $L$-structure defined by $x\cdot (a_1, \dots , a_n)= (a_2, \dots a_{n})$ is a locally elementary $L$-module. On the other hand the space $N$ of infinite sequences $(a_1, a_2, \dots )$ with the $L$-structure defined by $ (a_1, a_2, \dots)\cdot x = (0, a_1, a_2, \dots )$ is an enriched $L$-module.

\subsection{Standard constructions for $L$-modules}

(1) If $M$ and $N$ are respectively a locally elementary and an enriched $L$-module then Hom$(M,N)$ is an enriched $L$-module.

(2) The sub and quotient modules of a locally elementary (resp. enriched module) is locally elementary (resp. enriched). (The second assertion follows from Proposition \ref{p5.1}.)

(3) The tensor product of locally elementary modules is locally elementary.

(4) The adjoint representation of $L$ in any closed ideal $I$ is an enriched representation.

\vspace{3mm}\noindent {\bf Remarks.}

1. When dim$\, L/[L,L]<\infty$, any finite dimensional nilpotent $L$-module is elementary, but this is not true in general. In fact, if $L= L_0$ is an infinite dimensional abelian Lie algebra with basis $\{x_i\}$, the subspaces $I(k)$ spanned by the $x_i$, $i\geq k$, define an enriched structure in $L$. Then,   the finite dimensional nilpotent $L$-module $M = \mathbb Q a \oplus \mathbb Q b$   defined by 
$$x_i\cdot a = b \hspace{3mm}\mbox{and } x_i\cdot b = 0 $$
is not an elementary $L$-module.

2. Since for each elementary left $L$-module $Q$ there is an ideal $I_\alpha$ with $I_\alpha \cdot Q= 0$, the module $Q$ is naturally a $\widehat{UL_\alpha}$-module for some $\alpha$. Thus the representations of $L$ in enriched   $L$-modules naturally extend to representations of the algebra $\overline{UL}$. 

3. The surjections $\overline{UL}\to \widehat{UL_\alpha}/I_\alpha^n$ make $\overline{UL}$ into a right enriched $L$-module under right multiplication.

4. Let $N$ be a right enriched $L$-module. The space  of \emph{decomposable elements} of $N$ is by definition $
\overline{N\cdot L} = \varprojlim_\alpha  N_\alpha \cdot L$.
For instance for $N=L$ considered as a right enriched module over itself via the adjoint representation, we have $\overline{N\cdot L} = L^{(2)}$. 
 
 \subsection{Augmentations} If $N$ is a right enriched $L$-module, the   trivial $L$-modules, $\widetilde{N}_\tau:= N_\tau \otimes_{\overline{UL}} \mathbb Q= N_\tau/N_\tau \cdot L$ define the trivial right enriched $L$-module,   $\widetilde{N}= \varprojlim_\tau \widetilde{N_\tau}$   (i.e., an enriched vector space). The induced morphism
 $$\varepsilon_N : N\to \widetilde{N}$$
 is a morphism of enriched $L$-modules: $\varepsilon_N$ is the \emph{augmentation} for $N$. By construction, $\varepsilon_N : N\to \widetilde{N}$ depends naturally on $N$. 
 
 Let $M$ and $\widetilde{M}$ be the locally elementary $L$-modules defined by $M^\vee = N$ and $\widetilde{M}^\vee = \widetilde{N}$. Then   $\widetilde{M} = \{m\in M\, \vert\, L\cdot m= 0\}.$

\subsection{Quadratic $(\land V,d)$-modules and holonomy representations}

\vspace{3mm}\noindent {\bf Definition.} (i) A \emph{quadratic} $(\land V,d)$-module is a $(\land V,d)$-module of the form $\land V\otimes M$ in which $d: 1\otimes M\to V\otimes M$ and $M= \cup_{k\geq 0} M_k$, where
$$M_0 = M\cap \mbox{Ker}\, d \hspace{3mm}\mbox{and } M_{k+1} = M \cap d^{-1}(V\otimes M_k).$$
A \emph{morphism} of quadratic $(\land V,d)$-modules is a morphism of the form $id\otimes \varphi : \land V\otimes M\to \land V\otimes M'$. 

Thus quadratic $(\land V,d)$-modules are examples of the semifree $(\land Vd)$-modules defined in \cite{FHTI}

(ii) The \emph{holonomy representation} for a quadratic $(\land V,d)$-module $\land V\otimes M$ is the left representation of $L$ in $M$ given by
$$x\cdot m = -\sum <v_i, sx> m_i,$$
where $d(1\otimes m) = \sum v_i\otimes m_i.$
In other words, the action $L\otimes M\to M$ is the composition
$$\xymatrix{L\otimes M\ar[rr]^{s\otimes d} && sL\otimes V\otimes M \ar[rr]^{<\,>\otimes id_M} && \mathbb Q\otimes M = M.}$$

\begin{lem}
\label{l5.1}
\begin{enumerate}
\item[(i)] The correspondence $(\land V\otimes M,d) \leadsto M$ is a natural isomorphism between the category of quadratic $\land V$-modules and locally elementary $L$-modules.
\item[(ii)] A morphism $(id\otimes \varphi) : \land V\otimes M\to \land V\otimes M'$ is a quasi-isomorphism if and only if $\varphi : M\to M'$ is an isomorphism
\end{enumerate}
\end{lem}

\vspace{3mm}\noindent {\sl proof.} (i).  It is immediate that the holonomy representation makes $M$ into a locally elementary left $L$-module. In the reverse direction, suppose $M= \varinjlim_\sigma M_\sigma$ is a locally elementary left $L$-module. The classical Cartan-Eilenberg-Serre construction then has the form
$$C^*(L_{\alpha (\sigma)}, M_\sigma) = \land (sL_{\alpha (\sigma)})^\vee \otimes M_\sigma = \land V_{\alpha (\sigma)} \otimes M_\sigma$$
in which $d : M_\sigma \to V_{\alpha (\sigma)}\otimes M_\sigma$. Passing to direct limits constructs the quadratic $\land V$-module $\land V\otimes M = \varinjlim_\sigma \land V_{\alpha (\sigma)}\otimes M_\sigma$.  

\vspace{2mm}(ii). The wedge gradation of $\land V$ makes a quadratic $\land V$-module into a chain complex, and if $id\otimes \varphi$ is a quasi-isomorphism,
$$M\cap \mbox{ker}\, d = H^{[0]}(\land V\otimes M)  
\stackrel{\cong}{\longrightarrow}  
H^{[0]}(\land V\otimes M')= M'\cap \,\mbox{ker}\, d\,.
$$

Now suppose by induction that $\varphi$ restricts to an isomorphism $M_k\stackrel{\cong}{\longrightarrow} M'_k$. If $x\in M_{k+1}$ and $\varphi x= 0$ then $(id\otimes \varphi) dx= 0$. Since $dx\in V\otimes M_k$ it follows that $dx=0$. Thus $x\in M\cap \mbox{ker}\,d$, and so $x=0$.  A similar argument shows that $\varphi : M_{k+1} \to M'_{k+1}$ is surjective.

\hfill$\square$

Finally, as recalled in \S 4, associated with $L$ is the differential coalgebra $(\land sL, \partial)$ in which $\partial (sx\land sy) = (-1)^{deg\, sx} s[x,y]$. More generally, (eg. \cite[Chap 2]{RHTII}) associated with $M$ is the differential $(\land sL, \partial)$-comodule, $(\land sL\otimes M, \partial)$, characterized by
$$\partial (sx\otimes m)= (-1)^{deg\, sx} x\cdot m.$$
As recalled in \cite[Chap 10]{RHTII}, there are natural isomorphisms
$$\mbox{Tor}_p^{UL}(\mathbb Q, M) \cong H_{[p]}(\land sL\otimes M), \hspace{5mm} p\geq 0,$$
where $H_{[p]}$ denotes the subspace of $H(\land sL\otimes M)$ represented by cycles in $\land^psL\otimes M$.

 \subsection{Acyclic closures} 

Recall (\cite[Chap 3]{RHTII}) that the acyclic closure   is the special case, $\land V\otimes \land U$,  of a $\Lambda$-extension constructed inductively as follows. Write $V = \cup_n V_n$ with $V_0= V\cap \mbox{ker}\, d$ and $V_{n+1}= V\cap d^{-1}(\land^2V_n)$. Then $U = \cup_{n\geq 0} U_n$ and there is a degree $1$ isomorphism $p: U\stackrel{\cong}{\to} V$ restricting to isomorphisms $U_n \stackrel{\cong}{\to} V_n$. Thus this identifies $L= U^\vee$. Finally, the differential is determined by the conditions
$$du=pu, \hspace{3mm} u\in U_0\hspace{3mm}\mbox{and } (d-p) : U_{n+1}\to V_n\otimes \land U_n.$$
In particular the augmentation $\varepsilon_U : \land U\to \mathbb Q$, $U\to 0$ together with the unique augmentation in $\land V$ define a quasi-isomorphism,
$$\land V\otimes \land U \stackrel{\simeq}{\longrightarrow} \mathbb Q.$$

\emph{This construction identifies $\land V\otimes \land U$ as a quadratic $(\land V,d)$-module.} In this case the holonomy representation is a representation $\theta : L\to \mbox{Der}(\land U)$ of $L$ by derivations in $\land U$, and if $x\in L$ and $u\in U_n$ then 
\begin{eqnarray}\label{i9}
\theta (x)u \,\,+ <pu, x>\in V_{n-1}\otimes \land U_{n-1}.
\end{eqnarray}

 Now right multiplication by $L$ in $\overline{UL}$ makes (by definition) $\overline{UL}$ into a right enriched $L$-module. On the other hand (\ref{i8}) the holonomy representation in $\land U$ dualizes to make $(\land U)^\vee$ into a right enriched $L$-module, and hence into a $\overline{UL}$-module.
 
 \begin{Prop}
 \label{p5.2}
 Denote by $\varepsilon_U : \land U\to \mathbb Q$ the augmentation defined by $\varepsilon_U (\land^{\geq 1}U)= 0$. Then an isomorphism
 $$\eta_L : \overline{UL} \stackrel{\cong}{\longrightarrow} (\land U)^\vee $$
 of right $L$-modules is defined by $$<\Phi, \eta_L(a)> = (-1)^{deg\, a} \varepsilon_U (a\cdot \Phi).$$
 \end{Prop}

 \vspace{3mm}\noindent {\sl proof.}  First observe that for $x\in L$,
 $$\begin{array}{ll}
 <x\cdot \Phi, \eta_L(a)> &= (-1)^{deg \, a} \varepsilon_U (a\cdot x\cdot \Phi)\\
 &= (-1)^{deg \, x} <\Phi, \eta_L(ax)>.
 \end{array}$$
 Thus $\eta_L$ is a morphism of right $L$-modules. In \cite[Theorem 6]{RHTII} it is proved that $\eta_L$ is an isomorphism if dim$\, L_0/[L_0,L_0]<\infty$ and each $L_k$ is finite dimensional. Thus each $\eta_{L_\alpha}$ is an isomorphism. It follows that $\eta_L = \varprojlim_\alpha \eta_{L_\alpha}$ is an isomorphism. \hfill$\square$

  \subsection{Profree $L$-modules}

\vspace{3mm}\noindent {\bf Definition.}
Let $L$ be an enriched Lie algebra and $V$ an enriched vector space. The \emph{profree $L$-module freely generated by $V$} is the $L$-module 
$$V \, \widehat{\otimes}\, \overline{UL}.$$

\vspace{2mm} Profree $L$-modules satisfy the following property: If $M$ is an enriched $L$-module  every coherent map $f : V\to M$ extends in a unique way to a map of enriched $L$-module $V\widehat{\otimes}\, \overline{UL}\to M$.

In particular a morphism of enriched $L$-module of $M\widehat{\otimes}\overline{UL} $ to itself is the identity if it induces the identity on $M$.

\vspace{3mm} If $\land V\otimes \land U$ is the acyclic closure of the quadratic Sullivan model of $L$, then (Proposition \ref{p5.2})  $(\land U)^\vee \cong \overline{UL}$ is the profree module freely generated by $\mathbb Q$.

\subsection{Closure of an ideal}

\begin{lem}
\label{l5.2} If $S\subset L$ is any subspace then $S\cdot \overline{UL}$ is an ideal. If $S$ is a graded space of finite type, then $S\cdot \overline{UL}$ is closed.\end{lem}

 \vspace{3mm}\noindent {\sl proof.} Observe that
 $$[x\cdot \Phi, y]= x\cdot \Phi y\,, \hspace{1cm} x,y\in L,\, \Phi\in \overline{UL},$$
 and so $S\cdot \overline{UL}$ is an ideal. Since a graded space is closed if and only if each component in a given degree is closed, if $S$ has finite type it is sufficient to show that each $S^k\cdot \overline{UL}$ is closed. Since a finite sum of closed subspaces is closed we have only to show that $x\cdot \overline{UL}$ is closed for each $x\in L$.
 
 Denote $\rho_\alpha x$ by $x_\alpha$. Then the closure of $x\cdot \overline{UL}$ is the inverse limit
 $$\varprojlim_\alpha \rho_\alpha (x\cdot \overline{UL}) = \varprojlim_\alpha x_\alpha \cdot \widehat{UL_\alpha} = \mathbb Q x_\alpha \oplus \varprojlim_\alpha x_\alpha \cdot \widehat{J_\alpha} = \mathbb Q x_\alpha \oplus \varprojlim_{n,\alpha} (x_\alpha \cdot J_\alpha)/(x_\alpha\cdot J_\alpha^n).$$
 On the other hand, because each $J_\alpha/J_\alpha^n$ is finite dimensional, the surjections
 $J_\alpha /J_\alpha^n \to (x_\alpha \cdot J_\alpha)/x_\alpha \cdot J_\alpha^n)$ induce a surjection
 $$\overline{J} \to \varprojlim_{n,\alpha} x_\alpha\cdot J_\alpha /(x_\alpha \cdot J_\alpha^n).$$
 This factors through the map $\overline{J}\to x\cdot \overline{J}$, and therefore shows that
 $$x\cdot \overline{J} \to \varprojlim_\alpha \rho_\alpha (x\cdot \overline{J})$$
 is surjective. But this is the inclusion of $x\cdot \overline{J}$ in its closure, and so $x\cdot \overline{J}$ is closed.

 \hfill$\square$

\subsection{The $L/I$-module $I/I^{(2)}$}

Recall from \S 4.1 that a closed ideal $I\subset L$ decomposes $\land V$ as a $\Lambda$-extension $\land V = \land W\otimes \land Z$ in which the inclusion $(\land W, d_W) \to (\land V,d)$ and the surjection $(\land V,d)\to (\land Z, d_Z)$ are morphisms of quadratic Sullivan algebras corresponding respectively to the morphisms
$$L\to L/I=L_W \hspace{5mm}\mbox{and } I= L_Z\hookrightarrow L.$$

Thus filtering by the ideals $\land^{\geq k}W\otimes \land Z$ produces a spectral sequence converging from $(\land W\otimes H(\land Z,d_Z), \overline{d}) $ to $H(\land V,d)$. Note that since $(\land Z, d_Z)$ is quadratic
$$H(\land Z, d_Z) = \oplus_k H^{[k]}(\land Z, d_Z),$$
where $H^{[k]}(\land Z)$ is the subspace represented by classes having a representative in $\land^kZ$. Moreover, by definition
$$\overline{d} : 1\otimes H^{[k]}(\land Z) \to W\otimes H^{[k]}(\land Z).
$$
It follows that each $(\land W\otimes H^{[k]}(\land Z, \overline{d})$ is a quadratic $(\land W,d)$-module. Hence (\S 5.3) this defines a holonomy representation of $L_W= L/I$ in each $H^{[k]}(\land Z)$. 

\begin{Prop}
\label{p5.3} With the hypotheses and notation above, there is a natural isomorphism
$$s^{-1}(Z\cap \mbox{ker}\, d_Z)^\vee = s^{-1}H^{[1]}(\land Z)^\vee \stackrel{\cong}{\longrightarrow} I/I^{(2)}$$
of enriched $L/I$-modules, in which
\begin{enumerate}
\item[$\bullet$] $L/I$ acts in $H^{[1]}(\land Z)^\vee$ by the dual of the holonomy representation, and
\item[$\bullet$] $L/I$ acts in $I/I^{(2)}$ by the representation induced from the right adjoint representation of $L$ in $I$.
\end{enumerate}
\end{Prop}

\vspace{3mm}\noindent {\sl proof.}  By definition, $H^{[1]}(\land Z, d_Z) = Z\cap \mbox{ker}\, d_Z$, while from Lemma \ref{l4.1} we obtain
$$s^{-1}(Z\cap \mbox{ker}\, d_Z)^\vee = L_Z/L_Z^{(2)} \hspace{2mm} =I/I^{(2)}.$$
It remains to verify that this isomorphism is an isomorphism of $L_Z/L_Z^{(2)}$-modules.

Denote the holonomy representation of $L/I$ in $Z\cap \mbox{ker}\, d_Z$ by $\theta$. Then if $z\in Z\cap \mbox{ker}\, d_Z$, 
$$\overline{d}z= \sum\, w_i\otimes z_i$$
with $w_i\in W$ and $z_i\in Z\cap \mbox{ker}\, d_Z$. Thus (\S 5.3), for $y\in L$, $\theta (y)$ is determined by the equation
$$<\theta(y)z, sx> = - \sum_i <w_i, sy><z_i, sx> = - <dz, sx, sy>.$$
Now it follows from Proposition \ref{p5.1}(iv)  that
$$<\theta (y)z, sx> = <z, (-1)^{deg\, y} s[x,y]> = <z, \mbox{r ad}x(y)>.$$
\hfill$\square$

\newpage

\newpage
\part{Profree Lie algebras}

Free graded Lie algebras freely generated by a space $T$, and which we denote by  $\mathbb L_T$, play a key role in Lie algebra theory. Profree Lie algebras are the analogue of free Lie algebras in the category of complete enriched Lie algebras. More precisely:

\vspace{3mm}\noindent {\bf Definition} A \emph{profree Lie algebra} is a complete enriched Lie algebra of the form $\overline{\mathbb L}_T$ in which $(T, \{\rho_\alpha\})$ is a complete enriched vector space and the enriched structure in $\mathbb L_T$ is defined by the surjections
$$\mathbb L_T\to \mathbb L_{\rho_\alpha T} \to \mathbb L_{\rho_\alpha T} / \mathbb L^n_{\rho_\alpha T}.$$
We call $\overline{\mathbb L}_T$ the profree Lie algebra generated by $T$.

\vspace{3mm}\noindent {\bf Remark}; As observed in \S 1.1 the map $\mathbb L_T\to \overline{\mathbb L}_T$ is injective.

\section{Characterization of profree Lie algebras} 

\begin{Prop}
\label{p6.1} Let $(T, \{\rho_\alpha\})$ be a complete enriched vector space.
\begin{enumerate}
\item[(i)] Any coherent linear map $f : T\to E$ of degree $0$ into a complete enriched Lie algebra $E$ extends uniquely to a morphism
$$\overline{\mathbb L}_T \to E$$
of complete enriched Lie algebras.
\item[(ii)] $\overline{\mathbb L}_T = \overline{\mathbb L}^{(2)}_T \oplus T.$
\item[(iii)] If $T'$ is any closed direct summand of $\overline{\mathbb L}_T^{(2)}$ in $\overline{\mathbb L}_T$ then the inclusion extends to an isomorphism
$$\overline{\mathbb L}_{T'} \stackrel{\cong}{\longrightarrow} \overline{\mathbb L}_T$$
of complete enriched Lie algebras.
\item[(iv)] If $T$ is a graded vector space of finite type then the completion of the universal enveloping algebra $U\overline{\mathbb L}_T$ is given by 
$$\overline{U\overline{\mathbb L}_T} = \prod_{n\geq 0} \bigotimes^n(T)$$
where $\bigotimes^n$ denotes the $n^{th}$ tensor power.
\end{enumerate}
\end{Prop}

\vspace{3mm}\noindent {\sl proof.} (i) Denote by $E\to E_\beta$ the surjections in the enriched structure of $E$. Since $E_\beta$ is finite dimensional and nilpotent it follows that the extension of $f$ to a morphism $\mathbb L_T \to E$ factors to yield morphisms
$$\mathbb L_{\rho_{\alpha (\beta)}T} \to E_\beta.$$
Passing to inverse limits gives the morphism $\overline{\mathbb L}_T \to E$, and it is immediate that it is the unique extension of $f$. 

(ii) Since $\mathbb L_{\rho_\alpha T} = \rho_\alpha T\oplus \mathbb L^{2}_{\rho_\alpha T}$ it follows from Lemma \ref{l2.2} that
$$\overline{\mathbb L}_T^{(2)} = \varprojlim_\alpha \mathbb L_{\rho_\alpha} T / \rho_\alpha T = \overline{\mathbb L}_T/T.$$
This gives (ii).

(iii) If $T'$ is a second closed direct summand of $\overline{\mathbb L}_T^{(2)}$ in $\overline{\mathbb L}_T$ then (ii) provides morphisms
$$\alpha : \overline{\mathbb L}_{T'}\to \overline{\mathbb L}_T \hspace{4mm}\mbox{and } \beta : \overline{\mathbb L}_T\to \overline{\mathbb L}_{T'}$$
which induce isomorphisms $\overline{\mathbb L}_{T'} / \overline{\mathbb L}_{T'}^{(2)} \cong \overline{\mathbb L}_T / \overline{\mathbb L}_T^{(2)}$. Now filtering  by the ideals $\overline{\mathbb L}_{T'}^{(n)}$ and $\overline{\mathbb L}_T^{(n)}$ shows that $\alpha$ and $\beta$ are inverse isomorphisms.

(iv) Because $T$ has finite type, (cf \S1.2)
$$\overline{U\overline{\mathbb L}_T}= \widehat{U\mathbb L_T} = \widehat{\bigotimes (T)}= \prod_n \bigotimes^n(T)$$
where $\bigotimes (T)$ is the tensor algebra on $T$.
\hfill$\square$

\begin{Theorem}
\label{t1}
The following conditions are equivalent on a complete enriched Lie algebra, $L$.
\begin{enumerate}
\item[(i)] $L= \overline{\mathbb L}_T$ is profree.
\item[(ii)]   The quadratic model $\land V$ of $L$ satisfies
$$H(\land V) = \mathbb Q \oplus (V\cap \mbox{ker}\, d)$$
\end{enumerate}
\end{Theorem}

  \vspace{3mm} In this case division by $\overline{\mathbb L}_T^{(2)}$ converts the isomorphism $s\overline{\mathbb L}_T \stackrel{\cong}{\to} V^\vee$   to an isomorphism $sT\stackrel{\cong}{\longrightarrow} (V\cap \mbox{ker}\, d)^\vee$. 

\subsection{Quadratic Sullivan algebras} 

For any quadratic Sullivan algebra, $\land V$, recall (\S 4) that the standard filtration of $V$ is given by
$$V_0= V\cap \mbox{ker}\, d \hspace{5mm}\mbox{and } V_{n+1} = V\cap d^{-1}(\land^2V_n).$$
In particular,
$$H(\land V) = \oplus_k H^{[k]}(\land V),$$
where $H^{[k]}(\land V)$ is the image in $H(\land V)$ of the space of cycles in $\land^kV$. 
Our objective here is to prove (Lemma \ref{l6.2}) that
$$H^{[2]}(\land V)= 0 \Longleftrightarrow H(\land V) = \mathbb Q \oplus (V\cap \mbox{ker}\, d).$$

\begin{lem} 
\label{l6.1} Suppose $\psi: E\to L$ is a morphism of complete enriched Lie algebras.
\begin{enumerate}
\item[(i)] $\psi $ is surjective if and only if $\psi (2): E/E^{(2)}\to L/L^{(2)}$ is surjective.
\item[(ii)] Suppose the quadratic model, $\land V$, of $L$ satisfies $H^{[2]}(\land V)= 0$. Then the following assertions are equivalent:
\begin{enumerate}
\item[(a)] $\psi (2)$ is an isomorphism,
\item[(b)] $\psi(n)$ is an isomorphism, $n\geq 2$,
\item[(c)] $\psi$ is an isomorphism.
\end{enumerate}
\end{enumerate}
\end{lem}

\vspace{3mm}\noindent {\sl proof.} First recall that $\psi$ is the dual of   $\varphi :  V\to  W$ where $\varphi : \land V\to \land W$ is the corresponding morphism between the quadratic models of $L$ and $E$. In view of Lemma \ref{l4.1}, $\psi(n+2): E/E^{(n+2)}\to L/L^{(n+2)}$ is the dual of $\varphi_n: V_n\to W_n$.  

(i) We have only to show that if $\varphi_0$ is injective then $\varphi$ is injective. Assume by induction that $\varphi_n$ is injective. if $v\in V_{n+1}$ and $\varphi v= 0$ then since $dv\in \land ^2V_n$ and $\varphi_n$ is injective it follows that $dv= 0$. Thus $v\in V_0$ and since $\varphi_0$ is injective $v= 0$.

(ii) Suppose (a) is satisfied and assume by induction that $\varphi_n : V_n \stackrel{\cong}{\to} W_n$. If $w\in W_{n+1}$ then $dw$ is a cycle in $\land^2W_n$ and so $dw= \varphi\Phi$ where $\Phi$ is a cycle in $\land^2V_n$. Because $H^{[2]}(\land V)= 0$, it follows that $\Phi = dv$ for some $v\in V_{n+1}$. Thus $d(w-\varphi v)= 0$ and so $w-\varphi(v)\in W_0$. By hypothesis, $w-\varphi v = \varphi v_0$ for some $v_0\in V_0$, and so $W_{n+1} \subset \varphi (V_{n+1})$. On the other hand, by (i), $\varphi_{n+1} : V_{n+1}\to W_{n+1}$ is injective. This proves that (a) $\Rightarrow$ (b).

But then $\varphi = \varinjlim \varphi_n$ and so $\varphi$ is an isomorphism. Thus (b) $\Rightarrow$ (c). Finally, it is immediate that if $\varphi$ is an isomorphism then so is $\varphi_0$. \hfill$\square$

\begin{lem}
\label{l6.2}
If a morphism $\rho : \land V\to \land Z$ of quadratic Sullivan algebras restricts to a surjection $\rho_0 : V_0\to Z_0$, and if $H^{[2]}(\land V)= 0$, then
$$H^{[k]}(\land Z)= 0, \hspace{3mm} k\geq 2.$$
\end{lem}

\vspace{3mm}\noindent {\sl proof.} The proof is in two Steps.

\vspace{2mm}\noindent \emph{Step One. If, in addition, $\rho$ is surjective then $H^{[k]}(\land Z) = 0$, $k\geq 2$.}

Here we first show that $H^{[2]}(\land Z)= 0$. Since  $\rho$ is surjective, by \cite[Cor.3.4]{RHTII} $\rho$ extends to a quasi-isomorphism
$$\rho : \land V\otimes \land U \stackrel{\simeq}{\to} \land Z$$
in which $\land V\otimes \land U$ is a $\Lambda$-extension of $\land V$ and $d: \land U\to V\otimes \land U$. Since $\rho$ is surjective there is a quasi-isomorphism $\sigma : \land Z\to \land V\otimes \land U$ satisfying $\rho\circ \sigma = id$.  Because $d: \land U\to V\otimes \land U$ and $\land V\otimes \land U$ is a $\Lambda$-extension it is straightforward to verify that $\sigma$ may be constructed so that $\sigma: Z\to V\otimes \land U$. In particular, the decomposition $\land V\otimes \land U = \oplus_k \land^kV\otimes \land U$ induces a decomposition of its homology, and
$$\sigma: H^{[k]}(\land Z) \stackrel{\cong}{\to} H^{[k]} (\land V\otimes \land U), \hspace{5mm} k\geq 0.$$

On the other hand, because $\land V\otimes \land U$ is a $\Lambda$-extension, $\land U$ is the increasing union of the subspaces $(\land U)_q$ given by
$$(\land U)_0= \land U \cap d^{-1}(\land V\otimes 1)\hspace{5mm}\mbox{and } (\land U)_{q+1} = \land U \cap d^{-1}(\land V\otimes (\land U)_q).$$
Thus
$$H^{[k]}(\land V\otimes \land U) = \varinjlim_q H^{[k]}(\land V\otimes (\land U)_q).$$
But the differential in the quotients $\land V\otimes ((\land U)_{q+1}/(\land U)_q)$ is just $d\otimes id$, and it follows that $H^{[2]}(\land V\otimes \land U)= 0$.

It remains to show that $H^{[k]}(\land Z)= 0$, $k\geq 3$. Suppose $\Phi\in \land^kZ$ is a cycle. There is then a sequence $z_1, \dots , z_r$ of elements in $Z$ such that $dz_1= 0$, $dz_{i+1}\in \land^2(z_1, \dots , z_i)$, $\Phi\in \land (z_1, \dots , z_r)$ and such that division by $z_1, \dots , z_r$ maps $\Phi$ to zero. When $r=1$, $\Phi = z_1^k$ is a boundary. We use induction on $r$ and on $k$ to show that $\Phi$ is a boundary.

Observe first that division by $z_1$ gives a quadratic Sullivan algebra $(\land Z',d')$. By what we just proved, $H^{[2]}(\land Z',d')= 0$. Thus by induction on $r$, the image of $\Phi$ in $\land Z'$ is a boundary. It follows that for some $\Phi'\in \land^{k-1}Z$,
$$\Phi -d(1\otimes \Phi') = z_1\otimes \Phi'',$$
with $\Phi''\in \land Z$. In particular, in $\land Z'$, $\Phi''$ is a $d'$-cycle. By induction on $k$ in $\land Z'$, $\Phi'' = d'\Omega$ for some $\Omega \in \land ^{k-2}Z'$. Therefore $\Phi''= d\Omega + z_1\otimes \Psi$ for some $\Psi$.

If deg$\, z_1$ is odd then $z_1\otimes \Phi''= d(-z_1\otimes \Omega)$, and so $\Phi$ is a boundary. If deg$\, z_1$ is even then, since $H^{[2]}(\land Z)= 0$, we may choose $z_2$ so $dz_2= z_1^2$. Division by $z_2$ and $z_1^2$ then gives a quasi-isomorphism $\land Z\stackrel{\cong}{\longrightarrow} \left((\land z_1)/z_1^2\right) \otimes \land Z''$ and the same argument as above shows that $\Phi$ is a boundary.

\vspace{2mm}\noindent \emph{Step Two. Completion of the proof of Lemma \ref{l6.2}.}

  We define a sequence of surjective morphisms
$$\land V= \land V(1) \twoheadrightarrow \land V(2) \twoheadrightarrow \dots \twoheadrightarrow \land V(p)\twoheadrightarrow \dots$$
such that $\rho$ factors through each to yield morphisms
$$\rho(p) : \land V(p)\twoheadrightarrow \land Z.$$
In fact, if $\rho(p)$ has been defined let $\varphi_p: \land V(p)\twoheadrightarrow \land V(p+1)$ be obtained by division by $V(p)_0 \cap \mbox{ker} \rho(p).$

Then the kernels of the surjections $\land V \to \land V(p)$ form an increasing sequence of subspaces $K(p)\subset V$. Set $K = \cup_p K(p)$ and let
$$\varphi : \land V\to \land W$$
be the surjection obtained by dividing $V$ by $K$. Then $\varphi = \varinjlim_p \varphi_p$, and so   $\varphi$ and $\varphi_0$ are surjective. Thus by Step One   $H^{[k]}(\land W)= 0$, $k\geq 2$.

Now by construction,   $\rho$ factors as
$$\land V\stackrel{\varphi}{\to} \land W\stackrel{\gamma}{\to} \land Z.$$
Moreover, $\land W = \varinjlim_p \land V(p)$ and so $W_0 = \varinjlim_p V(p)_0$. Let $I(p)$ be the image of $V(p)_0$ in $V(p+1)_0$. Then also
$$W_0 = \varinjlim_p I(p).$$
Thus by construction, $\gamma   : W_0\stackrel{\cong}{\to} Z_0$. Now Lemma \ref{l6.1}(ii) asserts that $\gamma$ is an isomorphism, and hence
$$H^{[k]}(\land Z) = H^{[k]}(\land W)= 0, \hspace{3mm} k\geq 2.$$

\hfill$\square$

\subsection{Proof of Theorem \ref{t1}}

This is in two steps. We fix a closed summand, $T$, of $L^{(2)}$ in $L$  and we denote by $(\land V,d)$ the quadratic model of $L$. Then it follows from Lemma \ref{l4.1} that the isomorphism $L\stackrel{\cong}{\longrightarrow} V^\vee$ restricts to an isomorphism $T\stackrel{\cong}{\to} (V\cap \mbox{ker}\, d)^\vee$.

\vspace{3mm}\noindent {\emph{Step One.} If $L= \overline{\mathbb L}_T$ is profree   then $H(\land V)= \mathbb Q \oplus (V\cap \mbox{ker}\, d)$.

\vspace{2mm} In this case, by adjoining 	additional variables construct an inclusion $\lambda : \land V\to \land W$ of quadratic Sullivan algebras for which 
$$V_0\stackrel{\cong}{\longrightarrow} W_0\hspace{5mm}\mbox{and } H^{[2]}(\land W)= 0.$$
Dualizing $V\to W$ gives a surjection $L\stackrel{\rho}{\leftarrow} L_W$ of complete enriched Lie algebras. Moreover, since $V_0\stackrel{\cong}{\to} W_0$, $L_W = L_W^{(2)}\oplus T_W$ where $T_W$ is closed and 
$$\rho_T : T\stackrel{\cong}{\longleftarrow} T_W.$$

Here both $\rho_T$ and its inverse, $\sigma$, are coherent and so $\sigma$ extends to a morphism $\varphi : L\to L_W$. Since $\rho\circ \varphi \vert_T = id\vert_T$ it follows that $\rho\circ \varphi = id_L$. Now dualizing, $\varphi$,  gives a morphism
$$\psi : \land V\leftarrow \land W$$
such that $\psi\circ \lambda = id$.

On the other hand, the inclusion $\lambda$ defines a $\Lambda$-extension
$$\land V\otimes \land Z \stackrel{\cong}{\longrightarrow} \land W,$$
in which, if $Z\neq 0$ then some $z\in Z$ satisfies $dz\in \land^2V$. Thus
$$d\psi z= \psi dz= dz.$$
Hence $\psi z-z$ is a cycle in $W$. This gives
$$\psi z -z\in W_0= V_0$$
and so $z\in V$. This contradicts $z\in Z$ and it follows that $Z=0$ and $\psi$ is an isomorphism inverse to $\lambda$. In particular $H^{[2]}(\land V)= 0$ and now Lemma \ref{l6.2} implies that $H(\land V) = \mathbb Q \oplus (V\cap \mbox{ker}\, d).$

\vspace{3mm}\noindent {\emph{Step Two.}  If $H(\land V)= \mathbb Q \oplus (V\cap \mbox{ker}\, d)$ then $L$ is profree.

\vspace{2mm} In this case, by Proposition \ref{p6.1}, the inclusion $T\to L$ extends to a morphism
$$\varphi : \overline{\mathbb L}_T \to L.$$
On the one hand, Proposition \ref{p4.2} implies that the image of $\varphi$ is closed. On the other hand, since the image contains $T$, Lemma \ref{l3.2} implies that $\varphi$ is surjective. And, finally, Proposition \ref{p4.2} also asserts that ker$\, \varphi$ is closed.

Now (\S4.1), taking the quadratic models for the short exact sequence
$$0 \leftarrow L\leftarrow \overline{\mathbb L}_T\leftarrow \mbox{ker}\, \varphi\leftarrow 0,$$
yields the $\Lambda$-extension
$$\land V\to \land V\otimes \land Z\to \land Z$$
in which $\land V\otimes \land Z$ is the quadratic model of $\overline{\mathbb L}_T$. We show that $\varphi$ is an isomorphism by proving that $Z=0$.

But if $Z\neq 0$ then for some $z\neq 0$ in $Z$, $dz\in \land^2V$. Since by hypothesis $H^{[2]}(\land V)= 0$ this implies that $dz=dv$ some $v\in V$. But then $v-z$ is a cycle representing a class in $H^{[1]}(\land V\otimes \land Z)$. However, because $\varphi$ induces an isomorphism $\overline{\mathbb L}_T/ \overline{\mathbb L}_T^{(2)} \to L/L^{(2)}$ it follows that
$$V\cap \mbox{ker}\, d= H^{[1]}(\land V) \stackrel{\cong}{\to} H^{[1]}(\land V\otimes \land Z) = (V\oplus Z)\cap \mbox{ker}\, d.$$
Thus $z=0$ and therefore $Z=0$.
\hfill$\square$

\subsection{Properties of profree Lie algebras}

\vspace{3mm} We have from Theorem \ref{t1}:

\vspace{3mm}\noindent {\bf Corollary 1.} {\sl 
Suppose $E= R \oplus E^2$ is a free graded Lie algebra. If dim$\, R<\infty$ and $E$ is equipped with the unique (Proposition \ref{p3.3}) enriched structure then $\overline{E}$ is profree.
}

\vspace{3mm}\noindent {\sl proof.}  By Proposition \ref{p3.1}, $\overline{E}= R\oplus \overline{E}^{(2)}$. Now apply Theorem \ref{t1}.
 \hfill$\square$

\vspace{3mm}\noindent {\bf Corollary 2.} {\sl  Suppose $\overline{\mathbb L}_T$ is a profree Lie algebra. Then
\begin{enumerate}
\item[(i)]   If $S\subset T$ is any closed subspace then the closure $\overline{E}$ of the sub Lie algebra $E$ generated by $S$ satisfies $$\overline{E} = S\oplus \overline{E}^{(2)},$$ and $\overline{E}^{(2)} = \overline{E}\cap L^{(2)}$.
\item[(ii)] If $0\neq E$ is a solvable sub Lie algebra of  $L$, then $E$ is a free Lie algebra on a single generator.
\end{enumerate}}

\vspace{3mm}\noindent {\sl proof.} (i)  By Lemma \ref{l3.2}, $\overline{E}^{(2)}= \overline{E^2}\subset L^{(2)}$. Therefore $E= S\oplus  \overline{E}^{(2)}$. Since $\overline{E}$ is profree, $S$ generates a free Lie algebra.  

(ii) Denote by $E \supset \dots \supset E^{[k]}\supset \dots$ the sequence of ideals defined by $E^{[k+1]} = [E^{[k]}, E^{[k]}]$. By hypotheses, some $E^{[n+1]}= 0$. Then, by induction on $n$, we may assume $[E,E]$ is either zero or a free Lie algebra on a single generator. It follows that the closure, $\overline{E}\subset L$ satisfies dim$\, \overline{E}^{(2)}\leq 2$, since (Lemma \ref{l3.2}) $\overline{E}^{(2)} = \overline{E^2}$.

But since $E\neq 0$, Corollary 1 gives that $\overline{E}= S\oplus \overline{E}^{(2)}$ and $S$ generates a free Lie algebra. This implies that $\overline{E}$ is the free Lie algebra generated by a single element. In particular dim$\, \overline{E}\leq 2$ and so $E= \overline{E}$. 

\hfill$\square$

\vspace{3mm}\noindent {\bf Corollary 3.} {\sl If $E$ is a complete enriched Lie algebra then the inclusion of a direct summand $T$ of $E^{(2)}$ in $E$ extends to a surjection $\overline{\mathbb L}_T\to E$.}

\vspace{3mm}\noindent {\sl proof.} This follows by the same argument as at the start of the proof of Step Two in the Theorem \ref{t1}.
\hfill$\square$

\vspace{3mm}\noindent {\bf Corollary 4.} {\sl If $\psi : (L,{\mathcal I})\to (F, {\mathcal G})$ is a surjective morphism from a complete enriched Lie algebra to a profree Lie algebra then there is a morphism
$$\sigma : (F, {\mathcal G})\to (L, {\mathcal I})$$
such that $\psi\circ \sigma = id_F$.}

\vspace{3mm}\noindent {\sl proof.} Let $\varphi : \land V_F\to \land V$ be the injective morphism of the corresponding quadratic models which dualizes to $\psi$.   Write $V = \varphi (V_F)\oplus S$. Inverting $\varphi$ then defines an isomorphism $\gamma$ from the sub quadratic algebra $\land \varphi(F)$ to $\land V_F$ and clearly $\gamma\circ \varphi = id$. It remains to extend $\gamma$ to $S$ so that $\gamma (dw)= d\gamma (w)$, $w\in S$.

For this set $\gamma = 0$ in $S\cap V_0$. Now suppose $\gamma$ has been defined in $S\cap V_n$ and let $\{z_i\}$ be a basis of a direct summand of $S\cap V_n$ in $S\cap V_{n+1}$. Then $\gamma (dz_i)$ is a cycle in $\land^2V_F$. Hence $\gamma (dz_i)= dw_i$ for some $w_i\in V_F$. Extend $\gamma $ by setting $\gamma z_i = w_i$. \hfill$\square$

 \vspace{3mm} Now recall from \cite[Chap 9]{RHTII} that the category of a minimal Sullivan algebra, $\land V$, is the least $p$ (or $\infty$) such that $\land V$ is a homotopy retract of $\land V/\land^{>p}V$.
 
 \begin{Prop}
 \label{p6.2}
 (i) A complete enriched Lie algebra is profree if and only if its quadratic model satisfies cat$(\land V)= 1$.  
 
 (ii) Any closed sub Lie algebra of a profree Lie algebra is profree.
 \end{Prop}
 
 \vspace{3mm}\noindent {\sl proof.} (i) The condition cat$(\land V)= 1$ for a quadratic Sullivan algebra is equivalent to the condition $\mathbb Q \oplus (V\cap \mbox{ker}\, d) \stackrel{\simeq}{\to} \land V$. Thus (i) follows from Theorem \ref{t1}.
 
 (ii) Suppose $\land V\to \land Z$ is the morphism of quadratic models corresponding to an inclusion $E\to L$ of a closed sub Lie algebra in a complete enriched Lie algebra. Then \cite[Theorem 9.3]{RHTII} gives cat$(\land Z)\leq $ cat$(\land V)$. Thus (ii) follows from (i). 
  \hfill$\square$

\vspace{3mm}\noindent {\bf Remark.} Recall that any sub Lie algebra of a free Lie algebra is free (\cite{Sh}). By Proposition \ref{p6.2}(ii) the analogous statement for complete Lie algebras and closed sub algebras is also true. 

\vspace{3mm}\noindent {\bf Example.}
 Let $L= \overline{\mathbb L}_T$ be a profree Lie algebra. Then the sub Lie algebra generated by   $T$ is free.  The closed sub Lie algebras of $L$ are profree, but the sub Lie algebras are not necessarily free. For instance let $L= \widehat{\mathbb L (a,b)}$ and $E\subset L$ be the sub Lie algebra generated by the elements $a,b$ and $\omega = \sum_{n\geq 1} \mbox{ad}^n_a (b)$. Since $\omega = [a,\omega]- [a,b]$, $\omega\in E^2$, so that $E/E^2$ is generated by the classes of $a$ and $b$. If $E$ was free, then $E$ would be equal to $\mathbb L(a,b)$, which is impossible because $\omega \not\in \mathbb L(a,b)$.

 \subsection{The structure of a profree Lie algebra}
 
 \begin{Prop}
 \label{p6.3}
 Let $L$ be a complete enriched Lie algebra.
 \begin{enumerate}
 \item[(i)] Suppose $L/[L,L]$ has finite type. Then $L$ is profree if and only if the inclusion $T\to L$ of a direct summand of $[L,L]$ extends to an isomorphism
$$\widehat{\mathbb L_T} := \varprojlim_n {\mathbb L}_T/{\mathbb L}_T^n \stackrel{\cong}{\longrightarrow} L.$$
\item[(ii)]   If $L= \varprojlim_\sigma L_\sigma$ is any inverse limit of profree Lie algebras then $L$ is profree.
\end{enumerate}
 \end{Prop}
 
 \vspace{3mm}\noindent {\sl proof.} (i). Proposition \ref{p3.2} identifies $\widehat{\mathbb L_T}\stackrel{\cong}{\longrightarrow}  \overline{\mathbb L}_T$, since $\mathbb L_T$ has a unique enriched structure. Thus (i) follows from Theorem \ref{t1}.
 
 (ii) If $\land V$ and $\land V_\sigma$ are respectively the quadratic models of $L$ and $L_\sigma$ then
 $\land V = \varinjlim_\sigma \land V_\sigma$ and the map $V\cap \mbox{ker}\, d\to \varinjlim_\sigma V_\sigma \cap \mbox{ker}\, d$ is an isomorphism. It follows (since each $L_\sigma$ is profree) from Theorem \ref{t1} that $V\cap \mbox{ker}\, d\stackrel{\cong}{\longrightarrow} H^{\geq 1}(\land V)$, and hence (Theorem \ref{t1}) that $L$ is profree.
 
 \hfill$\square$

 \vspace{3mm} Now, suppose that $\land V$ is the quadratic model of a profree Lie algebra, $\overline{\mathbb L}_T$.   In particular
 $$\mathbb L_T = \oplus_{k\geq 1} T(k),$$
 where $T(k)$ is the linear span of the iterated Lie brackets of length $k$ in elements of $T$. In general,   however the subspaces $T(k)$ may not be closed. Nevertheless we do have an analogous structure for $\overline{\mathbb L}_T$.
 
 For this, denote the enriched structure of $T$ by $\{\rho_\sigma : T\to T(\sigma)\}$. Then by definition,
 $$\overline{\mathbb L}_T = \varprojlim_\sigma \varprojlim_n \mathbb L_{T(\sigma)}/ \mathbb L^n_{T(\sigma)} = \varprojlim_\sigma L(\sigma),$$
 where $L(\sigma) = \widehat{\mathbb L}_{T(\sigma)}$. Then write $\mathbb L_T = \oplus_{k\geq 1}T(k)$, where $T(k)$ is the span of the iterated Lie brackets of length 
 $k$ in elements of $T$.

 \begin{Prop}
 \label{p6.4}
Let $L= \overline{\mathbb L}_T$ be a profree Lie algebra.  Then, with the notation above,
 \begin{enumerate}
 \item[(i)] $L^{(k)} = \varprojlim_\sigma L(\sigma)^k$.
 \item[(ii)] $L^{(k)}= \overline{T(k)} \oplus L^{(k+1)}$.
  \item[(iii)] For each $k, \ell$, $[\overline{T(k)}, \overline{T(\ell)}] \subset \overline{T(k+\ell)}$. In particular, $\oplus_k \overline{T(k)}$ is a weighted Lie algebra with closure $L = \prod_k \overline{T(k)}$. 
  \end{enumerate}
 \end{Prop}
 
 \vspace{3mm}\noindent {\sl proof.} (i) Since $L= \varprojlim_\sigma L(\sigma)$ it follow that $V = \varinjlim_\sigma V_\sigma$ where $\land V_\sigma$ is the quadratic model of $L(\sigma)$. This gives
 $$V/V_k = \varinjlim_\sigma V(\sigma) / V(\sigma)_k,$$
 and (Lemma \ref{l4.1}) this dualizes to $L^{(k)}= \varprojlim L(\sigma)^{(k)}$. Since $L(\sigma) = \widehat{\mathbb L}_{T(\sigma)} $ we have $L(\sigma)^{(k)} = L(\sigma)^k$. 
 
 (ii) By definition, the enriched structure in $T(k)$ is provided by the surjections ${\rho_\sigma}\vert_{T(k)}$. But since $\rho_\sigma : T\to T(\sigma)$, it restricts to surjections
 $$\xymatrix{\rho_\sigma : T(k) \ar@{->>}[rr] && T(\sigma, k),}$$
 where $T(\sigma, k)$ is the span of the iterated Lie brackets of length $k$ in elements of $T(\sigma)$.  This gives
 $$\overline{T(k)} \stackrel{\cong}{\longrightarrow} \varprojlim_\sigma T(\sigma, k).$$
 
 On the other hand the inclusion of $T(k)$ in $\mathbb L_T^k$ induces a map $\overline{T(k)} \to L^{(k)}/L^{(k+1)}$.
 But this factors as
 $$\overline{T(k)} \stackrel{\cong}{\longrightarrow} \varprojlim_\sigma T(k, \sigma) = \varprojlim_\sigma L(\sigma)^{(k)}/L(\sigma)^{k+1)} = L^{(k)}/L^{(k+1)}.$$
 Thus this is an isomorphism and $L^{(k)} = \overline{T(k)} \oplus L^{(k+1)}.$
 
 (iii) This is immediate from $[T(k), T(\ell)]\subset T(k+\ell)$.
 
 \hfill$\square$

  \vspace{3mm} The next Proposition explains the difference between $L^2$ and $L^{(2)}$ for general profree Lie algebras.
  
 \begin{Prop}\label{p6.5} If $L=L_0$ is a profree Lie algebra   and dim$\, L/[L,L]=\infty$, then $L^2\neq L^{(2)}$. If in particular, $L/L^{(2)}$ is concentrated in degree $r$, then $L^2_{2r}\neq L^{(2)}_{2r}$.
 \end{Prop}

  \vspace{3mm}\noindent {\sl proof.}
  Indeed, let   $(\land W,d)$ be the quadratic model of $L$. We can decompose $W$ as an union $W = \cup_n W_n$ with $W_0 = W \cap \mbox{ker}\, d$,  and for $n>1$, $W_n= d^{-1}(\land^2W_{n-1})$. We denote by $Z_n$ a direct summand of $W_{n-1}$ in  $W_n$. Then by hypothesis $W_0$ is infinite dimensional and there is a quasi-isomorphism $\varphi : (\land W,d)\to (\mathbb Q \oplus W_0,0)$ that is the identity on $W_0$ and that maps each $Z_n$ to $0$.

For sake of simplicity, we write $V = W_0$ and $Z = Z_1$. By construction the isomorphism $L= (sW)^\vee$    induces isomorphisms
$$L/L^{(2)} \cong (sV)^\vee ,\hspace{1cm}\mbox{and } L^{(2)}/L^{(3)} \cong (sZ)^\vee.$$

  Suppose first $W$ countable and denote by $w_1, w_2, \dots$ a basis of $  V$. Since $d=0$ on $V$, $d: Z\to \land^2V$ is an isomorphism. We denote by $w_{ij}$, $i<j$ the basis of $Z$ defined by $d(w_{ij})= w_i \land w_j$.  Denote by $E$ the vector space of column matrices $X = (x_i)$ with only a finite number of nonzero $a_i$. Then the map $(a_i) \leadsto \sum a_i w_i$ defines an isomorphism $E\stackrel{\cong}{\to} V$. 
  
 Let represent an element $\varphi\in Z^\vee$ by the infinite dimensional antisymmetric matrix $M_\varphi$,$$(M_\varphi)_{ij}= \varphi (w_{ij}), \hspace{5mm} i<j.$$
  In a similar process, an element $f\in V^\vee$ can be represented by a column matrix $A_f$, with $(A_f)_i = f(w_i)$.  The vector space ker$\, f$ can then be identified with the sub vector space of $E$ formed by the column matrices $X$ satisfying $A_f^t\cdot X = 0$. (Here  $A^t$ denotes the line matrix  transposition of a column matrix $A$.)
  
  Note that when we have two column matrices $A$ and $B$, we can form the antisymmetric matrix $A\cdot B^t - B\cdot A^t$. Now remark that for $f, g\in V^\vee$, we have
  $$M_{[f,g]} = A_f\cdot B_g^t-B_g\cdot A_f^t.$$
  
  Let $\varphi_0\in Z^\vee$ be the particular element defined for $i<j$ by
  $$\varphi_0 (w_{ij}) = \left\{
  \begin{array}{ll} 1 & \mbox{if $w_{ij}= w_{2k+1, 2k+2}$, for some $k$}\\0 & \mbox{otherwise}\end{array}
  \right.
  $$
  The element $\varphi_0$ can be extended to all of $W$, by $\varphi_0(V)=0$ and $ \varphi_0(Z_n)= 0$, for $n>1$. By construction $\varphi_0\in L^{(2)}$. 
  The associated matrix is
  $$M_0 = \left( \begin{array}{cccc}
  B_0 & 0 & 0 & \dots\\
  0 & B_0 & 0 & \dots\\
  0 & 0 & B_0 & \dots \\
  \dots & \dots &\dots &\dots\end{array}
  \right) \hspace{15mm}\mbox{with } B_0 = \left(\begin{array}{cc} 0 & 1\\ -1 & 0\end{array}\right).$$
   
   Now consider a finite sum $\sum_{i=1}^n [f_i, g_i]$, with $f_i$ and $g_i\in V^\vee$. The associated matrix is $\sum_{i=1}^n M_{[f_i, g_i]}$. Then $$K = (\cap_{i=1}^n \mbox{ker}\, f_i) \cap (\cap_{i=1}^n \mbox{ker}\, g_i)$$
   is infinite dimensional, and so for some non-zero $X\in Z^\vee$,   $(\sum M_{[f_i,g_i]})\cdot X = 0$. Since $M_0\cdot X\neq 0$, it follows that $L^2\subset_{\neq} L^{(2)}$.   
   
   In the general case, let $E\subset V$ be a countable subvector space and let $\land T$ be the minimal model of $(\mathbb Q \oplus E,0)$. Then its homotopy Lie algebra $L_T$ is a retract of $L$, 
   $$\xymatrix{ 
   L_T \ar@/^/[r]^j & L\ar@/^/[l]^\rho .}$$
   Let $\varphi_0\in L_T^{(2)}$, not in $L_T^2$.Then $j(\varphi_0)\in L^{(2)}$ and not in $L^2$ because otherwise $\varphi_0= \rho(\varphi)$ would belong to $L_T^{2}$.

   \hfill$\square$

\begin{Prop}
\label{p6.6} Let $L= \overline{\mathbb L}_T$ be a profree Lie algebra in which $T= T_0$ and $sT_0$ is the dual of a vector space of countably infinite dimension. Then $\overline{\mathbb L}_T$   is not pronilpotent. \end{Prop}

\vspace{3mm}\noindent {\sl proof.} Denote the quadratic Sullivan model of $\overline{\mathbb L}_T$ by $(\land S,d)$. Then (Theorem \ref{t1})
$$sT_0 = (S\cap \mbox{ker}\,d)^\vee, \hspace{3mm} \mathbb Q \oplus (S\cap \mbox{ker}\, d) \stackrel{\cong}{\to} H(\land S), \hspace{3mm}\mbox{and } s\overline{\mathbb L}_T = S^\vee.$$
Thus $S\cap \mbox{ker}\, d$ has a basis of the form $v_i, w_i$, $i\in \mathbb N$. This in turn determines elements $x_i, y_i\in T_0$ by the conditions
 $$<v_i, sx_j>= \delta_{ij}= <w_i, sy_j> \hspace{5mm}\mbox{and } <v_i, sy_j>=0=<w_i, sx_j>.$$
 Now define $\omega\in \overline{\mathbb L}_T^{(2)}$ by
 $$\omega = \sum_{n\geq 1} ad^n_{x_n}y_n.$$
 We shall show by construction that $\omega \not\in \overline{\mathbb L}_T^2$ and then use Lemma \ref{l3.1}(iii) to prove that $\overline{\mathbb L}_T$ is not pronilpotent.
 
 Suppose then that $\omega\in \overline{\mathbb L}_T^{2}$, and write
 $$\omega = \sum_{i=1}^N [\omega'_i, \omega_i''].$$
 Since $sT = S^\vee$ the elements of $T$ are linear functions $S\to \mathbb Q$ of degree $-1$. Thus
 $$K = (\cap \mbox{ker}\, \omega_i') \cap (\cap \mbox{ker}\, \omega_i'')$$
 is a subspace of $S$ of finite codimension. Thus for some non-zero linear combination,
 $$a= \sum_{i=p}^{p+q} \lambda_i v_i,$$
 we have $a\in K$. We may, and do, assume that $\lambda_p= 1$.
 
 Now define a sequence $(a_r)_{r\geq 1}$ of elements of $S$ by the conditions
 $$da_1= aw_p, \hspace{3mm} da_r = aa_{r-1} \hspace{2mm} \mbox{if }r\geq 2, \hspace{3mm} \mbox{and } <a_r, sT_0>= 0, r\geq 1.$$
 It then follows from (\ref{i6}) that since $<a, s\omega_i'>= 0= <a, s\omega_i''>$, 
 $$<a_r , s\omega> = - \sum_{i=1}^N \, <da_r, s\omega_i', s\omega_i''> = 0, \hspace{3mm} r\geq 1.$$
 On the other hand, an iterated induction argument using (6) shows that
 $$<a_p, s\, ad_{x_k}^k y_k> = \left\{\begin{array}{ll} 1 & \mbox{if } k=p\\0 & \mbox{otherwise}.\end{array}\right.$$
 It follows that $<a_p, s\omega>= 1$, which is the desired contradiction.
 
 Finally, set 
 $$z_n = \sum_{k=1}^n ad_{x_n}^ky_k.$$
 Then $z_{n+1}-z_n \in \overline{\mathbb L}_T^{n+1}$. Thus the sequence $(z_n)$ is a coherent sequence in the tower
 $$\dots \to \overline{\mathbb L}_T/\overline{\mathbb L}_T^{n+1} \to \overline{\mathbb L}_T/\overline{\mathbb L}_T^n \to \overline{\mathbb L}_T/\overline{\mathbb L}_T^{n-1} \to \dots$$
 and therefore defines an element $z\in \widehat{\overline{\mathbb L}_T} = \varprojlim_n \overline{\mathbb L}_T/\overline{\mathbb L}_T^n$. Now Lemma \ref{l3.1} (iii) provides a retraction $\overline{\mathbb L}_T \stackrel{\varphi}{\rightarrow} \widehat{\overline{\mathbb L}_T} \stackrel{\psi}{\rightarrow} \overline{\mathbb L}_T$ in which $\psi z=\omega$. Since $\psi\circ \varphi = id$, $\psi (\varphi \omega ) = \omega$. Thus if $\psi$ is an isomorphism, $z=\omega$ and the images of $z$ and $\omega$ in $\overline{\mathbb L}_T/\overline{\mathbb L}_T^2$ coincide. But the image $z$ is in $\overline{\mathbb L}_T^2$ and $\omega \not\in \overline{\mathbb L}_T^2$. \hfill$\square$

\newpage

 \newpage
 \part{Topological spaces, Sullivan algebras, and homotopy Lie algebras}
 
 Recall that we identify the categories of topological spaces and simplicial sets via the singular simplex and Milnor realization functors (\cite{May}).

 \section{Sullivan models: spatial realizations and Sullivan completions}
 
 Sullivan's approach to rational homotopy theory is based on a two-way relation between topological spaces and his Sullivan algebras (\cite{S}, \cite[\S 1.6]{RHTII}), which we briefly review for the convenience of the reader. This relation is obtained from two contravariant adjoint functors
\begin{quote}
$A_{PL}$: simplicial sets $\leadsto$ cdga's \hspace{5mm}\mbox{and } $\langle\,,\,\rangle$: cdga's $\leadsto$ simplicial sets. 
\end{quote}
These are constructed as follows from a simplicial cdga $A_{PL}= \{(A_{PL})_n\}$ of acyclic algebras:
$$
A_{PL}(X) = Simpl (A, A_{PL})\hspace{5mm}\mbox{and } \langle A\rangle= Cdga (A, A_{PL}).
$$
Here the cdga structure of $A_{PL}(X)$ and the simplicial structure of $\langle A\rangle$ are inherited from the corresponding structures of $A_{PL}$. The adjointness of these functors follows from the equalities for all cdga's $A$ and simplicial sets $X$:
\begin{eqnarray}
\label{i10} 
Cdga (A, A_{PL}(X)) = \mbox{Hom}(A\times X, A_{PL})= Simpl (X, \langle A\rangle).
\end{eqnarray}
In particular, adjoint pairs will be denoted by
$$\varphi : A\to A_{PL}(X) \hspace{5mm}\mbox{and } \vert\varphi\vert : X\to \langle A\rangle.$$

Now specialize this construction from all cdga's to the category of Sullivan algebras, $\land V$. (For simplicity we denote $(\land V,d)$ simply by $\land V$ when we do not need to specify the differential.)

There is a natural equivalence relation (\emph{homotopy}) among cdga morphisms from a minimal Sullivan algebra (\cite{S}, \cite[\S 1.7]{RHTII}). It is straightforward from the definitions that (\ref{i10}) induces a natural bijection 
\begin{eqnarray}
\label{i11}
[\land V, A_{PL}(X)] \cong [X, \langle \land V\rangle]
\end{eqnarray}
of based homotopy classes of morphisms.

In particular, the two-way relation between topological spaces and Sullivan algebras is defined by the

\vspace{3mm}\noindent {\bf Definition.} \begin{enumerate}
\item[(i)] A \emph{(minimal) Sullivan model} of a connected space $X$ is a quasi-isomorphism
$$\varphi: \land V\stackrel{\simeq}{\longrightarrow} A_{PL}(X) $$ from a (minimal) Sullivan algebra.
\item[(ii)] The \emph{spatial realization} of a Sullivan algebra, $\land V$, is the connected based simplicial set $\langle\land V\rangle$. \end{enumerate}

\vspace{3mm}\noindent {\bf Remark.} The uniqueness (up to homotopy class of quasi-isomorphism) of a Sullivan model implies the uniqueness (up to homotopy class of homotopy equivalence) of a Sullivan  completion.

\vspace{3mm}
Now observe that any minimal Sullivan algebra, $\land V$, computes $\pi_*\langle\land V\rangle$ via a natural bijection
\begin{eqnarray}
\label{i12}
\pi_*\langle\land V\rangle \cong V^\vee.
\end{eqnarray}
In fact, if $\land W(n)$ is the minimal model of the $n$-sphere $S^n$ we have natural cdga quasi-isomorphisms
$$\xymatrix{A_{PL}(S^n) && \land W(n)\ar[ll]^\simeq \ar[rr]^\simeq && H(S^n) = \mathbb Q \oplus \mathbb Q a_n,}$$
$a_n$ denoting the orientation class of $S^n$. The bijections (\ref{i11}) convert these to bijections
$$\pi_n\langle\land V\rangle = [S^n, \langle\land V\rangle] = [\land V, H(S^n)] = (V^n)^\vee.$$
 These bijections naturally endow $\pi_*\langle\land V\rangle$ with the structure of a graded rational vector space. (Moreover (\cite[Theorem 1.4]{RHTII}) for $n\geq 2$, they are isomorphisms of abelian groups).
 
 \vspace{2mm}\noindent To summarize:
 \begin{enumerate}
 \item[(i)] \emph{If $\varphi : \land V\to A_{PL}(X)$ is a Sullivan model then $\land V$ computes the rational cohomology algebra $H(X)$ via
 $$H(\varphi) : H(\land V) \stackrel{\cong}{\longrightarrow} H(X).$$}
 \item[(ii)] \emph{If $\land V$ is any minimal Sullivan algebra then $\land V$ computes $\pi_*\langle\land V\rangle$ via the natural bijection
 $$\pi_*\langle\land V\rangle \cong V^\vee.$$}
 \item[(iii)] \emph{If $\varphi : \land V\to A_{PL}(X)$ is a minimal Sullivan algebra then the adjoint $\vert \varphi\vert : X\to \langle\land V\rangle$ induces 
 $$\pi_*(\vert\varphi\vert) : \pi_*(X)\to \pi_*\langle\land V\rangle$$
 connecting $\pi_*(X)$ to $\pi_*\langle\land V\rangle$.}
 \end{enumerate}
 
 \vspace{3mm}
 A more general connection between Sullivan algebras $\land V$ and its spatial realization $\langle\land V\rangle$ is provided by the morphism
 $$\iota_{\land V}: \land V\to A_{PL}\langle\land V\rangle$$
 adjoint to $id_{\langle\land V\rangle}$. In fact, if $\psi : \land V\to A_{PL}(X)$ is a morphism from an arbitrary Sullivan algebra then it is easy to verify that
 \begin{eqnarray}
 \label{i13}
 \psi = A_{PL}(\vert \psi\vert)\circ \iota_{\land V}.
 \end{eqnarray}
 
We also note that immediately from the definitions we obtain

 \begin{Prop}
 \label{p7.1}
 Suppose $\land V$ and $\land W$ are Sullivan algebras. Then
 \begin{enumerate}
  \item[(i)]  The natural map
 $$\langle\land V\otimes \land W\rangle \to \langle\land V\rangle \times \langle\land W\rangle$$
 is a homotopy equivalence
 \item[(ii)] If $\land Z$ is the minimal Sullivan model of $\land V\times_{\mathbb Q}\land W$ then the natural map
 $$\langle\land Z\rangle \to \langle\land V\rangle\vee \langle\land W\rangle$$
 is a homotopy equivalence.
 \item[(iii)] The diagram
 $$\xymatrix{
 \land V\otimes \land W \ar[d] \ar[r] & A_{PL}\langle\land V\rangle \otimes A_{PL}\langle\land W\rangle \ar[d] & A_{PL}(\langle\land V\rangle\times \langle\land W\rangle)\ar[d]\ar[l]\\\land V\times_{\mathbb Q}\land W \ar[r] & A_{PL}\langle\land V\rangle\times_{\mathbb Q}A_{PL}\langle\land W\rangle & A_{PL}(\langle\land V\rangle\vee \langle\land W\rangle)\ar[l]_\simeq}$$
 commutes. If either $H(\langle\land V\rangle)$ or $H(\langle\land W\rangle)$ is a graded vector space of finite type then $A_{PL}\langle\land V\rangle\otimes A_{PL}\langle\land W\rangle \longleftarrow A_{PL}(\langle\land V\rangle\times \langle\land W\rangle)$ is also a quasi-isomorphism.

 \end{enumerate}
 \end{Prop}

  \section{The homotopy Lie algebra   of a minimal Sullivan algebra}
 
 Associated with a minimal Sullivan algebra, $(\land V,d)$, is the quadratic Sullivan algebra, $(\land V,d_1)$ defined by: \emph{$d_1$ is the component in $\land^2V$ of $dv$}. This in turn determines (Proposition \ref{p4.1}) the complete enriched Lie algebra $L_V= (L_V)_{\geq 0}$ given by
 $$sL_V = V^\vee \hspace{4mm}\mbox{(as a graded vector space)}$$
 and 
 $$<v, s[x,y]> = (-1)^{1+ deg\, y} <d_1v, sx,sy>.$$
 
 \vspace{3mm}\noindent {\bf Definition.} $L_V$ is the \emph{homotopy Lie algebra} of $(\land V,d)$. 
 
 \vspace{3mm}\noindent {\bf Remarks.} 1. Filtering by the ideals $\land^{\geq p}V$ yields a spectral sequence, natural with respect to morphisms of minimal Sullivan algebras, and converging from $(\land V, d_1)$ to $H(\land V,d)$. 
 
 2.  If $\land V,d)= \varinjlim_\sigma (\land V_\sigma,d)$ then $L_V= \varprojlim_\sigma L_{V_\sigma}$.
 
 3. The enriched structure $L_V= \varprojlim_\alpha L_{V_\alpha}$ enables the extension of (\cite[Chapter 2]{RHTII}) to all homotopy Lie algebras, and with the same proofs. This removes the restriction dim$\, H^1(\land V)<\infty$, which is assumed in \cite{RHTII}.
 
 \begin{Prop}
 \label{p8.1} Suppose $\varphi, \psi : (\land V,d)\to (\land W,d)$ are homotopic morphisms between minimal Sullivan algebras. Then
 \begin{enumerate}
 \item[(i)] The induced morphisms $\varphi_1, \psi_1: (\land V, d_1)\to (\land W, d_1)$ coincide.
 \item[(ii)] The induced morphisms $L_V\leftarrow L_W$ coincide.
 \item[(iii)] $\varphi\vert_{V^1} = \psi\vert_{V^1}.$
 \end{enumerate}
 \end{Prop}
 
 \vspace{3mm}\noindent {\sl proof.} (i) Suppose $\Phi : (\land V,d) \to (\land W,d)\otimes \land (t,dt)$ is a homotopy from $\varphi$ to $\psi$ and denote by $(\varepsilon_0, \varepsilon_1): \land W\otimes \land (t,dt) \to \land W\times \land W$ the morphism sending $t$ to $0,1$. Then let $U= \mbox{ker}\,\varepsilon_0\cap \mbox{ker}\, \varepsilon_1$ and $A= \land^{\geq 1}W\otimes \land (t,dt) \oplus \mathbb Q$. This yields the quasi-isomorphisms
 $$\xymatrix{
 A &&A\otimes \land U\otimes \land dU \ar[ll]_\simeq^\rho \ar[rr]^\simeq_\gamma && \land W\otimes \land (t,dt)}$$
 in which $d: U\stackrel{\cong}{\to} dU$, $\rho (U)= 0$, and $\gamma\vert_U$ is the inclusion.
 
 Now lift $\Phi$ through $\gamma$ and then compose with $\rho$ to obtain a morphism
 $$\Psi : \land V\to A.$$
 Composing $\Psi$ with the inclusion $j : A\to \land W\otimes \land (t,dt)$ gives a second homotopy, and because of the choice of $U$ it is also a homotopy from $\varphi$ to $\psi$. Now filtering by the ideals $\land^{\geq j}V$ and $\land^{\geq j}W$ provides a homotopy
 $$\Gamma : (\land V, d_1) \to (\land W, d_1)\otimes \land (t,dt)$$
 from $\varphi_1$ to $\psi_1$. Moreover, division by $\land^{\geq 2}V$ and $\land^{\geq 2}W$ converts $\Gamma$ into a chain homotopy between
 $$\varphi_1\vert_V, \psi_1\vert_V : (V,0)\to (W,0).$$
 It follows that $\varphi_1\vert_V=\psi_1\vert_V$ and therefore that $\varphi_1=\psi_1$.
 
 (ii) This is immediate from (i).
 
 (iii) This follows because $\varphi\vert_{V^1}= \varphi_1\vert_{V^1}= \psi_1\vert_{V^1}= \psi\vert_{V^1}$. \hfill$\square$
 
 \vspace{3mm} Now, for any minimal Sullivan algebra $\land V$  the identifications (\ref{i12}) extend to natural bijections
 \begin{eqnarray}
 \label{i14}
 \pi_n\langle\land V\rangle= (V^n)^\vee= s(L_V)_{n-1}.
 \end{eqnarray}
 This converts the Lie brackets in $L_V$ to the Whitehead products in $\pi_*\langle\land V\rangle$:

 \begin{Prop}
 \label{p8.2} Suppose $(\land V,d)$ is a minimal Sullivan algebra. Then
 \begin{enumerate}
 \item[(i)] For $n\geq 2$, (\ref{i14}) is a linear isomorphism. In particular, $\pi_n\langle\land V\rangle$ is the rational vector space, $s(L_V)_{n-1}$.
 \item[(ii)] For $p,q \geq 2$, (\ref{i14}) identifies the Lie bracket $(L_V)_p \times (L_V)_q \stackrel{[\,,\,]}{\longrightarrow} (L_V)_{p+q}$ with the Whitehead product $\pi_{p+1}\langle\land V\rangle \times \,\pi_{q+1}\langle\land V\rangle \to \pi_{p+q+1}\langle\land V\rangle$. 
 \item[(iii)] The composite
 $$\xymatrix{G_L  && (L_V)_0 \ar[ll]^\cong_\exp  && \pi_1\langle\land V\rangle\ar[ll]^\cong}$$
 is an isomorphism of groups, where $G_L$ is the fundamental group of $L_V$ defined in \S 1.2.
 \item[(iv)] Finally (\ref{i14}) identifies the right adjoint  representation of $\pi_1\langle\land V\rangle$ in $\pi_n\langle\land V\rangle$ as
 $$\beta \bullet \exp \alpha = e^{-ad\,\alpha}(\beta), \hspace{5mm} \alpha \in (L_V)_0, \beta\in L_V.$$
 \end{enumerate}
 \end{Prop}

 \vspace{3mm}\noindent {\sl proof.} Recall that $(\land V,d)= \varinjlim_\alpha (\land V_\alpha,d)$ in which the $V_\alpha$ are finite dimensional subspaces of $V$. It follows that $\pi_*\langle\land V\rangle = \varprojlim_\alpha \langle\land V_\alpha\rangle$. Now Proposition \ref{p8.2} follows respectively from \cite[Theorem 1.4]{RHTII}, \cite[Proposition 13.16]{FHTI}, \cite[Theorem 2.4]{RHTII} and \cite[Theorem 2.5]{RHTII}. 
 \hfill$\square$
 
 \vspace{3mm}\noindent \emph{To summarize, a Sullivan algebra $\land V$ computes $\pi_*\langle\land V\rangle$ together with its Whitehead products, the group structure in $\pi_1\langle\land V\rangle$, and the action of $\pi_1\langle\land V\rangle $ on $\pi_{\geq 2}\langle\land V\rangle$.}

 \vspace{3mm}\noindent {\bf Example.} Denote by $\land V$ the minimal Sullivan model of $S^1\vee S^3$. Then $V^1= \mathbb Q v$ while $V^3 $ has a basis $w_0, w_1, \dots , w_n,\dots$ satisfying $dw_0= 0$ and $dw_{n+1} = vw_n$. Define $x\in (L_V)_0$ by $<v,sx>= 1$ and identifies $\pi_3\langle\land V\rangle = (V^3)^\vee$ as the space of infinite series  
 $$g(t)= \sum_{n\geq 0} c_n t^n$$ with $<w_n, g(t)> = c_n$. Then the action of $\pi_1\langle\land V\rangle$ in $\pi_3\langle\land V\rangle$ is given by
 $$g(t)\bullet \exp x= e^t\cdot g(t).$$

\vspace{3mm} 
 Finally, the connection between the properties of $H(\land V)$ and $L_V$ is illustrated in
 
\begin{Prop}
\label{p8.3} The following conditions on a minimal Sullivan algebra $(\land V,d)$ are equivalent:
\begin{enumerate}
\item[(i)] $L_V = \varprojlim_{n\in \mathbb N} L(n)$ in which each $L(n)$ is a finite dimensional nilpotent Lie algebra.
\item[(ii)] dim$\, V$ is at most countable.
\item[(iii)] dim$\, H(\land V,d)$ is at most countable.
\end{enumerate}
\end{Prop}

\vspace{3mm}\noindent {\sl proof.} (i) $\Rightarrow $ (ii) $\Rightarrow$ (iii). Let $\land V_n$ be the quadratic Sullivan model of $L(n)$. Then $(\land V, d_1) = \varinjlim_n (\land V_n, d_1)$ is at most countable. This proves (ii) and (ii) $\Rightarrow$ (iii) is obvious.

(iii) $\Rightarrow$ (i). Define subspaces $W_n\subset V$ by $W_0 = V\cap \mbox{ker}\, d$ and $W_{n+1} = V\cap d^{-1}(\land W_n)$. It follows from (iii) that each dim$\, W_n$ is at most countable. But because $(\land V,d)$ is a Sullivan algebra it follows that $V = \cup_n W_n$. and so dim$\, V$ is at most countable. Thus because $(\land V,d)$ is a Sullivan algebra $V= \cup_\alpha V_\alpha$ where dim$\, V_\alpha<\infty$ and $\land V_\alpha$ is preserved by $d$. Since dim$\, V$ is at most countable this implies that $V = \varinjlim_n V_n$ where $(\land V_n,d)$ is a sub Sullivan algebra and dim$\, V_n<\infty$;  Thus $L_V = \varprojlim_n L(n)$, $L(n)$ denoting the homotopy Lie algebra of $(\land V_n,d)$. 
\hfill$\square$

\vspace{3mm}\noindent {\bf Corollary.} If $(\land V,d)$ is the minimal Sullivan of a space $X$ and satisfies the conditions of the Proposition then $H(X)$ is a graded vector space of finite type.

   \subsection{Profree homotopy Lie algebras}
 
 Profree Lie algebras are characterized in Part II in terms of their quadratic Sullivan models. This can, however, be extended to general minimal Sullivan algebras (Proposition \ref{p8.4} below), and they arise, with their sub Lie algebras,  in a natural way as homotopy Lie algebras of wedges of spheres.
 
 \begin{Prop}
 \label{p8.4} Denote by $L$ and by $(\land V,d)$ the homotopy Lie algebra and associated quadratic Sullivan algebra of a minimal Sullivan algebra $(\land V,d)$. The following conditions are then equivalent
 \begin{enumerate}
\item[(i)] There is a quasi-isomorphism $(\land V,d) \stackrel{\simeq}{\to} \mathbb Q \oplus S$ with zero differential in $S$ and $S\cdot S = 0$.
\item[(ii)] The generating space $V$ can be chosen so that $d : V\to \land^2V$ and 
$$\mathbb Q \oplus (V\cap \mbox{ker}\, d) \stackrel{\cong}{\longrightarrow} H(\land V,d).$$
\item[(iii)] $\mathbb Q \oplus V\cap \mbox{ker}\, d_1 \stackrel{\cong}{\to} H(\land V, d_1)$.
\item[(iv)] $L= \overline{\mathbb L}_T$ is profree, with $sT= H^{\geq 1}(\land V)^\vee=S^\vee$ defining the enriched structure.
\end{enumerate}
If they hold then $(\land V,d) \cong (\land V, d_1)$  and, if $V$ is chosen to satisfy (iii), then $V\cap \mbox{ker}\, d = V\cap \mbox{ker}\, d_1$.
 \end{Prop}

\vspace{3mm}\noindent {\sl proof.} (i) $\Leftrightarrow$ (ii). First, let $S = S^{\geq 1}$ be any graded vector space. Then $\mathbb Q \oplus S$, with   zero differential and $S\cdot S = 0$ is a cdga, and its minimal model has the form 
$$(\land V,d) \stackrel{\simeq}{\longrightarrow} (\mathbb Q \oplus S,0).$$

On the other hand, a successive adjoining of variables $w$ to $S$ constructs a quadratic Sullivan algebra $(\land W, d_1)$ such that $H^{[1]}(\land W, d_1) = S$ and $H^{[2]}(\land W) = 0$. (Here as in Part II, $H^{[k]}(\land W, d_1)$ denotes the homology classes represented by cycles in $\land^kW$.) Now Lemma \ref{l6.2} asserts that $H^{[k]}(\land W)= 0$, $k\geq 2$. Thus division by $\land^{\geq 2}W$ and by a direct summand of $S$ in $W$ defines a quasi-isomorphism
$$(\land W, d_1) \stackrel{\simeq}{\to} \mathbb Q \oplus S.$$
Since minimal models are unique it follows that there is an isomorphism
$$(\land W, d_1) \cong (\land V,d).$$
Now choose $V$ so this isomorphism takes $W\stackrel{\cong}{\to} V$. 

Now, if (ii) holds,  then with the given choice of $V$, division by $\land^{\geq 2}V$ and by a direct summand of $V\cap \mbox{ker}\, d$ in $V$ defines a quasi-isomorphism 
$$(\land V,d) \stackrel{\simeq}{\to} \mathbb Q \oplus (V\cap \mbox{ker}\, d).$$
This gives (i).

(i) $\Leftrightarrow$ (iii). If (i) holds, 
 the observations above yield an isomorphism $(\land W, d_1) \cong (\land V,d)$. This induces an isomorphism $(\land W, d_1)\stackrel{\cong}{\to} (\land V, d_1)$, which for appropriate choice of $V$  preserves wedge degree.   It follows that with this choice $(\land V,d)$ is quadratic and $V\cap \mbox{ker}\, d= V\cap \mbox{ker}\, d_1$.  This proves (iii), and the same argument shows that (iii) $\Rightarrow$ (ii).
 
Finally, the assertion (iii) $\Leftrightarrow $ (iv) is Theorem \ref{t1}. \hfill$\square$

\vspace{3mm}\noindent {\bf Corollary 1.} If  $H^{\geq 1}(X)\neq 0$ then the homotopy Lie algebra $L_X$ is profree if and only if $X$ is formal and $H^{\geq 1}(X)\cdot H^{\geq 1}(X) \neq 0$. (A space $X$ is \emph{formal} if its minimal Sullivan model is also a minimal Sullivan model of $(H(X),0)$.)

\vspace{3mm}\noindent {\bf Corollary 2.} If cat$\, X=1$ then the homotopy Lie algebra, $L_X$ is profree.

\vspace{3mm}\noindent {\sl proof.} It follows from (\cite[Formula (9.1) and Theorem 9.1]{RHTII}) that if cat$\, X= 1$ then its minimal Sullivan model satisfies condition (i) of the Proposition.  \hfill$\square$

\vspace{3mm}\noindent {\bf Corollary 3.} If a connected co-H-space, $X$, satisfies $H^{\geq 1}(X)\neq 0$, then $L_X$ is profree.

\vspace{3mm}\noindent {\sl proof.} Let $U$ be a contractible open neighbourhood of a base point in $X$. Then $X\vee U$ and $U\vee X$ form an open cover for $X\vee X$
 
Since  X is a co-H-space, the diagonal map $\Delta: X\to X\times X$ factors up to homotopy through $X\vee X$,
$$\xymatrix{
X \ar[rrd]^\Delta \ar[rr]^f && X\vee X\ar[d]\\
&& X\times X.}$$
The open sets $f^{-1}(X\vee U)$ and $ f^{-1}(U\vee X)$ make then a covering of $X$ by contractible open sets, and 
 cat$\, X= 1$. Now the result follows from Corollary 2. \hfill $\square$

 \section{The Sullivan completion, $X_{\mathbb Q}$, of  a space $X$}
 
\emph{ Throughout \S 9 we consider a connected space, $X$ with a minimal Sullivan model}
 $$\varphi : (\land V,d)\stackrel{\simeq}{\longrightarrow} A_{PL}(X).$$

\vspace{3mm}\noindent {\bf Definition.} The adjoint map
 $$\vert\varphi\vert : X\to X_{\mathbb Q}:=\langle\land V\rangle,$$
  is the \emph{Sullivan completion} of $X$. 

\begin{Prop}
\label{p9.1}
\begin{enumerate}
\item[(i)] $\pi_*\vert \varphi\vert$ converts Whitehead products in $\pi_*(X)$ to constructions in $L_X$ as set out in Proposition \ref{p8.2}. 
\item[(ii)] 
 {\rm (\cite[Theorem 15.11]{RHTII})}. If $X$ is a simply connected space and $H(X)$ is a graded vector space of finite type then $\pi_*\vert \varphi\vert$ factors as
$$\pi_*(X)\to \pi_*(X)\otimes \mathbb Q \stackrel{\cong}{\longrightarrow} \pi_*(X_{\mathbb Q}).$$
\end{enumerate}
\end{Prop}

\vspace{3mm} The definition above depends on the choice for any space $X$ of a minimal Sullivan model $\land V$. However, since two minimal Sullivan models $\land V, \land V'\stackrel{\simeq}{\longrightarrow} A_{PL}(X)$ are connected by a unique homotopy class of isomorphisms, it follows that  for any minimal model $\land V'\stackrel{\simeq}{\longrightarrow} A_{PL}(X)$, the  homeomorphisms $\langle\land V'\rangle\stackrel{\cong}{\longrightarrow} X_{\mathbb Q}$ are connected by unique homotopy classes of isomorphisms. This identifies the vector spaces $V'$, the groups $\pi_*\langle\land V'\rangle$, and the enriched Lie algebras $L_{V'}$, and we write
$$V_X=V', \pi_*(X_{\mathbb Q})= \pi_*\langle\land V'\rangle, \hspace{3mm}\mbox{and } L_X= L_{V'}.$$

\vspace{3mm}\noindent {\bf Definition.}  \emph{$L_X$ is the homotopy Lie algebra of $X$}. 

\vspace{3mm} In particular, (\ref{i14}) specializes to   
\begin{eqnarray}
\label{i15}
\pi_*(X_{\mathbb Q})=V_X^\vee = sL_X.
\end{eqnarray}
Moreover, a continuous map $f : Y\to X$ has a unique homotopy class of Sullivan representatives $\varphi : \land V_Y\leftarrow \land V_X$ between the respective minimal Sullivan models. Now it follows from (\ref{i14}) that
\begin{eqnarray}
\label{i16}
\pi_*(f_{\mathbb Q}) = \pi_*\langle \varphi\rangle= sL_f.
\end{eqnarray}

 \vspace{3mm}This
  makes $\pi_*(X_{\mathbb Q})$ effectively computable in many cases, which in turn can provide information about $\pi_*(X)$. In particular,
  Proposition \ref{p8.2} provides an explicit description of the Whitehead products  in $\pi_*(X_{\mathbb Q})$ and of the group structure of $\pi_1(X_{\mathbb Q})$ in terms of the minimal model $(\land V,d)$ of $X$.   
  
  However $H(X_{\mathbb Q})$ remains mysterious, and it is unknown what conditions are necessary and or sufficient for $H^n(X_{\mathbb Q};\mathbb Z)$, $n\geq 1$, to be a rational vector space. We do have
 
 \begin{Prop}
 \label{p9.2} Suppose $\varphi : (\land V,d)\stackrel{\simeq}{\to} A_{PL}(X)$ is a minimal Sullivan model of a space $X$. Then  $(\land V,d)$ is a retract of the minimal Sullivan model of $X_{\mathbb Q}$, and so   $H(X)$ is a retract of $H(X_{\mathbb Q})$.
 
\end{Prop}

\vspace{3mm}\noindent {\sl proof.}  
Let $(\land W,d) \stackrel{\simeq}{\to} A_{PL}\langle\land V\rangle$ be the minimal Sullivan model. Since (\ref{i13}) $A_{PL}\vert\varphi\vert \circ \iota_{\land V}\sim \varphi$   the homotopy commutative diagram
$$\xymatrix{ && A_{PL}\langle\land V\rangle \ar[rr]^{A_{PL}\vert\varphi\vert} && A_{PL}(X)\\
\land V \ar[rr]\ar[rru]^{\iota_{\land V}}  && \land W\ar[u]^\simeq \ar[rr] && \land V \ar[u]^\simeq_\varphi}$$
  exhibits $(\land V,d)$ as a homotopy retract of $(\land W,d)$.   \hfill$\square$

\vspace{3mm}\noindent {\bf Corollary.} Any continuous map $f : Y\to X_{\mathbb Q}$ extends to a map $Y_{\mathbb Q}\to X_{\mathbb Q}.$

\vspace{3mm}\noindent {\sl proof.} Compose $f_{\mathbb Q}$ with the retraction $(X_{\mathbb Q})_{\mathbb Q}  \to X_{\mathbb Q}$. \hfill$\square$
 
 \vspace{3mm}  Now observe that the retraction
$$\xi : \land^{\geq 1} V\to V$$
with kernel $\land^{\geq 2}V$ satisfies $\xi\circ d= 0$ and so induces a map $H(\xi): H^{\geq 1}(\land V)\to V$. 
On the other hand, denote by $hur: \pi_*(X)\to H_*(X;\mathbb Z)$ the Hurewicz homomorphism, and recall that
$$\pi_*(\vert\varphi\vert) : \pi_*(X)\to \pi_*(X_{\mathbb Q})= V^\vee.$$

\begin{Prop}
\label{p9.3} (\cite[Proposition 1.19]{RHTII}) For $\alpha \in H^{\geq 1}(\land V)$ and $\beta\in \pi_*(X)$,
$$<H(\varphi)\alpha, hur\, \beta> = <H(\xi)(\alpha) , \pi_*(\vert\varphi\vert)\beta>;$$
\end{Prop}

 \vspace{3mm} Moreover, because of the naturality of the Hurewicz homomorphism, this follows from a commutative diagram valid for any minimal Sullivan algebra, $\land V$. Recall from \S 7 the natural morphism $\iota_{\land V} : \land V \to A_{PL}\langle\land V\rangle$ adjoint to $id_{\langle\land V\rangle}$. Proposition \ref{p8.2} specializes to (and follows from) the commutative diagram
 \begin{eqnarray}
 \label{i17}
 \xymatrix{\pi_*\langle\land V\rangle \ar[d]^{hur}\ar[rr]^{\cong} && V^\vee\ar[d]^{H(\xi)^\vee}\\
 H_*(\langle\land V\rangle;\mathbb Q) \ar[rr]^{H(\iota_{\land V})^\vee} && H^{\geq 1}(\land V)^\vee}
 \end{eqnarray}
 in which the upper isomorphism is the isomorphism (\ref{i12}). Note that the generating space $V$ can always be chosen so that $\xi: \mbox{ker}\, d \to V\cap\mbox{ker}\, d$. In this case in the diagram above $H^{\geq 1}(\land V)^\vee$ may be replaced by $(V\cap \mbox{ker}\, d)^\vee$.

  \subsection{The Sullivan completion of a free group, $G$}
 
 Suppose $b_1, \dots , b_s$ freely generate a free group $G$. The classifying space for $G$ is then
 $$X = S^1_{b_1} \vee \dots \vee S^1_{b_s},$$
 and $b_i$ is the homotopy class of the inclusion of $S^1_{b_i}$ in $X$. In particular, $\pi_*(X)= \pi_1(X)= G$. 
 
 On the other hand, (Corollary 2 to Proposition \ref{p8.4}), the minimal Sullivan model of $X$ satisfies the conditions of Proposition \ref{p8.4}. Therefore it is a quadratic Sullivan algebra $(\land V, d_1)$ and the inclusion
 $$j : V\cap \mbox{ker}\, d_1  {\longrightarrow} \land^{\geq 1}V$$
 induces an isomorphism $V\cap \mbox{ker}\, d_1\stackrel{\cong}{\to} H^1(\land V)= H^1(X)$.
 It follows that $V= V^1$.
 
 Now recall  that the quasi-isomorphism $\varphi : (\land V, d_1)\stackrel{\simeq}{\to} A_{PL}(X)$ induces via adjointness a continuous map
 $$\vert \varphi\vert : X\to \langle \land V\rangle = X_{\mathbb Q}.$$
 
 \vspace{3mm}\noindent {\bf Definition.} The group morphism $\pi_1\vert \varphi\vert : G= \pi_1(X) \to \pi_1(X_{\mathbb Q})$ is the \emph{Sullivan completion} of $G$.
 
 \vspace{3mm} We now recall three important constructions:
 \begin{enumerate}
 \item[$\bullet$] The \emph{homotopy Lie algebra}, $L_V$, of $(\land V, d_1)$, a complete enriched Lie algebra, defined (\S 8) via the quadratic differential, $d_1$. By construction
 $$sL_V = V^\vee$$
 as graded vector spaces.
 \item[$\bullet$] The \emph{group} $G_{L_V}$ (\S 1.2)  connected with $L_V$ by natural inverse bijections $$
 \xymatrix{L_V\ar@/^/[rr]^\exp && G_{L_V}\ar@/^/[ll]^\log.}$$
 \item[$\bullet$] A \emph{natural bijection of sets} (12)
 $$\tau_V: \pi_1\langle \land V\rangle \stackrel{\cong}{\longrightarrow} V^\vee.$$
 With $L_V$ and $V^\vee$ identified as above, (\cite[Theorem 2.4]{RHTII}) provides a commutative diagram,
 \begin{eqnarray}
 \label{star}
 \xymatrix{ G \ar[rr]^{\pi_1\vert \varphi\vert} \ar[rrd] && \pi_1\langle \land V\rangle \ar[d]^\cong_{\sigma_V}\ar[rrr]^{\tau_V}_\cong &&&  L_V\\&&
 G_{L_V},\ar[rrru]^\cong_{\log}
 }
 \end{eqnarray}
 in which $\sigma_V$ is a natural isomorphism of groups. We denote by $\log_G$ the composite
 $$\log_G : G\longrightarrow G_{L_V}\stackrel{\log}{\longrightarrow} L_V.$$
\end{enumerate}
 
  The constructions above combine to form the commutative diagram,
$$
  \xymatrix{
  G \ar[rrr]^{\tau_V\circ\pi_1\vert\varphi\vert}\ar[dd] &&& V^\vee \ar[rrr]^=\ar[dd]^{j^\vee} &&& sL_V\ar[dd]\\ & (1)&&& (2)\\
  G/G^2\otimes \mathbb Q \ar[rrr]_\cong &&& (V\cap \mbox{ker}\, d_1)^\vee \ar[rrr]_\cong &&& sL_V/L_V^2,}
$$
  which we now establish.
  
  First, the isomorphism $V\cap \mbox{ker}\, d_1 \stackrel{\cong}{\to} H^1(X)$ dualizes to an isomorphism $H_1(X;\mathbb Q)\stackrel{\cong}{\to} (V\cap \mbox{ker}\, d_1)^\vee$. The composite
  $$G \to V^\vee \stackrel{j^\vee}{\longrightarrow} (V\cap \mbox{ker}\, d_1)^\vee \stackrel{\cong}{\longrightarrow} H_1(X;\mathbb Q)$$
  is a set map $G\to H_1(X;\mathbb Q)$. It follows from Proposition \ref{p9.3} that this is $hur\otimes \mathbb Q$, where $hur$  is the classical Hurewicz homomorphism. Thus it induces an isomorphism $G/G^2\otimes \mathbb Q\stackrel{\cong}{\to} (V\cap \mbox{ker}\, d_1)^\vee$. This establishes $(1)$.
  
  Second, the quasi-isomorphism $\mathbb Q\oplus (V\cap \mbox{ker}\, d_1) \stackrel{\simeq}{\longrightarrow} (\land V, d_1)$ admits a cdga retraction $\rho : (\land V, d_1) \stackrel{\simeq}{\longrightarrow} \mathbb Q\oplus (V\cap \mbox{ker}\, d_1)$, where $V\cap \mbox{ker}\, d_1$ is given the trivial product. Thus Proposition \ref{p8.4} identifies $L_V$ as a profree Lie algebra $\overline{\mathbb L}_{T'}$, in which $\rho^\vee : (V\cap \mbox{ker}\, d_1)^\vee \stackrel{\cong}{\longrightarrow}  T'$. Since $j^\vee \circ \rho^\vee= id_{(V\cap \mbox{ker}\, d_1)^\vee}$ it follows that ker$\, j^\vee$ is a direct summand of $T'$. But for $v\in V\cap \mbox{ker}\, d_1$, and for $x,y\in L_V$,
  $$<v, j^\vee s[x,y]> = <d_1j (v), s[x,y]>= 0.$$
  Therefore $L_V^{(2)}\subset \mbox{ker}\, j^\vee$ and so $L^{(2)} = \mbox{ker}\, j^\vee$. This establishes (2).

  \vspace{3mm}  Next, let $\overline{J}$ denote the augmentation ideal of $\overline{UL_V}$ (\S 1.2) and recall that the right adjoint representation of $L_V$ extends to a representation of $\overline{J}$.
  
  \begin{Prop}
  \label{p9.4} (i) The elements $\log_Gb_i$ are the basis of a subspace $T\subset L_V$ for which $L_V= \overline{\mathbb L}_T$.  
  
  (ii) If $c,d\in G$ then
    $$\log_Gcd\in \log_Gc +\log_Gd + \frac{1}{2} [\log_Gc, \log_Gd] + [\log_Gc, \log_Gd]\cdot \overline{J},$$
  and
  $$\log_G(cdc^{-1}d^{-1}) \in [\log_Gc, \log_Gd] + [\log_Gc, \log_Gd]\cdot \overline{J}.$$
  \end{Prop}
  
  \vspace{3mm}\noindent {\sl proof.} (i)    follows from (18),   diagrams (1) and (2) above, and the commutative diagram (\cite[Theorem 2.2]{RHTII}),
  $$\xymatrix{
  G_{L_V} \ar[d]\ar[rr]^\log && L_V\ar[d]\\
  G_{L_V}/G_{L_V}^2 \ar[rr]^\log && L_V/L_V^{(2)},}$$
and diagrams (1) and (2) above.

(ii) A naturality argument reduces this to the case $G$ is freely generated by $c$ and $d$. In this case $T= \mathbb Q \log_Gc \oplus \mathbb Q\log_Gd$, $[T,T= \mathbb Q [\log_Gc, \log_Gd]$ and $L_V^{(3)} \subset [\log_Gc, \log_Gd]\cdot \overline{J}$. Thus (ii) follows immediately from (\cite[Proposition 2.4]{RHTII}).
  \hfill$\square$
  
  \subsection{Remarks}

\vspace{3mm}   1. Let $\varphi : \land V\to A_{PL}(X)$ and $\psi : \land W\to A_{PL}(Y)$ be minimal Sullivan models. Then cf (\ref{i11}) we have bijections
$$[Y, X_{\mathbb Q}] \cong [\land V, A_{PL}(Y)] \cong [\land V, \land W].$$

2. By \cite[Proposition 17.13]{FHTI} Sullivan completions induce a morphism
$$[\land V, \land W]\to [\langle\land W\rangle, \langle\land V\rangle].$$
This morphism is not always injective. For instance, if $Y$ is an infinite wedge of spheres $S^3$, and $X= S^3$, then $\land V= (\land u, 0)$ with deg$\, u = 3$, and $\land W\simeq (H(Y), 0)$. Therefore
$$[\land V, \land W] \cong W^3\, \hspace{3mm}\mbox{and } [\langle\land W\rangle, \langle\land V\rangle] \cong H^3\langle\land W\rangle =(\pi_3\langle\land W\rangle)^\vee= \left( (W^3)^\vee\right)^\vee.$$ Suppose $\alpha \in H^3\langle\land W\rangle$ is in the kernel of $H(\vert \psi\vert) : H(
\langle\land W\rangle)\to H(Y)$. This gives a map $g : Y_{\mathbb Q}\to X_{\mathbb Q}$ with $g\circ \vert\psi\vert \sim * \sim *\circ \vert\psi\vert$.

3. When $Y$ is a simply connected CW complex of finite type, then (\cite[Theorem 17.12]{FHTI}) $\vert\psi\vert : Y\to Y_{\mathbb Q}$ is an isomorphism in rational cohomology and classical obstruction theory shows that in that case the morphism $[\land V, \land W]\to [\langle\land W\rangle, \langle\land V\rangle]$ is a bijection.

\vspace{3mm} 4.  Let $X= \langle\land V\rangle$ be the spatial realization of the Sullivan quadratic model of the profree Lie algebra $\overline{\mathbb L}_T$ where $T=T_1$ is a countably infinite dimensional vector space. Then $H_3(X)\neq 0$ and $H_4(X)\neq 0$.

\vspace{2mm}  We only provide the proof for $H_3(X)$. The proof for $H_4(X)$ is similar. 
Let $\omega\in \pi_3(X)$ be a homotopy class corresponding by suspension to an element in $L_2$, but not in  $[L_1, L_1]$. Then $X$ contains a simply connected finite CW complex $Y$ such that $\omega$ is in the image of $\pi_3(Y)$. We denote this element   by $\omega_Y\in \pi_3(Y)$. It follows that   $\omega_Y$ is not decomposable, and by Lemma \ref{l9.1} below, the image $a_Y$ of $\omega_Y$ in $H_3(Y)$ is non zero.

For any finite sub CW complex $Z$ of $X$ containing $Y$, the element $\omega_Z$ induces also a non-zero element $a_Z$ by the Hurewicz map, and by naturalness of the Hurewicz map, the map $H_3(Y)\to H_3(Z)$ maps $a_Y$  to $a_Z$. Now $X$ is the union of the finite sub CW complexes $Z_\alpha$ containing $Y$. It follows that $H_3(X) = \varinjlim_\alpha H_3(Z_\alpha)$ and therefore the family $a_{Z_\alpha}$ induces a nonzero element in $H_3(X)$. 
\hfill$\square$

\begin{lem}
\label{l9.1} Let $L$ be the homotopy Lie algebra of a finite type simply connected CW complex, $Y$.
 \begin{enumerate}
 \item[(i)] If $L_2 \neq [L_1, L_1]$, then the Hurewicz map $\pi_3(Y)\to H_3(Y)$ is non zero. 
 \item[(ii)] If $L_3\neq [L_1, L_2]$, then the Hurewicz map $\pi_4(Y)\to H_4(Y)$ is non zero. 
\end{enumerate}
\end{lem}

\vspace{3mm}\noindent {\bf Proof.} Denote by $(\land V,d)$ the Sullivan minimal model of $Y$. Then $V$ is a finite type vector space. Therefore the brackets $[\,,\,] :  \land^2 L_1   \to L_2$, and $L_1\otimes L_2\to L_3$ are dual to the differentials $d : V^3\to \land^2V^2$, and $V^4\to V^2\otimes V^3$. Consider the situation (i); the proof for (ii) is similar. Since $L_2\neq [L_1, L_1]$, $d$ is not injective. This implies that $H^3(\land V) \to H^3(\land V/\land^{\geq 2}V)$   is non zero. But this last map is the dual of the Hurewicz map. Thus if $\omega\in \pi_3(Y) \otimes \mathbb Q = (V^3)^\vee$ with $\omega (v)\neq 0$, then hur $(\omega)\neq 0$. 
\hfill$\square$

\subsection{$\pi_1\vert\varphi\vert : \pi_1(X)\to \pi_1(X_{\mathbb Q})$}

Recall that for any group, $G$, $G^n$ denotes the subgroup generated by iterated commutators of length $n$, and that for any graded Lie algebra, $L$, $L^n$ is the span of the iterated commutators of length $n$. In particular a normal subgroup $Tor (G)$ is defined by
$$Tor(G) = \{ a\in G\, \vert\, \mbox{for each $n$ there is a $k$ such that }a^k\in G^n\}.$$
(Thus if $G$ is abelian, $Tor(G)$ is the subgroup of torsion elements.)

\begin{Prop}
\label{p9.5}
If $H_1(X;\mathbb Z)$ is   finitely generated   then
$$\mbox{ker}(\pi_1\vert\varphi\vert : \pi_1(X)\to \pi_1(X_{\mathbb Q})) = Tor (\pi_1(X)).$$
\end{Prop}

\vspace{3mm}\noindent {\sl proof.} For simplicity we denote $\pi_1(X)$ by $G$. Then the classifying map $X\to BG$ induces an isomorphism in $H_1(-;\mathbb Z)$ and, by a theorem (\cite[Theorem 5.2]{brown}) of H. Hopf, a surjection in $H_2(-;\mathbb Q)$. Therefore it induces an isomorphism in $H^1(-;\mathbb Q)$ and an injection in $H^2(-;\mathbb Q)$. In particular the corresponding morphism of minimal Sullivan models restricts to an isomorphism
$$(\land V^1,d) \stackrel{\cong}{\longleftarrow} (\land W^1,d).$$
This identifies $\pi_1(X_{\mathbb Q})$ with $\pi_1((BG)_{\mathbb Q})$ and reduces the Proposition to the case $X = BG$.

Denote by $L$ the component of degree $0$ in the homotopy Lie algebra of $BG$, and by $G_L$ the group $\pi_1((BG)_{\mathbb Q})$. Since (by hypothesis) $G/G^2$ is finitely generated it follows that each $G^n/G^{n+1}$ is finitely generated. In particular,
$$\mbox{ker}\, (G^{n-1}/G^n \to G^{n-1}/G^n \otimes_{\mathbb Z}\mathbb Q)$$
is a torsion group. On the other hand (\cite[Theorem 7.5]{RHTII}), $\log \circ \pi_1\vert\varphi\vert$ induces isomorphisms
$$\psi_n : G^{n-1}/G^n \otimes_{\mathbb Z}\mathbb Q \stackrel{\cong}{\longrightarrow} L^{n-1}/L^n.$$

Now suppose $a\in \mbox{ker}\, \pi_1\vert\varphi\vert$ and that $a^k\in G^{n-1}$. Let $b =$ image of $a^k$ in $G^{n-1}/G^n$. Then $\psi_nb$ is the image of $\log \pi_1\vert\varphi\vert a^k = k \log \pi_1\vert\varphi\vert a= 0$. Thus $b$ is in the kernel of $G^{n-1}/G^n \to G^{n-1}/G^n \otimes_{\mathbb Z} \mathbb Q$, and so
$$a^{\ell k}= b^\ell\in G^n.$$
It follows by induction that $a\in Tor (G)$.

In the reverse direction the same argument shows that if $a\in Tor (G)$, then for each $n$ and some $k$, 
$$k\log \pi_1\vert\varphi\vert a = \log \pi_1\vert\varphi\vert a^k\in L^n.$$
It follows that $\log \pi_1\vert\varphi\vert a\in \cap_n L^n = 0$, and thus $a\in \mbox{ker}\pi_1\vert\varphi\vert$. \hfill$\square$

\subsection{The radical and center of $\pi_*(X_{\mathbb Q})$}

As recalled in \S 6, the category, cat$\, (\land V)$, is the least $m$ (or $\infty$) for which $\land V$ is a homotopy retract of $\land V/\land^{>m}V$. Moreover (\cite[Theorem 9.2]{RHTII}),
$$\mbox{cat}\, (\land V) \leq \mbox{cat}\, X,$$
where cat$\, X$ is the classical Lusternik-Schnirelmann category of $X$. On the other hand, the \emph{radical}, rad$\, L_V$, of $L_V$ is the sum of the solvable ideals of $L_V$. Thus combining (\cite[Theorem 1]{Sdepth}) with (\cite[Theorem C]{Malcev}) yields 

\begin{Prop}
\label{p9.6} Suppose cat$\, (X)<\infty$. Then
\begin{enumerate}
\item[(i)] The sum of the solvable ideals, rad$\, L_V$, is finite dimensional and
$$\mbox{dim(rad}\, (L_V)_{even}) \leq \mbox{ cat}\, X.$$
\item[(ii)] The union of the normal solvable subgroups of $\pi_1(X_{\mathbb Q})$ is a normal solvable subgroup of solv length $\leq$ cat$\, (X)$.
\end{enumerate}
\end{Prop}

\vspace{3mm}\noindent {\sl proof.} (i) is immediate from the remarks above. For the proof of (ii) note that rad$(L_V)_0$ is a closed ideal since it is finite dimensional. Thus $\exp$(rad$(L_V)_0$) is a normal subgroup of $\pi_1(X_{\mathbb Q})$ whose quotient is $\exp$($(L_V)_0/\mbox{rad}(L_V)_0)$. Therefore, since the center of $(L_V)_0/ \mbox{rad}(L_V)_0$ is zero, $\exp$($L_V/\mbox{rad}(L_V)_0)$ has no center and therefore no normal solvable subgroup. \hfill$\square$

 \vspace{3mm}Now define $C_*(X_{\mathbb Q})\subset \pi_*(X_{\mathbb Q})$ to be the subgroup of elements $\alpha$ satisfying
 $$\mbox{Wh}(\alpha, \beta) = 0, \hspace{5mm} \beta\in \pi_*(X_{\mathbb Q}).$$
 It follows from Proposition \ref{p8.2} that $C_*(X_{\mathbb Q})$ corresponds under $\exp$ to the centre of the homotopy Lie algebra $L_V$. Thus Proposition \ref{p9.6} has the 
 
 \vspace{3mm}\noindent {\bf Corollary.} If cat$\, X<\infty$ then the group $C_*(X_{\mathbb Q})$ is a finite dimensional rational vector space.

  \subsection{The Gottlieb groups $G_*(X)$ and $G_*(X_{\mathbb Q})$}

The \emph{$n^{th}$ Gottlieb group}, $G_n(X)$, is the subgroup of $\pi_n(X)$ formed by the homotopy classes $[f]$ of  maps $f : S^n\to X$ such that $f\vee \mbox{id} : S^n\vee X\to X$ extends to the product $S^n\times X$. The Whitehead products of $\alpha \in G_*(X)$ and $\beta\in \pi_*(X)$ are zero. In particular, $G_*(X_{\mathbb Q})\subset C_*(X_{\mathbb Q})$. 

 More generally, if $g : Y\to X$ is a continuous map, then the $n^{th}$ Gottlieb group of $g$, $G_n(g)$ is the subgroup of $\pi_n(X)$ formed by the homotopy classes of maps $ f : S^n\to X$ such that the map $f\vee g : S^n\vee Y\to X$ extends to $S^n \times Y$. 
 
 Finally, recall from (\ref{i11}) that for any minimal Sullivan algebra $\land W$,
 $$\pi_n\langle\land W\rangle = [\land W, \mathbb Q \oplus \mathbb Q a_n]$$
 where $a_n$ is the fundamental cohomology class in $H^n(S^n)$. We define the \emph{Gottlieb group},
 $$G_n(\land W) \subset \pi_n\langle\land W\rangle$$ to be the subgroup formed by those elements represented by morphisms $\sigma : \land W\to \mathbb Q \oplus \mathbb Q a_n$ such that $\sigma \times id : \land W\to (\mathbb Q \oplus \mathbb Q a_n) \times_{\mathbb Q} \land W$ factors as
 $$\land W \to (\mathbb Q \oplus \mathbb Q a_n)\otimes \land W \to (\mathbb Q \oplus \mathbb Q a_n) \times_{\mathbb Q} \land W.$$
 Note here that
 $$(\mathbb Q \oplus \mathbb Q a_n)\times_{\mathbb Q}\land W = \mathbb Q a_n \oplus \land W$$
 and that $\sigma$ is just the map
 $$W = \land^{\geq 1}W/\land^{\geq 2}W \to \mathbb Q a_n.$$

\begin{Prop} \label{p9.7} 
\begin{enumerate}
\item[(i)] $\pi_*\vert\varphi\vert : G_*(X)\to G_*(X_{\mathbb Q})$.
\item[(ii)] $G_*(\land V)=  G_*(\vert\varphi\vert)=  G_*(X_{\mathbb Q})$. 
\item[(iii)] $G_*(X_{\mathbb Q}) \subset C_*(X_{\mathbb Q})$. 
\end{enumerate}
\end{Prop}

\vspace{3mm}\noindent {\sl proof.} (i) If $f: S^n \to X$ represents an element of $G_*(X)$ we may extend $f\vee id_X$ 
to $g : S^n\times X\to X$. Now Proposition \ref{p8.3}(iii) identifies $(S^n\times X)_{\mathbb Q}= S^n_{\mathbb Q}\times X_{\mathbb Q}$. Thus we obtain $g_{\mathbb Q} : S^n_{\mathbb Q} \times X_{\mathbb Q}\to X_{\mathbb Q}$. It follows from (\ref{i14}) that $g_{\mathbb Q}$ restricts to $\vert \varphi\vert \circ f : S^n \to X_{\mathbb Q}$. Thus $g_{\mathbb Q}: S^n\times X_{\mathbb Q}\to X_{\mathbb Q}$ extends $\vert\varphi\vert \circ f \vee id_{X_{\mathbb Q}}$.

(ii) We proceed in two steps.

\vspace{2mm}\noindent {\bf Step One. 
$G_*(\land V) = G_*\langle\land V\rangle$.
}

\vspace{2mm} Let $\land Z$ denote the minimal model of $S^n$; Then (\ref{i11})
$$\pi_n\langle\land V\rangle = \mbox{Hom}(V^n, \mathbb Q) = [\land V, \land Z]= [\land V, A_{PL}(S^n)].$$
Thus an element in $G_n(\land V)$ is represented by a morphism $\sigma : \land V\to \land Z$ for which $\sigma\times_{\mathbb Q} id : \land V\to \land Z\times_{\mathbb Q}\land V$ factors up to homotopy as
$$\xymatrix{\land V \ar[rr] \ar[rrd]_{\sigma \times_{\mathbb Q} id} && \land Z\otimes \land V\ar[d]\\
&& \land Z\times_{\mathbb Q}\land V.}$$
The identifications in Proposition \ref{p8.2} then convert this diagram into the homotopy commutative diagram
$$\xymatrix{
\langle\land V\rangle && \langle \land Z\rangle\times \langle\land V\rangle \ar[ll] && S^n \times 
\langle\land V\rangle\ar[ll]\\
&&\langle\land Z\rangle\vee \langle\land V\rangle \ar[llu]^{\langle\sigma \times_{\mathbb Q}id\rangle} \ar[u] && S^n\vee \langle\land V\rangle\ar[ll]\ar[u]}$$
in which the vertical arrows are the standard inclusions. In particular, the map $$S^n \to \langle\land Z\rangle \stackrel{\langle\sigma\rangle}{\longrightarrow} \langle\land V\rangle$$ represents $[\sigma]\in \pi_n(\land V)$.

In the reverse direction, suppose $\tau \in G_n\langle\land V\rangle$. Then we have the homotopy commutative diagram
$$
\xymatrix{
\langle\land V\rangle && S^n\times \langle\land V\rangle\ar[ll] \\
&& S^n\vee \langle\land V\rangle.\ar[llu]^{(\tau, id)}\ar[u]}$$
Let $\land W$ be the minimal Sullivan model of $\langle\land V\rangle$. Then the identifications of Proposition \ref{p8.2}, together with the proof of Proposition \ref{p9.2},  convert this to the homotopy commutative diagram
$$\xymatrix{\land V \ar[r] \ar[rd] & A_{PL}(S^n)\otimes A_{PL}\langle\land V\rangle \ar[d] && \land Z\otimes \land W\ar[ll]_-\simeq\ar[d] \ar[rr] &&\land Z\otimes \land V\ar[d]\\
 & A_{PL}(S^n)\times_{\mathbb Q} A_{PL}\langle\land V\rangle && \land Z\times_{\mathbb Q}\land W \ar[ll]_-\simeq \ar[rr] && \land Z\times_{\mathbb Q}\land V.}$$
 
 \vspace{3mm}\noindent {\bf Step Two. $G_n(\vert \varphi\vert) = G_n(X_{\mathbb Q}$)}.
 
 \vspace{1mm}   The injection   $G_n(X_{\mathbb Q})\subset G_n(\varphi_X)$ is a direct consequence of the following commutative diagram associated to an element $[f]\in G_n(X_{\mathbb Q})$,
$$\xymatrix{
S^n \vee X \ar[d] \ar[rr]^{id\vee \vert\varphi\vert} && S^n \vee X_{\mathbb Q} \ar[d]\ar[rr]^{f\vee id} && X_{\mathbb Q}\\
S^n \times X \ar[rr]^{id \times \vert\varphi\vert} && S^n \times X_{\mathbb Q}. \ar[rru]
}$$

For the reverse inclusion $G_n(\vert \varphi\vert)\subset G_n(X_{\mathbb Q})$ suppose $f: S^n \to X_{\mathbb Q}$ represents an element in $G_n(\vert \varphi\vert)$. Then $\vert\varphi\vert \vee f$ factors as
$$X\vee S^n\to X\times S^n \stackrel{g}{\to} X_{\mathbb Q}.$$
Restricting $g_{\mathbb Q}$ to $X_{\mathbb Q}\times S^n$, and denoting by $r: (X_{\mathbb Q})_{\mathbb Q}\to X_{\mathbb Q}$ the retraction (Proposition \ref{p8.1}(iii)), then gives the commutative diagram,
$$
\xymatrix{
X\vee S^n \ar[d]\ar[rr]&& X\times S^n \ar[d]\ar[rr]^g && X_{\mathbb Q}\\
X_{\mathbb Q}\vee S^n \ar[rr] &&X_{\mathbb Q}\times S^n \ar[rr]^{q_{\mathbb Q}} && (X_{\mathbb Q})_{\mathbb Q}.\ar[u]^r}$$
This gives $[f]\in G_n(X_{\mathbb Q})$.

(iii) This is immediate from the definitions and Proposition \ref{p9.6}(ii).
\hfill$\square$

 \vspace{3mm} Recall next (\cite[\S 3.3]{RHTII}) that a \emph{minimal Sullivan extension} is a morphism,
 $$\land W \to \land W\otimes \land V, \hspace{1cm} w\mapsto w\otimes 1,$$
 of Sullivan algebras in which $\land W$ and the quotient $(\land V, d_V)= \mathbb Q \otimes_{\land W}(\land W\otimes \land V)$ are minimal Sullivan algebras. However the differential $d$ in $\land W\otimes \land V$ may have a component $d_0$ of wedge degree $0$: $d_0 : W \leftarrow V$. The dual of $d_0$ is then a linear map 
 $$\delta :W^\vee \to V^\vee$$
 of degree $-1$. Recall   also that $V^\vee = \mbox{Hom}(V,\mathbb Q)\supset G_*(\land V)$.
 
 \begin{Prop}
 \label{p9.8}
 With the notation and hypotheses above $\delta : W^\vee\to V^\vee$ is injective, and 
 $$\mbox{Im}\, \delta \subset G_*(\land V).$$
 In particular Im$\, \delta$ suspends to a subspace in the center of $L_V$.
 \end{Prop}
 
 \vspace{3mm}\noindent {\sl proof.} Fix $n$ and choose a cdga retraction
 $$\xi : (\land W,d)\to (\mathbb Q \oplus d_0V^n, 0)$$
 with $\land^{\geq 2}W \subset \mbox{ker}\, \xi$. Then extend $\mathbb Q \oplus d_0V^n$ to an acyclic cdga
 $$\land U \otimes (\mathbb Q \oplus d_0V^n) \stackrel{\simeq}{\longrightarrow} \mathbb Q$$
 with $U = U^{\geq n}$, and
 $$d : U^n\stackrel{\cong}{\to} d_0V^n \hspace{5mm}\mbox{and } d : U^{> n}\to \land^{\geq 1}U\otimes d_0V^n.$$
 Now applying $-\otimes_{\land W}(\land W\oplus d_0V)$ defines a surjective quasi-isomorphism
 $$\rho : \land U\otimes (\mathbb Q \oplus d_0V^n)\otimes \land V \stackrel{\simeq}{\longrightarrow} (\land V, d_V).$$
 
 Thus $\rho$ has a left inverse $\sigma : (\land V, d_V)\to \land U \otimes (\mathbb Q\oplus d_0V^n)\otimes \land V$, and it is immediate from the definitions that $\sigma$ may be chosen so that $\sigma v = 1\otimes   v$, $v\in V^{<n}$ and 
 $$\sigma v = 1\otimes v -\alpha (d_0v),$$
 where $\alpha$ is inverse to the isomorphism $U^n \stackrel{\cong}{\longrightarrow} d_0V^n$. Division by $d_0V^n$ and by $U^{>n}$ then yields the sequence of morphisms
 $$(\land V, d_V) \to (\land U^n, 0)\otimes (\land V, d_V) \to (\land U^n, 0) \times_{\mathbb Q}(\land V, d_V).$$
 In particular, any morphism $\land U^n \to H(S^n)$ induces a sequence
 $$(\land V, d_V) \to H(S^n) \otimes (\land V, d_V) \to H(S^n)\times_{\mathbb Q}(\land V, d_V).$$
 In particular, the resulting morphism $(\land V,d_V) \to H(S^n)$ is an element of $G_n(\land V)$.
 
 Finally, this element of $G_n(\land V)$ is determined by the composite $V^n\to U^n\to H^n(S^n)$, and by construction the map $V^n\to U^n$ is given by $v\mapsto -\alpha (d_0v)$. This gives the commutative triangles
 $$
 \xymatrix{
 U^n \ar[rr]^\cong && d_0V^n\\V^n\ar[u] \ar[rru]^{d_0}} \hspace{1cm}\mbox{and } \xymatrix{ (U_n)^\vee \ar[d]  && (d_0V^n)^\vee \ar[ll] \ar[lld]^\delta\\
 (V^n)^\vee,}$$
 This identifies the elements in Im$\, \delta$ as elements in $G_n(\land V)$ and, since $V^n\to d_0V^n$ is surjective, $\delta$ is injective.
 
 \hfill$\square$

\begin{Prop}
\label{p9.9} If $\land V$ is a minimal Sullivan algebra and cat$(\land V)<\infty$, then 
\begin{enumerate}
\item[(i)] $G_{2n}(\land V)= 0$, $n\geq 1$, and
\item[(ii)] dim $G_*(\land V) \leq $ cat$(\land V)$.
\end{enumerate}
\end{Prop}

\vspace{3mm}\noindent {\sl proof.} (i) Let $B = \land a/(a^2)$, with deg$\, a= 2n$. We suppose that there exists a morphism  $\varphi : \land V \to B$ that is a representative of an element of $G_{2n}(\land V)$. If $\varphi \neq 0$ there is an element $z_0\in V^{2n}$ with $\varphi (z_0)= a$.

Now let $(\land W,d) :=  \mathbb Q \otimes_{\land V^{<2n}}(\land V,d)$. By (\cite[Proposition 9.3]{RHTII})  cat$\, (\land W)\leq$ cat$(\land V) <\infty$. We denote by $\overline{\varphi}: \land W\to B$ the map induced by $\varphi$. Since $\overline{\varphi}$ is also in $G_n(\land W)$, we get a morphism
$$\psi : \land W\to \land W\otimes B$$ with $\psi (w)-w\in \land W\otimes B^+$  and $\psi (z_0)= z_0 + a.$
By iteration this gives
$$\psi (q) : \land W \stackrel{\psi}{\longrightarrow} \land W\otimes B \stackrel{\psi\otimes id}{\longrightarrow} \land W\otimes B\otimes B \to \dots \to \land W\otimes B^{\otimes^q}.$$
In particular, $\psi(q) (z_0^q) -a^{\otimes^q} \in \land^+W\otimes B^{\otimes^q}$. By construction $d(z_0)= 0$. Since cat$\land W<\infty$,  there is $y\in \land W$ with $dy= z_0^q$. It follows that  $\psi(q)(z_0^q) $ is a boundary which is impossible. This shows that   $\varphi=0$. 

\vspace{2mm} (ii). Decompose $\land V$ as a Sullivan extension $\land V = \land U^{\leq n}\otimes \land Z^{\geq n}$ in which $G_n(\land V)\to G_n(\land Z)$ is injective and $<v,x>= 1$ for some $x\in G_n(\land Z)$ and some cycle $v\in Z^n$. Since (\cite[Proposition 9.3]{RHTII}) cat$(\land Z)\leq$ cat$(\land V)$, it is sufficient to show that dim$\, G_*(\land Z)\leq$ cat$(\land Z)$. In other words, without loss of generality we may suppose $V = V^{\geq n}$ and that $x\in G_n(\land V)$ and $v$ is a cycle in $V^n$. Since (i) implies that $n$ is odd we also have
$$H(S^n) = \land \widehat{v} \cong \land v.$$
Moreover, since $x\in G_n(\land V)$, by definition we obtain the sequence of morphisms
$$\land V \stackrel{\varphi}{\longrightarrow} \land \widehat{v}\otimes \land V \longrightarrow \land \widehat{v}\times_{\mathbb Q}\land V$$
induced by $x$. 

Now set $W = \{w\in V\, \vert, <w,x>= 0\}$. Then $W = \cup_{k\geq 1} W(k)$ defined by $W(1) = W \cap \mbox{ker}\, d$ and $W(k+1) = W \cap d^{-1}(\land v \otimes \land W(k))$. Suppose next, by induction on $k$ that $\land W(k)$ is preserved by $d$ and that 
$$\varphi w- w\in \widehat{v}\land v \otimes \land W(k), \hspace{1cm} w\in W(k).$$iC1-. 
Finally, fix a basis element $z$ in a direct summand of $W(k)$ in $W(k+1)$. Then by definition,
$$dz= v\Psi_1+\Psi_2\hspace{5mm}\mbox{and } \varphi z= \widehat{v}\land (v\Phi_1+\Phi_2) + z,$$
where $\Phi_1, \Phi_2, \Psi_1$ and $\Psi_2 \in \land W(k)$. This yields
$$d\varphi z= \widehat{v}\, v\, d\Phi_1 - \widehat{v}\, d\Phi_2 + dz$$
$$\varphi dz = (\widehat{v}+v)\varphi \Psi_1 + \varphi \Psi_2 = v\, \Psi_1+\Psi_2 +\widehat{v}\, \Psi_1 + \widehat{v}\, v\, \Gamma,$$
with $\Gamma \in \land W(k)$. Since $\varphi dz=d\varphi z$ and $dz= v\Psi_1+\Psi_2$, it follows that
$$\widehat{v}v\,d\Phi_1- \widehat{v}\, d\Phi_2 = \widehat{v}\, \Psi_1+ \widehat{v}v\Gamma.$$
In other words, $\Psi_1= -d\Phi_2$. 

Now extend the splitting to $W(k+1)$ by replacing $z$ by $z-v\Phi_2$. It is straightforward to check that our conditions are satisfied. Thus by induction we may write $\land V = \land v\otimes \land W$ with $\land W$ preserved by $d$. It follows that $G_*(\land V)= \mathbb Q x \oplus G_*(\land W)$. Since by (\cite[Corollary 9.3]{RHTII}), cat$(\land W)= $ cat$(\land V)-1$, (ii) follows by induction on cat$(\land V)$. \hfill$\square$

 \section{Wedges of topological spaces}

 \subsection{The minimal Sullivan model of $X\vee Y$}

If $X$ and $Y$ are connected spaces then the surjections $H(X\vee Y)\to H(X), H(Y)$ define an isomorphism $H(X\vee Y) \stackrel{\cong}{\to} H(X)\times_{\mathbb Q}H(Y)$. It follows that Sullivan models $\land W\stackrel{\simeq}{\to} A_{PL}(X)$ and $\land Q\stackrel{\simeq}{\to} A_{PL}(Y)$ extend to a quasi-isomorphism
$$\land W\times_{\mathbb Q}\land Q\stackrel{\simeq}{\longrightarrow} A_{PL}(X\vee Y).$$
This identifies a Sullivan model $\land T\stackrel{\simeq}{\to} \land W\times_{\mathbb Q}\land Q$ as a Sullivan model for $X\vee Y$.

On the other hand, for any two Sullivan algebras $\land W$ and $\land Q$ the surjection $\land W\otimes \land Q\to \land W\times_{\mathbb Q}\land Q$ extends to a quasi-isomorphism
$$\varphi : \land T:=\land W\otimes \land Q\otimes \land R\stackrel{\simeq}{\longrightarrow} \land W\times_{\mathbb Q}\land Q$$
from a minimal Sullivan extension. Since $\land W\otimes \land Q\to \land W\times_{\mathbb Q}\land Q$ is surjective we may choose $R\subset T$ so that $\varphi(R)= 0$. It follows that $\land T$ is a Sullivan algebra, minimal if $\land W$ and $\land Z$ are minimal, and that in any case the quotient $(\land R, \overline{d})$ is a minimal Sullivan algebra.

\begin{Prop}
\label{p10.1}
With the notation and hypotheses above, 
\begin{enumerate}
\item[(i)] With appropriate choice of generating space $R\subset \land R$,
$$R \cap \mbox{ker}\, \overline{d} \stackrel{\cong}{\longrightarrow} H^{\geq 1}(\land R, \overline{d}).$$
In particular (Proposition \ref{p8.4}) the homotopy Lie algebra $L_R$ is profree.
\item[(ii)] If $\land W$ and $\land Q$ are quadratic Sullivan algebras then $\land T$ is quadratic.
\item[(iii)] In general, the quadratic Sullivan algebra $(\land T, d_1)$ for $\land T$ is the Sullivan model of $(\land W,d_1)\times_{\mathbb Q} (\land Q, d_1)$. 
\end{enumerate}
\end{Prop}

 \vspace{3mm}\noindent {\bf Remark.}  This result is analogous to the fact that the usual fibre of the injection $X\vee Y\to X\times Y$ is the join of $\Omega X$ and $\Omega Y$ and thus a suspension. (But note that $(X\vee Y)_{\mathbb Q}$ may be different from $X_{\mathbb Q}\vee Y_{\mathbb Q}$.)

\vspace{3mm}\noindent {\sl proof.}
(i)  Let $\land W\otimes \land U_W$ and $\land Q\otimes \land U_Q$ denote the respective acyclic closures. 
These, together with $\varphi : \land T \stackrel{\simeq} \land W\times_{\mathbb Q} \land Q$, yield  the quasi-isomorphisms   
 $$\land R \stackrel{\simeq}{\longleftarrow} \land T \otimes_{\land W\otimes \land Q} (\land W\otimes \land U_W\otimes \land Q\otimes \land U_Q) \stackrel{\simeq}{\longrightarrow} A:= (\land W\times_{\mathbb Q}\land Q) \otimes \land U_W\otimes \land U_Q.$$
Now divide $A$ by the ideal generated by $W$ to obtain the short exact sequence
 $$0\to \land^{\geq 1}W\otimes \land U_W\otimes \land U_Q \to A \to  \land Q\otimes \land U_W\otimes \land U_Q \to 0.$$

 Next, decompose the differential in $\land W\otimes \land U_W$ in the form $d= d_1+d'$ with $d_1(W)\subset \land^2W$, $d_1(U_W)\subset W\otimes \land U_W$, $d'(W)\subset \land^{\geq 3}W$ and $d'(U_W)\subset \land^{\geq 2}W\otimes \land U_W$. Then $d_1$ is a differential and $(\land W\otimes \land U_W,d_1)$ is the acyclic closure of $(\land W,d_1)$. Choose a direct summand,  $S$, of $d_1(\land^{\geq 1}U_W) $ in $W\otimes \land U_W$. Then $I = (\land^{\geq 2}W\otimes \land U_W)\oplus S$ is acyclic for the differential $d_1$ and therefore also for the differential $d$. 
Thus $J = I\otimes \land U_Q$ is an acyclic ideal in $A$ and $A\stackrel{\cong}{\to} A/J$.  

On the other hand, consider the short exact sequence
$$0\to d(\land U_W)\otimes \land U_Q\to A/J\to \land Q\otimes \land U_W\otimes \land U_Q\to 0$$
in which $d(\land U_W)\otimes \land U_Q$ is an ideal with trivial multiplication and trivial differential. It follows from the long homology sequence that
$$\mathbb Q \oplus d(\land U_W)\otimes \land^{\geq 1}U_Q \stackrel{\simeq}{\longrightarrow} A/J.$$
Since $\land R\simeq A \simeq A/J$, this completes the proof of (i).

(ii) Assign  $\land W$ and $\land Q$ wedge degree as a second degree and assign $U_W$ and $U_Q$   second degree $0$. Then $(\land W\otimes \land U_W,d_1)$ and $(\land Q\otimes \land U_Q,d_1)$ are the respective acyclic closures of $(\land W,d_1)$ and $(\land Q, d_1)$, and $d_1$ increases the second degree by 1. Now $ {\varphi}$ and $ {T}$ 
 may be constructed so that $ {R}$ is equipped with a second gradation for which $d$ increases the second degree by one and $ {\varphi}$ is bihomogeneous of degree zero. 

The argument in the proof of (i) now yields a sequence of bihomogeneous quasi-isomorphisms connecting
$$\mathbb Q \oplus \left( d_1(\land^+U_W)\otimes \land^+U_Q\right)\simeq \land  {R}.$$
Thus $H^{\geq 1}(\land {R})$ is concentrated in second degree 1. Denote by $S\subset R$ the subspace of elements of second degree 1. Then $(\land S,d)$ is a quadratic Sullivan algebra. Moreover by (i), any cycle $\Phi\in \land^2S$ is a boundary, necessarily of an element $v\in S$ of wedge degree 1. Thus $H(\land S)$ has no homology of wedge degree 2. Now Lemma \ref{l6.2} implies that $S\cap \mbox{ker}\, d = H^{\geq 1}(\land S)$. But if $S\neq R$ then for some $z\in R$, $dz\in \land S$. This implies that for some $\Psi\in \land S$, $z-d\Psi$ is a cycle and therefore is in $S$. This gives $S= R$. Thus $\land T$ is quadratic and $\varphi$ preserves wedge degrees.

\vspace{2mm} (iii) Denote by $d= \sum_{i\geq 1}d_i$ the differentials in $\land W$ and $\land Q$, where $d_i$ increases the wedge degree by $i$. Then (ii) constructs a quasi-isomorphism
$$\varphi : (\land T, d_1)= (\land W\otimes \land Q\otimes \land R, d_1) \stackrel{\simeq}{\longrightarrow} (\land  W, d_1)\times_{\mathbb Q}(\land Q, d_1)$$
in which $\varphi$ is surjective and preserves wedge degree. In particular,
$$H(\mbox{Ker}\, \varphi, d_1)=0.$$

Now we deform $d_1$ to a differential $d= \sum_{i\geq 1}d_i$ in $\land T$ with the following properties
\begin{enumerate}
\item[$\bullet$] $d$ restricts to the original differentials in $\land W$ and $\land Q$.
\item[$\bullet$] $d_i : R\to \mbox{Ker}\, \varphi$, $i\geq 2$.
\end{enumerate}

Write $R = \cup_{\neq \geq 0}R_n$, with $R_0= R\cap d_1^{-1}(\land W\otimes \land Q)$ and for $n>0$, $$R_n = R\cap d_1^{-1} 
(\land W\otimes \land Q\otimes \land R_{n-1}).$$ 
We define the derivations $d_2, \dots , d_q, ...$ on $R$ by induction on $n$ with the property that $$d_i(R) \subset (\land^+R \otimes \land (W\otimes \land Q)) \, \cap \, \land^{i+1}(R\oplus W\oplus Q)$$
and $(d_1+ \sum_{j\geq 2}d_j)^2= 0$ on $R_n$. 

So suppose this is done on $R_{n-1}$ and $z\in R_n$. We construct the $d_p$ by induction on $p$. So suppose $d_2, \dots , d_{p-1}$ constructed such that for $2\leq q\leq p-1$, we have
$$(d_1d_q+d_2d_{q-1}+ \dots + d_qd_1)z= 0.$$
Then let $R(z) = (d_pd_1 + d_{p-1}d_2 + \dots + d_2d_{p-1})(z)$. Since $d_1R(z)= 0$ and $R(z)\in \mbox{Ker}\, \varphi$, there is an element $u\in \mbox{Ker}\, \varphi$ with $R(z) = -d_1(u)$. We define $d_p(x)= u$. It follows by induction on $p$ that $d_p$ increases wedge degree by $p$. Thus $d_pz=0$ if $p>n$. We set $d= \sum_{i\geq 1}d_i$.

It is immediate from the construction that $d^2=0$ and that the two required properties are satisfied. In particular $(\land T,d)$ is a Sullivan algebra and since $\varphi$ is a quasi-isomorphism with respect to $d_1$, 
$$\varphi : (\land T,d) \to (\land W,d)\times_{\mathbb Q}(\land Q,d)$$
is a Sullivan model. Finally, the associated quadratic Sullivan algebra is by (ii) the Sullivan model of $(\land W, d_1)\times_{\mathbb Q}(\land Q, d_1)$.

\hfill$\square$

\vspace{3mm}\noindent {\bf Example 1.} If $\land W$ and $\land Q$ are minimal Sullivan algebras then (Proposition \ref{p10.1}) the surjection $\land W\otimes \land Q \to \land W\times_{\mathbb Q}\land Q$ extends to a quasi-isomorphism
 $$\land T = \land W\otimes \land Q\otimes \land R \stackrel{\simeq}{\longrightarrow} \land W\times_{\mathbb Q}\land Q$$
 in which $\land T$ is a minimal Sullivan algebra and the homotopy Lie algebra, $L_R$, is profree. Thus $\langle\land R\rangle$ is rationally wedge-like.
 
 \vspace{2mm}{\bf 2.} If one of the spaces $X$ and $Y$ has rational homology of finite type then the homotopy fibre of $(X\vee Y)_{\mathbb Q}\to (X\times Y)_{\mathbb Q}$ is rationally wedge-like.
 
 In fact, in this case if $\land W$ and $\land Q$ are respectively the minimal Sullivan models of $X$ and $Y$ then the inclusion $\land W\otimes \land Q\to \land T$ is a Sullivan representative for the inclusion $i : X\vee Y\to X\times Y$. This identifies $i_{\mathbb Q}$ as the map $\langle\land W\otimes \land Q\rangle\leftarrow \langle\land T\rangle$ and it follows from (\cite[Proposition 17.9]{FHTI}) that its homotopy fibre is $\langle\land R\rangle$.

 \subsection{Free products and the homotopy Lie algebra of $X\vee Y$}

The category of enriched Lie algebras has free products. In fact the classical construction of the free product, $L \amalg L'$, of two graded Lie algebras  extends naturally to  enriched Lie algebras, $(L, \{I_\alpha\})$ and $(L', \{I'_\beta\})$, with the enriched structure   given by the surjections
$$\xi_{\alpha, \beta, n} : L\,\amalg\, L'\to L_\alpha\, \amalg\, L'_\beta \,/\, (L_\alpha \, \amalg\, L'_\beta)^n.$$
Its completion will be denoted by $L\, \widehat{\amalg}\, L'$, and almost by definition,
$$L\, \widehat{\amalg}\, L' \stackrel{\cong}{\longrightarrow} \overline{L}\, \widehat{\amalg}\, \overline{L'}.$$

\vspace{3mm}\noindent {\bf Definition.} $L\, \widehat{\amalg}\, L'$ is the \emph{free product} of $L$ and $L'$.

 \vspace{3mm}\noindent {\bf Remarks. 1.} As in the classical case $L\, \widehat{\amalg}\, L'$ is equipped with inclusions $\xymatrix{L, L'\ar@{^{(}->}[r] &  L\, \widehat{\amalg}\, L'}$, and here if $L$ and $L'$ are complete the images are
 closed sub Lie algebras. Moreover these inclusions extend to a \emph{linear} inclusion
 $$\xymatrix{L\oplus L'\ar@{^{(}->}[r] & L\, \widehat{\amalg}\, L'}$$
 into a closed subspace of the free product.
 
 \vspace{2mm} {\bf 2.} If $C\subset L$ and $C'\subset L'$ are closed subspaces, we denote by 
 $$C\oplus C'\subset C \, \widehat{\amalg}\, C' \subset L\, \widehat{\amalg}\, L'$$
 the closure of the subspace of iterated Lie brackets of sequences of elements alternating between $C$ and $C'$. 

\begin{lem}\label{l10.1} (i). With the notation above, any two morphisms $f : (L, \{I_\alpha\}) \to E$ and $g : (L', \{I_\beta '\})\to E$ into a complete enriched Lie algebra  extend  uniquely to a morphism
$$L\,\widehat{\amalg}\, L' \to E.$$
This identifies $L\, \widehat{\amalg}\, L'$ as a coproduct in the category of complete Lie algebras, and characterizes $L\,\widehat{\amalg}\,L'$ up to natural isomorphism.

(ii). Suppose a profree Lie algebra $L= \overline{\mathbb L}_T$, and   $T= T(1)\oplus T(2)$ with $T(1)$ and $T(2)$ closed subspaces. Then $\overline{\mathbb L}_{T(1)}$ and $\overline{\mathbb L}_{(T(2)}$ are   the closures of the Lie algebras generated by $T(1)$ and $T(2)$, and 
$$\overline{\mathbb L}_{T(1)} \, \widehat{\amalg}\, \overline{\mathbb L}_{T(2)} \stackrel{\cong}{\longrightarrow} \overline{\mathbb L}_T.$$
 
\end{lem}

\vspace{3mm}\noindent {\sl proof.} (i) Let $\{J_\gamma\}$ denote the enriched structure for $E$. Then for each $\gamma$ there are indices $\alpha (\gamma), \beta(\gamma)$ such that $f$ and $g$ factor to yield morphisms
$$f_\gamma : L_{\alpha (\gamma)} \to E_\gamma \hspace{5mm}\mbox{and } g_\gamma : L'_{\beta (\gamma)} \to E_\gamma.$$
Moreover, since $E_\gamma^{n(\gamma)} = 0$, some $n(\gamma)$, $f_\gamma \amalg g_\gamma$ factors to yield a morphism
$$h_\gamma : (L_{\alpha (\gamma)} \amalg L'_{\beta(\gamma)}) / (L_{\alpha (\gamma)} \amalg L'_{\beta (\gamma)})^n \longrightarrow E_\gamma.$$
Since $E$ is complete these extend to 
$$h : = \varprojlim_\gamma h_\gamma : L\,\widehat{\amalg}\, L'\to E.$$
  The uniqueness is immediate. 
  
  (ii) The inclusion $T(1)$, $T(2)\to T$ extend to a morphism $\overline{\mathbb L}_{T(1)} \, \widehat{\amalg}\, \overline{\mathbb L}_{T(2)} \to \overline{\mathbb L}_T$, with inverse given by the extension to $\overline{\mathbb L}_T$ of $T= T(1)\oplus T(2)\to \overline{\mathbb L}_{T(1)} \, \widehat{\amalg}\, \overline{\mathbb L}_{T(2)}$. \hfill $\square$

\begin{Prop}\label{p10.2} Denote by $L_X$ and $L_Y$ the homotopy Lie algebras of connected spaces $X$ and $Y$. 
 Then the inclusions $X,Y \to X\vee Y$ induce an isomorphism  
$$L_{X\vee Y} \stackrel{\cong}{\longrightarrow} L_X \,\widehat{\amalg}\, L_Y.$$
  In particular, the correspondence $X\mapsto L_X$ preserves coproducts.
\end{Prop}

Proposition \ref{p10.2} is an immediate consequence of the explicit description in Proposition \ref{p10.3} below of the homotopy Lie algebra of a fibre product $\land W\times_{\mathbb Q} \land Q$ of any two minimal Sullivan algebras.
 Recall from Proposition \ref{p10.1} that the  natural surjection $\land W\otimes \land Q \to \land W\times_{\mathbb Q}\land Q$ extends to a minimal Sullivan model   
$$\varphi : \land T:= \land W\otimes \land Q\otimes \land R \stackrel{\simeq}{\longrightarrow} \land W\times_{\mathbb Q}\land Q.$$  Denote   the corresponding homotopy Lie algebras (\S 8) by   $L_W, L_Q, L_T$ and $L_R$. In particular, this  identifies $L_T$ as the homotopy Lie algebra of a minimal Sullivan model of $\land W\times_{\mathbb Q}\land Q$.
The inclusions $\land W, \land Q \to \land W\times_{\mathbb Q}\land Q$ then lift to unique based homotopy classes of morphisms $\land W, \land Q\to \land T$, and thus induce morphisms $L_T\to L_W, L_Q$. These provide a surjection $L_T\to L_W\times L_Q$.

\begin{Prop}
\label{p10.3} With the hypotheses and notation above, 
\begin{enumerate}
\item[(i)] $\varphi$ induces a quasi-isomorphism $(\land T, d_1) \to (\land W, d_1)\times_{\mathbb Q}(\land Q, d_1)$ of associated quadratic Sullivan algebras.
\item[(ii)]  The surjection $L_T\to L_W\times L_Q$ factors as
$$\xymatrix{L_T\ar[rr]^\cong_\sigma && L_W\, \widehat{\amalg}\, L_Q\ar[rr] && L_W\times L_Q,}$$
and $\sigma$ is natural with respect to morphisms $\land W\to \land W'$ and $\land Q\to \land Q'$. 
\item[(iii)] If $f : L\to L'$ and $g : E\to E'$ are morphisms of complete enriched Lie algebras, then
$$f\, \widehat{\amalg}\, g: L\, \widehat{\amalg}\, E\to L'\, \widehat{\amalg}\, E'$$
is   surjective if and only if $f$ and $g$ are surjective.
\end{enumerate}
\end{Prop}

\vspace{3mm}\noindent {\bf Corollary.}   If $L$ and $L'$ are complete enriched Lie algebras then the kernel of 
$$L\, \widehat{\amalg}\, L' \to L\times L'$$
is profree. Moreover, if $L$ and $L'$ are both profree then $L\,\widehat{\amalg}\, L'$ is also profree.

\vspace{3mm}\noindent {\sl proof of the Corollary.}   Apply Proposition \ref{p10.1} to identify the kernel as  a profree Lie algebra $L_R$.    If now $L$ and $L'$ are profree then by Proposition \ref{p8.4} there are quasi-isomorphisms
$$\land W \stackrel{\simeq}{\longrightarrow} \mathbb Q \oplus S \hspace{5mm}\mbox{and } \land W' \stackrel{\simeq}{\longrightarrow} \mathbb Q \oplus S'$$
with $S\cdot S = 0= S'\cdot S'$, and vanishing differentials in $S$ and $S'$. Thus $\land W\times_{\mathbb Q}\land W' \simeq \mathbb Q\oplus S\oplus S'$, and its homotopy Lie algebra, $L\, \widehat{\amalg}\, L'$, is profree by Proposition \ref{p8.4}.

\hfill$\square$

\vspace{3mm} \vspace{3mm}\noindent {\sl proof of Proposition \ref{p10.3}.} (i) As recalled at the start of \S 8, $L_T, L_W$ and $L_Q$ depend only on the quadratic Sullivan algebras $(\land T, d_1)$, $(\land W, d_1)$ and $(\land Q, d_1)$. Moreover, by Proposition \ref{p10.1}, 
$$(\land T, d_1) \stackrel{\simeq}{\longrightarrow} (\land W, d_1)\times_{\mathbb Q}(\land Q, d_1).$$
Therefore for any quadratic Sullivan algebra $(\land V, d_1)$, the based homotopy classes of morphisms satisfy
$$
\renewcommand{\arraystretch}{1.4}
\begin{array}{ll}
[(\land V, d_1), (\land T, d_1)]_* & = [(\land V, d_1), (\land W, d_1)\times_{\mathbb Q}(\land Q, d_1)]*\\
& = [(\land V, d_1), (\land W, d_1)]_* \times [(\land V, d_1), (\land Q,d_1)]*.
\end{array}
\renewcommand{\arraystretch}{1}
$$

(ii) Proposition \ref{p4.1} translates this equality to give
$${\mathcal C}(L_T, L_V) = {\mathcal C}(L_W, L_V) \times {\mathcal C}(L_Q, L_V),$$
where ${\mathcal C}(-,-)$ denotes the set of morphisms in the category of enriched Lie algebras. This identifies $L_T= L_W\, \widehat{\amalg}\, L_Q$.
The naturality is a standard diagram chasing.

(iii) To show that if $f$ and $g$ are surjective then $f\, \widehat{\amalg}\, g$ is surjective, it   is sufficient to consider the case $g = id_E : E\stackrel{=}{\longrightarrow} E$. In this case we have the row-exact commutative diagram
$$\xymatrix{
0\ar[r] &  K\ar[rr]\ar[d]^h && L\, \widehat{\amalg}\, E\ar[d]^{f\, \widehat{\amalg}\, id_E} \ar[rr] &&   L\times E\ar[r]\ar[d]^{f\times id_E} & 0\\
0 \ar[r] & K' \ar[rr] && L'\, \widehat{\amalg}\, E \ar[rr] && L'\times E \ar[r] & 0,}$$
and so it is also sufficient to show that $h$ is surjective.

Now let $\land W$, $\land W'$ and $\land Q$ be, respectively the quadratic models of $L$, $L'$ and $E$, and let $\varphi : \land W'\to \land W$ be the morphism corresponding to $f$. Then the restriction $\varphi : W'\to W$ dualizes to $f$, and so $\varphi$ is injective. In particular, $\varphi$ extends to an inclusion
$$\varphi : \land W'\otimes \land U_{W'}\to \land W\otimes \land U_W$$
of the acyclic closures, restricting to an inclusion $U_{W'}\to U_W$.

On the other hand, the diagram above is the diagram of homotopy Lie algebras associated with a commutative diagram,
$$\xymatrix{\land W\otimes \land Q \ar[rr] && \land W\otimes \land Q\otimes \land R \ar[rr] && \land R\\
\land W'\otimes \land Q\ar[u]^{\varphi\otimes id} \ar[rr] && \land W'\otimes \land Q\otimes \land R' \ar[u]\ar[rr] && \land R'\ar[u]^\sigma,}$$
which identifies $L_\sigma$ with $h$.

Moreover, the computations in the proof of Proposition \ref{p10.1} identify
$$H(\sigma) = \varphi : \mathbb Q \oplus d(\land U_{W'}) \to \mathbb Q \oplus d(\land U_W).$$
It follows that $H(\sigma)$ is injective. But since $L_R$ and $L_{R'}$ are profree, Lemma \ref{l4.1} identifies
$$L_R/L_R^{(2)} = H^{\geq 1}(\land R)^\vee \hspace{5mm}\mbox{and } L_{R'}/L_{R'}^{(2)} = H^{\geq 1}(\land R')^\vee.$$
It follows that $L_\sigma : L_R\to L_{R'}$ induces a surjection $L_R/L_R^{(2)}\to L_{R'}/L_{R'}^{(2)}$. Now Lemma \ref{l6.1} (i) implies that $L_\sigma$ is surjective.

In the reverse direction note that the maps $id_L$ and $E\to 0$ define a retraction $L\, \widehat{\amalg}\, E\to L$. The corresponding retraction $L'\, \widehat{\amalg}\, E'\to L'$ fits in the commutative diagram,
$$
\xymatrix{
L\, \widehat{\amalg}\, E\ar[d]\ar[rr]^{f\, \widehat{\amalg}\, g} && L'\, \widehat{\amalg}\, E'\ar[d]\\
L \ar[rr]^f&& L'.}
$$
Therefore if $f\, \widehat{\amalg}\, g$ is surjective so is $f$ and, similarly, so is $g$.

\hfill$\square$

\vspace{3mm}\noindent {\bf Example.} Free products of completions of weighted Lie algebras
 
\vspace{1mm} Let $E$ and $E'$ be weighted enriched Lie algebras with enriched structures given by morphisms $E\to E_\alpha$ and $E'\to E_\beta'$.  
 Denote by $L$ and $L'$ the completions of $E$ and $E'$.  The weight decompositions in $E$ and $E'$ then extend  in the standard way to a weight decomposition in  $E\,\amalg\, E'$.  
 It then follows from \S2.1 that $L\, \widehat{\amalg}\, L'$ is the completion of $E\, \amalg\, E'$.

 \subsection{Limits and colimits of complete enriched Lie algebras}
 
 \begin{Prop}
 \label{p10.4} The category $\overline{\mathcal C}$ of complete enriched Lie algebras is complete and cocomplete : i.e., each diagram of complete enriched Lie algebras has a direct  and an inverse limit.\end{Prop}
 
 \vspace{3mm}\noindent {\sl proof.} The construction of the coproduct of two complete enriched Lie algebras  can be generalized to any family $\{L(i)\}_{i\in I}$ of complete enriched Lie algebras. In fact, denote by $A_i$ the family of ideals defining the enriched structure in $L(i)$.  For each finite subset $J\subset I$, we denote by $p_J : \amalg_{i\in I} L(i) \to \amalg_{i\in J} L(i)$ the projection defined by $p_J(L(i)) = 0$ for $i\not\in J$. Then for $J\subset I$, finite, for $\alpha_i\in A_i$ for $ i\in J$, and for an integer $n$, we denote by $K_{J,(\alpha_i),n}$ the kernel of the projection
 $$
 \xymatrix{
 \amalg_{i\in I}L(i)\ar[rr]^{p_J} && \amalg_{i\in J} L(i) \ar[rr] && \amalg_{i\in J} L(i)_{\alpha_i} /(\amalg_{i\in J}L(i)_{\alpha_i})^n.}$$
 These ideals define an enriched structure on $\amalg_{i\in I} L(i)$ and its completion $\overline{\amalg_{i\in I} L(i)}$ is the coproduct of the the $L(i)$ is the category of complete enriched Lie algebras.
 
 Now let $(L(i), f_j)$ be a diagram in $\overline{\mathcal C}$. Its colimit in the category of Lie algebras is the quotient of $\amalg_{i\in I} L(i)$ by an ideal $J$ associated to the morphisms $f_j$. Denote by $\iota: \amalg_{i \in I }L(i) \to \overline{\amalg_{i\in I}L(i)}$ the natural injection and by $\overline{J}$ the closure of the ideal generated by $\iota(J)$  in $\overline{\amalg L(i)}$. Then by the Corollary to Lemma \ref{l2.2}, the quotient $\overline{\amalg L(i)}/\overline{J}$ inherits a complete enriched structure. This is clearly   the colimit of the diagram $(L(i), f_j)$ in $\overline{\mathcal C}$.  
 
 On the other hand the usual product of the complete enriched Lie algebras $L(i)$ is complete for the family of the ideals which are the kernels of the surjections 
 $$\xymatrix{
 \prod_{i\in I}L(i)  \ar[rr]^{p_J}&& \prod_{i\in J} L(i)\ar[rr]^{\prod \rho(i)_{\alpha_i} } && L(i)_{\alpha_i}
 }$$
 associated to the finite subsets $J\subset I$ and the surjections $\rho(i)_{\alpha_i} : L(i)\to L(i)_{\alpha_i}.$
 
Now consider a diagram $(L(i), f_j)$ in $\overline{\mathcal C}$.   Its limit in the category of Lie algebras is the sub Lie algebra $S\subset \prod L(i)$ consisting of families $a_i\in L(i)$ such that for any morphism $f_j : L(i_0)\to L(i_1)$ we have $f_j(a_{i_0})= a_{i_1}$.   Since for any   morphism $f : L(1)\to L(2)$ of complete enriched Lie algebras, $f$ is a limit of morphisms $f_\alpha : L(1)_{\beta (\alpha)}\to L(2)_\alpha$, the sub Lie algebra $S$ is closed and is the   limit of the diagram in $\overline{\mathcal C}$.

\hfill$\square$

 \subsection{Wedges of spheres}

Suppose $\{S^{n_\alpha}\}_{\alpha \in \mathcal S}$ is a collection of spheres in which $\mathcal S$ is a linearly ordered set and each $n_\alpha \geq 1$. For each finite subset $\sigma \subset {\mathcal S}$ we write $\sigma = \{\sigma_1, \dots , \sigma_r\}$ with $\sigma_1<\dots <\sigma_r$, and set $\vert \sigma\vert = r$. If $\sigma \subset \tau$, so that $\vert \sigma\vert \leq \vert\tau \vert=q$, the inclusion defines an inclusion
$$j_{\sigma, \tau} : S^{n_{\sigma_1}}\vee \dots \vee S^{n_{\sigma_r}} \longrightarrow S^{n_{\tau_1}}\vee \dots \vee S^{n_{\tau_q}}.$$

\vspace{3mm}\noindent {\bf Definition.} A \emph{wedge of spheres} is a direct limit of the form

$$X=\vee_{\alpha} S^{n_\alpha} := \varinjlim_{\sigma_1<\dots <\sigma_r }\,\, S^{n_{\sigma_1}}\vee \dots \vee S^{n_{\sigma_r}}.$$

\vspace{3mm} \noindent {\bf Remark.} It is immediate from the definition that a wedge of spheres is a co-H-space. Thus (Corollary 3 to Proposition 8.4) its minimal Sullivan model, $\land V$, is quadratic and its homotopy Lie algebra is profree. Thus we may assume   $d: V\to \land^2V$ and write $V = S\oplus T$ with
$$S = V\cap \mbox{ker}\, d\stackrel{\cong}{\longrightarrow} H^{\geq 1}(\land V).$$ 
Now let ${\mathcal S}= \{\alpha\}$ index a linearly ordered basis $\{v_\alpha\}$ of $S$. This then defines a wedge 
$$X = \vee_\alpha S^{n_\alpha}$$
of spheres.

\vspace{3mm}

\begin{Prop}
\label{p10.5} A minimal Sullivan algebra, $\land V$, is the Sullivan model of a wedge of spheres, $X$, if and only if
\begin{enumerate}
\item[(i)] The homotopy Lie algebra, $L_V$, is profree, and
\item[(ii)] $H(\land V)$ is the dual of a graded vector space.
\end{enumerate}
Moreover, if $H(\land V)$ is a graded vector space of finite type then $\land V$ determines $X$ up to homeomorphism.
\end{Prop}

\vspace{3mm}\noindent {\sl proof.} As observed above a wedge of spheres, $X$, has a profree homotopy Lie algebra. Moreover $H(\land V)= H(X) = H_*(X;\mathbb Q)^\vee$.

In the reverse direction suppose $\land V$ satisfies (i) and (ii), so that $L_V$ is profree and $H^{\geq 1}(\land V)= T^\vee$. Then by Proposition \ref{p8.4},
$$\land V \stackrel{\simeq}{\longrightarrow} (\mathbb Q \oplus T^\vee, 0).$$
Choose a linearly ordered basis $\{x_i\}$ of $T$ and set $n_i = \mbox{deg}\, x_i$ and ${\mathcal S}_X= \{i\}$. 
Then set 
$$X = \vee_{i\in {\mathcal S}_X} S^{n_i}_i.$$
By construction, $H_{\geq 1}(X;\mathbb Q)= T.$ 

Now let $\land W$ be the minimal Sullivan model of $X$. Since $X$ is a wedge of spheres,
$$\land W \simeq (\mathbb Q \oplus H^{\geq 1}(X), 0)= (\mathbb Q \oplus T^\vee, 0) = (\mathbb Q \oplus H^{\geq 1}(\land V), 0) \simeq \land V.$$
Since $\land W$ and $\land V$ are minimal Sullivan algebras, it follows that $\land V \cong \land W$. Therefore $\land V$ is the minimal Sullivan model of $X$. 
 
 Finally, suppose $H(\land V)$ is a graded vector space of finite type, and that $\land V$ is the minimal Sullivan model of a second wedge of spheres $Y$. Then $H^{\geq 1}(Y)$ and $H^{\geq 1}(X)$ are isomorphic vetor spaces and so $H_{\geq 1}(Y;\mathbb Q) \cong H_{\geq 1}(X;\mathbb Q)$. Since $Y$ is a wedge of spheres
 $$Y = \vee_{j\in {\mathcal S}_Y} S^{n_j}$$
 and the number of indices $j\in {\mathcal S}_Y$ with $n_j=n$ is dim$\, H_n(Y;\mathbb Q)=$ dim$\, H_n(X;\mathbb Q)$. Therefore ${\mathcal S}_X$ and ${\mathcal S}_Y$ are isomorphic sets and $X \cong Y$. \hfill$\square$
 
 \vspace{3mm} These constructions provide an explicit description of the morphism (\S 9),
 $$\pi_*\vert \varphi\vert : \pi_*(Y) \to \pi_*(Y_{\mathbb Q})$$
 when $Y = \vee_\sigma S^{n_\sigma}$ is a wedge of spheres of finite type. When $Y$ is a wedge of circles, this extends diagrams (1) and (2) in \S 9.1.
 
 \begin{Prop}
 \label{p10.6}
 With the hypotheses and notation above, denote by $\sigma \in \pi_{n_\sigma}(S^{n_\sigma})$ the homotopy class represented by id$_{S^{n_\sigma}}$. Then the homotopy Lie algebra $L_Y$ of $Y$ can be written $\overline{\mathbb L}_T$ in which $\pi_*\vert \varphi\vert$ maps the set $\{\sigma\}$ bijectively to a basis of $T$.\end{Prop}
 
 \vspace{3mm}\noindent {\sl proof.} Here (Proposition \ref{p8.4}) $L_Y= \overline{\mathbb L}_T$ with $sT= (V\cap \mbox{ker}\, d_1)^\vee$. Because $Y$ is a wedge of finite type,
 $$(V\cap \mbox{ker}\, d_1)^\vee = H^{\geq 1}(Y)^\vee = H_{\geq 1}(Y;\mathbb Q).$$
 However, a basis of $H_{\geq 1}(Y;\mathbb Q)$ is provided by the images of the homotopy classes $\sigma$ under the Hurewicz homomorphism. Since the Hurewicz homomorphism is natural, the Proposition follows from Proposition \ref{p9.3}. \hfill$\square$

\vspace{3mm}\noindent {\bf Definition.} A space $X$ is \emph{rationally wedge-like} if it is homotopy equivalent to a spatial realization $\langle\land V\rangle$ with the homotopy Lie algebra of $\land V$ a profree Lie algebra.

\vspace{3mm}\noindent {\bf Remark.} It follows from Proposition \ref{p10.5} that if $\vee S^n_{\alpha}= X$ is a wedge of spheres and $\land V$ is its minimal Sullivan model then
$$(\vee S^n_{\alpha})_{\mathbb Q} = X_{\mathbb Q}= \langle \land V\rangle$$
is rationally wedge-like.

 \section{Sullivan rational spaces}
 
Let $\varphi_X : \land V\stackrel{\simeq}{\to} A_{PL}(X)$ be a minimal Sullivan model for a connected space, $X$, and let $L_X$ denote the homotopy Lie algebra of $\land V$. Recall (\S 9) that the adjoint map   
  $$\vert\varphi_X\vert : X\to  \langle\land V\rangle:= X_{\mathbb Q}$$
is called   the \emph{Sullivan completion} of $X$.  
  
  \vspace{3mm}\noindent {\bf Definition.} A path connected space, $X$,  is   a \emph{Sullivan rational space} if the Sullivan completion   $\vert\varphi_X\vert: X\to X_{\mathbb Q}$  is a homotopy equivalence.
  
  \vspace{3mm}\noindent {\bf Remarks. 1.} The condition that $\vert\varphi_X\vert$ is a homotopy equivalence implies that the minimal Sullivan model of $X$ directly computes both $H(X_{\mathbb Q})$ and $\pi_*(X_{\mathbb Q})$.
  
  {\bf 2.} If $X $ is Sullivan rational   then  $X_{\mathbb Q}\simeq (X_{\mathbb Q})_{\mathbb Q}$ and $X_{\mathbb Q}$ is Sullivan rational.

 \vspace{3mm} Recall that (\cite{Q})   Quillen's rationalization   for simply connected spaces assigns to each such space a map $Y\to \mathbb Q (Y)$ which induces isomorphisms
  $$H_{\geq 1}(Y)\otimes \mathbb Q \stackrel{\cong}{\to} H_{\geq 1}(\mathbb Q(Y)) \hspace{5mm}\mbox{and } \pi_*(Y)\otimes \mathbb Q \stackrel{\cong}{\to} \pi_*(\mathbb Q (Y)).$$
  In particular, $\mathbb Q (\mathbb Q (Y)) = \mathbb Q (Y)$.

 Sullivan's completion, $X_{\mathbb Q}$, is analogous to Quillen's rationalization, and also to the Bousfield-Kan completion $\mathbb Q_{\infty}(X)$ - cf \cite{BK}. Moreover, if $H(X)$ has finite type then $\mathbb Q_\infty (X)$ and $X_{\mathbb Q}$ are homotopy equivalent (\cite{BG}). But it is not always true that every $X_{\mathbb Q}$ is Sullivan rational. For instance (\cite{IM}), $H^2((S^1\vee S^1)_{\mathbb Q})$ has uncountable dimension, where as if $X_{\mathbb Q}$ is Sullivan rational by the Corollary to Theorem \ref{t2} below, $H(X_{\mathbb Q})$ has finite type.

 \vspace{3mm}\noindent {\bf Remarks} {\bf 1.} If $\varphi_X : \land V\to A_{PL}(X)$ is a minimal Sullivan model then $X$ is Sullivan rational if and only if each $$\pi_k\vert\varphi_X\vert : (V^k)^\vee \stackrel{\cong}{\longrightarrow} \pi_k(X).$$

{\bf 2.}  Let $Y$ be a Sullivan rational space and let $X$ be a connected space. Since $f_{\mathbb Q}\circ \varphi_X= \varphi_Y\circ f$, and $\varphi_Y$ is a homotopy equivalence, it follows that $f\mapsto f_{\mathbb Q}$ yields an injection
$  [X,Y]\to   [X_{\mathbb Q}, Y_{\mathbb Q}].$

\vspace{3mm} Now recall from \S 7 that if $\land V$ is a Sullivan algebra then adjoint to $id_{\langle\land V\rangle}$ is a morphism
$$\iota_{\land V} : \land V\to A_{PL}\langle\land V\rangle.$$

 \begin{lem} 
 \label{l11.1}Suppose for some minimal Sullivan algebra, $\land V$, that $\iota_{\land V} : \land V\to A_{PL}\langle\land V\rangle$ is a quasi-isomorphism. Then
 \begin{enumerate}
 \item[(i)] $\langle\land V\rangle$ is Sullivan rational.
 \item[(ii)]   A basis of $V$ is at most countable.
 \item[(iii)] $H(\land V)$ is a graded vector space of finite type.
 \end{enumerate}
 \end{lem}
 
 \vspace{3mm}\noindent {\sl proof.} (i)  By hypothesis $\iota_{\land V}$ is a Sullivan model for $
 \langle\land V\rangle$. Since its adjoint is the homotopy equivalence $id_{\langle\land V\rangle}$, $\langle\land V\rangle$ is Sullivan rational.
 
(ii). It follows from (17) that  the
  surjection $V^\vee \to (V\cap \mbox{ker}\, d)^\vee$ decomposes as the composite
 $$\xymatrix{V^\vee \ar[rr]^{\mbox{hur}} &&H_*\langle\land V\rangle \ar[rr] && (V\cap \mbox{ker}\, d)^\vee.}$$
 Thus $H_*\langle\land V\rangle\to (V\cap \mbox{ker}\, d)^\vee$ is surjective. Dualizing gives injections $(V^n\cap \mbox{ker}\, d)^{\vee\vee} \to H^n(\land V).$
 
 Now assume by induction that $V^{<n}$ has an at most countable basis and let $Z$ be the space of cycles in $(\land V)^n$. Then $Z\subset \land V^{<n}\oplus V^n$, and division by $V^n\cap \mbox{ker}\, d$ gives an exact sequence
 $$0\to V^n\cap \mbox{ker}\, d\to Z\to \land V^{<n}.$$
 Since $(V^n\cap \mbox{ker}\, d)^{\vee\vee}$ embeds in $H^n(\land V)$, it follows that dim$\, V^n\cap \mbox{ker}\, d<\infty$. 
 
 On the other hand, $V^n= \cup_{p\geq 0}V^n_p$ with $V^n_0= V^n \cap d^{-1}(\land V^{<n})$ and $V_{p+1}^n = V^n \cap d^{-1}(\land V^{<n}\otimes \land V^n_p)$. Since dim$\,V^n\cap \mbox{ker}\, d<\infty$, it follows by induction on $p$ that each $V^n_p$, and hence $V^n $ itself, has an at most countable basis.
 
 (iii) It follows from (i) that each $H^n(\land V)$ has an at most countable basis. But since $H^n(\land V)= H^n(\langle\land V\rangle)= H_n\langle\land V\rangle^\vee$, this implies that each dim$\, H^n(\land V)<\infty$. \hfill$\square$

\begin{lem}
\label{l11.2} Suppose $\land V$ is a minimal Sullivan algebra. If $V^1=0$ then the following  conditions are equivalent:
\begin{enumerate}
\item[(i)] $\iota_{\land V} : \land V\to A_{PL}\langle\land V\rangle$ is a quasi-isomorphism.
\item[(ii)]  Each dim$\, V^k<\infty$. 
\item[(iii)] Each dim$\, H^k(\land V)<\infty.$
\item[(iv)] Each dim$\, \pi_k<\land V>$ is finite.
\end{enumerate}
\end{lem}

\vspace{3mm}\noindent {\sl proof.} (i) $\Leftrightarrow$ (ii). This is \cite[Theorem 1.5]{RHTII}

(ii) $\Leftrightarrow $ (iii). This follows because $d: V^k\to \land V^{<k}$.

(ii) $\Leftrightarrow$ (iv) This is because $\pi_k = (V^k)^\vee$. \hfill$\square$

\begin{lem}
\label{l11.3}
Suppose $G$ is a discrete group and its classifying space, $BG$, is Sullivan rational. Then
\begin{enumerate}
\item[(i)] Each $G^n/G^{n+1}$ is a finite dimensional rational vector space.
\item[(ii)] $G= \varprojlim_n G/G^n$.
\end{enumerate}
\end{lem}

\vspace{3mm}\noindent {\sl proof.} Let $\land V$ be the minimal Sullivan model of $BG$. Then
$$(V^n)^\vee = \pi_n\langle\land V\rangle = \pi_n(BG) = \left\{
\begin{array}{ll} 0 & n\geq 2\\G & n=1\end{array},\right.$$
Thus $V = V^1$, and (cf Proposition 8.2)
$$G = \pi_1\langle\land V^1\rangle = G_L$$
where $L=L_0$ is the homotopy Lie algebra of $\land V$. Now the Lemma follows from Proposition \ref{p8.2}. \hfill$\square$
 
 \begin{Prop}
\label{p11.1} Suppose $\varphi : \land V\to A_{PL}(X)$ is a Sullivan model of a connected space, $X$. Then
\begin{enumerate}
\item[(i)] $H(\vert\varphi\vert) : H(X) \stackrel{\cong}{\longleftarrow} H(X_{\mathbb Q}) \, \Longleftrightarrow \, H(\iota_{\land V}) : H(\land V) \stackrel{\cong}{\rightarrow}  H(\langle\land V\rangle)$. 
\item[(ii)] When the conditions of (i) hold then $\langle\land V\rangle$ is Sullivan rational, $X_{\mathbb Q}\to (X_{\mathbb Q})_{\mathbb Q}$ is an isomorphism of simplicial sets, and $X_{\mathbb Q}$ is Sullivan rational.
\end{enumerate}
\end{Prop}

\vspace{3mm}\noindent {\sl proof.} (i) According to (\ref{i13}),
$$\varphi = A_{PL}(\vert \varphi\vert) \circ \iota_{\land V}.$$
Since $H(\varphi)$ is an isomorphism this establishes (i).

(ii) By definition, $X_{\mathbb Q} = \langle\land V\rangle$. Therefore (i) $\Rightarrow$ (ii) follows from Lemma \ref{l11.1}.   \hfill$\square$

 \vspace{3mm} Now suppose that
\begin{eqnarray}
\label{i18}\xymatrix{\land W \ar[rr]^\lambda \ar[d]^{\varphi_B} && \land W\otimes \land Z\ar[d]^{\varphi_X}\ar[rr]^\rho &&\land Z\ar[d]^{\varphi_F}\\
 A_{PL}(B)\ar[rr]_{A_{PL}(p)} && A_{PL}(X)\ar[rr]_{A_{PL}(j)} && A_{PL}(F)}\end{eqnarray}
 be a commutative cdga diagram, in which
 \begin{enumerate}
 \item[$\bullet$] The upper row is a $\Lambda$-extension decomposition of a minimal Sullivan algebra $\land V = \land W\otimes \land Z$.
 \item[$\bullet$] The lower row is an application of $A_{PL}(-)$ to a fibration $\xymatrix{ B& X\ar[l]^p & F\ar[l]^j}$ of path connected spaces.
 \end{enumerate}

 \begin{Theorem}
 \label{t2}
 Suppose in (\ref{i18}) that $F$ is Sullivan rational and that $H(\varphi_F)$ is an isomorphism. Then
 \begin{enumerate}
 \item[(i)] The holonomy representation of $\pi_1(B)$ in each $H^k(F)$ is nilpotent.
 \item[(ii)] $H(\varphi_B) $ is an isomorphism $\Leftrightarrow$ $H(\varphi_X)$ is an isomorphism.
 \item[(iii)] If $H(\varphi_B)$ and $H(\varphi_X)$ are isomorphisms then $B$ and $X$ are Sullivan rational.
 \end{enumerate}
 \end{Theorem}
 
 \vspace{3mm}\noindent {\sl proof.} (i) The adjoint map $\langle\varphi_B\rangle : B\to \langle\land W\rangle$ induces a morphism of fundamental groups, and (\cite[Theorem 4.2]{RHTII}) the isomorphism $H(\varphi_F): H^k(\land Z) \stackrel{\cong}{\to} H^k(F)$ is equivariant with respect to the two corresponding holonomy  representations. The first is locally nilpotent (\cite[Proposition 4.4]{RHTII}) and the second is the dual of a representation of $\pi_1(F)$ in $H_*(F)$. It follows that both representations are nilpotent.

 (ii) Since $\land Z$ is the minimal Sullivan model of a Sullivan rational space it follows from Lemma \ref{l11.1} that $H(\land Z)$ is a graded vector space of finite type. Thus, if $H(\varphi_B)$ is an isomorphism, (\cite[Theorem 5.1]{RHTII}) provides a Sullivan algebra $(\land V',d') = (\land W\otimes \land Z,d')$ and a commutative cdga diagram
 $$\xymatrix{
 \land W\ar[d]^{\varphi_B}_\simeq \ar[rr] && (\land W\otimes \land Z,d') \ar[rr]\ar[d]^{\varphi'}_\simeq && \land Z\ar[d]^{\varphi_F}_\simeq\\
 A_{PL}(B) \ar[rr]_{A_{PL}(p)}&&A_{PL}(X)\ar[rr] && A_{PL}(F).}$$
 It is straightforward to extend the inclusion of $\land W$ to a lift of $\varphi_X$ through $\varphi'$ and to conclude that $H(\varphi_X)$ is an isomorphism.
 
 On the other hand, suppose $H(\varphi_X)$ is an isomorphism. Then extend $\varphi_B$ to a quasi-isomorphism, $\varphi_B' : \land W\otimes \land W'\stackrel{\simeq}{\longrightarrow} A_{PL}(B)$, from a minimal $\Lambda$-extension. By \cite[Theorem 5.1]{RHTII} this extends to a commutative cdga diagram
 $$
 \xymatrix{
 \land W\otimes \land W' \ar[d]^\simeq \ar[rr] && \land W\otimes \land W'\otimes \land Z \ar[d]^\simeq \ar[rr] &&\land Z\ar[d]^{\varphi_F}\\
 A_{PL}(B) \ar[rr] && A_{PL}(X) \ar[rr] && A_{PL}(F).
 }$$
 Since $\varphi_X : \land W\otimes \land Z\to A_{PL}(X)$ is also a quasi-isomorphism it follows that $W'=0$ and $H(\varphi_B)$ is an isomorphism.
 
 (iii) In (\ref{i18}) we are now assuming that each of $H(\varphi_B)$, $H(\varphi_X)$ and $H(\varphi_F)$ are isomorphisms. On the other hand taking the adjoints of $\varphi_B$, $\varphi_X$ and $\varphi_F$ produces the commutative diagram
 $$
 \xymatrix{\langle\land W\rangle &&\langle\land W\otimes \land Z\rangle\ar[ll] && \langle\land Z\rangle\ar[ll]\\
 B\ar[u]^{\langle\varphi_B\rangle} && X\ar[ll]\ar[u]^{\langle\varphi_X\rangle} && F\ar[ll]\ar[u]^{\langle\varphi_F\rangle}.}$$
 According to (\cite[Proposition 17.9]{FHTI}) the upper row is also a fibration. Thus it follows from the 5-lemma that if any two of $\pi_*\langle\varphi_B\rangle$, $\pi_*\langle\varphi_X\rangle$ and $\pi_*\langle\varphi_F\rangle$ are isomorphisms so is the third. \hfill$\square$

\begin{Prop}
  \label{p11.2} If $\land V$ is a Sullivan algebra and each dim$\, V^l<\infty$ then
   $\iota_{\land V}  $
  is a quasi-isomorphism.
  \end{Prop}

  \vspace{3mm}\noindent {\sl proof.}  First observe that if deg$\, v= 1$ then
  $$[S^k, \langle\land v\rangle] = [\land v, A_{PL}(S^k)] = \left\{
  \begin{array}{ll}
  \mathbb Q & \mbox{if } k=1\\
  0 & \mbox{if } k>1.\end{array}\right.$$
  Thus $\pi_k\langle\land v\rangle = \mathbb Q$ if $k=1$ and $0$ if $k>1$. It follows that 
  $$H^k(\langle\land v\rangle) = \left\{\begin{array}{ll} \mathbb Q & \mbox{if } k=1\\
  0 & \mbox{if } k>1\end{array}\right.$$
  Moreover the identity of $\langle\land v\rangle$ is adjoint to $\iota_{\land v}$ and so $\iota_{\land v}$ is not null homotopic; Therefore
  $$H(\iota_{\land v}) : \land v \stackrel{\cong}{\longrightarrow} A_{PL}\langle\land v\rangle.$$ 
  
  Now we proceed by induction on dim $V^1$. If $V^1= 0$ the Proposition is proved in (\cite[Theorem 1.5]{RHTII}). Otherwise we have a $\Lambda$-extension decomposition of the form $\land V = \land v \otimes \land W$. In this case the Proposition follows from Theorem \ref{t2}. \hfill$\square$

\newpage

\newpage
\part{Enriched dgl's  and their homotopy theory}

\section{Enriched dgl's and semi-quadratic Sullivan algebras}

\noindent {\bf Definitions.}
\begin{enumerate}
\item[(1)] A \emph{differential graded Lie algebra} (dgl), $(L, \partial)$ is a graded Lie algebra $L= \{L_k\}_{k\geq 0}$ together with a differential $\partial : L_k\to L_{k-1}$ which is also a right derivation:
\begin{eqnarray}
\label{i19}\partial^2= 0 \hspace{5mm}\mbox{and } \partial [x,y] = (-1)^{deg\, y} [\partial x,y]+ [x,\partial y].
\end{eqnarray}
\item[(2)] An \emph{enriched dgl} (edgl) is a dgl, $(L, \partial)$,  together with an inverse system of morphisms $\rho_\alpha : (L, \partial)\to (L_\alpha, \partial_\alpha)$ such that
\begin{enumerate}
\item[$\bullet$] each $(L_\alpha, \partial_\alpha)$ is a finite dimensional nilpotent dgl, and
\item[$\bullet$] $L\stackrel{\cong}{\longrightarrow} \varprojlim_\alpha L_\alpha$. 
\end{enumerate}
In particular, \emph{$L$ is a complete enriched Lie algebra} and $\partial$ is a coherent linear map.
\item[(3)] An enriched dgl, $(L, \partial)$ is \emph{minimal} if $\partial : L\to L^{(2)}$. 
\item[(4)] A \emph{morphism of enriched dgl's} is a morphism of enriched Lie algebras which commutes with the differentials. In particular, morphisms are coherent. We denote the set of morphisms $(L, \partial)\to (L', \partial')$ by $\mbox{dgl}((L, \partial), (L', \partial'))$. 
\item[(5)] A \emph{dgl ideal} in an enriched dgl, $L$, is a closed differential ideal; in which case (Lemma \ref{l1.1}) $L/I$ is an enriched dgl.
\end{enumerate}

\vspace{3mm}\noindent {\bf Remarks.} \begin{enumerate}
\item[1)] For simplicity we will frequently denote an enriched dgl simply by $(L, \partial)$.
\item[2)] The convention that $\partial$ be a right derivation arises from the fact that in applications $L$ frequently acts from the right.
\end{enumerate}

\vspace{3mm} The 1-1 correspondence in \S 4 between complete enriched Lie algebras and quadratic Sullivan algebras extends to enriched dgl's, $(L, \partial)$ as follows. First, let $(\land V, d_1)$ be the quadratic Sullivan model of $L$. Recall from \S 4 that the inverse limit $sL = \varprojlim_\alpha sL_\alpha$ is the dual of a unique direct limit $V = \varinjlim_\alpha V_\alpha$, where $(\land V_\alpha, d_\alpha) = C^*(L_\alpha)$. Then
$$(\land V, d_1) = \varinjlim_\alpha \, (\land V_\alpha, d_{1,\alpha})$$
is by definition the quadratic model of $L$. 

On the other hand, each $\partial_\alpha$ is the dual (up to suspension) of a unique differential $d_{0,\alpha}$ in $V_\alpha$ and, by Proposition \ref{p2.1}, $\partial$ is the dual of the differential $d_0= \varinjlim_\alpha d_{0,\alpha}$ in $V$:
$$<d_0v, sx>= <v, s\partial x>, \hspace{1cm} v\in V, x\in L.$$
Finally, we extend $d_0$ to a derivation in $\land V_0$. Now a minimal Sullivan model $(\land W,d) \stackrel{\simeq}{\to} (\land V, d_0)$ induces an isomorphism $W\stackrel{\cong}{\to} H(V, d_0)$. Let $L_W$ denote the homotopy Lie algebra of $\land W$. Then  there are natural isomorphisms
\begin{eqnarray}
\label{i20}
\left.
\renewcommand{\arraystretch}{1.6}
\begin{array}{l}
s(L, \partial) = \varprojlim_\alpha s(L_\alpha, \partial_\alpha) = \left[\varinjlim_\alpha (V_\alpha, d_{0,\alpha}\right]^\vee = (V, d_0)^\vee\\
H(L, \partial) = s^{-1}H(V,d_0)^\vee= s^{-1}W^\vee = L_W,
\end{array}
\renewcommand{\arraystretch}{1}
\right\}
\end{eqnarray}
which in particular provide an isomorphism $H(L,\partial)\cong L_W$ of complete enriched Lie algebras. It also follows from (\ref{i19}) that $(d_0+d_1)^2= 0$. Thus Lemma \ref{l12.1}(i) immediately below shows that $(\land V, d_0+d_1)$ is a Sullivan algebra.  

\vspace{3mm}\noindent {\bf Definition.} (i) A \emph{semi-quadratic Sullivan algebra} is a Sullivan algebra of the form $(\land V,d=d_0+d_1)$ in which $d_0: V\to V$ and $d_1: V\to \land^2V$.

(ii) The semi-quadratic Sullivan algebra corresponding to an enriched dgl $L$ is its \emph{semi-quadratic Sullivan model}, and is denoted $\land V_L$.

(iii) A cdga morphism $\varphi: (\land V, d_0+d_1) \to (\land W, d_0+d_1)$ is \emph{bihomogeneous} if $\varphi : V\to W$.

\begin{lem}
\label{l12.1}   Let $(\land W,d_0+d)$ be a cdga in which $d_0 : W\to W$, $d : W\to \land^{\geq 2}W$ and $d^2=0$. Then $(\land W, d_0+d)$ is a Sullivan algebra if and only if $(\land W, d)$ is a Sullivan algebra.

\end{lem}

\vspace{3mm}\noindent {\sl proof.}   Suppose first that $(\land W, d_0+d)$ is a Sullivan algebra. Then $\land W = \varinjlim_n \land W_n$ in which $d_0+d : W_{n+1}\to \land W_n$. But this implies that $d : W_{n+1}\to \land W_n$ and so $(\land W, d)$ is a Sullivan algebra.

In the reverse direction, suppose $(\land W,d)$ is a Sullivan algebra. If $P\subset W$ is finite dimensional then $P\subset$ some $W'$ with dim$\, W'<\infty$ and $\land W'$ preserved by $d$. It follows that $\land (W'+ d_0W')$ is preserved by both $d$ and $d_0$. Therefore it is sufficient to prove that $(\land W, d_0+d)$ is a Sullivan algebra when dim$\, W<\infty$.

In this case we set $W_0= W\cap \mbox{ker}\, d$. Then $\land W$ decomposes as the $\Lambda$-extension
$$(\land W_0, d_0) \to (\land W_0 \otimes \land Z, d_0+d) \to (\land Z, \overline{d_0}+ \overline{d})$$
of cdga's and the result follows by induction on dim$\, W$. \hfill$\square$

\begin{Prop}
\label{p12.1} 
\begin{enumerate}
\item[(i)] The correspondence $(L, \partial)\leadsto (\land V, d_0+d_1)$ is a bijection between enriched dgl's and semi-quadratic Sullivan algebras.
\item[(ii)] This correspondence identifies the morphisms $f : (L, \partial) \to (L', \partial')$ with the bihomogeneous morphisms $\varphi : (\land V, d_0+d_1)\leftarrow (\land V', d_0'+d_1')$ via 
$$\mbox{suspension of }f= (\varphi\vert_{V'})^\vee.$$
In particular $\varphi : (\land V, d_1)\leftarrow (\land V', d_1')$ is the morphism of quadratic Sullivan algebras corresponding to $f : L\to L'$. 
\item[(iii)]   $(L, \partial)$ is minimal if and only if $d_0(V\cap \mbox{ker}\, d_1)= 0$.
\item[(iv)] If $I\subset L$ is a dgl ideal in an enriched dgl, $L$, then the morphisms $\land V_{L/I} \to \land V_L \to \land V_I$ decompose  $\land V_L$ as a Sullivan extension
$$\land V_L = \land V_{L/I}\otimes \land V_I.$$
Conversely, if an inclusion $\land W\to \land V_L$ of semi-quadratic Sullivan algebras decomposes $\land V_L$ as a Sullivan extension $\land W\otimes \land Z$ then the enriched dgl corresponding to $\land Z$ is a dgl ideal $I\subset L$ and $\land W = \land V_{L/I}$.  
\end{enumerate}
\end{Prop}

\vspace{3mm}\noindent {\sl proof.}

(i)   This is   immediate.

(ii) The correspondence associates to a morphism $f$ the map $\varphi  : V\leftarrow V'$, via the equation
$$<\varphi v', sx> = <v', sf(x)>.$$
In this case $\varphi \circ d_0= d_0'\circ \varphi$ and $\varphi \circ d_1= d_1'\circ \varphi$.  The bijection is a  straightforward consequence of the definitions and Propositions \ref{p2.1} and \ref{p4.1}. 

(iii) This follows because (Lemma \ref{l4.1}) the isomorphism $sL \cong V^\vee$ induces an isomorphism $s(L/L^{(2)})\cong (V\cap \mbox{ker}\, d_1)^\vee$.

(iv) This is immediate from the definitions and \S 4.1. \hfill$\square$

\vspace{3mm}\noindent {\bf Definition.} If $(\land V, d_0+d_1)$ is the semi-quadratic Sullivan algebra corresponding to an enriched dgl $(L, \partial)$ then it is the \emph{semi-quadratic Sullivan model} of $(L, \partial)$ and $(L, \partial)$ is the \emph{enriched dgl model} of $(\land V, d_0+d_1)$. 

\vspace{3mm}\noindent {\bf Remark.} In working with enriched dgl's and their semi-quadratic models it is important to note that
\begin{enumerate}
\item[1)] Not every morphism is bihomogeneous, and
\item[2)] The semi-quadratic Sullivan model of a minimal enriched dgl is a Sullivan algebra, but is not minimal unless $d_0=0$.
\end{enumerate}

  \vspace{3mm}
Finally, if $(\land V, d_0+d_1)$ is a semi-quadratic Sullivan algebra then it is immediate that $\land V(n) := \land V^{\leq n} \otimes \land d_0V^n$ is preserved by both $d_0$ and $d_1$. Thus
  $$(\land V, d_0+d_1) = \varinjlim_n \, (\land V(n), d_0+d_1).$$
  
  \begin{Prop}
  \label{p12.2} \begin{enumerate}
  \item[(i)] A bihomogeneous morphism $\varphi : (\land V, d_0+d_1) \to (\land W, d_0+d_1)$ of semi-quadratic Sullivan algebras restricts to morphisms
  $$\varphi(n) : (\land V(n), d_0+d_1) \to (\land W(n), d_0+d_1), \hspace{5mm} n\geq 0.$$
  \item[(ii)] If $\psi\sim \varphi : (\land V, d_0+d_1)\to (\land W, d_0+d_1)$ then 
  $$\varphi(n) \sim \psi(n), \hspace{5mm} n\geq 0.$$
  \item[(iii)] A bihomogeneous morphism $\varphi : (\land V, d_0+d_1) \to (\land W, d_0+d_1)$ is a quasi-isomorphism if and only if each $\varphi (n)$ is a quasi-isomorphism.
  \end{enumerate}
  \end{Prop}
  
  \vspace{3mm}\noindent {\sl proof.} (i) Since $\varphi$ preserves degrees it maps $\land V^{\leq n}$ to $\land W^{\leq n}$. But if $v\in V^n$ then
  $$\varphi (d_0v) = \varphi(dv-d_1v) = d\varphi v-\varphi d_1v= d_0\varphi v +d_1\varphi v -\varphi d_1v \in W(n),$$
  and so $\varphi : V(n)\to W(n)$.
  
  (ii) This follows in the same way as (i).
  
  (iii) Since $\varphi$ is a quasi-isomorphism of Sullivan algebras it has a homotopy inverse $\widehat{\varphi} : \land V\leftarrow \land W$. It follows from (i) and (ii) that $\widehat{\varphi}$ and $\varphi$ restrict to homotopy inverses between each $V(n)$ and $W(n)$. 
  
  \hfill$\square$

 \subsection{Homology of a semi-quadratic Sullivan algebra}
 
 \begin{lem}
 \label{l12.2} If $(\land W, d_0+d_1)$ is a semi-quadratic Sullivan algebra then
 $$(\land W, d_0+d_1) = \varinjlim_\alpha (\land W_\alpha, d_0+d_1)$$ 
  is the limit of sub semi-quadratic Sullivan algebras $(W_\alpha \subset W)$ for which each $H(W_\alpha, d_0) \to H(W, d_0)$ is injective.
  \end{lem}
  
 \vspace{3mm}\noindent {\sl proof.}  A sub cdga $(\land W_\beta, d_0+d)$ of $\land W$ will be called \emph{regular} if dim$\, W_\beta <\infty$ and if
 $$H(W_\beta, d_0)\to H(W, d_0)$$
 is injective. We show by induction on $n$ that any finite dimensional subspace $P\subset W^{\leq n}$ is contained in a regular cdga of the form $(\land V^{\leq n} \oplus d_0V^n).$
 
 Assuming this for some $n$, let $P\subset W^{\leq n+1}$ be finite dimensional. It follows from (i) that $P\subset Z= Z^{\leq n+1} + d_0(Z^{\leq n+1})$, where dim$\, Z<\infty$ and $\land Z$ is preserved by $d_0$ and $d_1$. Now choose a subspace $S\subset W^n$ so that $S\cap Z^n= 0$ and
 $$0\to S\stackrel{d_0}{\to} H^{n+1}(Z, d_0) \to H^{n+1}(W, d_0)$$
 is exact. By our induction hypothesis there is a regular cdga of the form $\land (V^{\leq n} \oplus d_0V^n)$ and for which 
  $$S\oplus Z^{\leq n}\subset V^{\leq n}.$$
Now, set $V^{n+1} = Z^{n+1}+ d_0V^n$. We show that $\land (V^{\leq n+1} \oplus d_0V^{n+1})$ is regular. Since $P\subset W^{\leq n+1}$ this will complete the proof.
  
  But since $\land Z$ and $\land (V^{\leq n} \oplus d_0V^n)$ are both preserved by $d_0$ and by $d_1$ it follows that so is $\land (V^{\leq n+1}\oplus d_0V^{n+1})$. Moreover, for $k\leq n$ the map $H^k(V^{\leq n+1}\oplus d_0V^{n+1})\to H^k(W)$ coincides with the injection $H^k(V^{\leq n} \oplus d_0V^n)\to H^k(W)$. Also, a cycle in $V^{n+1}$ is homologous to a cycle in $Z^{n+1}$. Thus if it is a boundary in $W^{n+1}$, it is in $d_0(S)$ and hence a boundary in $W^{n+1}$.  
  
  \hfill$\square$

\vspace{3mm}\noindent {\bf Corollary.} If $(L, \partial)$ is complete and dim$\, H(L)/ [H(L), H(L)]<\infty$ then $H(L)$ is pronilpotent.

\vspace{3mm}\noindent {\sl proof.} This follows directly from Proposition \ref{p3.1}.
\hfill$\square$

\begin{Prop}
\label{p12.3} The enriched structure in an enriched dgl can be represented by surjections $\rho_\beta : (L, \partial)\to (L_\beta, \partial_\beta)$ for which $(H(L), \{H(\rho_\beta)\})$ is a complete enriched Lie algebra.  
\end{Prop}

\vspace{3mm}\noindent {\sl proof.} An enriched structure in $(L, \partial)$ is represented by surjections $\rho_\beta : L\to L_\beta$ which are dual to inclusions $j_\beta : (V_\beta, d_0+d_1) \to (V, d_0+d_1)$. By Lemma \ref{l12.2} we may assume that each $H(j_\beta) : H(V_\beta, d_0)\to H(V, d_0)$ is injective. Then $H(\rho_\beta) = s^{-1}H(j_\beta)^\vee$ is surjective and $(H(L), \{H(\rho_\beta)\})$ is an enriched Lie algebra.
Since $L$ is complete, $(L, \partial) = \varprojlim_\beta  (L_\beta, \partial_\beta)$. Thus it follows from Proposition \ref{p1.1} that $H(L) = \varprojlim\, H(L_\beta)$ is   complete. \hfill$\square$

\begin{Prop}
\label{p12.4} Suppose $\varphi : (\land V, d_0+d_1) \leftarrow (\land V', d'_0+d'_1)$ is the bihomogeneous morphism corresponding to a morphism $f: (L, \partial) \to (L', \partial')$ of enriched dgl's. Then
\begin{enumerate}
\item[(i)]  $H(\varphi, d_0) : (\land H(V,d_0), H(d_1)) \leftarrow (\land H(V', d_0'), H(d_1))$ is the bihomogeneous morphism corresponding to $H(f) : (H(L), 0) \to (H(L'), 0)$. 
\item[(ii)] $H(f)$ is an isomorphism $\Leftrightarrow $ $H(\varphi, d_0)$ is an isomorphism.
\item[(iii)] If $H(\varphi, d_0)$ is an isomorphism then $H(\varphi, d_0+d_1)$ is an isomorphism.  
\end{enumerate}
\end{Prop}

\vspace{3mm}\noindent {\sl proof.} (i) It follows from (\ref{i20}) that $H(\varphi, d_0)^\vee$ is the suspension of $H(f): L_V\to L_{V'}$, where $L_V$ and $L_{V'}$ are the homotopy Lie algebras of $\land H(V,d_0)$ and $\land H(V', d'_0)$. 

(ii) This is immediate from (\ref{i20}).  

(iii) Filtering by wedge degree   produces a spectral sequence  converging   from $H(\land V, d_0)$   to $H(\land V, d_0+d_1)$, and this establishes (ii). \hfill$\square$

   \begin{Prop}
   \label{p12.5} Suppose $\varphi : (\land V, d_0+d_1) \to (\land W,D)$ is a quasi-isomorphism from a semi-quadratic Sullivan algebra to a minimal Sullivan algebra.
   \begin{enumerate}
   \item[(i)] Then $\varphi$ is surjective and restricts to surjective dga morphisms $$(\land^{\geq k}V, d_0+d_1) \stackrel{\simeq}{\longrightarrow} (\land^{\geq k}W,D).$$
   \item[(ii)]The morphism of wedge-filtration spectral sequences induced by $\varphi$ begins with the isomorphism of $E_1$-terms
   $$(\land H(V, d_0), H(d_1)) \stackrel{\cong}{\longrightarrow} (\land W, D_1),$$
   where $D_1$ is the quadratic component of $D$. 
   \item[(iii)] The isomorphism, (ii), of quadratic Sullivan algebras identifies the isomorphism 
   $ H(L, \partial) \cong L_W$ of (\ref{i20}) as an isomorphism of enriched Lie algebras. 
  \end{enumerate}
   \end{Prop}
   
 \vspace{3mm}\noindent {\sl proof.} (i) Since $\varphi$ is a quasi-isomorphism of Sullivan algebras it has a homotopy inverse, $\psi : (\land W,D)\stackrel{\simeq}{\longrightarrow} (\land V, d_0+d_1)$. Thus $\varphi\circ \psi$ is a quasi-isomorphism of minimal Sullivan algebras, and is therefore an isomorphism. In particular it restricts to dga automorphisms of each $\land^{\geq k}W$. It follows that $\varphi$ restricts to surjective dga morphisms $\land^{\geq k}V\to \land^{\geq k}W$, $k\geq 0$.

   (ii) Decompose $V$ as
   $$V = Z\oplus S\oplus d_0S,$$
   in which ker$\, d_0= Z\oplus d_0S$ and $d_0: S\stackrel{\cong}{\to} d_0S$. Filtration by wedge degree then shows that the inclusions $Z, S$ and $d_0S$ in $\land V$ extend to an isomorphism
   $$\chi : \land Z\otimes \land S\otimes \land d_0S \stackrel{\cong}{\longrightarrow} \land V$$
   of graded algebras. It follows that a differential, $d_Z$, in $\land Z$ is determined by the requirement that $\chi$ be a dga isomorphism. Since $\varphi$ is a quasi-isomorphism its restriction, $\chi_Z$, to $(\land Z, d_Z)$ is a quasi-isomorphism of minimal Sullivan algebras, and hence an isomorphism.
   
   On the other hand $\chi$ also induces a morphism of the spectral sequences induced by filtration by wedge degrees, and the $E_0$-term is an isomorphism
   $$(\land Z,0)\otimes (\land (S\oplus d_0S), d_0) \stackrel{\cong}{\longrightarrow} (\land V, d_0).$$ Thus the $E_1$- morphism is an isomorphism
   $$(\land Z, (d_Z)_1) \stackrel{\cong}{\longrightarrow} (\land H_0(V), H(d_1)). $$
   
   Finally, filtering $\land Z$ and $\land W$ by wedge degrees converts $\chi_Z$ into a morphism of spectral sequences for which the $E_1$-term is an isomorphism $(\land Z, (d_Z)_1) \stackrel{\cong}{\to} (\land W, D_1)$. This, combined with the isomorphism above, establishes (ii).
   
   (iii) is immediate from the definitions.
 
   \hfill$\square$

 \subsection{Free products}
 
 The differentials in enriched dgl's $(L', \partial')$ and $L'', \partial'')$ extend uniquely to a differential, $\partial,$ in the free product $L'\, \widehat{\amalg}\, L''$ (recalled in \S 10.2). Then $(L= L'\, \widehat{\amalg}\, L'', \partial)$ is an enriched dgl, the \emph{free product} of $(L', \partial')$ and $(L'', \partial'')$. 
 
 \begin{Prop}
 \label{p12.6} \begin{enumerate}
 \item[(i)] Any two morphisms, $\iota', \iota'': (L', \partial'), (L'', \partial'') \to (E, \partial_E)$ extend uniquely to a morphism
 $$(L'\, \widehat{\amalg}\, L'', \partial) \to (E, \partial_E).$$
 \item[(ii)] If $(\land W, d_0+d_1)$ and $(\land Q, d_0+d_1)$ are the semi-quadratic models of $(L', \partial')$ and $(L'', \partial'')$,   then there is a bihomogeneous quasi-isomorphism
 $$\xymatrix{
 (\land V, d_0+d_1) \ar[rr]^-\simeq && (\land W, d_0+d_1)\times_{\mathbb Q}(\land Q, d_0+d_1).}$$
 \item[(iii)]  Denote by $(L, \partial)$ the edgl corresponding to $(\land V, d_0+d_1)$. The edgl morphisms $L',L''\to L$ corresponding to the surjections $\land V\to \land W, \land Q$ then extend to an isomorphism
 $$(L',\partial')\, \widehat{\amalg}\, (L'', \partial'') \stackrel{\cong}{\longrightarrow} (L,\partial).$$
 \end{enumerate}
 \end{Prop}
 
 \vspace{3mm}\noindent {\sl proof.} (i) This follows from Lemma \ref{l10.1} since by construction the morphism $L'\, \widehat{\amalg}\, L''\to E$ is compatible with the differentials.
 
 (ii) The proof of Proposition \ref{p10.1} provides a quasi-isomorphism
 $$\varphi : (\land V, d_1) \to (\land W, d_1)\times_{\mathbb Q} (\land Q, d_1)$$
 in which $(\land V, d_1)$ is a quadratic Sullivan algebra and $\varphi : V\to W\oplus Q$. Let $L$ be the homotopy Lie algebra of $(\land V, d_1)$. Then (Proposition \ref{p10.3}) the surjections $\land V\to \land W, \land Q$ induce an isomorphism,
 $$L'\, \widehat{\amalg}\, L''\stackrel{\cong}{\longrightarrow} L,$$
 of enriched Lie algebras.
 
 In particular, $\partial'$ and $\partial''$ extend to a differential, $\partial$, in $L$. This in turn defines a differential $d_0$ in $\land V$ for which $(\land V, d_0+d_1)$ is a semi-quadratic Sullivan algebra and $\varphi\circ d_0= d_0\circ \varphi$. Filtering by degree $-$ wedge produces a morphism of spectral sequences converging from $H(\varphi, d_1)$ to $H(\varphi, d_0+d_1)$. Thus
 $$\varphi : (\land V, d_0+d_1) \to (\land W, d_0+d_1) \times_{\mathbb Q} (\land Q, d_0+d_1)$$
 is a quasi-isomorphism.

 (iii) This is immediate from the construction in (ii).
 \hfill$\square$

 \section{Profree dgl's}
 
 As in Part II, $\mathbb L_T$ denotes the free graded Lie algebra freely generated by a subspace $T= T_{\geq 0}$. At the start of Part II a complete enriched structure in $T$ is extended to a natural enriched Lie algebra structure in $\mathbb L_T$, and  for which $T$ is a closed subspace. Then, by definition, a \emph{profree Lie algebra} is a complete enriched Lie algebra of the form $L= \overline{\mathbb L}_T$.     In this case $L$ is the direct sum,
 $$L = T\oplus L^{(2)},$$
in which $T$ is closed in $L$. Finally Theorem 1, together with Lemma \ref{l6.2}, establishes the equivalence of the following three conditions on a complete enriched Lie algebra, $L$, and its quadratic model $(\land V,d_1)$:
\begin{eqnarray}
\label{i21}\left. \begin{array}{ll}
(i) & L = \overline{\mathbb L}_T\mbox{ is profree}.\\
(ii) & H(\land V, d_1) = \mathbb Q \oplus (V\,\cap  \mbox{ker}\, d_1).\\
(iii) & H^{[2]}(\land V, d_1)= 0.
 \end{array}\right\}
 \end{eqnarray}
 (Here $H(\land V, d_1)= \oplus_k H^{[k]}(\land V, d_1)$ is the decomposition induced by wedge degree.)

 \vspace{3mm}\noindent {\bf Definition.} A \emph{profree dgl} is a dgl $(L, \partial)$ in which $L= \overline{\mathbb L}_T$ is profree.
 
   \begin{Prop}
 \label{p13.1} 
 \begin{enumerate}
 \item[(i)] An enriched dgl $(L, \partial)$ admits a surjective quasi-isomorphism from a profree dgl $(\overline{\mathbb L}_T, \partial)$.
 \item[(ii)] If $\rho : (E, \partial)\to (F, \partial)$ is a surjective quasi-isomorphism of enriched dgl's, then any morphism $g : (\overline{\mathbb L}_T, \partial)\to (F, \partial)$ from a profree dgl lifts to provide a commutative diagram
 $$
 \xymatrix{ && (E, \partial)\ar[d]^\rho\\
 (\overline{\mathbb L}_T, \partial) \ar[rr]^g \ar[rru]^h&& (F, \partial).}$$
 \item[(iii)] Any profree dgl $(L, \partial)$ is the free product $(L, \partial)= (L', \partial') \, \widehat{\amalg}\, (L'', \partial'')$ of a profree dgl satisfying $H(L')=0$ and a minimal profree dgl $(L'', \partial'')$.
 \end{enumerate}
 \end{Prop}

 \vspace{3mm}\noindent {\sl proof.} (i). Let $(\land V, d_0+d_1)$ be the semi-quadratic model of $(L, \partial)$. We construct by induction on $n$ semi-quadratic Sullivan models of the form $(\land V\otimes \land W(n), d_0+d_1)$ with $W(0)\subset W(1) \subset \dot W(n)\subset W(n+1)\dots$, and such that
 \begin{enumerate}
 \item[(1)] The induced map $H^{[2]}(\land V\otimes \land W(n),d_1)\to H^{[2]}(\land V\otimes \land W(n+1), d_1)$ is zero.
 \item[(2)] The induced map $H(V\oplus W(n), d_0)\to H(V\oplus W(n+1), d_0)$ is an isomorphism.
 \end{enumerate}
 Here $H^{[2]}$ denotes as before the component of the cohomology   of wedge degree 2. 
 Suppose $\land V\otimes \land W(n)$ has been defined. Then let $Z= s^{-1}H^{[2]}(\land V\otimes \land W(n),d_1)$, $W(n+1) = W(n) \oplus Z\oplus d_0(Z)$, and for $x=s^{-1}\Omega\in Z$, we choose a $d_1$-cycle $\omega$ with $[\omega]=\Omega$, and we define $d_1(x)= \omega$ and $d_1d_0x = -d_0\omega$.
 
 Denote $W = \cup_n W(n)$. Since $H^{[2]}=0$, by Lemma \ref{l6.2},  $(\land V\otimes \land W,d_0+d_1)$ is the semi-quadratic model of a profree dgl $(\overline{\mathbb L}_T,\partial)$, and the injection $\land V\to \land V\otimes \land W$ corresponds to a surjective quasi-isomorphism $(\overline{\mathbb L}_T, \partial)\to (L, \partial)$. 
 
 (ii) Since $T$ is closed and $\rho$ is a surjective coherent map there is a coherent linear map $\alpha : T\to E$ such
 that $\rho\circ \alpha = id$. Moreover, since $\rho$ is a coherent surjective quasi-isomorphism, $H(\mbox{ker}\, \rho)= 0$. Thus there is a coherent linear map $\tau : \mbox{ker}\, \rho \cap \mbox{ker}\, \partial\to \mbox{ker}\, \rho$ such that $\partial \circ \tau = id$. 
 
 Now assume $h$ has been defined in $T_{<k}$, and therefore also in $\overline{\mathbb L}_{T_{\leq k}}$.  Then $\sigma \circ \partial_T -\partial \circ \alpha : T_k\to \mbox{ker}\,\rho \cap \mbox{ker}\, \partial$. Extend $h$ to $T_k$ by setting
 $$h x= \alpha x + \tau (h\partial_T-\partial \alpha)(x).$$
 It is immediate that the extension to $\overline{\mathbb L}_{T_{\leq n}}$ is a morphism and that $\rho\circ h= g$.

 (iii). Denote by $\overline{\partial}$ the quotient differential in $L/L^{(2)}$. Then (cf. Proposition \ref{p2.1})
 $$L/L^{(2)} = S \oplus \overline{\partial}S \oplus H,$$
 where $S, \overline{\partial}S$ and $H$ are closed subspaces, $\overline{\partial} : S \stackrel{\cong}{\to} \overline{\partial}S$, and $\overline{\partial}(H)= 0$. Now the isomorphism $T \cong L/L^{(2)}$ from a closed direct summand of $L^{(2)}$ yields a decomposition
 $$T = P\oplus \partial P \oplus Q,$$
 in which $\partial : P\stackrel{\cong}{\to} \partial P$ and $\partial : Q \to L^{(2)}$. It follows that $L= \overline{\mathbb L}_{(T\oplus \partial P\oplus Q)}$.
 
 Moreover, filtering the Lie bracket length provides a spectral sequence converging from $\overline{\mathbb L}_{Q}$ to $H(L)$. This implies that division by $T$ and $\partial T$ yields a quasi-isomorphism
 $$\rho : (L, \partial) \stackrel{\simeq}{\to} (\overline{\mathbb L}_Q, \partial_Q)$$
 where $\partial_Q$ is the quotient differential. 
 
 Finally, by (ii) $\rho$ admits a splitting $\sigma : (\overline{\mathbb L}_Q, \partial_Q) \to (L, \partial)$ and it is immediate from the construction that
 $$id \, \widehat{\amalg}\, \sigma : \overline{\mathbb L}_{P\oplus \partial P)} \, \widehat{\amalg}\, (\overline{\mathbb L}_Q, \partial_Q) \stackrel{\cong}{\longrightarrow} (L, \partial).$$
 Since $\partial : Q \to \mathbb L_Q^{(2)}$, this completes the proof. 
 
 \hfill$\square$

 \vspace{3mm} Given a profree dgl, $(\overline{\mathbb L}_T, \partial)$, division by $\overline{\mathbb L}_T^{(2)}$ is a coherent surjection 
 $$\rho : (\overline{\mathbb L}_T, \partial) \to (T, \overline{\partial}).$$
 
 \begin{Prop}
 \label{p13.2} Suppose $(\overline{\mathbb L}_T, \partial)$ is a profree dgl with semi-quadratic Sullivan model $(\land V, d_0+d_1)$, and let $\xi : \land^{\geq 1}V\to V$ denote the surjection (\S 9).
 \begin{enumerate}
 \item[(i)] Any coherent chain map $(T, \overline{\partial})\to (L, \partial)$ into an enriched complete dgl extends uniquely to a morphism $(\overline{\mathbb L}_T, \partial) \to (L, \partial)$.
 \item[(ii)] $H(V\cap \mbox{ker}\, d_1, d_0) \stackrel{\cong}{\longrightarrow} H^{\geq 1}(\land V, d_0+d_1)$.
 \item[(iii)] The isomorphism (\ref{i20}), $s(\overline{\mathbb L}_T, \partial) \stackrel{\cong}{\longrightarrow} (V, d_0)^\vee$, extends to the commutative diagram
 $$\xymatrix{
 s(\overline{\mathbb L}_T, \partial) \ar[rr]^\cong\ar[d]_{s\rho} && (V, d_0)^\vee\ar[d]\ar[rrd]^{\xi^\vee}\\
 s(T, \overline{\partial}) \ar[rr]^\cong && (V\cap \mbox{ker}\, d_1, d_0)^\vee\ar[rr]^\simeq && (\land^{\geq 1}V, d_0+d_1)^\vee.
 }$$
 \item[(iv)] $(\overline{\mathbb L}_T, \partial)$ is minimal if and only if $\overline{\partial}= 0$, or equivalently, if $d_0\vert_{V\cap \mbox{ker}\, d_1} = 0$.
 \item[(v)] $H(\overline{\mathbb L}_T, \partial)= 0$ if and only if $\overline{\mathbb L}_T^{(2)}$ has a closed direct summand of the form $T= T(1) \oplus \partial (T(1))$ where $\partial : T(1) \stackrel{\cong}{\longrightarrow} \partial T(1)$. 
 \end{enumerate}
 \end{Prop}

 \vspace{3mm}\noindent {\bf Corollary.} If $(\overline{\mathbb L}_T, \partial)$ is minimal then $sT \cong H^{\geq 1}(\land V, d_0+d_1)^\vee$.

 \vspace{3mm}\noindent {\sl proof of Proposition \ref{p13.2}.} (i) This is an immediate consequence of the definition of a profree dgl.
 
 (ii) Filtration of $\land V$ by degree $-$ wedge degree results in a spectral sequence converging to $H(\land V, d_0+d_1)$. Its $E_1$-term and $E_2$-term are respectively $H(\land V, d_1)$ and $(\mathbb Q \oplus (V\cap \mbox{ker}\, d_1), d_0)$. This proves (ii).
 
 (iii) By Lemma \ref{l4.1}, the inclusion $V\cap \mbox{ker}\, d_1 \to V$ is dual to $\overline{\mathbb L}_T \to \overline{\mathbb L}_T/ \overline{\mathbb L}_T^{(2)}$, since $(\land V, d_1)$ is the quadratic model of $\overline{\mathbb L}_T$. This gives (iii).
 
 (iv) By definition, $(\overline{\mathbb L}_T, \partial)$ is minimal if and only if $\partial : \overline{\mathbb L}_T \to \overline{\mathbb L}_T^{(2)}$. This gives (iv). 
 
 (v) If $H(\overline{\mathbb L}_T, \partial) = 0$ it follows from the isomorphism $s(\overline{\mathbb L}_T, \partial) \cong (V, d_0)^\vee$ that $H(V, d_0)= 0$. This then implies that $(\land V, d_0+d_1)$ is a contractible Sullivan algebra, and therefore has zero homology in positive degrees. Now (ii) and (iii) combine to give $H(T, \overline{\partial}) = 0$.
 
 Write $T = T(1) \oplus \overline{\partial} T(1)$ and observe that it follows that
 $$\overline{\mathbb L}_T = \overline{\mathbb L}_{T(1)}\, \widehat{\amalg}\, \overline{\mathbb L}_{\overline{\partial}T(1)}.$$
 Now it follows from the Corollary to Lemma \ref{l3.1} that
 $$\overline{\mathbb L}_{T(1)} \,\widehat{\amalg}\, \overline{\mathbb L}_{\partial T(1)} \to \overline{\mathbb L}_T$$
 is an isomorphism, and this is clearly an isomorphism of dgl's.
 
 \hfill$\square$

 \begin{Prop}
 \label{p13.3}
 Let $f : (L, \partial)\to (L', \partial')$ be a morphism of   profree dgl's with corresponding bihomogeneous morphism
 $$\varphi : (\land V, d_0+d_1) \longleftarrow (\land V', d'_0+d'_1)$$
 of their semi-quadratic Sullivan algebras. Then the following three conditions are equivalent
 \begin{enumerate}
 \item[(i)] $H(\overline{f})$ is an isomorphism
 \item[(ii)] $H(f)$ is an isomorphism
 \item[(iii)] $H(\varphi)$ is an isomorphism.
 \end{enumerate}
 \end{Prop}

 \vspace{3mm}\noindent {\sl proof.} It follows from the diagram in Proposition \ref{p13.2}(iii) that $H(\overline{f})$ is an isomorphism if and only if $H(\varphi)$ is an isomorphism. Moreover, by Proposition \ref{p13.2}(iii), $H(f)$ is an isomorphism if and only if $H(\varphi, d_0)$ is an isomorphism. Finally, $\varphi : (\land V', d_0+d_1) \to (\land V, d_0+d_1)$ induces a morphism of minimal Sullivan algebras of the form
 $$\psi : (\land H(V', d_0), \widehat{d}) \longrightarrow (\land H(V, d_0), \widehat{d}).$$
 Here the linear map induced from $\psi$ by division by $\land^{\geq 2}$ is $H(\varphi, d_0)$.  Therefore,
 $$\renewcommand{\arraystretch}{1.3}
 \begin{array}{ll}
 H(\varphi, d_0) \mbox{ is an isomorphism } & \Longleftrightarrow \psi \mbox{ is an isomorphism}\\
 & \Longleftrightarrow H(\psi) \mbox{ is an isomorphism }\\
 & \Longleftrightarrow H(\varphi) \mbox{ is an isomorphism.}
 \end{array}
 \renewcommand{\arraystretch}{1}
 $$
 
 \hfill$\square$
 
 \vspace{3mm}\noindent {\bf Corollary.} If $L$ and $L'$ are minimal then the conditions of Proposition \ref{p13.3} are equivalent to each of the conditions
 $$ (i) f \mbox{ is an isomorphism} \hspace{5mm} (ii) \varphi \mbox{ is an isomorphism.}$$
 
 \vspace{3mm}\noindent {\sl proof. } If $L$ and $L'$ are minimal then $\overline{f}= H(\overline{f})$ is an isomorphism. Thus $f$ maps a closed direct summand, $T$, of $L^{(2)}$ isomorphically onto a direct summand, $T'$, of $(L')^{(2)}$. Since $f$ is coherent, $T'$ is closed and $L'$ is the closure of $\mathbb L_{T'}$. It follows that $f$ is an isomorphism. Since, in particular, $\varphi : (\land V', d_1)\to (\land V, d_1)$ is the morphism of quadratic Sullivan algebras corresponding to $f$ it follows that $\varphi$ is also an isomorphism. \hfill$\square$

 \section{Profree extensions and homotopy}
 
 \subsection{Profree extensions}
 
 \vspace{3mm}\noindent {\bf Definition.} 
 \begin{enumerate}
 \item[1.] A \emph{profree extension} is an inclusion of the form $i : (L, \partial) \to (L\, \widehat{\amalg}\, \overline{\mathbb L}_S, \partial)$ of enriched dgl's.
 \item[2.] The \emph{quotient} of a profree extension is the profree dgl $(\overline{\mathbb L}_S, \overline{\partial})$ obtained by dividing by the ideal generated by $L$. The profree extension is \emph{minimal} if the quotient is a minimal profree dgl.
 \item[3.] The \emph{secondary gradation} in $L\, \widehat{\amalg}\, \overline{\mathbb L}_S$ is that uniquely determined by assigning to an iterated bracket of Lie elements $x_i\in L$ and $y_i\in  S$ the secondary degree, $\sum \mbox{deg}\, y_i$.
 \item[4.] The \emph{spectral sequence} of the profree extension is the spectral sequence obtained by the decreasing filtration by secondary degree.
 \item[5.] Two profree extensions $L\to L\, \widehat{\amalg}\, \overline{\mathbb L}_S$, $L\, \widehat{\amalg}\, \overline{\mathbb L}_T$ define a single profree extension $L\, \widehat{\amalg}\, \overline{\mathbb L}_{S\oplus T}$ and $$L\, \widehat{\amalg}\, \overline{\mathbb L}_{S\oplus T} = L\, \widehat{\amalg}\,\overline{\mathbb L}_S\, \widehat{\amalg}\, \overline{\mathbb L}_T =  L\, \widehat{\amalg}\,\overline{\mathbb L}_T\, \widehat{\amalg}\, \overline{\mathbb L}_S.$$ 
 \end{enumerate}
 
 \vspace{3mm} Profree extensions have the five important properties listed in Proposition \ref{p14.1} below. The first three follow immediately from the definitions. The fourth follows from the same argument as the proof of Proposition \ref{p13.1}(ii). The fifth follows from the definition (\S 10.2) of free products.
 
 \begin{Prop}
 \label{p14.1} Let $(L, \partial) \to (L\, \widehat{\amalg}\, \overline{\mathbb L}_S, \partial)$ be a profree extension.
 \begin{enumerate}
 \item[(i)] The spectral sequence of the extension converges from $H(L)\, \widehat{\amalg}\, \overline{\mathbb L}_S$ to $H(L\, \widehat{\amalg}\, \overline{\mathbb L}_S)$.
 \item[(ii)] A morphism $f: (L, \partial)\to (L', \partial')$ of enriched dgl's determines a unique differential in $L'\, \widehat{\amalg}\, \overline{\mathbb L}_S$ for which $f\, \widehat{\amalg}\, id_{\overline{\mathbb L}_S}$ is a morphism of profree extensions.
 \item[(iii)] If $f$ is a quasi-isomorphism so is $f\, \widehat{\amalg}\, id_{\overline{\mathbb L}_S}$.
 \item[(iv)] In a commutative diagram of enriched dgl morphisms,
 $$\xymatrix{L\ar[rr]^f \ar[d]^\iota && E\ar[d]^\rho\\
 L\, \widehat{\amalg}\, \overline{\mathbb L}_S \ar[rr]^g && F,}$$
 suppose that $\iota$ is the inclusion of a profree extension and that $\rho$ is a surjective quasi-isomorphism. Then there is a morphism
 $$h : L\, \widehat{\amalg}\, \overline{\mathbb L}_S \longrightarrow E$$
 such that $h\circ \iota = f$ and $g\circ h=g$.
 \item[(v)] The natural map $L\, \widehat{\amalg}\, \overline{\mathbb L}_S \to \varprojlim_n L/L^{(n)}\, \widehat{\amalg}\, \overline{\mathbb L}_S$ is an isomorphism.
 \end{enumerate}
 \end{Prop}

\begin{Prop}
\label{p14.2}
Suppose $\varphi : (E, \partial_E)\to (F,\partial_F)$ is a morphism of   enriched dgl's. Then
\begin{enumerate}
\item[(i)] $\varphi$ factors as
$$\xymatrix{\varphi : (E, \partial_E) \ar[rr]^\iota && (E\, \widehat{\amalg}\, \overline{\mathbb L}_S, \partial) \ar[rr]^\rho && (F, \partial_F)}$$
in which $\iota$ is the inclusion of a profree extension and $\rho$ is a surjective quasi-isomorphism.
\item[(ii)] $\varphi$ also factors as
$$\xymatrix{\varphi : (E, \partial_E)\ar[rr]^\lambda &&(E\, \widehat{\amalg}\, \overline{\mathbb L}_T, \partial') \ar[rr]^\gamma && (F, \partial_F)}$$
in which $\lambda$ is the inclusion of a minimal profree extension and $\gamma$ is a quasi-isomorphism.
\end{enumerate}
\end{Prop}

\vspace{3mm}\noindent {\sl proof.} (i) Proposition \ref{p13.1} gives a surjective quasi-isomorphism $\sigma : (\overline{\mathbb L}_{S(0)}, \partial) \stackrel{\simeq}{\longrightarrow} (F, \partial)$. By Proposition \ref{p12.6} (i) this extends to the morphism
$$\xymatrix{\rho(0) = \varphi \, \widehat{\amalg}\, \sigma : (E, \partial_E)\, \widehat{\amalg}\, (\overline{\mathbb L}_{S(0)}, \partial) \ar[rr] && (F, \partial_F).}$$
By construction, both $\varphi\, \widehat{\amalg}\, \sigma$ and $H(\varphi\, \widehat{\amalg}\, \sigma)$ are surjective.

Next we construct $S(1) = \oplus_{n\geq 1}S(1)_n$ and a sequence of profree extensions
$$E \, \widehat{\amalg}\, \overline{\mathbb L}_{S(0)} \, \widehat{\amalg}\, L(n) = \left[ E\, \widehat{\amalg}\, \overline{\mathbb L}_{S(0)}\, \widehat{\amalg}\, L(n-1)\right] \, \widehat{\amalg}\, \overline{\mathbb L}_{S(1)_n}.$$
These will be equipped with morphisms
$$\rho(n) : E\, \widehat{\amalg}\, \overline{\mathbb L}_{S(0)}\, \widehat{\amalg}\, L(n) \longrightarrow F$$
with $\rho(n)$ extending $\rho(n-1)$ and $H_{<n}(\rho(n))$ an isomorphism. Then $\overline{\mathbb L}_{S(1)}= \varinjlim_n L(n)$ and the $\rho(n)$ define a quasi-isomorphism
$$\xymatrix{\rho : E\, \widehat{\amalg}\, \overline{\mathbb L}_S \ar[rr]^\simeq && F}$$
with $S= S(0)\oplus S(1)$. 

The construction is by induction on $n$. If $\rho(n-1)$ is constructed we may, since $\rho(n-1)$ is coherent, write
$$\mbox{ker}\, \rho(n-1)\cap \mbox{ker}\, \partial = C \oplus \partial (\mbox{ker}\, \rho(n-1))$$
in which $C = C_{\geq n-1}$ is a closed subspace. Set $S(1)_n = sC_{n-1}$ and define $\partial\vert_{S(1)_n}$ to be the isomorphism $S(1)_n \stackrel{\cong}{\longrightarrow} C_{n-1}$. Then extend  $\rho(n-1)$ to $\rho(n)$ by setting $\rho(n)( S(1)_n) = 0$.

(ii) Denote by $\rho : (E\, \widehat{\amalg}\, \overline{\mathbb L}_S, \partial)\to F$ the construction of (i), and denote by $\partial_S$ the quotient differential in $\overline{\mathbb L}_S$. By Proposition \ref{p13.1}(iii) we may write
$$(\overline{\mathbb L}_S, \partial_S) = (\overline{\mathbb L}_T, \partial_T) \, \widehat{\amalg}\, \overline{\mathbb L}_{(U\oplus \partial_SU)}$$
with $\partial_S : U\stackrel{\cong}{\to} \partial_SU$ and $(\overline{\mathbb L}_T, \partial_T)$ a minimal profree dgl. Filtering by the closed ideals $E^{(n)}\, \widehat{\amalg}\, \overline{\mathbb L}_S$   shows   (Proposition \ref{p14.1} (v)) that
$$E\, \widehat{\amalg}\, (\overline{\mathbb L}_T \, \widehat{\amalg}\, \overline{\mathbb L}_{(U\oplus \partial U)})\longrightarrow E\, \widehat{\amalg}\, \overline{\mathbb L}_S$$
is an isomorphism of dgl's.

This provides the commutative diagram of dgl morphisms
$$
\xymatrix{
E \ar[d]\ar[rr] && 
(E\,\widehat{\amalg}\, \overline{\mathbb L}_{U\oplus \partial U})\, \widehat{\amalg}\, \overline{\mathbb L}_T 
\ar[d]^\simeq_\psi 
\ar[rr]^-\cong && E\, \widehat{\amalg}\, \overline{\mathbb L}_S \ar[rr]^\simeq && F
\\
 E\, \widehat{\amalg}\, \overline{\mathbb L}_T \ar[rr]_=\ar[rru]^h && E\,\widehat{\amalg}\, \overline{\mathbb L}_T
}$$
in which the spectral sequence shows that $\psi$ is a quasi-isomorphism. Thus $h$ exists by Proposition \ref{p14.1}(iv). In this case $h$ must also be a quasi-isomorphism, and the composite
$$\xymatrix{E\ar[rr] && E\, \widehat{\amalg}\, \overline{\mathbb L}_T \ar[rr]^\simeq && F}$$
is the factorization of $\varphi$.

\hfill$\square$

\vspace{3mm}

Next suppose that 
$$\xymatrix{ E \ar[d]_{\gamma_E} \ar[rr] && E\, \widehat{\amalg}\, \overline{\mathbb L}_S \ar[rr]^{\alpha_E}_\simeq && F\ar[d]^{\gamma_F}\\
E' \ar[rr] && E'\widehat{\amalg}\, \overline{\mathbb L}_{S'} \ar[rr]^{\alpha_{E'}}_\simeq && F'}$$
is a commutative diagram of morphisms of enriched dgl's. Then extend the lower row to
$$\xymatrix{ E' \ar[rr] && E'\,\widehat{\amalg}\, \overline{\mathbb L}_{P\oplus \partial P} \, \widehat{\amalg}\, \overline{\mathbb L}_{S'} \ar[rr]_\simeq^{\alpha}&& F'
}$$ by setting $P= F'_{\geq 1}$, $\partial : P \stackrel{\cong}{\longrightarrow} \partial P$, and requiring $\alpha\vert_{P}$ to be the inclusion.

\begin{Prop}
\label{p14.3}
With the notation and hypotheses above there is a morphism $\gamma : E\, \widehat{\amalg}\, \overline{\mathbb L}_S \to E'\, \widehat{\amalg}\, \overline{\mathbb L}_{P\oplus \partial P} \, \widehat{\amalg}\, \overline{\mathbb L}_{S'}$ for which the diagram,
$$\xymatrix{ E \ar[d]_{\gamma_E} \ar[rr] && E\, \widehat{\amalg}\, \overline{\mathbb L}_S \ar[d]^\gamma \ar[rr]^{\alpha_E}_\simeq && F\ar[d]^{\gamma_F}\\
E' \ar[rr] && E'\, \widehat{\amalg}\, \overline{\mathbb L}_{P\oplus \partial P} \, \widehat{\amalg}\, \overline{\mathbb L}_{S'} \ar[rr]^\alpha_\simeq && F',}$$
commutes. \end{Prop}

\vspace{3mm}\noindent {\sl proof.} Assume by induction that $\gamma $ is defined in $S_{<n}$, and let $x$ be a basis element in $S_n$. Then $\partial x\in E\, \widehat{\amalg}\, \overline{\mathbb L}_{S_{<n}}$ and 
$$\alpha \gamma \, \partial x = \partial \gamma_F\alpha_E x.$$
Since $\gamma \partial x$ is a cycle and $\alpha$ is a quasi-isomorphism it follows that 
$$\gamma \partial x= \partial y, \hspace{5mm}\mbox{some } y\in E'\, \widehat{\amalg}\, \overline{\mathbb L}_{P\oplus \partial P} \, \widehat{\amalg}\, \overline{\mathbb L}_{S'}.$$
In particular, $\partial \alpha y= \gamma_F\alpha_E\partial x= \partial \gamma_F\alpha_E x$. Since $\alpha$ is a quasi-isomorphism
$$\alpha y - \gamma_F\alpha_E x= \alpha z + \partial w$$
for some cycle $z$ and some element $w\in F'$. But necessarily deg$\, w\geq 1$, so that $w= \alpha (u)$ for some $u\in P$. Thus
$$\alpha (y-z-\partial u) = \gamma_F\alpha_Ex,$$
and we may extend $\gamma$ by setting $\gamma x= y-z-\partial u$.
\hfill$\square$

\begin{Prop}
\label{p14.4} Suppose
$$\xymatrix{
L \ar[d]^\varphi \ar[rr] && L\, \widehat{\amalg}\, \overline{\mathbb L}_S \ar[d]^\psi \ar[rr] && \overline{\mathbb L}_S\ar[d]^\chi\\
L' \ar[rr] && L'\, \widehat{\amalg}\, \overline{\mathbb L}_{S'} \ar[rr] && \overline{\mathbb L}_{S'}
}$$
is a commutative diagram connecting two profree extensions. If any two of $\varphi, \psi$ and $\chi$ are quasi-isomorphisms, then so is the third.

\end{Prop}

\vspace{3mm}\noindent {\sl proof.} We proceed in a number of Steps.

\vspace{2mm}\noindent \emph{Step One. If $\overline{\mathbb L}_S = \overline{\mathbb L}_{S'}$ and $\psi = \varphi \, \widehat{\amalg}\, id_{\overline{\mathbb L}_S}$, then $\psi$ is a quasi-isomorphism if and only if $\varphi$ is a quasi-isomorphism.}

\vspace{1mm} If $\varphi$ is a quasi-isomorphism it follows from the associated spectral sequence that $\varphi \, \widehat{\amalg}\, id_{\overline{\mathbb L}_S}$ is a quasi-isomorphism.

In the reverse direction, suppose $\varphi \, \widehat{\amalg}\, id_{\overline{\mathbb L}_S}$ is a quasi-isomorphism. Factor $\varphi$ as (Proposition \ref{p14.2})
$$L \to L\, \widehat{\amalg}\, \overline{\mathbb L}_T \stackrel{\simeq}{\longrightarrow} L'$$
in which $L\, \widehat{\amalg}\, \overline{\mathbb L}_T$ is a minimal profree extension. This then extends to
$$L\, \widehat{\amalg}\, \overline{\mathbb L}_S \to L\, \widehat{\amalg}\, \overline{\mathbb L}_T \, \widehat{\amalg}\, \overline{\mathbb L}_S \to L'\, \widehat{\amalg}\, \overline{\mathbb L}_S,$$
in which the second morphism is a quasi-isomorphism by the argument above. Since the hypothesis $L\, \widehat{\amalg}\, \overline{\mathbb L}_S \to L'\, \widehat{\amalg}\, \overline{\mathbb L}_S$ is also a quasi-isomorphism, it follows that
$$L\, \widehat{\amalg}\, \overline{\mathbb L}_S \to L\, \widehat{\amalg}\, \overline{\mathbb L}_S \, \widehat{\amalg}\, \overline{\mathbb L}_T = L\, \widehat{\amalg}\, \overline{\mathbb L}_T \, \widehat{\amalg}\, \overline{\mathbb L}_S$$
is also a quasi-isomorphism.

If $T = T_{\geq n}$ a simple argument then shows that $T_n \subset \partial (T_{n+1})+ I $ where $I$ is the ideal generated by $L\, \widehat{\amalg}\, \overline{\mathbb L}_S$. This contradicts   the minimality of the extension $L\to L\,\widehat{\amalg}\, \overline{\mathbb L}_T$.

\vspace{2mm}\noindent \emph{Step Two. If $\varphi$ is a quasi-isomorphism then $\psi$ is a quasi-isomorphism if and only if $\chi$ is a quasi-isomorphism.}

\vspace{1mm} Proposition \ref{p14.2} provides a quasi-isomorphism $L\, \widehat{\amalg}\, \overline{\mathbb L}_T\stackrel{\simeq}{\to} 0$ from a minimal profree extension. If $\varphi$ is a quasi-isomorphism it follows from Step One that $L\, \widehat{\amalg}\, \overline{\mathbb L}_T \stackrel{\simeq}{\longrightarrow} L'\, \widehat{\amalg}\, \overline{\mathbb L}_T$ and so $H(L'\, \widehat{\amalg}\, \overline{\mathbb L}_T)= 0$. This implies that the inclusions
$$L\, \widehat{\amalg}\, \overline{\mathbb L}_T \longrightarrow L\, \widehat{\amalg}\overline{\mathbb L}_T \, \widehat{\amalg}\, \overline{\mathbb L}_S \hspace{5mm}\mbox{and }\hspace{2mm} L'\, \widehat{\amalg}\, \overline{\mathbb L}_T \longrightarrow L'\, \widehat{\amalg}\, \overline{\mathbb L}_T\, \widehat{\amalg}\, \overline{\mathbb L}_{S'}$$
extend to isomorphisms
$$(L\, \widehat{\amalg}\, \overline{\mathbb L}_T, \partial) \, \widehat{\amalg}\, (\overline{\mathbb L}_S, \partial_S) \stackrel{\cong}{\longrightarrow} L\, \widehat{\amalg}\, \overline{\mathbb L}_T\, \widehat{\amalg}\, \overline{\mathbb L}_S = L\, \widehat{\amalg}\, \overline{\mathbb L}_S\,\widehat{\amalg}\, \overline{\mathbb L}_T$$
and
$$(L'\, \widehat{\amalg}\, \overline{\mathbb L}_T, \partial) \, \widehat{\amalg}\, (\overline{\mathbb L}_{S'}, \partial_{S'}) \stackrel{\cong}{\longrightarrow} L'\, \widehat{\amalg}\, \overline{\mathbb L}_T\, \widehat{\amalg}\, \overline{\mathbb L}_{S'} = L'\, \widehat{\amalg}\, \overline{\mathbb L}_{S'}\,\widehat{\amalg}\, \overline{\mathbb L}_T$$

Now by Step One, $\psi$ is a quasi-isomorphism if and only if $\psi\, \widehat{\amalg}\, id_{\overline{\mathbb L}_T}$ is a quasi-isomorphism. Since $H(L\, \widehat{\amalg}\, \overline{\mathbb L}_T)= 0= H(L'\, \widehat{\amalg}\, \overline{\mathbb L}_T)$ the isomorphisms above show that this holds if and only if $\chi$ is an isomorphism.

\vspace{2mm}\noindent \emph{Step Three. Completion of the proof.}

\vspace{1mm} Here we may assume that $\psi$ and $\chi$ are quasi-isomorphisms, and we have to show that $\varphi$ is a quasi-isomorphism. Now Proposition \ref{p14.2} extends $\varphi$ to a quasi-isomorphism
$$\alpha : L\, \widehat{\amalg}\, \overline{\mathbb L}_U \stackrel{\simeq}{\longrightarrow} L'$$
from a minimal profree extension. This gives the commutative diagram,
$$\xymatrix{
L\, \widehat{\amalg}\overline{\mathbb L}_U \, \widehat{\amalg}\, \overline{\mathbb L}_S\ar[d]_\beta && L\, \widehat{\amalg} \overline{\mathbb L}_S\ar[ll]_{\gamma \, \widehat{\amalg}\, id_{\overline{\mathbb L}_S}} \ar[dll]^\psi\\
L'\, \widehat{\amalg}\, \overline{\mathbb L}_{S'},
}$$
in which $\gamma$ is the inclusion $L\to L\, \widehat{\amalg}\, \overline{\mathbb L}_U$. Since $\alpha$ and $\chi$ are quasi-isomorphisms, it follows from Step Two that $\beta$ is a quasi-isomorphism. Since $\psi$ is a quasi-isomorphism, so is $\gamma \, \widehat{\amalg}\, id_{\overline{\mathbb L}_S}$. By Step One, $\gamma $ is also a quasi-isomorphism. It is straightforward to see that, since $L\, \widehat{\amalg} \, \overline{\mathbb L}_U$ is a minimal extension, this implies that $U= 0$ and hence that $\varphi$ is a quasi-isomorphism.

\hfill$\square$

 \begin{Prop}
 \label{p14.5} Any quasi-isomorphism $f : (L, \partial) \to (L', \partial)$ of enriched dgl's extends to a surjective quasi-isomorphism
 $$\widehat{f} : (L, \partial) \, \widehat{\amalg}\, (\overline{\mathbb L}_T, \partial) \to (L', \partial)$$
 from a profree extension for which $H(\overline{\mathbb L}_T,  {\partial})= 0$.
 \end{Prop}
 
 \vspace{3mm}\noindent {\sl proof.} Let $P$ be a complete enriched vector space isomorphic with $L'_{\geq 1}$ and let $T = P\oplus \partial P$ with $\partial : P\stackrel{\cong}{\to} \partial P$. Use the inclusion $P\to L'$  to extend $f$ to a morphism
 $$\widehat{f} : L\, \widehat{\amalg}\, \overline{\mathbb L}_T \to L'.$$
 Since (cf Proposition \ref{p13.3}) $H(\overline{\mathbb L}_T)= 0$ it follows that the inclusion $i : L\to L\, \widehat{\amalg}\, \overline{\mathbb L}_T$ is a quasi-isomorphism.
 
 Finally, by construction $\widehat{f}$ is surjective in degrees $\geq 1$. Since $H(f)$ is an isomorphism an element $x'\in L'_0$ has the form $x'= f(x) + \partial y$ for some $x\in L$ and $y\in L_1'$. Since $y\in$ Im$\, \widehat{f}$ it follows that $\widehat{f}$ is surjective in degree $0$. \hfill$\square$

 \vspace{3mm} Finally, if $(L, \partial)\to (L\, \widehat{\amalg}\, \overline{\mathbb L}_S, \partial)$ is a profree extension then division by the ideal $I$ generated by $L$ is a surjection
 $$(L\, \widehat{\amalg}\, \overline{\mathbb L}_S, \partial) \to (\overline{\mathbb L}_S, \partial_S)$$
 of enriched dgl's. The corresponding morphism of semi-quadratic Sullivan algebras is an inclusion $\land W \leftarrow \land V_S$. It decomposes $\land W$ as a Sullivan extension
 $$\land W = \land V_{\overline{\mathbb L}_S}\otimes \land Z.$$
 In particular, this gives the short exact sequence,
 $$0\to Z^\vee \to W^\vee \to V_S^\vee\to 0,$$
 and this establishes (cf. Proposition \ref{p12.1} (iv))
 
 \begin{Prop}
 \label{p14.6} With the hypotheses and notation above, the quotient $\land Z$ is the semi-quadratic Sullivan model of the ideal $I$.
 \end{Prop}

 \subsection{Cylinder objects and homotopy}

 \vspace{3mm} Important examples of a profree extension are the cylinder objects defined next, which we use to define homotopy.

 \vspace{3mm}\noindent {\bf Definition.}
 \begin{enumerate}
 \item[1.] A \emph{cylinder object} for an enriched dgl, $(L, \partial)$ is a commutative diagram of enriched dgl morphisms:
 $$\xymatrix{
 (L\, \widehat{\amalg}\, L, \partial) \ar[rr]^i \ar[rd]_p && (L\, \widehat{\amalg}\, L\, \widehat{\amalg}\, \overline{\mathbb L}_S, \partial)\ar[ld]_q\\
 & (L, \partial)
 }$$
 with the following properties:
 \begin{enumerate}
 \item[$\bullet$] $i$ is a profree extension
 \item[$\bullet$] $p$ is the identity in $L\, \widehat{\amalg}\, 0$ and in $0\, \widehat{\amalg}\, L$.
 \item[$\bullet$] $q$ is a quasi-isomorphism and $q(S)= 0$.
 \end{enumerate}
 \item[2.] Two morphisms $f_1, f_2 : (L, \partial) \to (L', \partial')$ of enriched dgl's are \emph{homotopic} if $f_1\, \widehat{\amalg}\, f_2 : L\, \widehat{\amalg}\, L\to L'$ extend to a morphism (called the \emph{homotopy}),
 $$L\, \widehat{\amalg}\, L\widehat{\amalg}\, \overline{\mathbb L}_S \to L', $$
 from a cylinder object.
 \end{enumerate}
 
 \begin{Prop}
 \label{p14.7} Any enriched dgl $(L, \partial)$ has a cylinder object.
 \end{Prop}

 \vspace{3mm}\noindent {\sl proof.} We first establish the Proposition when $L= \overline{\mathbb L}_T$ is profree, and in this case we set $S = sT$. For convenience we denote $\overline{\mathbb L}_{T_{\leq n}}$ by $L(n)$ and for $x\in L$ denote by $x_\ell$ and $x_r$ the corresponding elements in $L\, \widehat{\amalg}\, 0$ and $0\, \widehat{\amalg}\, L$. Now suppose by induction we have constructed $\partial$ and $q$ in $L(n)\, \widehat{\amalg}\, L(n)\, \widehat{\amalg}\, \overline{\mathbb L}_{S_{\leq n+1}}$. 
 
 Then for $x\in T_{n+1}$ we have $q(\partial x_\ell - \partial x_r)= 0$. Since $q$ is a surjective quasi-isomorphism it follows, exactly as in the proof of Proposition \ref{p13.1}(ii),  that for some coherent map $\omega : L(n) \to L(n)\, \widehat{\amalg}\, L(n)\, \widehat{\amalg}\, \overline{\mathbb L}_{S_{\leq n +1}}$, we have
 $$\partial \omega (x) = \partial x_\ell -\partial x_r\hspace{5mm}\mbox{and } q\circ \omega= 0.$$
 Extend $q$ and $\partial$ to $L(n+1)\, \widehat{\amalg}\, L(n+1)\, \widehat{\amalg}\, \overline{\mathbb L}_{S_{\leq n+2}}$ by setting
 $$\partial \omega (x)=   x_\ell -  x_r -\omega\hspace{5mm}\mbox{and } q(x)= 0.$$
 
 Now for a general enriched dgl, $(L, \partial)$, Proposition \ref{p13.1}   provides a surjective quasi-isomorphism $f: (\overline{\mathbb L}, \partial) \to (L, \partial)$ from a profree dgl. Now it follows from Proposition \ref{p14.1}   that $f$ transforms the diagram for $\overline{\mathbb L}$ to a diagram
 $$
 \xymatrix{
 L \, \widehat{\amalg}\, L \ar[rr]^i\ar[d]_p && L\, \widehat{\amalg}\, L\, \widehat{\amalg}\, \overline{\mathbb L}_S\ar[lld]_q\\L}$$
 satisfying the conditions of the Proposition.
 
 \hfill$\square$

 \begin{Prop}
 \label{p14.8}
 \begin{enumerate}
 \item[(i)] Homotopy is an equivalence relation, and is independent of the choice of cylinder object.
 \item[(ii)]   If $f_1\sim f_1: (L, \partial)\to (L', \partial)$, and $g_1\sim g_2 : (L',\partial)\to (L'', \partial)$ then $g_1\circ f_1\sim g_2\sim f_2$.  
 \end{enumerate}
 \end{Prop}
 
 \vspace{3mm}\noindent {\sl proof.} Suppose $f : L\to L'$ is a morphism of enriched dgl's and that $q: L\, \widehat{\amalg}\, L\widehat{\amalg}\, \overline{\mathbb L}_S \to L$ and $q': L'\, \widehat{\amalg}\, L'\, \widehat{\amalg}\, \overline{\mathbb L}_{S'}\to L'$ are cylinder objects. Then Proposition \ref{p14.1}(iv) yields the commutative diagram
 $$\xymatrix{
 L\, \widehat{\amalg}\, L \ar[rr]^{f\, \widehat{\amalg}\, f} \ar[d] && L'\, \widehat{\amalg}\, L'\, \widehat{\amalg}\, \overline{\mathbb L}_{S'}\ar[d]^{q'}\\
 L\, \widehat{\amalg}\, L\, \widehat{\amalg}\, \overline{\mathbb L}_S \ar[rr]_{f\circ q}\ar@{-->}[rru]^{\widehat{h}} && L'}$$ which extends $f\, \widehat{\amalg}\, f$ to a morphism of the cylinder objects.
 
 (i) The commutative diagram above, applied to $id_L$ provides a quasi-isomorphism connecting any two cylinder objects for an enriched dgl. Now it follows in the standard way that homotopy is independent of the choice of cylinder object.
 
 On the other hand, the automorphism of $L\, \widehat{\amalg}\, L$ exchanging the two factors extends
 to the cylinder objects, and it follows that homotopy is symmetric. Moreover, denote by $L(1)$, $L(2)$ and $L(3)$ three copies of an enriched dgl, $L$, and let
 $$L(1)\, \widehat{\amalg}\, L(2)\, \widehat{\amalg}\, \overline{\mathbb L}_{S(1)}\stackrel{\simeq}{\longrightarrow} L \hspace{3mm}\mbox{and } L(2)\, \widehat{\amalg}\, L(3)\, \widehat{\amalg}\, \overline{\mathbb L}_{S(2)}\stackrel{\simeq}{\longrightarrow} L$$
 be cylinder objects. Then Proposition \ref{p14.1}(iv) yields the commutative diagram,
 $$
 \xymatrix{
 L(1)\, \widehat{\amalg}\, L(2)\,\widehat{\amalg} \, L(3)\,  \widehat{\amalg}\, \overline{\mathbb L}_{S(1)}\,  \widehat{\amalg}\, \overline{\mathbb L}_{S(2)}\ar[rr]^-\simeq && L\\
L(1)\, \widehat{\amalg}\, L(3) \,  \widehat{\amalg}\, \overline{\mathbb L}_{S}.\ar[u]\ar[rru]}$$
It follows that homotopy is transitive.
 
 (ii) If $g_1\sim g_2$ via $\widehat{g} : L'\, \widehat{\amalg}\, L'\, \widehat{\amalg}\, \overline{\mathbb L}_{S'} \to L''$ then composing $g$ with the morphism $\widehat{h}$ in the diagram above yields a homotopy $g_1\circ f\sim g_2\circ f$. On the other hand composing $g_i$ with the homotopy $L\, \widehat{\amalg}\, L\, \widehat{\amalg}\, \overline{\mathbb L}_S \to L'$ from $f_1$ to $f_2$ yields a homotopy $g_i\circ f_1\sim g_i\circ f_2$. \hfill$\square$
 
\vspace{3mm}\noindent {\bf Definition.} The set of morphisms equivalent via homotopy to a morphism $f : (L, \partial) \to (L', \partial')$ of enriched dgl's is its \emph{homotopy class} and is denoted by $[f]$. 

The set of homotopy classes of these morphisms is denoted by $[L, L']$.

\begin{Prop}
\label{p14.9}
\begin{enumerate}
\item[(i)] If $f_1\sim f_2: (L, \partial) \to (L', \partial')$ then $H(f_1)= H(f_2)$. Moreover, if $L=\overline{\mathbb L}_T$ and $L'= \overline{\mathbb L}_{T'}$ then the quotient maps $\overline{f_1}, \overline{f_2} : (T, \overline{\partial})\to (T', \overline{\partial'}$ are homotopic.
\item[(ii)] The bijection $dgl(L\, \widehat{\amalg}\, L',L'') \stackrel{\cong}{\longrightarrow} dgl(L, L'') \times dgl (L', L'')$ (cf. Proposition \ref{p12.6}) induces a bijection 
$$[L\, \widehat{\amalg}\, L', L''] \stackrel{\cong}{\longrightarrow} [L, L''] \times [L', L''].$$
\item[(iii)] If $(L, \partial)$ is a profree dgl and $h : (L', \partial') \to (L'', \partial'')$ is a quasi-isomorphism of enriched dgl's, then composition with $h$ induces a bijection
$$[L, L'] \stackrel{\cong}{\longrightarrow} [L, L''].$$
\item[(iv)] A morphism between profree dgl's is a quasi-isomorphism if and only if it is a homotopy equivalence. Moreover, a homotopy equivalence between minimal profree dgl's is an isomorphism.
\end{enumerate}
\end{Prop}

\vspace{3mm}\noindent {\sl proof.}  (i) The left and right inclusions of $L$ in a cylinder object $L\am L\am \overline{\mathbb L}_S$ are both left inverses to $q: L\am L\am \overline{\mathbb L}_S \to L$. They therefore induce isomorphisms in homology which coincide.  

Then let $\psi : \overline{\mathbb L}_{T'}\, \widehat{\amalg}\, \overline{\mathbb L}_{T'} \, \widehat{\amalg}\, \overline{\mathbb L}_S \to \overline{\mathbb L}_T$ be a homotopy. Division on both sides by $(\,)^{(2)}$ yields a linear map
  $\overline{\psi} : T'\oplus T' \oplus sT' \to T$, and 
   $$\overline{\partial} \overline{\psi} sx) = \overline{f_1}(\overline{\partial}x) - \overline{f_2}(\overline{\partial}x). $$

(ii) This is immediate from the definition.

(iii) We show first that $[L, L']\to [L, L'']$ is injective. For this, suppose $f_1, f_2 : L\to L'$ and that $h\circ f_1\sim h\circ f_2$. Then use
 Proposition \ref{p14.5} to extend $h$ to a surjective quasi-isomorphism,
 $$\widehat{h} : L'\, \widehat{\amalg}\, \overline{\mathbb L}_T \stackrel{\simeq}{\to} L'',$$
 in which $\partial : T\to \overline{\mathbb L}_T$ and $H(\overline{\mathbb L}_T, \partial)= 0$. Next, let $i : L\, \widehat{\amalg}\, L \to L\, \widehat{\amalg}\, L\, \widehat{\amalg}\, \overline{\mathbb L}_S$ be a cylinder object for $L$. Then combining the hypothesis $h\circ f_1\sim h\circ f_2$ with Proposition \ref{p14.1}(iv) yields the commutative diagram
 $$
 \xymatrix{
 L\, \widehat{\amalg}\, L \ar[d]^i\ar[rr]^{f_1\, \widehat{\amalg}\, f_2} && L'\, \widehat{\amalg}\, \overline{\mathbb L}_T\ar[d]^{\widehat{h}}_\simeq\\
 L\, \widehat{\amalg}\, L\, \widehat{\amalg}\, \overline{\mathbb L}_S \ar[rr]_q\ar[rru]^{\widehat{q}} && L''.}$$
 Finally, division by the ideal generated by $\overline{\mathbb L}_T$ converts $\widehat{q}$ to a morphism $L\, \widehat{\amalg}\, L\, \widehat{\amalg}\, \overline{\mathbb L}_S\to L'$. This shows that $f_1\sim f_2$. 
 
 To show that composition with $h$ induces a surjection, let $g : L\to L''$ be a morphism. As above, Propositions \ref{p14.1}(iv) and \ref{p14.5} yield a commutative diagram
 $$\xymatrix{
 && L'\am \overline{\mathbb L}_T\ar[d]^{\widehat{h}}\\
 L \ar[rru]^{\widehat{g}} \ar[rr]^g&& L''.}$$
 On the other hand, division by the ideal generated by $\overline{\mathbb L}_T$ yields a morphism $p : L'\am \overline{\mathbb L}_T \to L'$ which is a right inverse to the inclusion $j : L'\to L'\am \overline{\mathbb L}_T$.
 
 In particular, $p$ is a surjective quasi-isomorphism, and since $(p\circ j\circ p)\circ \widehat{g}= p\circ \widehat{g}$ it follows from the injectivity established above that
 $$j\circ p\circ \widehat{g}\sim \widehat{g}.$$
 But this gives
 $$g= \widehat{h}\circ \widehat{g} \sim \widehat{h}\circ j\circ p\circ \widehat{g}= h\circ p\circ \widehat{g}.$$
 In particular $g$ is homotopic to a morphism in the image of composition with $h$.
 
 (iv) The first assertion follows immediately from   (i) and (iii) above. The second then follows immediately from Proposition \ref{p13.3}. \hfill$\square$

 \subsection{Sullivan homotopy and semi-quadratic Sullivan algebras}
 
 A \emph{Sullivan cylinder object} for an augmented cdga, $\varepsilon : A\to \mathbb Q$, is a decomposition of $(id, id) : A\to A\times_{\mathbb Q}A$ as
 $$A \to A\otimes \land (S^{\geq 0}, dS) \stackrel{\rho}{\to} A\times_{\mathbb Q}A$$
 in which $d : S\stackrel{\cong}{\longrightarrow} dS$ and $\rho$ is surjective. In particular the morphism $(\varepsilon_0, \varepsilon_1) : A\otimes \land (t,dt) \to A$ defined by $\varepsilon_0(t)= 0$, $\varepsilon_1(t)= 1$, is a Sullivan cylinder object. It is straightforward to check that two Sullivan cylinder objects for $A$ are connected by a quasi-isomorphism, $\varphi$, extending $id_A$ and satisfying $\rho'\circ \varphi = \rho$. 
 
 Two augmentation preserving morphisms $\varphi_0, \varphi_1 : A'\to A$ are \emph{homotopic} $(\varphi_0\sim \varphi_1$) if $(\varphi_0, \varphi_1) : A'\to A\times_{\mathbb Q}A$ lifts to a morphism (called the \emph{homotopy}), 
 $$\varphi : A'\to A\otimes \land (S, dS),$$
 to a Sullivan cylinder object such that $\rho\circ \varphi = (\varphi_0, \varphi_1)$. The remarks above show that the condition is independent of the choice of cylinder object.
 
 Now recall that  a morphism
 $$\varphi : (\land V', d_0+d_1)\to (\land V, d_0+d_1),$$
 of semi-quadratic Sullivan algebras is bihomogeneous if $\varphi : V'\to V$.
 
 \begin{Prop}
 \label{p14.10} Let $\land V$ and $\land V'$ be the semi-quadratic models of  profree dgl's $L$ and $L'$. 
 \begin{enumerate}
 \item[(i)] The correspondence $dgl(L',L) \mapsto cdga (\land V, \land V')$ induces a bijection
 $$[L',L] \stackrel{\cong}{\longrightarrow} [\land V, \land V'].$$
 \item[(ii)] A morphism $f : L'\to L$ is a quasi-isomorphism if and only if the corresponding morphism $\land V\to \land V'$ is a quasi-isomorphism.
 \end{enumerate}
 \end{Prop}
 
 \vspace{3mm}\noindent {\sl proof.} (i). For simplicity, we denote $d_0+d_1$ simply by $d$. Then (i) follows immediately from the following two assertions:
 \begin{enumerate}
 \item[(a)] Any morphism $\varphi : \land V'\to \land V$ is homotopic to a bihomogeneous morphism.
 \item[(b)] If $\varphi_0, \varphi_1 : \land V'\to \land V$ are bihomogeneous morphisms corresponding to enriched dgl morphisms $f_0, f_1 : L'\leftarrow L$ then 
 $$\varphi_0\sim \varphi_1 \Longleftrightarrow f_0\sim f_1.$$
 \end{enumerate}

 \vspace{3mm}\noindent {\sl proof of 
 (a)}. It is sufficient to consider the case that $\land V' = \land W\otimes \land Z$ in which $d: Z\to \land W$ and $\varphi\vert_{\land W}$ is bihomogeneous, and then to construct
 $$\psi : \land W\otimes \land Z \to \land V\otimes \land (t,dt)$$
 so that $\psi\vert_{\land W}= \varphi$, $\varepsilon_0\circ \psi = \varphi$, and $\varepsilon_1\circ \varphi : Z\to V$.
 
 Fix a basis of $Z$ and let $z$ be an element of that basis. Since $\land V'$ is semi-quadratic and $\varphi_{\vert \land W}$ is bihomogeneous, it follows that $\varphi (dz)$ is a cycle in $\land^{\leq 2}V$. In particular the component of $\varphi(dz)$ in $\land^2V$ is a $d_1$-cycle. Since $L$ is profree it follows from Theorem 1   that that component is a $d_1$-boundary : for some $v_1\in V$:
 $$\varphi (dz) -d_1v_1 \in V.$$
 
 This gives $d\varphi (z) - dv_1\in V$. Thus $d\varphi(z)-dv_1$ is a cycle in $V$ that is a boundary in $\land V$. Since $L$ is profree it follows from Proposition \ref{p13.2} that for some $v_2\in V\cap \mbox{ker}\, d_1$, $$d(\varphi z-v_1-v_2) = 0.$$
 But, again by Proposition \ref{p12.2}, any cycle in $\land V$ is homologous to some cycle $v_3\in V$, and so
 $$\varphi z = v_1+v_2+v_3 + d\Phi.$$
 Define $\psi : \land W\otimes \land Z\to \land V\otimes \land (t,dt)$ by
 $$\psi w=   w \hspace{3mm}\mbox{and } \psi (z) = \varphi (z) - d(t\Phi).$$

 \vspace{3mm}\noindent {\sl proof of (b)}. Let $L\, \widehat{\amalg}\, L \am \overline{\mathbb L}_S $ be a cylinder object for $L$. The inclusions $i_0, i_1 : L\to L\am L\am \overline{\mathbb L}_S$ and the quasi-isomorphism $L\am L\am \overline{\mathbb L}_S\to L$ correspond to bihomogeneous morphisms
 $$\xymatrix{ \land V \ar[rr]^\lambda && \land R \ar[rr]^-{\rho= (\rho_0, \rho_1)} && \land V\times_{\mathbb Q}\land V}$$
 of semi-quadratic Sullivan algebras. Here $\land V$ is a retract of $\land R$, $\rho$ is surjective and (Proposition \ref{p13.3}) $\lambda$ is a quasi-isomorphism. It follows that $\land R = \land V \otimes \land (S\oplus dS)$ is a cylinder object for $\land V$.
 
 Now if $f_0\sim f_1 : L\to L'$, the homotopy $L\am L\am \overline{\mathbb L}_S\to L'$ corresponds to a bihomogeneous morphism $\varphi : \land V'\to \land R$. It is immediate from the definitions that $\varphi$ is a homotopy from $\varphi_0$ to $\varphi_1$. On the other hand, suppose $\varphi : \land V'\to \land R$ is a (not necessarily bihomogeneous) homotopy from $\varphi_0$ to $\varphi_1$. Since $\varphi_0$ and $\varphi_1$ are bihomogeneous it is then straightforward to check that a bihomogeneous homotopy, $\widehat{\varphi} : \land V'\to \land R$ is defined by
 $$\widehat{\varphi} : V'\to R \hspace{3mm}\mbox{and } \varphi-\widehat{\varphi} : V'\to \land^{\geq 2}R.$$
 This then corresponds to a homotopy $f_0\sim f_1.$
 
 (ii). This follows immediately from the Corollary to Proposition \ref{p13.3}.
 \hfill$\square$
 
 \vspace{3mm}\noindent {\bf Corollary.} A morphism $\varphi : \land V\to \land V'$ is a quasi-isomorphism if and only if it corresponds to a homotopy equivalence $f: L\to L'$.

 \subsection{The   profree dgl model of an enriched dgl}

 \noindent {\bf Definition.} A (minimal) profree dgl model of an enriched  dgl $(L, \partial)$ is a quasi-isomorphism
 $$f : (\overline{\mathbb L}_T, \partial) \stackrel{\simeq}{\longrightarrow} (L, \partial)$$
 from a (minimal) profree dgl.
 
 \begin{Theorem}
 \label{t3}
 \begin{enumerate}
 \item[(i)] Each   enriched dgl, $(L, \partial)$, has a minimal profree dgl model.
 \item[(ii)] Two   profree dgl models $f(1), f(2) : \overline{\mathbb L}_{T(1)}, \overline{\mathbb L}_{T(2)} \stackrel{\simeq}{\to} L$ determine a unique homotopy class of homotopy equivalences
 $$g : \overline{\mathbb L}_{T(1)} \to \overline{\mathbb L}_{T(2)}$$
 such that $f(2)\circ g\sim f(1)$. If the models are minimal, then $g$ is an isomorphism.
 \item[(iii)] More generally, suppose $f(1) : \overline{\mathbb L}_{T(1)} \to L(1)$ and $f(2) : \overline{\mathbb L}_{T(2)} \to L(2)$ are   profree dgl models. Then any morphism
 $$h : L(1)\to L(2)$$
 of enriched dgl's determines a unique homotopy class of morphisms $\widehat{h} : \overline{\mathbb L}_{T(1)} \to \overline{\mathbb L}_{T(2)}$ such that 
 $$f(2)\circ \widehat{h}\sim h\circ f(1).$$
 \end{enumerate}
 \end{Theorem}
 
 \vspace{3mm}\noindent {\sl proof.} (i) Proposition \ref{p13.1} provides a quasi-isomorphism $$\overline{\mathbb L}_T \to L$$
 from a minimal profree dgl.  
 
 (ii) By Proposition \ref{p14.9}(iii) composition with $f$ yields a bijection $$[\overline{\mathbb L}_{T(1)}, \overline{\mathbb L}_{T(2)}] \stackrel{\cong}{\longrightarrow} [ \overline{\mathbb L}_{T(1)}, L(2)].$$ This proves the existence and uniqueness, up to homotopy, of a homotopy equivalence, $g$.  Proposition \ref{p14.9}(iv) asserts that if the models are minimal then $g$ is an isomorphism.  
 
 (iii) This follows from Proposition \ref{p14.9} in exactly the same way as (ii). \hfill$\square$

  \newpage
  \part{The profree dgl model of a topological space}

As described in \S 9, Sullivan's rational completion $X\to X_{\mathbb Q}$ is constructed from the minimal model of $A_{PL}(X)$. 
Here we provide an alternative, but equivalent, characterization of $X_{\mathbb Q}$ in terms of minimal profree dgl's. Namely, to any rational cdga $(A,d)$ with $H^0(A)= \mathbb Q$ we associate a unique (up to isomorphism) minimal profree dgl $L= (\overline{\mathbb L}_T, \partial)$ together with a dga quasi-isomorphism
$$\land V_L \stackrel{\simeq}{\longrightarrow} (A,d)$$
from its semi-quadratic Sullivan model. Thus $\land V_L$ is a (usually non-minimal) Sullivan model of $(A,d)$. In particular, in the case of $A_{PL}(X)$, $$\langle\land V_L\rangle \simeq X_{\mathbb Q}  
 \hspace{3mm}\mbox{and } H(\land V_L)\cong H(X).$$

\vspace{3mm}\noindent {\bf Remark.} Sullivan algebras $\land V$ are free commutative graded algebras equipped with a differential. Profree Lie algebras are the completions (\S 6) of free graded Lie algebras, and profree dgl's are profree Lie algebras equipped with a differential.

  \section{The   profree dgl model of a cdga}
  
Let (\S 12)    $\land V_L$ denote  the semi-quadratic model of an enriched dgl, $(L, \partial)$. We will usually suppress the $\partial$ from the notation and refer simply to an enriched dgl, $L$.

\vspace{3mm}\noindent {\bf Definition.} A \emph{(profree) enriched dgl model} of a cdga $A$  is a (profree) enriched dgl, $L$,  together with a quasi-isomorphism
$$\land V_L\stackrel{\simeq}{\longrightarrow} A.$$

 If $L'$ is an enriched dgl model of a cdga $A'$ then a morphism $\ell : L\to L'$ is a \emph{dgl representative} of a homomorphism $\alpha : A'\to A$ if the morphism $\lambda (\ell) : \land V_{L'}\to \land V_L$ corresponding to $\ell$ makes the diagram
 $$\xymatrix{
 \land V_{L'} \ar[rr]^{\lambda(\ell)} \ar[d]_\simeq && \land V_L\ar[d]^\simeq\\
 A'\ar[rr]_\alpha && A}$$
 homotopy commutative.
 
 In this case we also say that \emph{$A$ is a cdga model of $L$} and that $\alpha$ is a \emph{cdga representative of $\ell$}.

\vspace{3mm} The main step in establishing the introduction above is then Theorem \ref{t4}, below, which depends crucially on \S 14.

\begin{Theorem}
\label{t4}
\begin{enumerate}
\item[(i)] Any cdga, $A$ with $H^0(A) = \mathbb Q$ has a minimal profree dgl model, $L$.
\item[(ii)] If $L$ and $L'$ are both   profree dgl models of $A$ then there is a unique homotopy class of homotopy equivalences, $f :  (L', \partial) \stackrel{\simeq}{\to} (L, \partial)$ such that the diagram
$$\xymatrix{\land V_L \ar[rr]^\simeq \ar[d]_\simeq && A\\
\land V_{L'} \ar[rru]^\simeq}$$
homotopy commutes. If $L$ and $L'$ are minimal then $f$ is an isomorphism.
\item[(iii)] Suppose $L$ and $L'$ are profree dgl models of cdga's $A$ and $A'$. Then any cdga morphism $\varphi : A\to A'$ determines a unique homotopy class of morphisms $f : (L',\partial) \to (L, \partial)$ such that the corresponding cdga diagram
$$\xymatrix{\land V_L \ar[rr]\ar[d]^\psi && A\ar[d]^\varphi\\
\land V_{L'} \ar[rr] && A'}$$
homotopy commutes. Moreover,   $\varphi$ is a quasi-isomorphism if and only if  $f$ is a homotopy equivalence.
\end{enumerate}
\end{Theorem}

 \vspace{3mm}\noindent {\sl proof.} When (i) is established, (ii) and (iii) follow immediately from \S 14.3, \S 14.4 and Proposition \ref{p12.1}. 
 
 For the proof of (i), first observe that if there is a quasi-isomorphism $\varphi : (\land W, d_0+d_1) \to A$ from  a semi-quadratic Sullivan algebra then there is also a quasi-isomorphism of the form $\land V_L\stackrel{\simeq}{\to} A$, in which $(L, \partial)$ is a minimal profree dgl. In fact (Proposition \ref{p12.1}), $(\land W, d_0+d_1)$ is the semi-quadratic Sullivan model of an enriched dgl, $E$.  Proposition \ref{p14.2} then yields a quasi-isomorphism
 $$L = \overline{\mathbb L}_T \stackrel{\simeq}{\longrightarrow} E$$
 from a minimal profree dgl. By Proposition \ref{p12.1}, this yields a quasi-isomorphism
 $$(\land W, d_0+d_1) \stackrel{\simeq}{\longrightarrow} \land V_L$$
 of Sullivan algebras. Since such morphisms have homotopy inverses we obtain a quasi-isomorphism
 $$\land V_L \stackrel{\simeq}{\longrightarrow} (\land W, d_0+d_1) \stackrel{\simeq}{\longrightarrow} A.$$
 It remains to construct a quasi-isomorphism of the form $(\land W, d_0+d_1) \to (A,d)$. Since $H^0(A)= \mathbb Q$ we may without loss of generality suppose that $A= \mathbb Q \oplus A^{\geq 1}$. 
 
 To start, let $Z\subset A^{\geq 1}$ be a subspace of cycles for which $Z \stackrel{\cong}{\longrightarrow} H^{\geq 1}(A)$. Then set $W_0= Z$ and define a quadratic Sullivan algebra $(\land W, d_1)$ with $W = \oplus_{n\geq 0} W_n$ by requiring that
 $$0\to W_{n+1}\stackrel{d_1}{\longrightarrow} \land^2W_{\leq n}\cap \mbox{ker}\, d_1$$
 be exact. It then follows from the construction and from Lemma \ref{l6.2} that
 \begin{eqnarray}
 \label{i22}
 \mathbb Q \oplus W_0 \stackrel{\cong}{\longrightarrow} H(\land W, d_1).
 \end{eqnarray}
 We denote $\oplus_{i\leq n} W_i$ by $W_{\leq n}$. 
 
Next we construct $d_0$ in $\land W$ and a morphism $\varphi : (\land W, d_0+d_1)\to A$ by induction on $n$, and  so that $d_0 : W_n  \to W_{<n}$ and
 \begin{enumerate}
 \item[(i)] $(d_0+d_1)^2= 0$, and 
 \item[(ii)] $\varphi : (\land W_{\leq n}, d_0+d_1) \to A$ is a dga morphism.
 \end{enumerate}
 
 To begin set $d_0=0$ in $W_0$ and let $\varphi : W_0=Z\hookrightarrow \mbox{ker}\, d_A$ be the inclusion. Now suppose $\varphi$ and $d_0$ are constructed in $W_{<n}$. Then filter $\land^2W_{<n}$ by the subspaces
 $$F_k= \oplus_{i+j=k, i,j<n} W_i \land W_j.$$
 If $w\in W_{\leq n}$ and $d_1w\in F_{n-2}$ then $w\in W_{<n}$ and $\varphi w$ and $d_0w$ are defined. On the other hand, for any $w\in W_{\leq n}$, $d_1w\in F_{2(n-1)}$. It is therefore sufficient to define $\varphi_0w$ and $d_0w$ if $w\in W_n$ and for some $k$,
 \begin{enumerate}
 \item[$\bullet$] $d_1w\in F_k$, and
 \item[$\bullet$] $\varphi_0$ and $d_0$ are already defined in $W_{\leq n} \cap d_1^{-1}(F_{k-1})$.
 \end{enumerate}
 
By hypothesis, $d_0$ is defined in $W_{<n}$. Since $d_1d_0d_1w= -d_0d_1^2w= 0$, $d_0d_1w = d_1w_1$ for some $w_1\in W_{\leq n}$. But since $d_0 : W_p\to W_{p-1}$ it follows that $d_0d_1w\in F_{k-1}$. Thus $d_1w_1\in F_{k-1}$, and so by hypothesis $\varphi_0w_1$ and $d_0w_1$ are defined. In particular
 $$d_1d_0w_1 = -d_0d_1w_1= -d_0^2d_1w = 0.$$
 It follows from (\ref{i22}) that $d_0w_1\in W_0$. On the other hand,
 $$\varphi (d_0w) = -\varphi (d_0+d_1) (d_1w-w_1) = d_A\varphi (w_1-d_1w).$$
 Since $d_0w\in W_0$ it follows that $d_0w = 0$.
 
 This in turn implies that $(d_0+d_1)(d_1w-w_1) = 0$ and so
 $$\varphi (d_1w-w_1) = \varphi z+ d_A\Phi$$
 for some $z\in W_0$. Set
 $$d_0w= -w_1-z \hspace{5mm}\mbox{and } \varphi w = \Phi$$
 to close the induction and complete the construction.
 
 Finally, since $W_0\subset \mbox{ker} (d_0+d_1)$ and $\varphi$ induces an isomorphism $W_0\stackrel{\cong}{\to} H(A)$ it follows that $H(\varphi)$ is surjective. On the other hand, given a homology class in $H(\land W, d_0+d_1)$ choose a representative cycle in $\land^{\leq k}W$ with $k$ minimized. If $k\geq 2$ the component in $\land^kW$ is a $d_1$-cycle and therefore of the form $d_1\Phi$ with $\Phi \in \land^{k-1}W$. Subtracting $(d_0+d_1)\Phi$ then produces a representative cycle in $\land^{\leq k-1}W$. It follows that this cycle is an element $w\in W$ in which $d_0w= 0$ and $d_1w=0$. But this implies that $w\in W_0$ and thus $H(\varphi)[w]\neq 0$ unless $w=0$. 
 \hfill$\square$
 
 \vspace{3mm}\noindent {\bf Corollary.} Suppose a cdga $A$ admits a quasi-isomorphism from a quadratic Sullivan algebra $(\land V, d_1)$ with homotopy Lie algebra $L_V$. Then $(L_V, 0)$ is the minimal dgl model for $A$.
 
 \vspace{3mm}\noindent {\sl proof.} In this case $(\land V, d_1)$ is a semi-quadratic model for $A$, and $(L_V,0)$ is the corresponding dgl. \hfill$\square$

 \section{The profree dgl model of a topological space}
 
 The (profree) enriched dgl model of a topological space, $X$, is a (profree) enriched dgl model of $A_{PL}(X)$. More precisely,
 
 \vspace{3mm}\noindent {\bf Definition.} A \emph{(profree) enriched dgl model} of a topological space $X$ is a (profree) enriched dgl, $(L(X), \partial)$, together with a cdga quasi-isomorphism,
 $$(\land V_{L(X)}, d_0+d_1) \stackrel{\simeq}{\longrightarrow} A_{PL}(X)$$  from its semi-quadratic Sullivan model.
 
 \subsection{Connecting spaces, dgl models, and Sullivan models}
 
 \vspace{3mm} Suppose $f : Y\to X$ is a continuous map between connected spaces. Then denote by 
 $$\varphi_f : (\land V_Y, d_0+d_1) \longleftarrow (\land V_X, d_0+d_1)$$
 a bihomogeneous morphism of semi-quadratic Sullivan algebras induced from $A_{PL}(f) : A_{PL}(Y) \leftarrow A_{PL}(X)$. Then $\varphi_f$ extends to a homotopy commutative diagram
 \begin{eqnarray}
 \label{i23}
 \xymatrix{
 (\land V_Y, d_0+d_1) && (\land V_X, d_0+d_1)\ar[ll]^{\varphi_f}\\
 (\land W_Y, d) \ar[u]^\simeq && (\land W_X,d) \ar[u]^\simeq\ar[ll]^{\psi_f}}
 \end{eqnarray}
 in which $\land W_Y$ and $\land W_X$ are minimal Sullivan models of $Y$ and $X$ and $\psi_f$ is a Sullivan representative of $f$. In particular this produces the homotopy commutative diagram
 \begin{eqnarray}
 \label{i24}
 \xymatrix{
 \langle \land V_Y\rangle \ar[d]^\simeq \ar[rr]^{\langle \varphi_f\rangle} && \langle \land V_X\rangle\ar[d]^\simeq\\
 Y_{\mathbb Q} \ar[rr]_{f_{\mathbb Q}} && X_{\mathbb Q}
 }
 \end{eqnarray}
 
 On the other hand, a dgl representative $\ell_f$ for $f$ is the morphism
 $$\xymatrix{(L(Y), \partial) \ar[rr]^{\ell_f} && (L(X), \partial)}$$
 defined by the commutative diagram
 \begin{eqnarray}
 \label{i25}
 \xymatrix{
 (L(Y), \partial) \ar[d]^= \ar[rr]^{\ell_f} && (L(X), \partial)\ar[d]^=\\
 s(V_Y, d_0)^\vee \ar[rr]^{s\varphi_f^\vee} && s(V_X, d_0)^\vee.
 }
 \end{eqnarray}
 As recalled in (\ref{i20}), passage to homology converts (\ref{i25}) to the commutative diagram
 \begin{eqnarray}
 \label{i26}
 \xymatrix{
 H(L(Y), \partial) \ar[d]^\cong \ar[rr]^{H(\ell_f)} && H(L(X), \partial)\ar[d]^\cong \\
 L_Y \ar[rr]^{L_f} && L_X}
 \end{eqnarray}
 in which $L_f$ is the morphism of homotopy Lie algebras induced by $f$. This identifies
 \begin{eqnarray}
 \label{i27}
 \pi_*(f_{\mathbb Q}) = sL_f.
 \end{eqnarray}
 
 Finally, combining Proposition \ref{p13.1} and \ref{p14.1} converts $\ell_f : L(Y) \to L(X)$ to an equivalent minimal profree extension
 $$\xymatrix{  (\overline{\mathbb L}_Q, \partial) \ar[rr]^{\ell_f} && (\overline{\mathbb L}_Q\, \widehat{\amalg}\, \overline{\mathbb L}_T, \partial)}$$
 in which $\overline{\mathbb L}_Q$ is also a minimal profree dgl. Thus we may assume that $\land V_Y$ and $\land V_X$ are the corresponding semi-quadratic Sullivan algebras. But then Proposition \ref{p13.2} provides the commutative diagram
 $$
 \xymatrix{ (Q,0) \ar[rr]\ar[d]^\simeq && (Q\oplus T, \overline{\partial}) \ar[d]^\simeq \ar[rr] && (T,0)\\
 (\land^{\geq 1} W_Y)^\vee \ar[rr] && (\land^{\geq 1} W_X,d)^\vee.}$$
 In particular, if $H(Y)$ and $H(X)$ are graded vector spaces of finite type this identifies
 $$\to Q_k\to H_n(Q\oplus T,\partial) \to T_k \stackrel{H(\partial)}{\rightarrow} Q_{k-1}\to $$
 with the long exact homology sequence
 $$\to H_k(Y;\mathbb Q) \to H_k(X;\mathbb Q) \to H_k(X,Y;\mathbb Q) \to H_{k-1}(Y;\mathbb Q)\to $$ for the pair $(X,Y)$.

\vspace{3mm}\noindent {\bf Remark.} Suppose $(\overline{\mathbb L}_T, \partial)$ is a minimal profree dgl model of $X$. Then Proposition \ref{p13.2} gives
$$sT \cong H^{\geq 1}(X)^\vee.$$
Thus $H(X)$ is a graded vector space of finite type if and only if $T$ is a graded vector space of finite type, and then
$$sT\cong H_{\geq 1}(X;\mathbb Q).$$

\subsection{The Hurewicz homomorphism}

We continue to use the notation above  and, as in \S 9, we   denote the Hurewicz homomorphism by
$$hur : \pi_*(X) \to H_*(X;\mathbb Z)\,.$$
Proposition \ref{p9.3} provides the commutative diagram
$$
\xymatrix{\pi_*(X) \ar[d]\ar[rr] && \pi_*(X_{\mathbb Q})\ar[d] \ar[r]^= & W^\vee\ar[d]^{H(\xi)^\vee}\\
H_*(X;\mathbb Z) \ar[rr] && H(X)^\vee \ar[r]^= & H(\land W)^\vee,}$$
where $\xi : \land^{\geq 1}W\to W$ is division by $\land^{\geq 2}W$. Moreover the discussion above replaces $W^\vee$ by $H(V_L, d_0)^\vee$ and $H(\land W)^\vee$ by $H(\land V_L, d_0+d_1)^\vee$. Finally, Proposition \ref{p13.2} then converts this diagram to the commutative diagram
$$
\xymatrix{\pi_*(X)\ar[d]^{hur}\ar[rr] && sH(\overline{\mathbb L}_T, \partial)\ar[d]^{sH(\rho)}\\
H_*(X;\mathbb Z) \ar[rr] && sH(T, \overline{\partial}).}$$

In summary, this identifies $sH(\rho)$ as the \emph{Sullivan rationalization} of the Hurewicz homomorphism.

\begin{Prop} 
\label{p16.1}Let $X$ be a simply connected space whose   homotopy Lie algebra $L=L_X$ is profree. Then the Hurewicz map induces an injective map
$$hur : L^{(2)}/ L^2 \to H_*(X_{\mathbb Q}).$$
\end{Prop}

\vspace{3mm}\noindent {\sl proof.} Let $\alpha $ be an element in $L^{(2)}_{q-1}$, and write $\alpha = [g]$, with $g : S^q\to X_{\mathbb Q}$. Suppose $hur (\alpha)= 0$. Then there is a  simply connected finite CW complex $K$ of dimension $<q$ and a map $f: K\to X_{\mathbb Q}$ such that $\alpha \in f_*\pi_q(K)$ (\cite{Milnor}).

Denote by $(L_K, \partial)$ the minimal cdgl model of $K$. Thus $L_K = \widehat{\mathbb L}(W)$ with $W= W^{<q-1}$. Denote by $\varphi : (L_K, \partial) \to (L_X, 0)$ be a cdgl representative of $f$. Since  all the elements in $L_X$ are cycles, $\varphi (W)$ is a vector space of cycles. On the other hand by construction there is a cycle $\beta\in (L_K)_{q-1}$ such that $\varphi (\beta)$ is homologous to $\alpha$, thus equal to $\alpha$. But for degree reasons $\beta$ is a decomposable element in $L_K$ and so $\alpha \in L^2$. This gives the result. \hfill$\square$

 \subsection{Wedges of spheres and disks}

Recall from \S 10.4    that if $X= \vee_\alpha S^{n_\alpha}$ is a wedge of spheres then its minimal 
 Sullivan model   is a quadratic Sullivan algebra, $(\land V, d_1)$ and that $(\land V, d_1)  \stackrel{\simeq}{\to} (\mathbb Q\oplus S,0)$ with $S\cdot S= 0$.

More generally, (Theorem \ref{t1}) for any graded vector space $S= S^{\geq 1}$ the minimal Sullivan model $\land V_S\stackrel{\simeq}{\to} (\mathbb Q\oplus S,0)$ is quadratic. Moreover, its homotopy Lie algebra is a profree Lie algebra $\overline{\mathbb L}_T$ and the identity $s\overline{\mathbb L}_T = V_S^\vee$ restricts to an isomorphism $sT = S^\vee$.  Thus $(\overline{\mathbb L}_T,0)$ is a profree dgl model of $\land V_S$,
and if $\land V_S$ is the model of a wedge of spheres, $X$,  then $sT = H^{\geq 1}(X)^\vee$.

On the other hand, the wedge $\vee_\alpha D^{n_\alpha +1}$ of disks is contractible. Thus the zero Lie algebra is a dgl model of $\vee_\alpha D^{n_\alpha + 1}$. Therefore if $(\overline{\mathbb L}_T,0)$ is a profree dgl model of $\vee_\alpha S^{n_\alpha + 1}$ then
$$(\overline{\mathbb L}_T\, \widehat{\amalg}\, \overline{\mathbb L}_{sT}, \partial), \hspace{1cm}\partial sx= x,$$
is a profree dgl model of $\vee_\alpha D^{n_\alpha +1}$. Furthermore, the inclusion
$$(\overline{\mathbb L}_T,0) \to (\overline{\mathbb L}_T\, \widehat{\amalg}\, \overline{\mathbb L}_{sT}, \partial)$$
is a dgl representative of the inclusion $i : \vee_\alpha S^{n_\alpha}\to \vee_\alpha D^{n_\alpha + 1}.$

 \subsection{The quotient space, $X/Y$}
 
 Again let $\varphi : Y\to X$ be a map of connected spaces, but assume  as well that $\varphi$ is the inclusion of a sub-simplicial set, or of a sub CW complex, or of a deformation retract of an open set. In each of these case the maps
 $$\xymatrix{ Y\ar[rr]^\varphi && X \ar[rr]^\pi && X/Y}$$
 provide morphisms
 $$A_{PL}(Y, y_0) \longleftarrow A_{PL}(X, y_0) \longleftarrow A_{PL}(X/Y, \overline{y_0}).$$
 Moreover, these morphism then induce a long exact homology sequence.
 
 Now, as above, denote by $\overline{\mathbb L}_{S_Y}$ and $\overline{\mathbb L}_{S_X}$ profree dgl models of $Y$ and $X$. A profree representative of $\varphi$ is then a morphism $f : \overline{\mathbb L}_{S_Y} \to \overline{\mathbb L}_{S_X}$. Moreover (Proposition \ref{p14.2}) $f$ factors as
 $$f: \xymatrix{\overline{\mathbb L}_{S_Y} \ar[rr]^\lambda && \overline{\mathbb L}_{S_Y}\, \widehat{\amalg}\, \overline{\mathbb L}_T \ar[rr]^\gamma_\simeq && \overline{\mathbb L}_X,}$$
 in which $\lambda$ is a profree extension and $\gamma $ is a quasi-isomorphism. In particular, $\lambda$ is a profree representation of $\varphi$.
 
 \begin{Prop}
 \label{p16.2} With the hypotheses above, the surjection
 $$\rho : \overline{\mathbb L}_{S_Y}\, \widehat{\amalg}\, \overline{\mathbb L}_T \longrightarrow \overline{\mathbb L}_T$$
 is a profree representative of the map $\pi : X\to X/Y$. In particular, $\overline{\mathbb L}_T$ is a profree dgl model of $X/Y$.
 \end{Prop}
 
 \vspace{3mm}\noindent {\sl proof.} Denote by 
 $$\psi : \land V_X \to \land V_Y$$
 the bihomogeneous morphism corresponding to $\lambda : \overline{\mathbb L}_{S_Y} \to \overline{\mathbb L}_{S_Y}\,\widehat{\amalg}\, \overline{\mathbb L}_T$. Since $\lambda $ is a profree representative of $\varphi$, we obtain a homotopy commutative diagram
 $$
 \xymatrix{
 \land V_X \ar[d]^\psi \ar[rr]^\simeq && A_{PL}(X)\ar[d]^{A_{PL}(\varphi)}\\
 \land V_Y\ar[rr]^\simeq && A_{PL}(Y).}$$
 Our hypothesis on $\varphi$ implies that $A_{PL}(\varphi)$ is surjective. Therefore (\cite[Proposition 3.7]{RHTII}) we may assume that this diagram is commutative.
 
 On the other hand a long exact homology sequence argument show that $A_{PL}(\pi)$ restricts to a quasi-isomorphism
 $$A_{PL}(\pi) : \mbox{ker}\, A_{PL}(\varphi) \stackrel{\simeq}{\longrightarrow} A_{PL}(X/Y).$$
 Moreover, since $\lambda$ is an inclusion it follows that $\psi$ is surjective. Therefore a second long exact homology sequence argument shows that ker$\, \psi \stackrel{\simeq}{\longrightarrow} \mbox{ker}\, A_{PL}(\varphi)$. 
 This yields the quasi-isomorphism
 \begin{eqnarray}
 \label{i28} \mbox{ker}\, \psi \stackrel{\simeq}{\longrightarrow} A_{PL}(X/Y).
 \end{eqnarray}
 
 Finally, let
 $$\land V_Y \stackrel{\psi}{\longleftarrow} \land V_X \stackrel{\eta}{\longleftarrow} \land V$$
 be the sequence of bihomogeneous morphisms in which $\eta$ corresponds to the surjection $\overline{\mathbb L}_S\, \widehat{\amalg}\, \overline{\mathbb L}_T \to \overline{\mathbb L}_T$. Proposition \ref{p13.2}(iii) then provides the commutative diagram
 $$
 \xymatrix{
 sS_Y \ar[d]^\simeq \ar[rr] && S_Y\oplus T\ar[d]^\simeq \ar[rr] && T\ar[d]^\simeq\\
 (\land^{\geq 1} V_Y)^\vee \ar[rr] && (\land^{\geq 1}V_X)^\vee \ar[rr] && (\land V)^\vee.}$$
 A final long exact homology sequence argument then shows that $$(\mbox{ker}\, \psi)^\vee = (\land^{\geq 1}V_X)^\vee / (\land^{\geq 1}V_Y)^\vee \stackrel{\simeq}{\longrightarrow} (\land V)^\vee.$$
 It follows that $\eta$ is in fact a quasi-isomorphism
 $$\mbox{ker}\, \psi \stackrel{\simeq}{\longleftarrow} \land V.$$
 
 Thus we obtain the commutative diagram,
 $$\xymatrix{\land V_X \ar[rr]^\simeq && A_{PL}(X)\\
 \land V \ar[u]^\eta \ar[rr]^\simeq && A_{PL}(X/Y)\ar[u]^{A_{PL}(\varphi)}.}$$
 This identifies $\overline{\mathbb L}_T$ as the profree dgl model of $X/Y$ and $\rho$ as a profree representative of $\pi$. 
 
 \hfill$\square$

\subsection{Topological and dgl cofibres}

Suppose $f_0, f_1 : Y\to X_0, X_1$ are maps between connected spaces. Identifying $y\times \{0\}$ with $f_0(y)$ and $y\times \{1\}$ with $f_1(y)$ produces the space $X_0\cup_{f_0}\, Y\times I\, \cup_{f_1} X_1$. Then
$$X_0\cup_Y X_1 := (X_0\cup_{f_0}Y\times I \cup_{f_1} X_1)/ (Y\times \{1/2\})$$
is the \emph{cofibre} of $f_0$ and $f_1$. Here we have the commutative diagram
\begin{eqnarray}
\label{i29}
\xymatrix{Y\times \{1/2\} \ar[d] \ar[rr] && X_1\cup_{f_1} (Y\times [1/2, 1])\ar[d]\\
X_0\cup_{f_0}(Y\times [0, 1/2]) \ar[rr] && X_0\cup_Y X_1,}
\end{eqnarray}
in which $X_0 \simeq X_0\cup_{f_0}(Y\times [0, 1/2])$ and $X_1\simeq X_1\cup_{f_1}(Y\times [1/2, 1])$. 

In parallel, suppose $\ell (0), \ell(1) : L\to L(0), L(1)$ are morphisms of enriched dgl's. These factor (Proposition \ref{p14.2}) through minimal profree extensions as
$$L\to L\, \widehat{\amalg}\, \overline{\mathbb L}_{T(0)} \stackrel{\simeq}{\longrightarrow} L(0) \hspace{1cm}\mbox{and } L\to L\, \widehat{\amalg}\, \overline{\mathbb L}_{T(1)} \stackrel{\simeq}{\longrightarrow} L(1).$$
The corresponding \emph{cofibre square} for $\ell (0)$ and $\ell(1)$ is the commutative diagram,
\begin{eqnarray}
\label{i30}
\xymatrix{ L \ar[d] \ar[rr] && L\,\widehat{\amalg}\, \overline{\mathbb L}_{T(1)}\ar[d]\ar[rr]^\simeq && L(1)\ar[d]\\
L\, \widehat{\amalg}\, \overline{\mathbb L}_{T(0)} \ar[rr] &&
L\, \widehat{\amalg}\, \overline{\mathbb L}_{T(1)\oplus T(0)}\ar[rr]^\simeq && L(1)\, \widehat{\amalg}\, \overline{\mathbb L}_{T(0)},
}
\end{eqnarray}
and $L\, \widehat{\amalg}\, \overline{\mathbb L}_{T(0)\oplus T(1)}$ is a \emph{dgl cofibre} for $\ell(0)$ and $\ell(1)$.

\begin{Prop}
\label{p16.3}
With the hypotheses above suppose that $\ell(0)$ and $\ell(1)$ are enriched dgl representatives for $f_0$ and $f_1$.  Then \begin{enumerate}
\item[(i)] The dgl cofibre $L\, \widehat{\amalg}\, \overline{\mathbb L}_{T(0)\oplus T(1)}$ is an enriched dgl model for $X_0\cup_YX_1$.
\item[(ii)] The morphisms in (\ref{i30}) are dgl representatives for the maps in (\ref{i29}).
\end{enumerate} \end{Prop}

 As a first Corollary we obtain the following extension of Proposition \ref{p10.2}:
 
 \vspace{3mm}\noindent {\bf Corollary}. If $\overline{\mathbb L}_V$ and $\overline{\mathbb L}_W$ are profree dgl models for $X$ and $Y$, then a profree dgl model for $X\vee Y$ is given by $\overline{\mathbb L}_V\, \widehat{\amalg}\, \overline{\mathbb L}_W$.

 \vspace{3mm}  The proof of Proposition \ref{p16.3} is deferred to \S 16.7. First, \S 16.6 establishes the cdga model of a general dgl cofibre.

\subsection{The cdga model of a dgl cofibre}

First suppose that two pairs $\ell(0), \ell(1) : L\to L(0), L(1)$ and $\ell'(0), \ell'(1) : L'\to L(0)', L(1)'$ are connected in a commutative diagram by morphisms
$$\beta, \beta(0), \beta(1) : L, L(0), L(1) \to L', L(0)', L(1)'.$$
Proposition \ref{p14.3} then provides diagrams of morphisms
$$\xymatrix{
L\ar[d]^\beta \ar[rr] && L\, \widehat{\amalg}\, \overline{\mathbb L}_{T(i)} \ar[d]^{\gamma (i)}\ar[rr]^\simeq && L(i)\ar[d]^{\beta (i)}\\
L' \ar[rr] && L'\, \widehat{\amalg}\, \overline{\mathbb L}_{P(i)\oplus \partial P(i) }\, \widehat{\amalg} \, \overline{\mathbb L}_{T'(i)}\ar[rr]^\simeq && L'(i),}
$$
and in which $\partial : P(i)\stackrel{\cong}{\longrightarrow} \partial P(i)$.

Division by $P(i)$ and by $\partial P(i)$ then converts $\gamma (i)$ into a morphism
$$\alpha (i) : L\, \widehat{\amalg}\, \overline{\mathbb L}_{T(i)} \to L'\, \widehat{\amalg}\, \overline{\mathbb L}_{T(i)'}$$ extending $\beta$. In particular, the $\alpha (i)$ connect the respective cofibre squares in a commutative diagram.

\begin{Prop}
\label{p16.4}
With the notation above suppose $\beta, \beta(0)$ and $\beta(1)$ are quasi-isomorphisms. Then $\alpha (0)$ and $\alpha(1)$ are quasi-isomorphisms and induce isomorphisms
$$\overline{\mathbb L}_{T(i)} \stackrel{\cong}{\longrightarrow} \overline{\mathbb L}_{T(i)'}.$$
\end{Prop}

\vspace{3mm}\noindent {\sl proof.} It is immediate from the diagrams above that $\gamma (0)$ and $\gamma (1)$ are quasi-isomorphisms. Hence so are $\alpha (0)$ and $\alpha (1)$: Thus it follows from Proposition \ref{p14.4} that $\overline{\mathbb L}_{T(i)}\to \overline{\mathbb L}_{T(i)'}$ is a quasi-isomorphism.
Since these are minimal profree dgl's the Corollary to Proposition \ref{p13.3} shows that this is an isomorphism.
\hfill$\square$

\vspace{3mm}\noindent {\bf Corollary.} If $L\to L'$ is an isomorphism then both $\alpha (0)$ and $\alpha (1)$ are isomorphisms. In particular, the cofibre square of $\ell(0), \ell(1)$ is uniquely determined up to isomorphism.

\vspace{3mm} Next suppose again that (\ref{i30}) is the cofibre diagram for morphisms $\ell(0), \ell(1) : L\to L(0), L(1)$. Then denote the corresponding semi-quadratic models for $L, L\, \widehat{\amalg}\, \overline{\mathbb L}_{T(0)}, L\, \widehat{\amalg}\overline{\mathbb L}_{T(1)}$ and $L\, \widehat{\amalg}\, \overline{\mathbb L}_{T(0)\oplus T(1)}$ respectively by $\land V_L$, $\land V(0), \land V(1)$ and $\land V$.

\begin{Prop}
\label{p16.5} With the hypotheses   above, the morphism $$\varphi: \land V\to \land V(0) \times_{\land V_L} \land V(1)$$ resulting from (\ref{i30}) is a quasi-isomorphism.
\end{Prop}

\vspace{3mm}\noindent {\sl proof.} First observe that the morphisms in (\ref{i30}) are injections and the duals of the morphisms $V\to V(0),V(1)\to V_L$. Thus these linear maps are surjections. It follows that $H(\land V(0))\times H(\land V(1)), H(\land V_L)$, and $H(\land V(0)\times_{\land V_L} \land V(1))$ are connected by a long exact homology sequence.  On the other hand, Proposition \ref{p13.1}, (i) and (iii), provide a quasi-isomorphism
$$\psi : \overline{\mathbb L}_T \stackrel{\simeq}{\longrightarrow} L$$
from a minimal profree dgl, $\overline{\mathbb L}_T$.  The discussion above, together with Proposition \ref{p16.4}, connects the cofibre square for $\ell(0), \ell(1)\circ \psi : \overline{\mathbb L}_T \to L(0), L(1)$ by quasi-isomorphisms to the cofibre square for $\ell(0), \ell(1)$. The corresponding morphisms of semi-quadratic Sullivan algebras are therefore (Proposition \ref{p12.4}) also quasi-isomorphisms. Furthermore, a long exact homology sequence argument then shows that the induced map of fibre products from $\land V(0) \times_{\land V_L}\land V(1)$ is also a quasi-isomorphism. To summarize, \emph{this reduces Proposition \ref{p16.5} to the case $L$ is a minimal profree dgl, $\overline{\mathbb L}_T$.}

Finally, since $L\, \widehat{\amalg}\, \overline{\mathbb L}_{T(0)}$ and $L\, \widehat{\amalg}\, \overline{\mathbb L}_{T(1)}$ are minimal profree extensions, Proposition \ref{p13.2}(iii) identifies the desuspension of the dual of $H(\varphi)$ as the isomorphism
$$\frac{T\oplus T(0) \oplus T\oplus T(1)}{T} \longrightarrow T\oplus T(0)\oplus T(1)$$
in which the inclusion of $T$ in $T\oplus T$ is the map $x\mapsto (x,-x)$. \hfill$\square$

\subsection{Proof of Proposition \ref{p16.3}.} 

Without loss of generality we may replace $X_0, Y, X_1$ respectively by
$$X_0' = X_0\cup_{f_0} (Y\times [0, 1/2]), \hspace{7mm} Y' = Y\times \{1/2\}, \hspace{6mm}\mbox{and } X_1'= X_1\cup_{f_1}(Y\times [1/2, 1]).$$
Then we may choose Sullivan models $\land V_{X_0}, \land V_Y$ and $\land V_{X_1}$ for $X_0', Y'$ and $X_1'$ and surjective morphisms
$$\psi_0, \psi_1 : \land V_{X_0}, \land V_{X_1}\to \land V_Y$$
so that $\psi_0$ and $\psi_1$ correspond to the inclusions $Y'\to X_0', X_1'$. It then follows that $\land V_{X_0}\times_{\land V_Y} \land V_{X_1}$ is quasi-isomorphic to $A_{PL}(X_0\cup_YX_1).$

On the other hand, since $\ell(0), \ell(1) : L\to L(0), L(1)$ are dgl representatives of the inclusions $Y\to X_0, X_1$ they are also dgl representatives for the inclusions $Y'\to X_0', X_1'$. Thus by definition we have homotopy commutative diagrams
$$
\xymatrix{
\land V(0) \ar[d]\ar[rr]^\simeq_{\sigma_0} && \land V_{X_0}\ar[d]^{\psi_0}\\
\land V_L\ar[rr]_\simeq && \land V_Y}
\hspace{15mm}\mbox{and } 
 \xymatrix{
\land V(1) \ar[d]\ar[rr]^\simeq_{\sigma_1} && \land V_{X_1}\ar[d]^{\psi_1}\\
\land V_L\ar[rr]_\simeq && \land V_Y.}
$$
Moreover, since $\psi_0$ and $\psi_1$ are surjective we may choose $\sigma_0$ and $\sigma_1$ to make the diagrams commutative. A long exact homology sequence argument then shows that $\sigma_0$ and $\sigma_1$ induce a cdga quasi-isomorphism
$$\land V(0) \times_{\land V_L} \land V(1) \stackrel{\simeq}{\longrightarrow }\land V_{X_0}\times_{\land V_Y} \land V_{X_1} \simeq A_{PL}(X_0 \cup_Y X_1).$$
Thus Proposition \ref{p16.3} follows from Proposition \ref{p16.5}.

\hfill$\square$

 \section{Cell attachments}

\subsection{Topological  cell attachments}

\vspace{3mm}\noindent {\bf Definition.} A (topological) \emph{attaching map} is a continuous map
$$g : \vee_\alpha S^{n_\alpha}\to X$$
from a wedge of connected spheres to a connected space, $X$. Then the map $$\iota_X : X\to X\cup_g \vee_\alpha D^{n_\alpha + 1}$$
is a (topological) \emph{cell attachment}.

 \vspace{3mm} Given such a cell attachment observe that $X\cup_g \vee_\alpha D^{n_\alpha +1}$ is the cofibre of the pair of maps,
 $$g, i : \vee_\alpha S^{n_\alpha} \to X, \vspace{3mm}\vee_\alpha D^{n_\alpha + 1},$$
and
 $$\xymatrix{
 \vee_\alpha S^{n_\alpha} \ar[d]^i \ar[rr]^g&& X\ar[d]^{i_X}\\
 \vee_\alpha D^{n_\alpha + 1} \ar[rr] && X\cup_g \vee_\alpha D^{n_\alpha + 1}
 }$$
 is the corresponding commutative diagram (\ref{i29}).

Now denote by $(\overline{\mathbb L}_T,0)\to (\overline{\mathbb L}_{T\, \widehat{\amalg}\, sT}, \partial sx= x)$ the dgl representative of  $\vee_\alpha S^{n_\alpha}\to \vee_\alpha D^{n_\alpha +1}$ constructed in \S 16.3, and let $(L, \partial_L)$ be a dgl model of $X$. Then (\S 16) the dgl diagram (\ref{i30}) corresponding to the diagram above is given by 
$$
\xymatrix{
(\overline{\mathbb L}_T, 0) \ar[d]\ar[rr]^{\ell_g} && L\ar[d]^{j_X}\\
\overline{\mathbb L}_T\, \widehat{\amalg}\, \overline{\mathbb L}_{sT} \ar[rr] && L\, \widehat{\amalg}\, \overline{\mathbb L}_{sT}.}$$
Thus Proposition \ref{p16.3} specializes to

 \begin{Prop}
 \label{p17.1} With the hypotheses and notation above, 
 $$(L\, \widehat{\amalg}\, \overline{\mathbb L}_{sT}, \partial), \hspace{1cm} \partial sx= \ell_g(x),$$
 is an enriched dgl model for $X\cup_g \vee_\alpha D^{n_\alpha + 1}$. Furthermore, the inclusions
 $$\ell_g : (\overline{\mathbb L}_T,0) \to L \hspace{5mm}\mbox{and }\hspace{3mm} j_X: L\to L\, \widehat{\amalg}\, \overline{\mathbb L}_{sT}$$
 are respectively dgl representatives of $g$ and of $i_X$. 
 \end{Prop}

 \vspace{3mm}\noindent {\bf Example.} \emph{Attachment of 2-cells to a finite wedge of circles}
 
 \vspace{2mm} Suppose $g : Y= \vee_{i=1}^r S^1_{a_i} \to X= \vee_{j=1}^s S^1_{b_j}$ is an attaching map with corresponding cell attachment
 $$i_X : X\to X\cup_g \vee_i D^2_{a_i}.$$
 Then $g$ determines and is determined by a group homomorphism
 $$\omega : G_a\to G_b$$
 between the corresponding fundamental groups. Both groups are free and freely generated respectively by element $a_1, \dots , a_r$ and $b_1, \dots , b_s$ corresponding to the circles.
 
 In particular, dgl models for $Y$ and $X$ respectively have the form - \S 16.3 - $(L_Y,0)$ and $(L_X,0)$ and we have the commutative diagram,
 $$
 \xymatrix{
 (L_Y,0) \ar[rr]^{\ell_g} \ar[d]^= && (L_X,0\ar[d]^=\\
 sV_Y^\vee \ar[rr]_{s\psi^\vee} && sV_X^\vee,}$$
 where $\psi : (\land V_Y,d_1) \longleftarrow (\land V_X, d_1)$ is the bihomogeneous morphism of the quadratic Sullivan algebra models of $Y$ and $X$.
 
 On the other hand, by Proposition \ref{p9.4}
 $$L_Y = \overline{\mathbb L}_{T(a)} \hspace{5mm}\mbox{and } L_X= \overline{\mathbb L}_{T(b)}$$
 in which $T(a)$ has as a basis the elements $\log_{G_a}(a_i)$ and $T(b)$ has as a basis the elements $\log_{G_b}(b_j)$. Moreover, it is immediate from the definitions (\S 9.1) that the diagram
 $$
 \xymatrix{
 G_a \ar[d]^\omega\ar[rr] && G_{L_Y} \ar[d] \ar[rr]^\log && L_Y\ar[d]^{\ell_g}\\
 G_b \ar[rr] && G_{L_X} \ar[rr]^\log && L_X}$$
 commutes. Therefore
 $$\ell_g (\log_{G_a}(a_i)) = \log_{G_b}(\omega a_i).$$
 Thus it follows from Proposition \ref{p17.1} that a dgl representative of $i_X$ is the morphism
 $$L_X\to (L_X\, \widehat{\amalg}\, \overline{\mathbb L}_{sT(a)}, \partial)\hspace{5mm}\mbox{with } \partial s\log_{G_a} a_i = \log_{G_b}(\omega a_i).$$

 \subsection{dgl cell attachments}
 
 Proposition \ref{p17.1} permits the generalization of topological cell attachments to general dgl cell attachments.

 \vspace{3mm}\noindent {\bf Definition.}  A \emph{dgl attaching map} is a morphism
 $$\gamma : (\overline{\mathbb L}_T,0) \to (L, \partial_L)$$
 of enriched dgl's, and the corresponding \emph{dgl cell attachment} is the inclusion
 $$j: L\to L\, \widehat{\amalg}\, \overline{\mathbb L}_{sT}, \hspace{5mm} \partial sx= \gamma (x).$$

 \vspace{3mm}\noindent {\bf Remark.} If $\gamma : \overline{\mathbb L}_T\to L$ is a dgl attachment then, since $\gamma$ is coherent, both $\gamma (\overline{\mathbb L}_T)$ and $\gamma (T)$ are closed subspaces of $L$.

 \subsection{Sullivan cell attachments}

 Sullivan algebras also admit a definition of cell attachments and this is equivalent to the definition of dgl attachments (Proposition \ref{p17.3} below). This, however, provides an explicit minimal model representation of a cell attachment (Proposition \ref{p17.4}). 
 
 In fact, as
 recalled in \S 16.3, Theorem 1 provides for each graded vector space $S= S^{\geq 1}$ a quadratic Sullivan model,
 $$(\land V_S, d_1)\to (\mathbb Q \oplus S, 0)$$
 in which $S\cdot S= 0$. Moreover, (\S 16.3)  $(\land V_S, d_1)$ is the semi-quadratic model of the profree dgl,  $(\overline{\mathbb L}_T,0)$, with $sT= S^\vee$. 
 
 \vspace{3mm}\noindent {\bf Definition.} A \emph{Sullivan  attaching map} is a morphism
 $$\rho : (\land W,d)\to (\mathbb Q \oplus S,0)$$
 from a minimal Sullivan algebra, where $S= S^{\geq 1}$ is an arbitrary such graded vector space. If $\rho$ is surjective then it is a \emph{surjective attaching map}.
 
 \vspace{3mm} Now suppose $\gamma : (\overline{\mathbb L}_T, 0) \to (L, \partial_L)$ is a dgl attaching map. Then (Proposition \ref{p12.1}) $\gamma$ induces a morphism
 $$\psi :\land V_L \to \land V_S$$
 between the semi-quadratic Sullivan algebras corresponding to $(L, \partial_L)$ and $(\overline{\mathbb L}_T,0)$. By definition, $\psi$ is a Sullivan representative for $\gamma$. A minimal Sullivan model $\land W\stackrel{\simeq}{\to} \land V_L$ has a surjective homotopy inverse $\land V_L\stackrel{\simeq}{\to} \land W$. This (with (\ref{i20})) then determines the commutative 
 diagrams
 \begin{eqnarray}
 \label{i31}
 \xymatrix{\land V_L\ar@{->>}[rr]^\simeq \ar[d]^\psi && \land W\ar[d]^\rho\\
 \land V_S\ar@{->>}[rr]^\simeq_{\rho_S} && \mathbb Q\oplus S,} \hspace{1cm}\mbox{ and } 
 \xymatrix{ sT \ar[rr]^\cong\ar[d]_{sH(\gamma)} && S^\vee\ar[d]^{\rho^\vee}\\
 sH(L) \ar[rr]^= && W^\vee,}
 \end{eqnarray}
 in which the first horizontal morphisms are surjective quasi-isomorphisms. In particular, there follows
 
 \begin{Prop}
 \label{p17.2} The map $H(\gamma) : T\to H(L)$ is injective if and only if $\rho : \land W\to \mathbb Q\oplus S$ is surjective. \end{Prop}
 
 \vspace{3mm}\noindent {\bf Definition.} $\rho : \land W\to \mathbb Q\oplus S$ is the Sullivan attaching map corresponding to $\gamma$.
 
 \vspace{3mm}\noindent {\bf Remark.} Of course $\rho$ lifts through $\rho_S$ to provide a morphism $\land W\to \land V_S$, and this as well as $\psi$ are Sullivan representatives of $\gamma$. We shall, however, abuse notation and refer to $\rho$ as a Sullivan representative of $\gamma$.   
 
 \begin{Prop}
 \label{p17.3}
 The correspondence $\gamma \mapsto \rho$ described above defines a bijection
 $$[\overline{\mathbb L}_T, L]\stackrel{\cong}{\longrightarrow} [\land W, \mathbb Q\oplus S]$$
 from the homotopy classes of dgl attaching maps to the homotopy classes of Sullivan  attaching maps.
 \end{Prop}
 
 \vspace{3mm}\noindent {\sl proof.} Since $\land V_L \stackrel{\simeq}{\longrightarrow} \land W$ is a surjective quasi-isomorphism, it follows that any morphism $\land W\to (\mathbb Q\oplus S)$ lifts to a morphism $\land V_L\to \land V_S$. Now Proposition \ref{p14.10} shows that this determines a unique homotopy class of dgl morphisms $(\overline{\mathbb L}_T,0)\to (L, \partial_L)$. These by definition induce morphisms $\land W\to \mathbb Q \oplus S$ homotopic to the original one. \hfill$\square$
 
 \vspace{3mm}\noindent{\bf Remark.} The relationships among attaching maps 
 $$g : \vee_\alpha S^{n_\alpha}\to X, \hspace{6mm} \gamma : (\overline{\mathbb L}_T,0)\to (L,\partial_L) \hspace{5mm}\mbox{and }\hspace{3mm} \rho : \land W\to \mathbb Q\oplus S$$
 can be summarized as follows:
 \begin{enumerate}
 \item[(i)] If $\gamma$ is obtained from $g$ then $\gamma = \ell_g$ is a dgl representative of $g$.
 \item[(ii)] If $\rho$ is obtained from $\gamma$ then $\rho$ is a Sullivan representative of $\gamma$.
 \item[(iii)] If $\rho$ is obtained from $\gamma$ which was obtained from $g$ then $\rho$ is a Sullivan representative of the continuous map $g$.
 \end{enumerate}
 
 \vspace{3mm} Again, let $\rho : \land W\to \mathbb Q\oplus S$ be the Sullivan   attaching map corresponding to a dgl attaching map
 $$\gamma : (\overline{\mathbb L}_T,0)\to (L, \partial_L).$$
 Then $\rho$ factors over the surjection $\land W\to \mathbb Q\oplus W$ to define a linear map
 $$\overline{\rho} : W\to S.$$
 This in turn defines two cdga's:
 \begin{enumerate}
 \item[$\bullet$] $(A,d) := \mathbb Q \oplus \mbox{ker}\, \overline{\rho} \oplus \land^{\geq 2} W\subset \land W$, and
 \item[$\bullet$] $(\mathbb Q \oplus s\,\mbox{coker}\, \overline{\rho}, 0)$ with $s\,\mbox{coker}\,\overline{\rho}\cdot s\,\mbox{coker}\, \overline{\rho} = 0$.
 \end{enumerate}

 \begin{Prop}
 \label{p17.4}
Denote by
 $\land V_L \longleftarrow \land \widehat{V}$ 
 the morphism of semi-quadratic Sullivan algebras corresponding to $j : L\to L\, \widehat{\amalg}\, \overline{\mathbb L}_{sT}$. Then   there is a quasi-isomorphism,
 $$\widehat{\varphi} : \land \widehat{V} \stackrel{\simeq}{\longrightarrow} A\times_{\mathbb Q} (\mathbb Q \oplus s\, \mbox{coker}\, \overline{\rho}),$$
 which makes the diagram
 \begin{eqnarray}
 \label{i32}\xymatrix{\land \widehat{V}\ar[d]^\simeq \ar[rrrr] &&&&\land V_L\ar[d]^\simeq\\
 A\times_{\mathbb Q}(\mathbb Q\oplus s\, \mbox{coker}\, \overline{\rho}) \ar[rr] && A\ar[rr] && \land W}\end{eqnarray}
 homotopy commutative.
 \end{Prop}

 \vspace{3mm}\noindent {\sl proof.} 
We shall use without future reference the following basic fact. If $C(0), C(1)\to C$ are morphisms of chain complexes and if the first is surjective then the long exact homology sequence for $C(0\times_C C(1)$ computes $H(C(0)\times_C C(1))$. 
 
Next, let
 $ \land V(0)  $
 denote   the semi-quadratic Sullivan model  of   $(\overline{\mathbb L}_{T\oplus sT}, \partial)$. Thus (Proposition \ref{p12.1}) the dgl morphisms
 $$\overline{\mathbb L}_T \to   \overline{\mathbb L}_{T\oplus sT} \hspace{5mm}\mbox{and }  \hspace{3mm}  \gamma \,\overline{\amalg}\, id: \overline{\mathbb L}_T\, \widehat{\amalg}\, \overline{\mathbb L}_{sT} \to L\, \widehat{\amalg}\, \overline{\mathbb L}_{sT}$$
 induce morphisms
 $$ \land V(0) \to \land V_S \hspace{5mm}\mbox{and} \,\, \land  \widehat{V}\to  \land V(0).$$
 In particular (\S 16.3) the quasi-isomorphism $\land V_S\stackrel{\simeq}{\to} \mathbb Q\oplus S$ extends to a quasi-isomorphism
 $$
 \land V(0) \stackrel{\cong}{\longrightarrow}                                                                                                                    
  \mathbb Q \oplus S\oplus U, \hspace{5mm} d: U\stackrel{\cong}{\to} S.$$

 Now (Proposition \ref{p14.2}) factor $\gamma: \overline{\mathbb L}_T\to L$ as
 $$\overline{\mathbb L}_T \to \overline{\mathbb L}_{T\oplus T(1)} \stackrel{\simeq}{\longrightarrow} L.$$
 This yields the commutative diagram
$$
 \xymatrix{
 \overline{\mathbb L}_T \ar[d] \ar[rr] && \overline{\mathbb L}_{T\oplus T(1)} \ar[d]\ar[rr]^\simeq && L\ar[d]\\
 \overline{\mathbb L}_{T\oplus sT} \ar[rr] && \overline{\mathbb L}_{T\oplus T(1)\oplus sT} \ar[rr]^\simeq && L\, \widehat{\amalg}\, \overline{\mathbb L}_{sT}
 }
$$
 In view of the Remarks above and Proposition \ref{p16.5} it follows that $$\land  \widehat{V}\to \land V(0)\times_{\land V_S}\land V_L$$
 is a quasi-isomorphism.
 
 Now the surjections
 $$\land V(0) \stackrel{\simeq}{\to} \mathbb Q \oplus S\oplus U, \hspace{4mm}\land V_S \stackrel{\simeq}{\to} \mathbb Q  \oplus S, \hspace{3mm}\mbox{and the surjection } \land V_L \stackrel{\simeq}{\to} \land W$$
 yield a quasi-isomorphism
 $$\land  \widehat{V}\stackrel{\simeq}{\longrightarrow} (S\oplus U)\times_S\land W .$$
 Write $S = \mbox{Im}\, \psi \oplus C$ and denote $U_C = d(C)$. Then a straightforward check shows that the inclusion
 $$(U_C,0)\times A \to (S\oplus U) \times_S \land W$$
 is a cdga quasi-isomorphism. Since $U_C \cong s\, \mbox{coker}\, \overline{\rho}$, this completes the proof.
  \hfill$\square$

 \vspace{3mm}\noindent {\bf Corollary.} Suppose $g : \vee_\alpha S^{n_\alpha}\to X$ is a topological attaching map and that
 $$j_X : L(X)\to L(X)\, \widehat{\amalg}\, \overline{\mathbb L}_{sT}$$
 is a dgl representative of $\iota_X : X\to X\cup_g \vee_\alpha D^{n_\alpha + 1}$. Then the diagram  
 $$\xymatrix{
 \land \widehat{V} \ar[d]^\simeq \ar[rr] && \land V_L  \ar[d]^\simeq\\
 A\times_{\mathbb Q}(\mathbb Q \oplus \mbox{coker}\, \overline{\rho}) \ar[rr] && \land W}$$
 identifies $A\times_{\mathbb Q}(\mathbb Q \oplus \mbox{coker}\, \overline{\rho})\to \land W$ as a cdga representative of $i_X$. 
 
 \vspace{3mm}\noindent {\sl proof.} Since $j_X$ is a dgl representative of $i_X$ it follows by definition - see (\ref{i23}) - that
 $$\xymatrix{
 \land \widehat{V} \ar[d]\ar[rrr] &&& \land V_L\ar[d]\\
 A_{PL}(X\cup_g \vee_\alpha D^{n_\alpha +1}) \ar[rrr]^-{A_{PL}(i_X)} &&& A_{PL}(X)}$$
 homotopy commutes. Thus the Corollary follows from Proposition \ref{p16.4} \hfill$\square$

 \subsection{Surjective Sullivan attaching maps}

A Sullivan attaching map $\rho : \land W\to \mathbb Q\oplus S$ is surjective   if and only if coker$\, \overline{\rho}=0$, i.e. if and only if the diagram (\ref{i32}) of Proposition \ref{p17.4} reduces to the homotopy commutative diagram,
 $$\xymatrix{
 \land V \ar[d]^\simeq \ar[drr]^\lambda \ar[rr] && \land V_L\ar[d]^\simeq\\
 A\ar[rr] && \land W,}$$
 in which $\land \widehat{V}$ is replaced by its minimal Sullivan model $\land V$. Since (Proposition \ref{p12.5}) $\land V_L\to \land W$ is surjective we may choose $\land V\to \land V_L$ so that the diagram in fact commutes. 
 
 Thus Proposition \ref{p17.4} reduces to
 
 \begin{Prop}
 \label{p17.5} Suppose $\rho : \land W\to \mathbb Q\oplus S$ is surjective.
 \begin{enumerate}
 \item[(i)] The morphism
 $$\lambda : \land V\to \land W$$
 is a Sullivan representative of $j: L\to L\, \widehat{\amalg}\, \overline{\mathbb L}_{sT}$.
 \item[(ii)] Let $Q = s^{-1}S $ and consider the cdga $(\land W\oplus Q,D)$ in which
 $$Dw= dw + s^{-1}\rho w$$
 and $W\cdot Q= 0$. Then, although $(\land W,d)\to (\land W\oplus Q,D)$ is not always a cdga morphism, the composite
 $$A \to \land W\to \land W \oplus Q$$
 is a cdga quasi-isomorphism $(A,d) \stackrel{\simeq}{\to} (\land W\oplus Q,D)$. In particular the composite
 $$(\land V,d) \to (\land W\oplus Q,D)$$
 is a Sullivan model.
 \end{enumerate}
 \end{Prop}
 
 Now if $\rho$ is surjective if follows that $H^1(\lambda)$ is injective. Therefore $\lambda$ decomposes $\land W$ as a $\Lambda$-extension
 $$\lambda : \land V\to \land V\otimes \land Z\stackrel{\simeq}{\to} \land W$$
 in which $Z= Z^{\geq 1}$ and $\land Z$ is a minimal Sullivan algebra.
 
 This in turn provides the diagram
 $$\xymatrix{\langle\land V\rangle &&\langle\land V\otimes \land Z\rangle\ar[ll] && \langle\land Z
 \rangle\ar[ll]\\
 &&\langle\land W\rangle,\ar[llu]_{\langle\lambda\rangle} \ar[u]^\simeq}$$
 in which the horizontal row is a fibration (\cite[Proposition 17.9]{FHTI}). This identifies $\langle\land Z\rangle$ as the homotopy fibre of $\langle\lambda\rangle$.
 
 On the other hand, since $\land \widehat{V}\to \land V_L$ is by definition a Sullivan representative of $L\to L\, 
 \widehat{\amalg}\, \overline{\mathbb L}_{sT}$ it follows that $\land V\to \land V_L$ is also such a Sullivan representative . In particular, if $L \to L\, \widehat{\amalg}\, \overline{\mathbb L}_{sT}$ is obtained from a topological cell attachment $i_X : X\to \cup_g \vee_\alpha D^{n_\alpha +1}$, it follows that $\lambda$ is a Sullivan representative of $i_X$. This establishes
 
 \begin{Prop}
 \label{p17.6} With the hypotheses and notation above
 $$\langle\lambda \rangle = (i_X)_{\mathbb Q}$$
 and $\langle\land Z\rangle$ is the homotopy fibre of $(i_X)_{\mathbb Q}.$
 \end{Prop}

 \section{The profree dgl model and the Sullivan model of a CW complex}
 
 Fix a CW complex $X = \cup_n X_n$ in which $X_n$ denotes the $n$-skeleton and $X_0= \{ pt\}$. Then we may write
 $$X_{n+1} = X_n \cup_{g_n} \vee_\alpha D_\alpha^{n+1},$$
 where $g_n : \vee_\alpha S^n_\alpha \to X_n$ is the attaching map. By Proposition \ref{p17.1}, if $L_{X_n}$ is a profree dgl model for $X_n$ then $X_{n+1}$ has a profree dgl model of the form
 $$L_{X_n} \to (L_{X_n}\, \widehat{\amalg}\, \overline{\mathbb L}_{sT_{n-1}}, \partial)$$
 in which $L_{g_n} : \overline{\mathbb L}_{T_{n-1}} \to L_{X_n}$ is a profree representative for $g_n$.   Now the obvious induction produces a profree dgl $(\overline{\mathbb L}_T, \partial)$ in which $\overline{\mathbb L}_{T_{<n}}$ is a profree dgl model for $X_n$ and the inclusion $X_n\to X_{n+1}$ is represented by
 $\overline{\mathbb L}_{T_{<n}}\to \overline{\mathbb L}_{T_{\leq n}}$.

 \begin{Theorem}
 \label{t5}
 Let $X$ be a CW complex as above. Then
 \begin{enumerate}
 \item[(i)] $L= \overline{\mathbb L}_T$ is a profree dgl model of $X$.
 \item[(ii)] The isomorphisms $sH(T_{<n}, \overline{\partial}) \stackrel{\cong}{\longrightarrow} H^{\geq 1}(X_n)^\vee$ induce an isomorphism
 $$sH(T, \overline{\partial}) \stackrel{\cong}{\longrightarrow} \left[H^{\geq 1}(X)\right]^\vee.$$
 \end{enumerate}
 \end{Theorem}

 \vspace{3mm}\noindent {\sl proof.} Denote by $L(n)$ and by $\land V(n)$ the dgl $\overline{\mathbb L}_{T_{<n}}$ and its semi-quadratic Sullivan model. The inclusion $\overline{\mathbb L}_{T_{<n}} \to \overline{\mathbb L}_{T_{\leq n}}$ is then the dual (up to suspension) of a surjection $V(n+1)\to V(n)$. This yields the sequence of surjections
 $$\dots \leftarrow \land V(n)\leftarrow \land V(n+1) \leftarrow \dots$$
 of semi-quadratic Sullivan algebras. 
Moreover, for $k\leq n$,
$$V^k(n) \stackrel{\cong}{\longrightarrow} V^k(n+1).$$
It follows that
$$ \varinjlim_n \land V(n) = \varinjlim_n \land V^{\leq n}(n) = \land \left( \varinjlim_n V^{\leq n}(n)\right).$$

Now set $V = \varinjlim_n V^{\leq n}(n)$. Then the obstructions above imply that $\land V$ is a semi-quadratic Sullivan algebra and in fact is the semi-quadratic Sullivan model of $\overline{\mathbb L}_T$. Moreover, by definition $L(n)\to L(n+1)$ is a profree representative of the inclusion $X_n\to X_{n+1}$. Thus
$$\xymatrix{
\land V(n) \ar[rr]^{\simeq} && A_{PL}(X_n)\\
\land V(n+1) \ar[u]\ar[rr]^\simeq && A_{PL}(X_{n+1})\ar[u]}$$
is a homotopy commutative diagram. Now the same degree by degree argument shows that $\land V\stackrel{\simeq}{\to} A_{PL}(X)$. This completes the proof. \hfill$\square$

 \begin{Prop}
 \label{p18.1}
 If $X$ is a connected CW complex all of whose cells are odd dimensional then
 \begin{enumerate}
 \item[(i)] The minimal profree model of $X$ has the form $(\overline{\mathbb L}_T,0)$ with $T= T_{even}$.
 \item[(ii)] The minimal Sullivan model of $X$ is also   the minimal Sullivan model of a wedge of odd spheres.
 \end{enumerate}
 \end{Prop}
 
 \vspace{3mm}\noindent {\sl proof.} (i) follows directly from Theorem \ref{t5}, while (ii) follows directly from Proposition \ref{p10.6}.
 \hfill$\square$

 \subsection{Cell attachments via a wedge of spheres of finite type}
 
 Suppose $Y = \vee S^{n_\alpha}$ is a wedge of   connected spheres with only finitely many spheres in each dimension. Let $(\overline{\mathbb L}_T, 0)$ be the minimal profree dgl model of $Y$ (\S 16.3) and note that its quadratic Sullivan model is then the minimal Sullivan model of $Y$. We show now that the composite
 $$\pi_*(Y) \to \pi_*(Y_{\mathbb Q}) \stackrel{\cong}{\to} s\overline{\mathbb L}_T \to   sT$$
maps a set of generators, $\alpha_i$, of $\pi_*(Y)$ onto a basis of $sT$. (The central isomorphism is (\ref{i25}) in \S 16).
 
 In fact, let $\alpha_1, \dots , \alpha_k \in \pi_n(Y)$ be represented by the $n$-spheres in $Y$. Then $H^n(Y; \mathbb Q)$ has a basis represented by the fundamental cohomology classes of the spheres, and it follows that the map $\pi_*(Y) \to H_*(Y;\mathbb Q)$ maps the $\alpha_i$ to the dual of a basis $a_i$ for $H^n(Y;\mathbb Q)$. By Proposition \ref{p10.6}, the $a_i$ are a basis of $sT$. 
 
 Therefore a morphism from $\overline{\mathbb L}_T$ is uniquely determined by its composite with $\pi_*(Y) \to \pi_*(Y_{\mathbb Q})$. In particular, if $g: Y\to X$ is any map, then from the commutative diagram (cf (\ref{i16})
 $$\xymatrix{
 \pi_*(Y) \ar[rr]\ar[d]_{\pi_*(g)} && \pi_*(Y_{\mathbb Q})\ar[d]^{L_g}\\
 \pi_*(X) \ar[rr] && \pi_*(X_{\mathbb Q})
 }$$
  we obtain that $\pi_*(g)$ uniquely determines $L_g$. Therefore the profree representative of $X \to X\cup_g \vee D^{n_\alpha+1}$ is given by
 $$L_X \to L_X\, \widehat{\amalg}\, \overline{\mathbb L}_{sT}, \hspace{1cm} \partial sx= \pi_*(g)x,$$
 where $\pi_*(g)(x)$ is regarded as an element in $\pi_*(X_{\mathbb Q})$.
 
 In particular, if $X$ is a connected CW complex of finite type then the differential $\partial$ in $\overline{\mathbb L}_T$ in Theorem \ref{t5} can be constructed directly from the morphisms $\pi_*(g_\alpha): \pi_*(S^{n_\alpha}) \to  \pi_*(X^n)$ where the $g_\alpha$ are the attaching maps of $X$.

\newpage
\part{Inert attaching maps}

\section{dgl vs Sullivan  attaching maps}

Recall (\S 17) that a topological attaching map and the induced cell attachment are continuous maps of connected spaces 
$$g : \vee_\alpha S^{n_\alpha} \to X \hspace{5mm}\mbox{and }\hspace{3mm} i_X : X\to X\cup_g \vee_\alpha D^{n_\alpha +1}.$$
While Im$\, \pi_*(g)\subset $ Ker$\,\pi_*(\iota_X)$ this construction will often create additional elements in $\pi_*(X\cup_g \vee_\alpha D^{n_\alpha +1})$ and it is far from understood which conditions will assure that
$$\pi_*(i_X) : \pi_*(X) \to \pi_*(X\cup_g\vee_\alpha D^{n_\alpha +1})$$
is surjective.

By contrast this condition is much better understood after passage to Sullivan's rational completions: consistent with \cite{Anick} and \cite{HL} we introduce the

\vspace{3mm}\noindent {\bf Definition.} The attaching map, $g$, is \emph{rationally inert} if 
$$\pi_*(i_X)_{\mathbb Q} : \pi_*(X_{\mathbb Q}) \to \pi_*((X\cup_g \vee_\alpha D^{n_\alpha + 1})_{\mathbb Q})$$
is surjective.

\vspace{3mm}   This can be rephrased in terms of the dgl representatives (Proposition \ref{p17.1})
$$\ell_g : (\overline{\mathbb L}_T,0)\to (L, \partial_L) \hspace{5mm}\mbox{and } \hspace{3mm} j_X: (L, \partial_L)\to (L\, \widehat{\amalg}\, \overline{\mathbb L}_{sT}, \partial)$$
of $g$ and $i_X$. Since (\ref{i27}) $\pi_*(i_X)_{\mathbb Q}= sH(j_X)$ it follows that $g$ is rationally inert if and only if $H(j_X)$ is surjective. Therefore we can generalize the definition above to the category of enriched dgl's:

 \vspace{3mm}\noindent {\bf Definition.} A dgl attaching map $\gamma :(\overline{\mathbb L}_T, 0)\to (L, \partial_L)$ is \emph{inert} if
 $$H(j) : H(L) \to H(L\, \widehat{\amalg}\, \overline{\mathbb L}_{sT})$$
 is surjective.
 
 \vspace{2mm}\noindent 
 Thus  
 $$g : \vee_\alpha S^{n_\alpha}\to X \mbox{ is rationally inert } \Leftrightarrow \mbox{a dgl representative of $g$ is inert}.$$
 
 \vspace{3mm}On the other hand, (\S 17.3) a Sullivan attaching map is a morphism
 $$\rho : (\land W,d) \to (\mathbb Q \oplus S,0)$$
 from a minimal Sullivan algebra. Proposition \ref{p17.3} establishes a bijection between the homotopy classes of Sullivan attaching maps and the homotopy classes of dgl attaching maps $\gamma : (\overline{\mathbb L}_T,0)\to (L, \partial)$, in which (\ref{i31})
 \begin{enumerate}
 \item[$\bullet$] $\land W$ is the minimal Sullivan model of $(L, \partial)$, thereby identifying  
 $$sH(L) = W^\vee,$$
  \item[$\bullet$] $sT=S^\vee$ and $sH(\gamma) = \rho^\vee$,
 \end{enumerate}
 and
 \begin{enumerate}
 \item[$\bullet$] $sH(\gamma):  sT \to sH(L)$.

 \end{enumerate}
 
 \vspace{3mm}\noindent {\bf Definition.}
 Such a pair $(\gamma, \rho)$ of attaching maps is called \emph{equivalent}.
 
 \vspace{3mm}
 Now, as in \S 17.3, with $\rho$ we associate the cdga
 $$(A,d):= \mathbb Q \oplus \mbox{ker}\, \rho\vert_W \oplus \land^{\geq 2}W.$$
 But here we also introduce the minimal Sullivan model of $A$:
 $$\land V\stackrel{\simeq}{\longrightarrow} A,$$
 which combined with the inclusion $A\to \land W$ defines  morphisms
 $$\lambda : \land V\to \land W \hspace{5mm}\mbox{and }\hspace{3mm} L_\lambda:L_V\leftarrow L_W,$$
 $L_\lambda$ denoting the corresponding morphism of homotopy Lie algebras. Then we introduce the 
 
 \vspace{3mm}\noindent {\bf Definition.} A Sullivan attaching map $\rho : \land W\to \mathbb Q\oplus S$ is \emph{inert} if both (i) $\rho$ is surjective and (ii) $L_\lambda$ is surjective.
 
 \begin{Prop}
 \label{p19.1}
Suppose $\gamma : (\overline{\mathbb L}_T,0) \to (L, \partial)$ and $\rho : \land W\to \mathbb Q\oplus S$ are equivalent   attaching maps. Then
\begin{enumerate}
\item[(i)] $H(\gamma) : T\to H(L)$ is injective $\Leftrightarrow$ $\rho$ is surjective.
\item[(ii)]Then $\mbox{$\gamma$ is inert $\Longleftrightarrow$ $\rho$ is inert}.$ 
\item[(iii)] If $\gamma$ is a dgl representative of $g: \vee S^{n_\alpha} \to X$ then $\gamma$ is inert if and only if $g$ is rationally inert.
\item[(iv)] If $\gamma $ is inert then $H(\gamma) : T\to H(L, \partial)$ is injective.
\end{enumerate}
 \end{Prop}
 
 \vspace{3mm}\noindent {\sl proof.} (i) This is Proposition \ref{p17.2}.  (ii)  We rely on diagram (\ref{i32}) in Proposition \ref{p17.4}, which introduces the cdga $A\times_{\mathbb Q} (\mathbb Q\oplus s\,\mbox{coker}\, \overline{\rho})$. The homotopy Lie algebra of the minimal Sullivan model of this cdga has (Proposition \ref{p12.6}) the form $L_V\, \widehat{\amalg}\, E$, in which $E=0$ if and only if $\rho$ is surjective.
 
 On the other hand, we can now apply (\ref{i20}) to convert (\ref{i32}) to the commutative diagram
 $$
 \xymatrix{
 H(L\, \widehat{\amalg}\, \overline{\mathbb L}_{sT}, \partial) && H(L, \partial) \ar[ll]^{H(j)}\\
 L_V\, \widehat{\amalg}\, E\ar[u]_\cong & L_V\ar[l] & L_W\ar[l]^{L_\lambda}\ar[u]^\cong.}$$
 It follows that $H(j)$ is surjective if and only if $E=0$ and $L_\lambda$ is injective. 
 
 (iii) As observed above this follows because $\pi_*(\iota_X)_{\mathbb Q}$ is surjective if and only if $H(L)\to H(L\, \widehat{\amalg}\, \overline{\mathbb L}_{sT})$ is surjective.
 
 (iv) This follows from (i).
 
 \hfill$\square$
 
 \vspace{3mm}\noindent {\bf Corollary.} If $\rho : \land W\to \mathbb Q\oplus S$ is the Sullivan representative of an attaching map $g : \vee_\alpha S^{n_\alpha}\to X$, then $g$ is rationally inert if and only if $\rho$ is inert.
 
 \vspace{3mm}\noindent {\bf Definition.} If $\gamma : (\overline{\mathbb L}_T,0)\to (L, \partial)$ is an inert dgl attaching map then $H(\gamma)(T)$ is an \emph{inert subspace} of $H(L)$ and its non zero elements are \emph{inert elements}.
 
 \section{Characterisation of inert dgl attaching maps}

 Given a morphism $\gamma : \overline{\mathbb L}_T \to E$ of complete enriched Lie algebras, we denote by $I$ the closed ideal generated by $\gamma (\overline{\mathbb L}_T) $ and by $p: E\to E/I$ the quotient map. Then let $\lambda: C\to E$ be the inclusion of a closed subspace for which $p\circ \lambda : C\stackrel{\cong}{\to} E/I$. This (\S 10.2) defines a commutative diagram
 $$\xymatrix{ \overline{\mathbb L}_T\, \widehat{\amalg}\, C \ar[rr]^{\gamma \, \widehat{\amalg}\, \lambda} \ar[rd]_{0\, \widehat{\amalg}\, p\lambda} && E\ar[ld]^p\\
 & E/I.}$$

 \vspace{3mm}\noindent {\bf Definition.} $\gamma (\overline{\mathbb L}_T)$ \emph{freely generates} $I$ if $\gamma \, \widehat{\amalg}\, \lambda$ is an isomorphism.
 
 \vspace{3mm}Now consider a dgl attaching map $\gamma : (\overline{\mathbb L}_T,0) \to (L, \partial_L)$.  Then $\gamma$ determines the profree extension
 $$j : (L, \partial_L) \to (L\, \widehat{\amalg}\, \overline{\mathbb L}_{sT}, \partial_L+\partial_\gamma), \hspace{1cm} \partial_\gamma (sx) = \gamma x.$$ Denote by $I$ the closed ideal in $H(L)$ generated by $H(\gamma)( T)$.
 
 \begin{Theorem}
 \label{t6} With the hypotheses above, the following conditions are equivalent:
 \begin{enumerate}
 \item[(i)] $\gamma$ is inert.
 \item[(ii)] ker$\, H(j)$ is freely generated by $H(\gamma) (\overline{\mathbb L}_T)$.
 \item[(iii)] $H(\gamma): \overline{\mathbb L}_T \to H(L)$ is inert.
 \end{enumerate}
 In particular, in this case $I= \mbox{ker}\, H(j).$\end{Theorem}
 
 \vspace{3mm}\noindent {\sl proof.}
  \emph{Step One, the quasi-isomorphism
 $$(L, \partial_L)\stackrel{\cong}{\longrightarrow} (\overline{\mathbb L}_T \, \widehat{\amalg}\, L\, \widehat{\amalg}\, \overline{\mathbb L}_{sT}, \partial),$$
 where $\partial\vert_{\overline{\mathbb L}_T} = 0$, $\partial\vert_L= \partial_L$, and $\partial\vert_{sT}(sx) = \gamma x-x$.}
 
 \vspace{2mm} For convenience we write $\partial\vert_{\overline{\mathbb L}_{sT}}= \partial_\gamma + \partial'$, where $\partial_\gamma (sx) = \gamma x$, $\partial'(sx) = -x$, and $\partial_\gamma$ and $\partial'$ vanish on $\overline{\mathbb L}_T$ and $L$. 
 
Let $(\overline{\mathbb L}_T\, \widehat{\amalg}\, \overline{\mathbb L}_{sT})(n)$ denote the subspace of iterated Lie brackets of length $\leq n$ in $T$ and $sT$. Then the increasing filtration by the closed subspaces $L\, \widehat{\amalg}\, (\overline{\mathbb L}_T\, \widehat{\amalg}\, \overline{\mathbb L}_{sT})(n)$ provides a convergent spectral sequence whose first term is $(L\, \widehat{\amalg}\, (\overline{\mathbb L}_T\, \widehat{\amalg}\, \overline{\mathbb L}_{sT}), \partial_L+\partial')$.  
  It follows that
 $$(L, \partial_L )\stackrel{\cong}{\longrightarrow} (L\, \widehat{\amalg}\, \overline{\mathbb L}_T\, \widehat{\amalg}\, \overline{\mathbb L}_{sT}, \partial)$$
 is a quasi-isomorphism. The reverse quasi-isomorphism is then
 $$id_L \, \widehat{\amalg}\, \gamma \, \widehat{\amalg}\, 0 : L\, \widehat{\amalg}\, \overline{\mathbb L}_T\, \widehat{\amalg}\, \overline{\mathbb L}_{sT} \stackrel{\simeq}{\to} L.$$
 Note that this gives the commutative diagram
 $$
 \xymatrix{ (L, \partial_L) \ar[rr]^\simeq\ar[rrd]_j&& (\overline{\mathbb L}_T\, \widehat{\amalg}\, L\, \widehat{\amalg}\, \overline{\mathbb L}_{sT}, \partial)\ar[d]^q\\
 && (L\, \widehat{\amalg}\, \overline{\mathbb L}_{sT}, \partial_L+\partial_\gamma),}\hspace{1cm} q\vert_{\overline{\mathbb L}_T}= 0.$$
 This identifies $H(q)$ with $H(j)$.

 \vspace{2mm}\noindent \emph{Step Two. If $\gamma$ is inert then $I=$ ker$\, H(j)$ and is freely generated by $H(\gamma)(\overline{\mathbb L}_T)$.}
 
 \vspace{2mm} By hypothesis, $H(q)$ is surjective. Since $q$ is also surjective there is a coherent linear map
 $$\sigma : (L\, \widehat{\amalg}\, \overline{\mathbb L}_{sT}, \partial_L+\partial_\gamma) \to (\overline{\mathbb L}_T\, \widehat{\amalg}\, L\, \widehat{\amalg}\, \overline{\mathbb L}_{sT}, \partial)$$ such that $q\circ \sigma= id$. This in turn defines a chain map
 $$id_{\overline{\mathbb L}_T}\, \widehat{\amalg}\, \sigma : (\overline{\mathbb L}_T\, \widehat{\amalg}\, L\, \widehat{\amalg}\, \overline{\mathbb L}_{sT}, \partial_L+\partial_\gamma) \to (\overline{\mathbb L}_T\, \widehat{\amalg}\, L\,\widehat{\amalg}\, \overline{\mathbb L}_{sT}, \partial).$$
 Assign to iterated Lie brackets of $n$ elements $x_i\in T$ with elements in $L\, \widehat{\amalg}\, \overline{\mathbb L}_{sT}$ the filtration degree $n$. This defines an increasing filtration of both dgl's, and the chain map of the associated graded chain complexes, induced by $id_{\overline{\mathbb L}_{sT}} \, \widehat{\amalg}\, \sigma$   is just the identity. It follows that $id_{\overline{\mathbb L}_T}\, \widehat{\amalg}\, \sigma$ is an isomorphism of chain complexes. In particular it induces an isomorphism and a commutative diagram
 $$
 \xymatrix{\overline{\mathbb L}_T\, \widehat{\amalg}\, H(L\, \widehat{\amalg}\,\overline{\mathbb L}_{sT}, \partial_L+\partial_\gamma) \ar[rd]_{0\, \widehat{\amalg}\, id}\ar[rr]^\cong && H(L)\ar[ld]^{H(j)}\\
 & H(L\, \widehat{\amalg}\,\overline{\mathbb L}_{sT}, \partial_L+\partial_\gamma).}$$
 This identifies ker$\, H(j)$ as the ideal $I$ generated by $H(\gamma)(\overline{\mathbb L}_T)$ and shows that the ideal is freely generated.
 
 \vspace{2mm}\noindent \emph{Step Three. If ker$\, H(j)$ is freely generated by $H(\gamma)(\overline{\mathbb L}_T)$ then $\gamma$ and $H(\gamma$) are inert.}
 
 \vspace{2mm} Recall that $j$ is the profree extension
 $$j : (L,\partial_L)\to (L\,\widehat{\amalg}\, \overline{\mathbb L}_{sT}, \partial_L+\partial_\gamma).$$
 Thus $j$ induces a morphism of   spectral sequences, where $L$ has the trivial spectral sequence and $L\, \widehat{\amalg}\, \overline{\mathbb L}_{sT}$ has the spectral sequence induced from the secondary gradation. At the $E_1$ level this gives the morphism
 $$H(L) \to (H(L)\, \widehat{\amalg}\, \overline{\mathbb L}_{sT},H(\partial_\gamma)).$$

 Since $H(L)$ is freely generated by $H(\gamma)(\overline{\mathbb L}_T)$ we may use $H(\gamma)$ to write
 $$H(L) = \overline{\mathbb L}_T\, \widehat{\amalg}\, C$$
 where $C= H(L)/I$ and $I$ is the ideal generated by $H(\gamma)(\overline{\mathbb L}_T)$. Thus
 $$H(L)\, \widehat{\amalg}\, \overline{\mathbb L}_{sT} = \overline{\mathbb L}_T\, \widehat{\amalg}\, \overline{\mathbb L}_{sT}\, \widehat{\amalg}\, C, \hspace{1cm}\mbox{ and } H(\partial_\gamma)sx= x.$$
 It follows that
 $$H(H(L)\, \widehat{\amalg}\, \overline{\mathbb L}_{sT},H(\partial_\gamma))= C$$
 and so $H(L)\to H(H(L)\, \widehat{\amalg}\, \overline{\mathbb L}_{sT}, H(\partial_\gamma))$ is surjective. This shows that $H(\gamma)$ is inert. Moreover, if $H(\gamma)$ is inert, this spectral sequence trivially collapses and so
 $$H(j) : H(L)\to H(L\, \widehat{\amalg}\, \overline{\mathbb L}_{sT})$$
 is also surjective, and $\gamma$ is inert. 
 
 \hfill$\square$

 \vspace{3mm}\noindent {\bf Corollary.} Suppose $\gamma : (\overline{\mathbb L}_T,0) \to (L, \partial)$ and $\rho : (\land W,d) \to (\mathbb Q \oplus S,0) $ are equivalent attaching maps. Then
 \begin{enumerate}
 \item[(i)] $H(\gamma): (\overline{\mathbb L}_T,0) \to (H(L),0)$ and $\rho : (\land W, d_1)\to (\mathbb Q\oplus S,0)$ are also equivalent.
 \item[(ii)] If any one of the four attching maps is inert then all are.
 \end{enumerate}
 
 \vspace{3mm}\noindent {\sl proof.}  (i) If $(\land W_L, d_0+d_1)$ is the semi-quadratic Sullivan algebra corresponding to $(L, \partial)$, then \S 17.3 provides the commutative diagram
 $$\xymatrix{
 (\land W_L, d_0+d_1) \ar[d]^\simeq\ar[rr]^\varphi && (\land V_S, d_1)\ar[d]^\simeq\\
 \land W\ar[rr]_\rho && \mathbb Q\oplus S.}$$

 On the other hand (Proposition \ref{p12.4}(i)), $H(\gamma) : (\overline{\mathbb L}_T,0) \to (H(L), 0)$ corresponds to the morphism
 $$(\land (H(W_L, d_0), H(d_1)) \to (\land V_S, d_1).$$
 This is just $(\land W,d_1)\to (\land V_S, d_1)$, which is equivalent to $$\rho (\land W,d_1)\to \mathbb Q\oplus S.$$

(ii) For convenience denote $\rho: (\land W,d_1)\to (\mathbb Q\oplus S,0)$ simply by $\rho (1)$. Then combining Proposition \ref{p19.1} with Theorem \ref{t6} gives
$$\rho(i) \mbox{ inert} \Longleftrightarrow H(\gamma) \mbox{ inert } \Longleftrightarrow \gamma \mbox{ inert } \Longleftrightarrow \rho \mbox{ inert}.$$
 \hfill$\square$

\vspace{3mm}\noindent {\bf Remark.} Let $f: X \to Y$ be a continuous map from a wedge of spheres $X$. By Theorem 6, when $f$ is inert the homotopy of the cofibre $Z$ depends only on the map $\pi_*(X)\to \pi_*(Y)$. This is not the case in general.

\vspace{3mm} As an example denote by $f: S^2 \to \mathbb CP^3$ the injection of the $2$-skeleton in $\mathbb CP^3$. The homotopy cofiber   is $Z=\mathbb S^4 \vee S^6$ with homotopy Lie algebra $L_Z= \mathbb L (x_3, x_5) $. Now let $g : S^2\to S^3\times \mathbb  CP^\infty$ denote also the injection of the $2$-skeleton. The maps $\pi_*(f)\otimes \mathbb Q$ and $\pi_*(g)\otimes \mathbb Q$ coincide. The homotopy cofiber of $g$ is a formal space $Z'$ whose cohomology is
$$\land (x,y,z,t)/(x^3-y^2, zt, tx-yz, ty-x^2b) , deg\, x= 4, deg\, y= 6, deg\, z= 7, deg\, t=9.$$
The homotopy Lie algebras $L_Z$ and $L_{Z'}$ are clearly  different.

\section{Characterization of inert Sullivan attaching maps}

In this section we fix the following notation, recalling (\ref{i31}) that $\overline{\rho}^\vee : S^\vee \to W^\vee$ is the suspension of $H(\gamma) : \overline{\mathbb L}_T\to H(L)$.

\begin{eqnarray}
\label{i33}
\left.
\mbox{
\begin{minipage}{14cm}\noindent 
\begin{enumerate}
\item[$\bullet$]   $\rho : \land W\to \mathbb Q\oplus S$ is a surjective cdga morphism from a minimal Sullivan 
 algebra, where $S= S^{\geq 1}$ is an arbitrary graded vector space and $S\cdot S= 0 = d(S)$.  
\item[$\bullet$]
 $\overline{\rho}: W= \land^{\geq 1}W/\land^{\geq 2}W \to S$ is the surjection induced from $\rho$. 
 \item[$\bullet$]   $L_W$ is the homotopy Lie algebra of $\land W$ and $T\subset L_W$ is the closed subspace defined by $sT= \overline{\rho}^\vee (S^\vee)$. 
\item[$\bullet$]   $(A,d) = \mathbb Q\oplus \mbox{ker}\, \overline{\rho}\oplus \land^{\geq 2} W\subset \land W$.
\item[$\bullet$]  $\varphi : \land V\stackrel{\simeq}{\longrightarrow} A$ is a minimal Sullivan model.
\item[$\bullet$] $\lambda: \land V\to \land W$ is the composite $\land V\stackrel{\simeq}{\to} A\to \land W$, and 
$$L_\lambda : L_V \longleftarrow L_W$$
is the induced morphism of homotopy Lie algebras.
\item[$\bullet$] If $\rho$ is inert then the image of $\overline{\rho}^\vee : S^\vee \to W^\vee$ is an \emph{inert subspace} $sT$ of $W^\vee = \pi_*\langle\land W\rangle$, and $T\subset L_W$ is an \emph{inert subspace} of $L_W$
 \end{enumerate}
 \end{minipage}
 } \right\}
\end{eqnarray}

\vspace{3mm}
  Division by $\mathbb Q$ converts                                                                                                        $\rho$ to a surjective linear map $\widehat{\rho} : \land W\to S$. This yields  the row-exact commutative diagram
\begin{eqnarray}
\label{i34}
\xymatrix{
0 \ar[r] & A \ar[r] & \land W \ar[r]^{\widehat{\rho}} & S \ar[r] & 0\\
& & \land V.\ar[u]_\lambda\ar[ul]_\simeq^\varphi}
\end{eqnarray}
Since $\widehat{\rho}$ is surjective, $H^1(\lambda)$ is injective. Therefore (\cite[Theorem 3.1]{RHTII}) $\lambda$ factors as
\begin{eqnarray}
\label{i35}
\lambda: \xymatrix{\land V\ar[rr]^\eta && \land V\otimes \land Z \ar[rr]^\xi_\simeq && \land W}, \hspace{5mm} \eta (v) = v\otimes 1,
\end{eqnarray}
in which $\eta$ is a minimal $\Lambda$-extension, $\land V\otimes \land Z$ is a (not necessarily minimal) Sullivan algebra, and $\xi$ is a surjective quasi-isomorphism. Now, as usual, let     $(\land V\otimes \land U,d)$ be the acyclic closure of $\land V$. Applying $-\otimes_{\land V}\land V\otimes \land U$ to the short exact sequence in (\ref{i34}) and to the quasi-isomorphism $\xi$ yields the isomorphisms
\begin{eqnarray}
\label{i36}
H(\land Z) \cong H(\land W\otimes_{\land V}(\land V\otimes \land U)) \cong \mathbb Q \oplus (S\otimes \land U).\end{eqnarray}

 \vspace{3mm}\noindent {\bf Remark.} The isomorphism (\ref{i36}) 
  is a generalisation of the following result of \cite{FT}: Let $F$ be the homotopy fibre of a map $X\to X':=X\cup_gD^{n+1}$ with $X$ simply connected and $n\geq 1$. The boundary $S^n = \partial D^{n+1}$ then lifts to a map $S^n\to F$. Denote by $u\in H_n(F;\mathbb Q)$ the image of the fundamental class of $S^n$. The holonomy action of $H_*(\Omega X';\mathbb Q)$ on $H_*(F;\mathbb Q)$ then yields an isomorphism $H_{\geq 1}(F;\mathbb Q) \cong (\mathbb Q u)\otimes H_*(\Omega X';\mathbb Q)$.

 \begin{Theorem}
 \label{t7}
With the hypotheses and notation of (\ref{i33}) and (\ref{i35}), 
the following conditions are equivalent:
 \begin{enumerate}
 \item[(i)] $\rho : \land W\to \mathbb Q\oplus S$ is inert.
 \item[(ii)] $\land V\otimes \land Z$ is a minimal Sullivan algebra and $\xi$ is an isomorphism.
 \item[(iii)] For appropriate choice of generating space, $Z$, there is an isomorphism of graded algebras $\mathbb Q \oplus (Z\cap \mbox{ker}\, \overline{d}) \stackrel{\cong}{\longrightarrow} H(\land Z, \overline{d}),$ i.e. the homotopy Lie algebra, $L_Z$, of $\land Z$ is profree
 \end{enumerate}
 \end{Theorem}

 \vspace{3mm}\noindent {\sl proof.}  
 (i $\Longleftrightarrow$ (ii). By Theorem \ref{t6}, $\rho$ is inert if and only if   $\lambda$ induces an inclusion,
 $$V\to \land^{\geq 1}W/\land^{\geq 2}W.$$
 But this map is an inclusion if and only if $\lambda : \land V\to \land W$ is a minimal $\Lambda$-extension,  i.e., if and only if $\land V\otimes \land Z$ is a minimal Sullivan algebra. And in this case, $\xi$ is a quasi-isomorphism of minimal Sullivan algebras, and hence is an isomorphism.

 (ii) $\Longrightarrow$ (iii) Identify $\land V\otimes\land Z = \land W$ via the isomorphism $\xi$.  Then denote by $(\land V\otimes \land U,d)$ the acyclic closure of $\land V$. As in (\ref{i36}), applying $-\otimes_{\land V}\land V\otimes \land U$ to the quasi-isomorphism $\xi$ yields a quasi-isomorphism   
 $$\tau : (\land W\otimes \land U,d) \stackrel{\simeq}{\to} (\land Z, \overline{d}).$$
We show that $\tau$ has a right inverse
$$\sigma : (\land Z, \overline{d})\to (\land W\otimes \land U,d)$$
in which $\sigma : Z\to \land^{\geq 1}W\otimes \land U$ : 

In fact, write $Z= \cup_k Z_k$ with $Z_{k+1}= Z\cap d^{-1}(\land V\otimes \land Z_k)$. Then assume by induction that $\sigma$ has been constructed in some $Z_k$. This implies that for $z\in Z_{k+1}$
$$\sigma (\overline{d}z)\in \sigma (\land^{\geq 2}Z_k) \subset \land^{\geq 2}W\otimes \land U.$$
In particular, $\sigma (\overline{d}z)$ is a cycle in $A^{\geq 1}\otimes \land U$ which is quasi-isomorphic to $\land^{\geq 1}V\otimes \land U$. Since $H^{\geq 1}(\land V\otimes \land U)=0$,  $\sigma (\overline{d}z) = d\Phi$ for some $\Phi \in \land W\otimes \land U$. 

Next write $\Phi = \Phi_1+\Phi_2$ with $\Phi_1\in \land^{\geq 1}W\otimes \land U$ and $\Phi_2\in \land U$. Since $d : U\to \land^{\geq 1}V\otimes \land U$, and $\land W$ is minimal,   $d\Phi_1\in \land^{\geq 2}W\otimes \land U$, while also $d\Phi = \sigma (\overline{d}z)\in \land^{\geq 2}W\otimes \land U$. Thus $d\Phi_2\in \land^{\geq 2}W\otimes \land U$. But as with any acyclic closure, if $\Phi_2\neq 0$ then $d\Phi_2$ has a non-zero component in $V\otimes \land U$. Thus $\Phi_2= 0$ and $\sigma (\overline{d}z)= d\Phi_1$. Extend $\sigma$ to $Z_{k+1}$ by setting $\sigma (z) = \Phi_1$ for a basis element $z$ of a direct summand of $Z_k$ in $Z_{k+1}$.

Finally, apply $-\otimes_{\land V}\land V\otimes \land U$ to the row exact diagram in (\ref{i34}). This yields the commutative diagram     
\begin{eqnarray}
\label{i37}
 \xymatrix{
  0\ar[rr] && 
\mathbb Q 
 \ar[d]^\alpha_\simeq \ar[rr] 
 && \land Z\ar[rr]\ar[d]^\simeq_\sigma 
 && \land^{\geq 1}Z \ar[rr]\ar[d]_{\overline{\sigma}}^\simeq 
 && 0
\\
  0 \ar[rr] && 
A\otimes \land U 
 \ar[rr] 
 && \land W\otimes \land U\ar[rr]^{\widehat{\rho}\otimes id} 
 && S\otimes \land U \ar[rr] 
 && 0,
 }\end{eqnarray}
 in which $\alpha$ is a quasi-isomorphism because $\land V \stackrel{\simeq}{\to} A$. It follows that $\overline{\sigma}$ is a quasi-isomorphism in which $\overline{\sigma}(\land^{\geq 2}Z)= 0$. In particular, the injection $Z\cap \mbox{ker}\, \overline{d}\to \land^{\geq 1}Z$ is a quasi-isomorphism, and   $\mathbb Q \oplus (Z\cap \mbox{ker}\, \overline{d})\stackrel{\simeq}{\longrightarrow} \land Z$. 
 
 \vspace{2mm} (iii) $\Longrightarrow$ (ii) Here, even if $\land V\otimes \land Z$ is not a minimal Sullivan algebra, we have from (\ref{i34}) a row-exact commutative diagram of the form
\begin{eqnarray}
\label{i38}
\xymatrix{ 0 \ar[rr] && B\ar[d]\ar[rr]&& \land V\otimes \land Z\ar[d]^\gamma_\simeq \ar[rr]^{\widehat{\rho}\circ \xi} && S\ar[d]^=\ar[rr] && 0\\
0 \ar[rr] && A \ar[rr] && \land W \ar[rr]^{\widehat{\rho}} && S \ar[rr] && 0.
}
\end{eqnarray}
Since $\gamma$ is surjective, $B\to A$ is a quasi-isomorphism. Therefore so is the inclusion $\land V\to B$. 
 
As observed in \cite[Proposition 17.9]{RHTII}, spatial realization provides a fibration,
$$\xymatrix{ 0 & \langle \land V\rangle \ar[l] & \langle \land V\otimes \land Z\rangle \ar[l] & \langle \land Z\rangle \ar[l] & 0\ar[l]}$$
which identifies the sequence $0\leftarrow \pi_*\langle \land V\rangle \leftarrow \pi_*\langle \land V\otimes \land Z\rangle \leftarrow \pi_*\langle \land Z\rangle\leftarrow 0$ with the sequence $0\leftarrow V^\vee \leftarrow W^\vee \leftarrow Z^\vee \leftarrow 0$. The second sequence is short exact if and only if $\land V\otimes \land Z$ is a minimal Sullivan algebra, while the first is short exact if and only if the connecting homomorphism $\partial : \pi_*<\land V\rangle \to \pi_{*-1}\langle \land Z\rangle$ is zero.

However, a classical homotopy theory argument shows that if $\partial \alpha \neq 0$ then the Whitehead products $Wh(\partial \alpha,\beta)$, $\beta \in \pi_*\langle \land Z\rangle$, are all zero. If we identify $\pi_*\langle \land Z\rangle= sL_Z$ then by Proposition \ref{p8.2}, $\partial \alpha$ is in the center of the profree Lie algebra $L_Z$. This implies that $L_Z= \overline{\mathbb L}_x$ and dim$\, H^{\geq 1}(\land Z)= 1$.

On the other hand, it follows from (\ref{i38}) that $H^{\geq 1}(\land Z)\cong S\otimes \land U$.   
Thus if dim$\, H^{\geq 1}(\land Z)= 1$ then $U= 0$, $\land V\otimes \land Z= \land Z$  and $\land V\otimes \land Z$ is minimal. \hfill$\square$

 \vspace{3mm}\noindent {\bf Corollary 1.} If $\rho :\land W\to \mathbb Q \oplus S$ is inert and the center of $L_Z$ is non zero then $V= 0$, $W=Z$ and $L_Z$ is a free Lie algebra on a single generator.

   \vspace{3mm}\noindent     {\bf Corollary 2.}
 The following conditions on the homotopy Lie algebra, $L_W$, of a minimal Sullivan algebra, $\land W$, are equivalent.
 \begin{enumerate}
 \item[(i)] $L_W$ is profree.
 \item[(ii)] Each closed direct summand, $T$, of $L_W^{(2)}$ is inert.
 \item[(iii)] Some closed direct summand, $T$, of $L_W^{(2)}$ is inert.
 \end{enumerate}

 \vspace{3mm}\noindent {\sl proof.}  
 (i) $\Rightarrow$ (ii) By Proposition \ref{p8.4} we may assume $(\land W,d)$ is quadratic and that $W\cap \mbox{ker}\, d \stackrel{\cong}{\longrightarrow} H^{\geq 1}(\land W)$. Moreover, by Theorem \ref{t1} and the preceding Remarks, the inclusion of $T$ in $L_W$ is dual to a surjection $\rho : W\to W\cap \mbox{ker}\, d$. In particular, for this Sullivan attaching map $H(A)= \mathbb Q$, $V=0$, and $T$ is inert.

 (ii) $\Rightarrow$ (iii) is immediate.
 
 (iii) $\Rightarrow$ (i). By Theorem \ref{t7} we obtain from the inclusion $T\subset L_W$ and the corresponding Sullivan attaching map $\rho :\land W\to W\cap \mbox{ker}\, d$, the decomposition $$\xymatrix{\land V\otimes \land Z \ar[rr]^\cong_\xi && \land W.}$$
 By Lemma \ref{l4.1}, $L_W$ is the closure of the sub Lie algebra generated by $T$. Since $T\subset L_Z$ it follows that $L_W=L_Z$, $V= 0$ and   $L_W$ is profree. 
 
  \hfill$\square$

 \vspace{3mm}\noindent {\bf Example.} Suppose
 $$A\to A\oplus \mathbb Q y, \hspace{5mm} A^{\geq 1}\cdot y= 0$$
 is a cdga morphism defined as follows:
 \begin{enumerate}
 \item[$\bullet$] $A$ has a quadratic Sullivan model $\varphi : \land V\stackrel{\simeq}{\to} A$.
 \item[$\bullet$] $dy= \varphi x$ where $x\in (\land V)^m$ and $[x]\in H(\land V)$ is non zero.
 \item[$\bullet$] $A^0= \mathbb Q$, $H^m(\land V)= \mathbb Q x$ and $H^{>m}(A)= 0$.
 \end{enumerate}
 
 \begin{Prop}
 \label{p21.1} With the hypotheses above let $\varphi_W : \land W\to A\oplus \mathbb Q y$ be a minimal Sullivan model. Then the composite
 $$\rho : \land W\to A\oplus \mathbb Q y\to \mathbb Q \oplus \mathbb Q y$$
 is inert.
 \end{Prop}
 
 \vspace{3mm}\noindent {\sl proof.}  Since $A\to A\oplus \mathbb Q y$ is an isomorphism in homology in degrees $<m$ and is surjective in homology in degree $m$, the induced Sullivan representative $\land V\to \land W$ satisfies $V^{<m}\to W^{<m}$ is injective.
 
 Now for any minimal Sullivan algebra $(\land Z,d)$ denote by $\rho_Z$ the restriction $Z^\vee \to (Z\cap \mbox{ker}\, d)^\vee$. When the differential $d$ is quadratic, then (Lemma \ref{l4.1}), the sequence 
 $$0\to L_Z^{(2)} \to L_Z \stackrel{\rho_Z}{\longrightarrow} s (Z\cap \mbox{ker}\, d)^\vee\to 0$$
is exact.  By hypothesis, in the situation here $V^k \cap \mbox{ker}\, d= 0$, $k\geq m$. It follows that $L_V$ is generated as an enriched Lie algebra by elements in degrees $<m-1$. 
 
 Since $L_W\to L_V$ is surjective in degrees $<m$, this map is surjective.
 \hfill$\square$

\vspace{3mm}\noindent {\bf Corollary.} Let $X = Y\cup_{S^{m-1}} D^m$ be a connected   CW complex whose minimal Sullivan model $\land V$ satisfies $V= V^1$. If $m\geq 2$ and dim$\, Y<m$, the attaching map $S^{m-1}\to Y$ is rationally inert.

 \section{A third characterization of inert attaching maps} 
 
 Suppose (\S 19) that 
 $$(\overline{\mathbb L}_T, 0) \to (L, \partial_L) \hspace{5mm}\mbox{ and } (\land W,d)\to \mathbb Q\oplus S$$
 are equivalent attaching maps. Then (Corollary to Theorem \ref{t6}),
 $$\overline{\mathbb L}_T\to H(L) \hspace{5mm}\mbox{ and } (\land W, d_1)\to \mathbb Q \oplus S$$
 are also equivalent, and if any of these attaching maps are inert then all four are inert. Therefore in characterizing  inert attaching maps it is sufficient to consider two equivalent attaching maps of the form 
 $$\gamma : \overline{\mathbb L}_T \to (L,0) \hspace{5mm}\mbox{ and } \rho : (\land W,d_1) \to \mathbb Q\oplus S.$$
\emph{Henceforth in this section we fix such an equivalent pair $(\gamma, \rho)$. }

\vspace{3mm}
Then consistent with the notation in \S 19, \S 20, and \S 21,
\begin{eqnarray}\left.
\label{i39}
\begin{array}{ll}
\bullet & \mbox{ $L= L_W$ is the homotopy Lie algebra of $\land W$,}\\
\bullet & \mbox{ $s\gamma = \rho^\vee$, and   }\\
\bullet &\mbox{ $I\subset L$ is the closed ideal generated by $\gamma (T)$, so that $s I \subset W^\vee$.} 
\end{array}
\right\}
\end{eqnarray}

Now recall (\ref{i35}) that a quasi-isomorphism $\varphi : \land V\stackrel{\simeq}{\to} A= \mbox{ker}\, \widehat{\rho}$     extends to a minimal $\Lambda$-extension
$$\lambda: \xymatrix{\land V\ar[rr]^\eta&& \land V\otimes \land Z \ar[rr]^{\xi}_\simeq &&\land W}
$$in which $\xi\vert_{\land V}$ is the composite $\xymatrix{\land V  \ar[r]^\simeq & A \ar@{^{(}->}[r] & \land W}$. Moreover (\S 21) $\land V\otimes \land Z$ is a Sullivan algebra,  $\xi$ is a surjective quasi-isomorphism and (Theorem \ref{t7}) $\xi$ is an isomorphism if and only if $\rho$ is inert.

\begin{Prop}
\label{p22.1} If $(\gamma, \rho)$ is an inert pair, then $\varphi$ can be chosen so that $\land V$ is quadratic and $\lambda : V\to W.$
\end{Prop}

\vspace{3mm}\noindent {\sl proof.} Since $\rho : \land W\to \mathbb Q\oplus S$ preserves wedge degrees it follows that $A=\oplus_k A\cap \land^kW$. Therefore the morphism $\varphi$ can be chosen so that $V= \oplus_k V(k)$ with $\varphi : V(k)\to A\cap \land^kW$. Now if $\rho$ is inert then $\lambda (V)\cap \land^{\geq 2}W= 0$. Thus in this case $\lambda : V\to W$ and $(\land V,d)$ is quadratic. \hfill$\square$

\begin{Theorem} 
\label{t8} With the hypotheses and notation above and in (\ref{i33}),  the following conditions on an equivalent pair $(\gamma, \rho)$ are equivalent:
\begin{enumerate}
\item[(i)] $(\gamma, \rho)$ is an inert pair.
\item[(ii)]   The closed ideal $I$ in $L$ generated by $\gamma (T)$ is profree, and the right adjoint representation of $L/I$ in $I/I^{(2)}$ induces an isomorphism of $\overline{UL/I}$-modules
$$T\,\widehat{\otimes}\, \overline{UL/I} \stackrel{\cong}{\longrightarrow} I/I^{(2)},$$
extending $T \to I\to I/I^{(2)}$.
\end{enumerate}
\end{Theorem}

 \vspace{3mm}\noindent {\bf Corollary.}  If $\gamma : \mathbb Q x \to L$ is an inert attaching map, and if $y \in L$ is not a scalar multiple of $x$ or of $[x,x]$, then $[x, y]\neq 0$. 

\vspace{3mm}\noindent {\sl proof.} Let $I$ be the closed ideal generated by $x$. Then (Proposition \ref{p6.2}) $I$ is profree. If $y\in I$ and if $y\not\in \mathbb Q x + \mathbb Q [x,x]$ then   the sub Lie algebra generated by $x$ and $y$ is free, and $[x,y]\neq 0$.

On the other hand when $y\not\in I$, let $\overline{y} = [y] \in L/I$. Since $I/I^{(2)}$ is isomorphic to $x\otimes \overline{U(L_W/I)}$ via the adjoint representation, $x\cdot\overline{y}   = \overline{[{x}, y]}$ is non zero. 
\hfill$\square$

\vspace{3mm}
The proof of Theorem 8 occupies the next two subsections. In \S 22.1 we establish a Proposition which is core of the proof. Then \S 22.2 applies that Proposition to complete the proof of the Theorem.

\subsection{A comparison of holonomy representations  for inert pairs $(\gamma, \rho)$}

 \emph{In this sub section we fix an arbitrary inert equivalent pair $(\gamma, \rho)$.}  Thus   $\xi : \land V\otimes \land Z\stackrel{\cong}{\to} \land W$ and 
 $$0\leftarrow L_V\leftarrow L_W \leftarrow L_Z \leftarrow 0$$
 is a short exact sequence.
 
 Now the relation $<\lambda V, s\gamma (T)>= <\rho \lambda V, sT>= 0$ implies that $\gamma : T\to L_Z$. Since $L_Z$ is a closed ideal this implies that 
 $$I\subset L_Z.$$
 
 On the other hand, since $\rho (\land^{\geq 1}V) \subset \rho (A^{\geq 1})= 0$, $\rho$ factors over division by $\land^{\geq 1}V$ to define $\rho_Z : \land Z\to S\oplus \mathbb Q$. Thus restriction of $\rho_Z$ to $Z\cap \mbox{ker}\, d_Z$ and composition of $\gamma$ with the surjection $L_Z\to L_Z/L_Z^{(2)}$ yield linear maps
 $$\overline{\rho_Z} : Z\cap \mbox{ker}\, d_Z \to S \hspace{5mm}\mbox{and } \overline{\gamma} : T\to L_Z/L_Z^{(2)}.$$
 Note that by construction, $\overline{\gamma}$ is coherent. Now denote by $T\widehat{\otimes}\overline{UL_V}$ the completion of the enriched vector space $T\otimes \overline{UL_V}$ (\S 2).

\begin{Prop} 
\label{p22.2} Since $(\gamma, \rho)$ is an inert equivalent pair then there is a commutative diagram
$$
\xymatrix{
L_Z/L_Z^{(2)} \ar[rr]^\chi_\cong && T\,\widehat{\otimes}\, \overline{UL_V} \\
& T\ar[lu]^{\overline{\gamma}} \ar[ru]_{id \otimes 1}
}$$
in which $\chi$ is an isomorphism of $\overline{UL_V}$-modules, and
\begin{enumerate}
\item[$\bullet$] The representation in $L_Z/L_Z^{(2)}$ is induced   by the right adjoint representation of $\overline{UL_V}$ in $L_Z$. 
\item[$\bullet$] The representation in $T\widehat{\otimes} \overline{UL_V}$ is multiplication from the right.
\end{enumerate}
\end{Prop}

\vspace{3mm}\noindent {\sl proof.} The proof is in several Steps. Throughout, $\varepsilon_V$ and $\varepsilon_U$ denote the augmentations in $\land V$ and $\land U$, with $\varepsilon_U(U)= 0$. 

\vspace{2mm}\noindent \emph{ Step One: The quasi-isomorphism
$$ \varepsilon_V\otimes \rho_Z : (\land V\otimes (Z\cap \mbox{ker}\, d_Z), d) \stackrel{\simeq}{\longrightarrow} (S,0).$$}

\vspace{2mm} Since $\land W$ is a quadratic Sullivan algebra it follows that
$$d: 1\otimes Z\cap \mbox{ker}\, d_Z \to \land^2V\oplus (V\otimes Z\cap \mbox{ker}\, d_Z).$$
Thus division by $\land V\otimes 1 $ defines (\S 5.3) a quadratic $\land V$-module $(\land V\otimes (Z\cap \mbox{ker}\, d_Z),d)$. In particular,
$$\varepsilon_V \otimes \rho_Z : \land V\otimes Z\cap \mbox{ker}\, d_Z \to S$$
is a morphism of chain complexes.

On the other hand, the row-exact commutative diagram
\begin{eqnarray}
\label{i40}
\xymatrix{ 0 \ar[r] & (\land V,d) \ar[d]^\simeq\ar[r] & (\land V\otimes \land Z,d) \ar[d]^\cong \ar[r] & (\land V\otimes \land^{\geq 1}Z, \overline{d}) \ar[d]^\simeq \ar[r] & 0\\
0 \ar[r]& A\ar[r] & \land W\ar[r] & S \ar[r] & 0}
\end{eqnarray}
shows that $\rho: \land V\otimes \land^{\geq 1}Z\to S$ is a quasi-isomorphism. But since $\rho$ is inert, $L_Z$ is profree (Theorem \ref{t7}). Therefore (Theorem \ref{t1}) $Z\cap \mbox{ker}\, d_Z \to (\land^{\geq 1}Z, d_Z)$ is a quasi-isomorphism. Because $\land W$ is quadratic, it follows that $\land V \otimes (Z\cap \mbox{ker}\, d_Z)$ is preserved by $d$, and 
$$\xymatrix{\land V\otimes Z\cap \mbox{ker}\, d_Z \ar[rr]^\simeq && \land V\otimes \land^{\geq 1}Z  \ar[rr]^\simeq_\rho && S},$$
which completes the proof of Step One.

\vspace{3mm}\noindent {\emph{Step Two. The quasi-isomorphism
$$\land V\otimes Z\cap \mbox{ker}\, d_Z\stackrel{\simeq}{\longrightarrow} \land V\otimes \land U\otimes S.$$}

\vspace{2mm} Here $\land V\otimes \land U$ is the acyclic closure of $\land V$, and (\S 5.4) is a quadratic $\land V$-module. Thus  $\land V\otimes \land U\otimes S $ is also   a quadratic $\land V$-module, and
$$\varepsilon_V \otimes \varepsilon_U\otimes id : \land V\otimes \land U\otimes S \stackrel{=}{\longrightarrow} S $$
is a quasi-isomorphism.
Since quadratic $\land V$-modules are semifree the quasi-isomorphism $\varepsilon_V\otimes \rho_Z$ of Step One lifts through   $\varepsilon_V\otimes\varepsilon_U\otimes id$ to yield a commutative diagram of $\land V$-modules
$$\xymatrix{
\land V\otimes Z\cap\mbox{ker}\, d_Z \ar[rr]^\psi_\simeq
\ar[rd]_{\varepsilon_V \otimes \rho_Z}^\simeq &&\land V\otimes \land U\otimes S\ar[ld]_\simeq^{\varepsilon_V\otimes \varepsilon_U\otimes id}\\& S}$$

Now assign to any quadratic $\land V$-module a secondary gradation with $\land^kV\otimes -$ having secondary degree $k$. Then, if $S$ is assigned secondary degree zero, $\varepsilon\otimes \rho_Z$  and $\varepsilon_V\otimes \varepsilon_U\otimes id$ preserve the secondary degree, while the differentials increase it by 1. It follows that $\psi$ can be constructed to be a morphism of $\land V$-modules preserving secondary degree. This implies that $\psi$ is a quasi-isomorphism of quadratic $\land V$-modules.

\vspace{3mm}
\noindent \emph{Step Three. Completion of the proof of Proposition \ref{p22.1}.}

\vspace{2mm}It follows from Lemma \ref{l5.1} that $\psi$ is an isomorphism and that the diagram in Step Two restricts to a commutative diagram,
$$
\xymatrix{Z\cap \mbox{ker}\,d_Z\ar[rr]^\cong\ar[dr]_{\rho_Z} && \land U\otimes S\ar[ld]^{\varepsilon_U\otimes id}\\ & S}$$ in which the horizontal isomorphism is an isomorphism of $\overline{UL_V}$-modules. Dualizing (Theorem \ref{t1} and Proposition \ref{p5.2}) yields the commutative diagram,  
$$
\xymatrix{L_Z/L_Z^{(2)}=  s^{-1}(Z\cap \mbox{ker}\, d_Z)^\vee && T\widehat{\otimes} \overline{UL_V} \ar[ll]^-\cong\\& T\ar[lu]^{\overline{\gamma}} \ar[ru]_{id\otimes 1}}$$
in which   $\overline{UL_V}$ acts by right multiplication in $T\,\widehat{\otimes}\, \overline{UL_V}$.
 Finally, it follows from Proposition \ref{p5.3} that the representation of $\overline{UL_V}$ in $L_Z/L_Z^{(2)}$ is obtained from the right adjoint representation of $L$ in $L_Z$. 
 \hfill$\square$
 
  \subsection{Proof of Theorem \ref{t8}}

 (i) $\Rightarrow$ (ii). Since $(\gamma, \rho)$ is an inert pair, Proposition \ref{p22.2} shows that $\overline{\gamma}(T)$ generates $L_Z/L_Z^{(2)}$. This implies that $\overline{\gamma}(I) = L_Z/L_Z^{(2)}$, since $\overline{\gamma}(I)$ is a closed ideal. It follows that $I+L_Z^{(2)} = L_Z$ and therefore (Lemma \ref{l4.1}(ii)) that 
 $$I = L_Z \hspace{5mm}\mbox{and } I/I^{(2)} = L_Z/L_Z^{(2)}.$$
 Since (Theorem \ref{t7}) $L_Z$ is profree so is $I$, and so (ii) follows from Proposition \ref{p22.2}.
 
 (ii) $\Rightarrow$ (i). Except for the final step, the proof parallels  Proposition \ref{p22.2}. In fact, since $I\subset L$ is a closed ideal, it decomposes $\land W$ as a minimal $\Lambda$-extension 
 $$\land W = \land V_I\otimes \land Z_I$$
 in which $V_I = \{ w\in W\, \vert \, <w, sI>= 0\}.$ In particular,  $I= L_{Z_I}$ and, since 
 $$<\rho V_I, sT> = <V_I, s\gamma (T)> = 0,$$
 it follows that $\land V_I \subset A$. Moreover, by hypothesis, $I$ is profree.

 Now let $\land V_I\otimes \land U_I$ be the acyclic closure of $\land V_I$. 
 Then, exactly as in the proof of Proposition \ref{p22.2}, we obtain from $\rho$ the commutative diagram, 
 $$\xymatrix{Z_I\cap \mbox{ker}\, d_{Z_I} \ar[rr]^\nu\ar[rd]_{\rho_I} && \land U_I\otimes S\ar[ld]^{\varepsilon_{U_I}\otimes id}\\ & S,}$$
 in which $\nu$ is a morphism of $\overline{UL_{V_I}}$-modules. But because we have not assumed that $\land V_I\to A$ is a quasi-isomorphism we may not yet conclude that $\nu$ is an isomorphism.
 
 However, exactly as in the proof in Proposition \ref{p22.2}, the diagram dualizes to
 
 $$
 \xymatrix{ I/I^{(2)}   && T\widehat{\otimes}\overline{UL_{V_I}}\ar[ll]_{s\nu^\vee}\\
 & T \ar[lu]^{\widetilde{\gamma}} \ar[ru]_{id \otimes 1}
 },$$
 in which $\widetilde{\gamma}$ is the composite $T\stackrel{\gamma}{\to} I \to I/I^{(2)}$ and $s\nu^\vee$ is a morphism of the $\overline{UL_{V_I}}$-modules defined in (ii). Thus by hypothesis, $s\nu^\vee$ is an isomorphism. Therefore $\nu$ is an isomorphism, and it follows (because $I$ is profree) that $\land V_I\otimes \land^{\geq 1}Z_I \stackrel{\simeq}{\to} S$. Now it follows from the analogue of (\ref{i40}) in Proposition \ref{p22.2} that $\land V_I\stackrel{\simeq}{\to} A$.
 
 But by definition, $\varphi : \land V\stackrel{\simeq}{\to} A$ is a minimal Sullivan model. Therefore we can choose $\varphi$ so that it factors as
 $$\land V\stackrel{\cong}{\to} \land V_I\to A.$$
 In particular, $\varphi$ extends to an isomorphism $\land V\otimes \land Z_I\stackrel{\cong}{\to} \land W$ and, by Theorem \ref{t7}, this implies that $(\gamma, \rho)$ is inert.
   \hfill$\square$

\section{Properties of inert attaching maps}

\subsection{Pairs of dgl attaching maps}

Suppose (\S 17.2) that $\gamma : (\overline{\mathbb L}_T,0) \to (L, \partial_L)$ is a dgl attaching map, and that $$T = T(1) \oplus T(2)$$
decomposes $T$ as a direct sum of closed subspaces. This yields (\S 19) dgl attaching maps
$$\gamma (1) : (\overline{\mathbb L}_{T(1)},0) \to (L, \partial_L) \hspace{1cm}\mbox{and } \gamma(2) : (\overline{\mathbb L}_{T(2)},0) \to (L\, \widehat{\amalg}\, \overline{\mathbb L}_{sT(1)}, \partial_L+\partial_{\gamma(1)}).$$

\vspace{3mm}
\begin{Prop} 
\label{p23.1} With the hypotheses and notation above,
\begin{quote}
{$\gamma$ is inert $\Longleftrightarrow$ \, $\gamma(1)$ and $\gamma(2)$ are both inert}.\end{quote}
\end{Prop}

\vspace{3mm}\noindent {\sl proof.} Denote by
$$j(1) : L\to L\, \widehat{\amalg}\, \overline{\mathbb L}_{sT(1)} \hspace{5mm}\mbox{ and } j(2) : L\, \widehat{\amalg}\, \overline{\mathbb L}_{sT(1)} \to  L\, \widehat{\amalg}\, \overline{\mathbb L}_{sT(1)} \, \widehat{\amalg}\, \overline{\mathbb L}_{sT(2)} = L\, \widehat{\amalg}\, \overline{\mathbb L}_{sT}$$ the inclusions corresponding to $\gamma (1)$ and $\gamma (2)$. By definition, $j: L\to L\, \widehat{\amalg}\, \overline{\mathbb L}_{s(T)}$ is given by
$$j = j(2)\circ j(1).$$
In particular, if $\gamma(1)$ and $\gamma(2)$ are inert, then $H(j) = H(j(2)) \circ H(j(1))$ is surjective and by definition $\gamma $ inert.

In the reverse direction, suppose $\gamma $ is inert. Since $H(j)$ is surjective, $H(j(2))$ is also surjective and therefore $\gamma(2)$ is inert. To show that $\gamma (1)$ is inert we may assume (Theorem 6) that the differential in $L$ is zero. In this case, since $\gamma$ is inert Theorem \ref{t6}  asserts that ker$\, H(j)$ is freely generated by $\gamma (\overline{\mathbb L}_T)$. In other words, if $I$ is the closed ideal generated by $\gamma (T)$ we have a commutative diagram
$$\xymatrix{
\overline{\mathbb L}_T\, \widehat{\amalg}\, C \ar[rr]_\cong^{\gamma \, \widehat{\amalg}\, \lambda} \ar[rd]_{0 \, \widehat{\amalg}\, p\lambda} && L\ar[ld]^p\\
& L/I.}$$

Now write $\overline{\mathbb L}_T \, \widehat{\amalg}\, C = \overline{\mathbb L}_{T(1)}\, \widehat{\amalg}\, (\overline{\mathbb L}_{T(2)}\, \widehat{\amalg}\, C)$. Then this diagram shows that $I(1)$ is freely generated by $\gamma (T(1))$. Therefore (Theorem \ref{t6}) $\gamma (1)$ is inert. \hfill$\square$

\subsection{Inertness of elements in a profree Lie algebra}

\begin{Theorem}
\label{t10}
Suppose $\overline{\mathbb L}_T$ is a profree Lie algebra in which $T= T_{even}$. Then any non-zero morphism
$$\gamma : \overline{\mathbb L}_{\mathbb Q x}\to \overline{\mathbb L}_T$$
is an inert attaching map.
\end{Theorem}

\vspace{3mm}\noindent {\sl proof.} We accomplish the proof in three steps.

\vspace{2mm}\noindent \emph{Step One. Reduction to the case $T= T_0$}.

 Suppose $(\land W,d)$ is the quadratic model of a complete enriched Lie algebra, $L$. While wedge degree provides a second gradation for the algebra $\land W$, in general $(\land W,d)$ with this gradation is not a dga. However, if $L= L_{even}$ then $W= W^{odd}$ and, with the wedge gradation $(\land W,d)$ is a quadratic Sullivan algebra.
 
 For convenience, in this case we denote by $\widetilde{W}$ the graded vector space defined by
 \begin{enumerate}
 \item[$\bullet$] $\widetilde{W}= W$ as a vector space, but
 \item[$\bullet$] $\widetilde{W}= (\widetilde{W})^1$.
 \end{enumerate}
  The gradation in the cdga $(\land \widetilde{W},d)$ is then the wedge degree. Now denote the homotopy Lie algebra of $(\land \widetilde{W},d)$  by $\widetilde{L}$. By construction, $W^\vee \subset \widetilde{W}^\vee$ and it follows that
 $$sL = W^\vee \subset \widetilde{W}^\vee = s\widetilde{L}.$$
 This identifies $L$ as a closed sub Lie algebra of $\widetilde {L}$ in which the inclusion "forgets" the degrees in $L$.
 
 \begin{Prop}
 \label{p23.2}
 With the hypotheses and notation above
 \begin{enumerate}
 \item[(i)] $L$ is profree if and only if $\widetilde{L}$ is profree.
 \item[(ii)] If $\gamma : T\to L$ is a coherent linear map and $\widetilde{\gamma} : T\to \widetilde{L}$ is the composition of $\gamma $ with the inclusion $L\subset \widetilde{L}$ then
 \begin{quote} $\gamma$ is inert  $\Longleftrightarrow$ $\widetilde{\gamma}$ is inert.\end{quote}
 \end{enumerate}
 \end{Prop}
 
 \vspace{3mm}\noindent {\sl proof.} (i) According to Theorem \ref{t1}, $L$ (respectively $\widetilde{L})$ is profree if and only if $W\cap \mbox{ker}\, d\stackrel{\cong}{\longrightarrow} H(\land^{\geq 1}W)$, 
  respectively $\widetilde{W}\cap \mbox{ker}\, d\stackrel{\cong}{\longrightarrow} H(\land^{\geq 1} \widetilde{W},d)$. 
  But by definition the (non-degree preserving) isomorphism $id : W\stackrel{\cong}{\longrightarrow}\widetilde{W}$ extends to a dga isomorphism $(\land W,d) \stackrel{\cong}{\rightarrow} (\land \widetilde{W},d)$. This implies (i).

 (ii) The Sullivan attaching maps corresponding to $\gamma$ and $\widetilde{\gamma} $ are morphisms
 $$\land W \stackrel{\simeq}{\longrightarrow} \mathbb Q \oplus S \hspace{5mm}\mbox{and } \land \widetilde{W}\to \mathbb Q\oplus S$$
   which coincide if we forget degrees. Thus if either $\gamma $ or $\widetilde{\gamma}$ is surjective, both are, and 
 $$\mbox{ker}\, \rho \cap W = \{w\in W\, \vert\, <w, sT>= 0\}\hspace{5mm}\mbox{ and ker}\, \widetilde{\rho}\cap \widetilde{W} = \{w\in \widetilde{W}\, \vert\, <w, sT>= 0\}.$$
 
 Now let $A= \mathbb Q\oplus \mbox{ker}\, \rho\cap W\oplus \land^{\geq 2}W$ and let $\varphi : \land V \stackrel{\cong}{\to} (A,d)$ b a minimal Sullivan model. Forgetting degrees identifies $\land V$ as a minimal Sullivan model of $\widetilde{A} = \mathbb Q \oplus \mbox{ker}\, \widetilde{\rho} \cap \widetilde{W} \oplus \land^{\geq 2} \widetilde{W}$. Moreover, composing $\varphi$ with the inclusion $A\to \land W$ yields
 $$\lambda : \land V\to \land W$$
 which is identified then with $\lambda : \land V\to \land \widetilde{W}$. In particular, 
 $$\lambda (V) \cap \land^{\geq 2}W= 0 \Longleftrightarrow \lambda (V) \cap \land^{\geq 2} \widetilde{W} = 0,.$$
This establishes (ii) and completes the proof of Step One.
 
 \hfill$\square$

\vspace{3mm}\noindent \emph{Step Two. Reduction to the case $T=T_0$ and dim$\, T<\infty$.}

\vspace{1mm}By definition (\S 6), $T$ is equipped with a family of surjections $\rho_\alpha : T\to T_\alpha$ onto finite dimensional vector spaces and
$$\overline{\mathbb L}_T = \varprojlim_\alpha \overline{\mathbb L}_{T_\alpha}.$$
Therefore, if $x\in \overline{\mathbb L}_T$ is non-zero then some $x_\alpha= \overline{\mathbb L}_{\rho_\alpha}x$ is non-zero. On the other hand, (\S 2), $T= \mbox{ker}\, \rho_\alpha\oplus Q$ with $\rho_\alpha : Q\stackrel{\cong}{\to} T_\alpha$. It is immediate that $\overline{\mathbb L}_{\mbox{ker}\, \rho_\alpha}\to \overline{\mathbb L}_T$ is inert and therefore by Proposition \ref{p23.1}, $x$ is inert if and only if $\rho_\alpha x$ is inert.

\vspace{3mm}\noindent  \emph{Step Three. Proof  when $T= T_0$ and dim$\, T<\infty$.}

\vspace{2mm} By definition, $\gamma$ is inert if and only if
$$\overline{\mathbb L}_T \to H(\overline{\mathbb L}_T\, \widehat{\amalg}\,  {\mathbb L}_{\mathbb Q y})$$
is surjective where the differential in $\overline{\mathbb L}_T \, \widehat{\amalg}\, \mathbb L_{\mathbb Q y}$ is defined by $\partial y = x$. Now fix a basis $x_1, \dots , x_n$ of $T$ and observe that $\overline{\mathbb L}_T \, \widehat{\amalg}\, \mathbb L_{\mathbb Q y}$ is the completion of the free graded Lie algebra,
$$L= \mathbb L (x_1, \dots ,x_n, y).$$
Denote by $L(k)$ the subspace of $L$ of iterated Lie brackets of length $k$ in the $x_i$ and $y$. Thus an element in $(\overline{\mathbb L}_T \, \widehat{\amalg}\, \mathbb L_{\mathbb Q y})_p$ is an infinite sum of the form $z= \sum_{k\geq m} z_k$, with each $z_k \in L(k)_p$. For convenience, given a fixed $m$, we denote the subspace of such elements by $L(\geq m)_p$. In particular, for some $r\geq 1$, 
$$\partial y = \sum_{k\geq r} a_k$$
with $a_k \in L(k)_0$ and $a_r\neq 0$. 

Now consider the dgl $(L, \delta)$ with $\delta y= a_r$. For this dgl we can modify the degrees without changing $L\to H(L,\delta)$, as follows: simply assign deg $2$ to each $x_i$ and degree $2r+1$ to $y$. It follows from (\cite{HL}, Theorem 3.12) that with this new gradation, $L\to H(L, \delta)$ is surjective. This then remains true with the original gradation:
$$L_0\to H(L, \delta)$$
is also surjective and so $H_{\geq 1}(L, \delta)= 0$.

Next, let $\omega$ be a cycle of degree $p\geq 1$ in $\overline{\mathbb L}_T\, \widehat{\amalg}\, \mathbb L_{\mathbb Q y}$ and write
$$\omega= \sum_{k\geq m} \omega_k, \hspace{4mm}\mbox{with } \omega_k\in (L(k)_p,$$
and $\omega_m\neq 0$. Then $\delta \omega_m \in L(m+r-1)_{p-1}$ and
$$\partial \omega -\delta \omega_m \in L(\geq m+r)_{p-1}.$$
It follows that $\delta \omega_m= 0$ and so $\omega_m$ is a $\delta$-boundary:
$$\omega_m = \delta \Phi_{m-r+1}, \hspace{3mm}\mbox{some } \Phi_{m-r+1}\in L(m-r+1)_{p+1}.$$
 But this in turn implies that $$\omega - \partial \Phi_{m-r+1} \in L(\geq m+1).$$
Iterating this process yields a sequence $\Phi_{k-r+1}\in L(k-r+1)_{p+1}$ with
$$\partial (\sum_k\Phi_{k-r+1}) = \omega.$$
\hfill$\square$

\vspace{3mm}\noindent {\bf Corollary.} Suppose $a,b\in \overline{\mathbb L}_T$ and $T = T^{even}$. If $a\in \overline{\mathbb L}_T^{(2)}$, $b\in T$ and $a$ is not in the closed ideal generated by $b$, then $b$ represents an inert element in $\overline{\mathbb L}_T/I$, where $I$ is the closed ideal generated by $a$.

\vspace{3mm}\noindent {\sl proof.} It follows from the hypothesis that $\overline{\mathbb L}_T = \mathbb L_{\mathbb Q b}\, \widehat{\amalg}\, \overline{\mathbb L}_{T'}.$
Thus if $J$ is the closed ideal generated by $b$, then $\overline{\mathbb L}_T/J\cong \overline{\mathbb L}_{T'}$. Hence (Theorem 10) $a$ represents an inert element in $\overline{\mathbb L}_T/J$. Similarly, $b$ is inert in $\overline{\mathbb L}_T$. Thus (Proposition \ref{p22.1}) $\mathbb Q a \oplus \mathbb Q b$ is an inert subspace of $\overline{\mathbb L}_T$. Now (again by Proposition \ref{p22.1}) $b$ represents an inert element in $\overline{\mathbb L}_T/(a).$

\hfill$\square$

\subsection{Free products}

\begin{Prop} \label{p23.3} Suppose $Q\subset L$ and $Q'\subset L'$ are inert subspaces of complete enriched Lie algebras. Then 
$$Q\oplus Q'\subset L\, \widehat{\amalg}\, L'$$
is inert. \\
In particular,   $Q\subset L\,\widehat{\amalg}\, L'$ is inert.\end{Prop}

\vspace{3mm}\noindent {\sl proof.} By hypothesis we have dgl's $(L\, \widehat{\amalg}\, \overline{\mathbb L}_N,\partial)$ and $(H(L' \, \widehat{\amalg}\, \overline{\mathbb L}_{N'}, \partial')$ with $\partial :N\stackrel{\cong}{\to} Q$ and $\partial' : N'\stackrel{\cong}{\to} Q'$
 and
$$H(L)\to H'L\, \widehat{\amalg}\, \overline{\mathbb L}_N) \hspace{3mm}\mbox{ and } H(L')\to H(L'\, \widehat{\amalg}\, \overline{\mathbb L}_{N'})$$
both surjective.

 Since $H(L\, \widehat{\amalg}\, L')\cong H(L)\, \widehat{\amalg}\, H(L')$ this is a direct consequence of the following commutative diagram in which $\psi$ is  surjective (Proposition \ref{p10.1}(iii)),

$$\xymatrix{H(L\, \widehat{\amalg}\, L') \ar[d]^\cong \ar[rr] && H(L \, \widehat{\amalg}\, L' \, \widehat{\amalg}\, \overline{\mathbb L}_N \, \widehat{\amalg}\, \overline{\mathbb L}_{N'})\ar[d]^\cong\\
H(L) \, \widehat{\amalg}\, H(L') \ar[rr]^-\psi && H(L \, \widehat{\amalg}\, \overline{\mathbb L}_N, \partial) \, \widehat{\amalg}\, H(L' \, \widehat{\amalg}\, \overline{\mathbb L}_{N'}, \partial').}$$

   \hfill$\square$

 \subsection{Semi-direct products and weighted Lie algebras}

 A closed ideal, $I$, in a complete enriched Lie algebra, $L$, determines (\S 3) the weighted Lie algebra $g(I):= \oplus_{n\geq 1} I^{(n)}/I^{(n+1)}$. Then $g(I)$ is an enriched Lie algebra with completion
 $$\overline{g(I)} = \prod_{n\geq 1} I^{(n)}/I^{(n+1)}.$$
 Moreover, the right adjoint representation of $L$ induces representations of $L/I$ in each $I^{(n)}/I^{(n+1)}$. The resulting semi-direct product $$g(I,L):= g(I) \, \widehat{\times}\, L/I$$ inherits a natural enriched structure with completion the semi-direct product $\overline{g(I)}\, \widehat{\times}\, L/I$.
 
\vspace{2mm} \emph{Henceforth in this section we fix such a pair $I\subset L$ and assume further that 
 $$I = \overline{\mathbb L}_T$$
 is a profree Lie algebra.}  Then Proposition \ref{p6.4} provides closed subspaces $\overline{T(k)}\subset I^{(k)}$, $k\geq 1$,  satisfying
$$T= T(1), \hspace{3mm} I^{(k)} = \overline{T(k)} \oplus I^{(k+1)},\hspace{3mm}\mbox{and } [\overline{T(k)}, \overline{T(\ell)}] \subset \overline{T(k+\ell)}.$$
This identifies $\overline{T(k)} = I^{(k)}/I^{(k+1)}$, and $g(I) = \oplus_{k\geq 1} \overline{T(k)}$. Moreover (Proposition \ref{p6.4}(iii))
$$\overline{g(I)} = \prod_{k\geq 1} \overline{T(k)} = \overline{\mathbb L}_T = I.$$
In particular, the construction above defines the semi-direct product
$$I\, \widehat{\times}\, L/I = \overline{g(I)} \, \widehat{\times}\, L/I.$$

Finally suppose, in addition to the notation and hypotheses above, that $Q$ is a closed subspace of $T$  and  which generates the ideal $I= \overline{\mathbb L}_T$. The inclusion extends to a morphism $\psi : \overline{\mathbb L}_Q \rightarrow \overline{\mathbb L}_T$ and this extends to the morphism
$$\varphi := \psi \, \widehat{\amalg}\, id_{L/T} : \overline{\mathbb L}_Q\, \widehat{\amalg}\, L/I \to \overline{\mathbb L}_T \, \widehat{\times}\,L/I.$$

\begin{Prop}
\label{p23.4}
With the hypotheses and notation above the following conditions are equivalent :
\begin{enumerate}
\item[(i)] $Q$ is inert in $L$
\item[(ii)] $Q$ is inert in $\overline{\mathbb L}_T \, \widehat{\times}\, L/I$.
\item[(iii)] $\varphi$ is an isomorphism.
\end{enumerate}
\end{Prop}

\vspace{3mm}\noindent {\sl proof.} (i) $\Leftrightarrow$ (ii)  By Theorem \ref{t8}, $Q$ is inert in $L$ if and only if the right adjoint representation of $L/I$ in $I/I^{(2)}$ induces an isomorphism
$$Q \, \widehat{\otimes}\, \overline{UL/I} \longrightarrow I/I^{(2)}.$$
But by definition this is the right representation of $L/I$ in $\overline{\mathbb L}_T/\overline{\mathbb L}_T^{(2)}$ resulting from the semi-direct product, and so (i) $\Leftrightarrow$ (ii) follows from Theorem \ref{t8}. 

(ii) $\Leftrightarrow$ (iii) First observe (Corollary 2 to Theorem \ref{t7}) that $Q$ is inert in $\overline{\mathbb L}_Q$. Therefore (Proposition \ref{p23.3}) $Q$ is inert in $\overline{\mathbb L}_Q \, \widehat{\amalg}\, L/I$. The kernel of the surjection
$$\overline{\mathbb L}_Q \, \widehat{\amalg}\, L/I \to L/I$$
is the ideal $J$ generated by $Q$. Since $Q$ is inert in the free product, and since it generates the ideal $J$, it follows from Theorem \ref{t8} that the right representation of $L/I$ in $J/J^{(2)}$ defines an isomorphism
$$Q \, \widehat{\otimes }\, \overline{UL/I} \stackrel{\cong}{\longrightarrow} J/J^{(2)}.$$

On the other hand, since $\varphi$ is a morphism of enriched Lie algebras it restricts to a morphism $J\to I$. This in turn induces a morphism
$$\chi : J/J^{(2)} \to I/I^{(2)}$$
of $L/I$-modules, in which $\chi_{\vert Q} = id_Q$. In view of the observations above, $\chi$ may be identified with the morphism
$$Q \, \widehat{\otimes}\, \overline{UL/I} \to I/I^{(2)}$$
induced by the right representation of $L/I$ in $I/I^{(2)}$. Thus $Q$ is inert in $L$ if and only if $\chi$ is an isomorphism.

Finally, since $J$ is the ideal generated by $Q$ it follows from Theorem 8 that $J$ is profree. Thus $\varphi : J\to I$ is a morphism of profree Lie algebras, and is therefore an isomorphism if and only if it induces an isomorphism $J/J^{(2)} \stackrel{\cong}{\longrightarrow} I/I^{(2)}$. 
\hfill$\square$

 \newpage

\newpage
\part{Applications in topology}

\section{Characterization of rationally inert topological attaching maps}

Recall (\S 19) that a topological attaching map
$$g: \vee_\alpha S^{n_\alpha} \to X$$
is \emph{rationally inert} if
$$\pi_*(\iota_X)_{\mathbb Q} : \pi_*(X_{\mathbb Q})\to \pi_*(X\cup_g\vee_\alpha D^{n_\alpha +1})_{\mathbb Q}$$
is surjective. Moreover, if $\lambda : \land V\to \land W$ is a Sullivan representative of $\iota_X : X\to X\cup_g D^{n_\alpha +1}$, then (\ref{i16}) 
$$\pi_*(\iota_X)_{\mathbb Q} = \pi_*\langle \lambda\rangle = sL_\lambda.$$
It follows that $g$ is rationally inert if and only if $L_\lambda$ is surjective.

\vspace{3mm} Next, let $L_X$ denote the homotopy Lie algebra of $X$, and  recall that associated with $g$ are two Lie algebra morphisms:
\begin{enumerate}
\item[$\bullet$] A dgl representative $\ell_g: (\overline{\mathbb L}_T,0) \to (L, \partial)$ of $g$.
\item[$\bullet$] The desuspension $L_g : \overline{\mathbb L}_T \to L_X$ of $\pi_*(g_{\mathbb Q})$.
\end{enumerate}

It follows from (\ref{i27}) that $L_g = H(\ell_g)$. Thus Proposition \ref{p19.1} and Theorem \ref{t6} combine to yield the equivalence of the three assertions:
\begin{enumerate}
\item[$\bullet$] $g$ is rationally inert.
\item[$\bullet$] $\ell_g$ is inert.
\item[$\bullet$] $L_g$ is inert.
\end{enumerate}

\vspace{3mm} On the other hand recall from (\ref{i35}) the decomposition
$$\land V \stackrel{\eta}{\longrightarrow} \land V\otimes \land Z \stackrel{\simeq}{\longrightarrow} \land W$$
of the minimal Sullivan model, $\land W$, of $X$. This in turn (\cite[Proposition 17.9]{FHTI}) provides a fibration
\begin{eqnarray}
\label{i41}
\xymatrix{ \langle \land V\rangle &&\langle \land V\otimes \land W\rangle \ar[ll]_{(\iota_X)_{\mathbb Q}} && \langle \land Z\rangle \ar[ll].}
\end{eqnarray}
By Proposition  \ref{p17.6}, if $\rho : \land W\to \mathbb Q \oplus S$ is surjective then $\langle \land Z\rangle$ is homotopy equivalent to the homotopy fibre of $(\iota_X)_{\mathbb Q}$. 

\begin{Theorem}
\label{t9}
With the hypotheses and notation above, if $\rho$ is surjective then $$
\mbox{$g$ is rationally inert $\Longleftrightarrow$ $\langle \land Z\rangle$ is rationally wedge-like.}$$
\end{Theorem}

\vspace{3mm}\noindent {\sl proof.} If $g$ is rationally inert then (Theorem \ref{t7}) $\land V\otimes \land Z \stackrel{\cong}{\longrightarrow} \land W$ and the homotopy Lie algebra, $L_Z$, is profree. By definition (\S 10.4) $\langle \land Z\rangle$ is rationally wedge-like.

In the reverse direction we need to show that if $g$ is not rationally inert then $\langle \land Z\rangle$ is not rationally wedge-like. The proof is essentially the same as the proof that (iii) $\Rightarrow $ (ii) in Theorem \ref{t7}. If $g$ is not rationally inert then $\pi_*(\iota_X)_{\mathbb Q}$ is not surjective. Thus in the long exact homotopy sequence for (\ref{i41}) the connecting homomorphism
$$\partial : \pi_*\langle \land V\rangle \to \pi_{*-1}\langle \land Z\rangle$$
is not trivial. In particular for some $0\neq y\in \pi_{*-1}\langle \land Z\rangle$, $y= \partial x$. A standard argument from classical homotopy theory shows that Whitehead products 
$$Wh(\partial x,z), \hspace{3mm} z\in \pi_*\langle \land Z\rangle$$
all vanish.

On the other hand, if $\langle \land Z\rangle$ is rationally wedge-like then we would have $\langle \land Z\rangle \simeq \langle \land Q\rangle$ with the homotopy Lie algebra $L_Q$ of $\land Q$ a profree Lie algebra. Thus $y$ would correspond to an element $y'\in \pi_*\langle \land Q\rangle$ and all the Whitehead products $y'*z'$, $z'\in \pi_*\langle \land Q\rangle$ would vanish. Now it would follow from Proposition \ref{p8.2} that $y'$ would correspond to a non-zero element in the centre of $L_Q$. Since $L_Q$ is profree this would imply that $L_Q$ is a free Lie algebra on a single generator.

Finally, since the isomorphism $\pi_*\langle \land Q\rangle \cong \pi_*\langle \land Z\rangle$   preserves Whitehead products, it would also follow that the homotopy Lie algebra $L_Z$ is a free Lie algebra on a single generator, so that necessarily dim$\, H^1(\land Z)\leq 1$. In this case the proof of (iii) $\Rightarrow$ (ii) in Theorem \ref{t7} shows that $g$ is rationally inert. This   establishes the Theorem.\hfill$\square$

 \section{Poincar\'e duality complexes}  We say a CW complex $Y= X\cup_fD^{n+1}$ is a \emph{rational Poincar\'e duality complex} if $H(Y)$ is a Poincar\'e duality algebra and the top class is in the image of $H(Y,X)$. In this case it follows that $H^{\leq n}(X) \stackrel{\cong}{\to} H(X)$.  
 
\begin{Theorem} \label{t11} If $Y = X\cup_fD^{n+1}$ is a rational Poincar\'e duality complex and the algebra $H(Y)$ requires at least two generators,  then $f$ is rationally inert.\end{Theorem}
 
 \vspace{3mm} This Theorem was first established for simply connected spaces in \cite[Theorem 5.1]{HL}.
 
 \vspace{3mm} 
 \noindent {\sl proof.} Recall from \S 19 that $f$ is rationally inert if
 $$\pi_*(i_{\mathbb Q}) : \pi_*(X_{\mathbb Q}) \to \pi_*(X\cup_f D^{n+1})_{\mathbb Q}$$
 is surjective. Now let $\lambda : \land V\to \land W$ be a Sullivan representative of $i : X\to X\cup_f D^{n+1}$ and then decompose $\land W$ as a $\Lambda$-extension
 $$\xymatrix{\land V \ar[rr]^\eta&& \land V\otimes \land Z\ar[rr]^\xi && \land Z}$$
 as described in (\ref{i35}). Then according to Theorem \ref{t7} $f$ is rationally inert if and only if the homotopy Lie algebra, $L_Z$, is profree, i.e. (Theorem 1)
 $$\land Z \simeq (\mathbb Q\oplus S, 0)$$
 with $S\cdot S= 0$. 
 
Now let $\land V\otimes \land U$ be the acyclic closure of $\land V$. Then we have
 $$\land W\otimes\land U:= \land W\otimes_{\land V} (\land V\otimes \land U)\simeq \land Z.$$
Therefore we have only to prove that 
$$\land W\otimes \land U \simeq (\mathbb Q \oplus S,0)$$
with $S\cdot S= 0$. 
 
 For this, we first establish the following notation. Let $P$ be a direct summand
  in $(\land V)^{n+1}$ of $(\land V)^{n+1}\cap \mbox{ker}\, d$. Then division by $P$ and by $(\land V)^{>n+1}$ defines a surjective quasi-isomorphism $\land V\stackrel{\simeq}{\to} A$, and 
 $$A^{n+1}= A^{n+1} \cap \mbox{Im}\, d \oplus \mathbb Q \omega,$$
 where $\omega$ is a cycle representing the top cohomology class of $Y$. 
 
 In the same way, let $P'$ be a direct summand in $(\land W)^n$ of $(\land W)^n\cap \mbox{ker}\, d$. Then division by $P'$ and $(\land W)^{>n}$ defines a surjective quasi-isomorphism $\land W\stackrel{\simeq}{\to} B$. Moreover, by construction $\varphi$ factorizes in a commutative diagram,
 $$\xymatrix{\land V \ar[d]^\simeq\ar[rr]^\varphi &&\land W\ar[d]^\simeq\\
 A\ar[rr]^{\overline{\varphi}} && B,}$$
 in which $\overline{\varphi}\omega= 0$. 

 Now consider the commutative cdga diagram
 $$
 \xymatrix{
 A \ar[rr]^{\overline{\varphi}}\ar[rrd]_j && B\\
 && A\oplus \mathbb Q t\ar[u]^\psi_\simeq}$$
 in which $dt= \omega$, $t\cdot A^{\geq 1}= 0$, $t^2=0$, and $\psi t = 0$. Here $\psi$ is a quasi-isomorphism and $j$ is a cdga representative of $\varphi$. This provides a cdga quasi-isomorphism
 $$\land W\otimes \land U \simeq (A\oplus \mathbb Q t)\otimes_A (A\otimes_{\land V}(\land V\otimes \land U)) = (A\oplus \mathbb Q t)\otimes \land U.$$
 Thus from the short exact sequence of chain complexes
 $$0\to A\otimes \land U \to (A\oplus \mathbb Qt)\otimes \land U\to \mathbb Q t\otimes \land U\to 0,$$
 we deduce a linear isomorphism
 $$H^{\geq 1}((A\oplus \mathbb Qt)\otimes \land U ) \stackrel{\cong}{\longrightarrow} \mathbb Qt\otimes \land U $$
  of graded vector spaces. It remains to show that we can lift a basis of $\mathbb Q t \otimes \land U$ to cycles $\Phi_i\in (A\oplus \mathbb Q t)\otimes \land U$ such that $\Phi_i\cdot \Phi_j= 0$. 
 
Before going further we first eliminate two special cases. First if $V^1= 0$ the argument of (\cite[\S 5]{HL}) shows that $(A\oplus \mathbb Qt)\otimes \land U$ is rationally wedge-like, and so $f$ is rationally inert. (Note that in \cite{HL} it is assumed that $X$ is simply connected; however the proof of this assertion relies only on the fact that $V^1= 0$.) Secondly, if $n=1$ then $X \simeq_{\mathbb Q} S^1_1\vee \dots \vee S^1_{2q}$ and so $Y$ is rationally equivalent to an oriented Riemann surface. In this case Theorem \ref{t12} is established in \cite{Artin}.

 Thus we may now assume that $n\geq 2$ and that $A^1$ contains a non-zero cycle $v$. Since $H(A)$ is a Poincar\'e duality algebra there is a cycle $w\in A^n$ such that $wv=\omega$. The first step for the proof is then
 
\begin{lem} \label{l25.1} With the hypotheses and notation above, $A^{n+1}\otimes \land U \subset d(A^n\otimes \land U)$.\end{lem}
 
 \vspace{3mm}\noindent {\sl proof:} Choose $\overline{v}\in U^0$ so that $d\overline{v}=v$. Since $\land V$ is a minimal Sullivan algebra, $V$ is the union of an increasing sequence of subspaces $V(0)\subset \dots \subset V(q)\subset \dots$ in which $V(0)= \mathbb Q v$ and $d: V(q+1)\to \land V(q)$. It follows that $U$ is the union of an increasing sequence of subspaces $U(0)\subset \dots \subset U(q)\subset \dots$ in which $U(0)= \mathbb Q \overline{v}$ and $$d: U(q+1)\to A^{\geq 1}\otimes \land U(q).$$
 We show by induction on $q$ that
$$
 A^{n+1}\otimes \land U(q) \subset d(A^n\otimes \land U(q))
$$
 
 First note that any $z\in A^{n+1}$ has the form $z= dy + \lambda wv$, some $\lambda \in \mathbb Q$. Thus
 $$z\otimes 1 = d(y\otimes 1) + (-1)^n d(\lambda w\otimes \overline{v}) \in d(A^n\otimes \land U(0)).$$
 Then for $r\geq 1$,
 $$z\otimes \overline{v}^r = d(y\otimes \overline{v}^r + \frac{(-1)^n\lambda}{r+1} w \otimes \overline{v}^{r+1}) - (-1)^n  ryv\otimes \overline{v}^{r-1}.$$
 It follows by induction on $r$ that $A^{n+1}\otimes \land U(0)\subset d(A^n\otimes \land U(0)).$
 
 Now fix a direct summand, $T$, of $U(q)$ in $U(q+1)$, and assume by induction that for some $s$,
 $$A^{n+1} \otimes \land U(q)\otimes \land^{\leq s}T \subset d(A^n\otimes \land U(q)\otimes \land ^{\leq s} T).$$ Then write $\Phi\in A^{n+1}\otimes \land U(q)\otimes \land^{\leq s+1}T$ as
 $\Phi = \sum \Phi_i \otimes \Psi_i$ with $\Phi_i \in A^{n+1}\otimes \land U(q)$ and $\Psi_i\in \land^{\leq s+1}T$. By the hypothesis $\Phi_i = d\Omega_i$ with $\Omega_i \in A^n\otimes \land U(q)$. Therefore
 $$\sum \Phi_i\otimes \Psi_i = d(\sum \Omega_i\otimes \Psi_i) -(-1)^n \sum \Omega_i \land d\Psi_i.$$
 The first term is in $d(A^n\otimes \land U(q)\otimes \land^{\leq s+1}T)$. On the other hand, $d\Psi_i \in A^{\geq 1} \otimes \land U(q)\land^{\leq s}T$ and so the second term is in $ A^{n+1}\otimes \land U(q)\otimes \land^{\leq s} T$. By hypothesis, the second term is contained in $d(A^n\otimes \land U(q)\otimes \land^{\leq s}T)$. This closes the induction. \hfill$\square$
 
 \vspace{3mm}Now, let $\Phi\in \land U$. Then 
 $$t-(-1)^n w\overline{v} \in (A\oplus \mathbb Qt)\otimes \land U$$
 is a cycle, and 
 $$d((t-(-1)^nw\overline{v})\Phi) = -w\overline{v} \, d\Phi \in A^{n+1}\otimes \land U.$$
 By Lemma \ref{l24.1}, $w\overline{v}\, d\Phi = d\Psi$ for some $\Psi \in A^n\otimes \land U$. Thus $(t-(-1)^n)w\overline{v}) \Phi + \Psi$ is a cycle projecting to $t\otimes \Phi$ in $\mathbb Q t\otimes \land U$. Thus such cycles map to a basis of $\mathbb Q t\otimes \land U$. But because $n\geq 2$, $2n>n+1$ and so the product of any two of those cycles is zero. This shows that $f$ is rationally inert. \hfill$\square$

 \section{Whitehead's problem}

  \vspace{3mm} A famous unsolved problem of JHC Whitehead \cite{W} asks: is a subcomplex of an aspherical two-dimensional CW complex aspherical ? As observed by Anick \cite{An} it is sufficient to consider the case that both subcomplexes share the same $1$-skeleton and base point. The problem then reduces to the question: If $X$ is a finite 2-dimensional connected CW complex and $X\cup \left( \vee_{k=1}^p D_k^2\right)$ is aspherical, is $X$ aspherical ?
 
 In \cite{An} Anick provides a positive answer to an analogous question for simply connected rational spaces. Here we have a positive answer for   Sullivan spatial completions of   connected spaces.
 
\begin{Theorem}\label{t12} If $X$ is a connected CW complex and $(X\cup \vee_{k=1}^p D^2)_{\mathbb Q}$ is aspherical then $X_{\mathbb Q}$ is aspherical. \end{Theorem}
 
 \vspace{3mm}\noindent {\sl proof:} The obvious induction reduces the statement to the case $p=1$. Then, since $\pi_*((X\cup_fD^2)_{\mathbb Q}) \cong V^\vee$ as sets where $\land V$ is the minimal Sullivan model of $X\cup_fD^2$, our hypothesis simply implies that $V = V^1$. Let $\varphi : (\land V,d)\to (\land W,d)$ be a Sullivan representative for the inclusion $i : X \to X\cup_fD^2$.  Since $H^1(X\cup_fD^2)\to H^1(X)$ is injective, it follows that $\varphi$ is injective and so $\land W$ decomposes as $\land V\otimes \land Z$, with $Z= Z^{\geq 1}$. In particular $ f $ is rationally inert. Moreover, it follows from Theorem \ref{t7} together with (36) that 
 $$H^{\geq 1}(\land Z, \overline{d}) \cong \mathbb Q b\otimes \land U,$$
 where deg$\, b= 1$ and $\land V\otimes \land U$ is the acyclic closure of $\land V$. Since $V = V^1$, $U= U^0$ and $H^{\geq 1}(\land Z, \overline{d}) = H^1(\land Z, \overline{d})$. This in turn implies $Z= Z^1$ and $X_{\mathbb Q}$ is aspherical.
 
 \hfill $\square$

\section{Homotopy Lie algebra of a 2-cone}
 
 \vspace{5mm} A $2$-cone is a space $Y= \vee_j S^{n_j}\cup_g \vee_i D^{k_i+1}$ constructed from an attaching map $g: \vee_i S^{k_i} \to X=\vee_j S^{n_j}$. We denote by
 $$\iota_X : X = \vee_jS^{n_j}\to X\cup_g \vee_i D^{k_i+1}= Y$$
 the inclusion.
 
 \begin{Prop}
 \label{p27.1}
 Suppose $\iota_X : X\to Y$ is a two cone, and denote by $E$ the image of $L_X$ in $L_Y$. Assume that $H^{\geq 1}(g)=0$ and that $\pi_*(g_{\mathbb Q})$ is injective. Then there is a retraction $r : L_Y\to E$ whose kernel is a profree Lie algebra $K$, i.e. there is a short exact sequence
 $$0\to K\to L_Y\to E\to 0.$$
 In particular, the attaching map $g$ is rationally inert if and only if $K= 0$, and otherwise $L_Y$ contains a profree Lie algebra. 
 \end{Prop}

 \vspace{3mm}\noindent {\sl proof.}  Let $(\land W,d)$ be a quadratic model of $X$ and let $\rho : W\to S$ the map associated to $g$. Then (Proposition \ref{p17.5}) a cdga model of $Y$ is given by the cdga
 $$(\land W\oplus s^{-1}S,D).$$
 
 We equip $\land W\oplus s^{-1}S$ with a new gradation $\land W\oplus s^{-1}S = \oplus_{k\geq 0} (\land W\oplus s^{-1}S)^{(k)}$ with $W = W^{(1)}$ and $s^{-1}S= (s^{-1}S)^{(2)}$. Then $D$ increases this gradation by $1$. Therefore the cohomology of $\land W\oplus s^{-1}S$ inherits a new gradation
 $$H^*(\land W\oplus s^{-1}S,D) = \oplus_{p\geq 0} H^{(p)}.$$
 The short exact sequence
 $$0\to (s^{-1}S,0)\to (\land W\oplus s^{-1}S,D)\to (\land W,d)\to 0$$
 shows that $H^{(p)}= 0$ for $p>2$. On the other hand, clearly, $H^{(0)}= \mathbb Q$ and, since $X$ is a wedge of spheres, $H^{(1)}= H^{\geq 1}(X)$. 

Now, let $\varphi :(\land V,d)\stackrel{\simeq}{\longrightarrow} (\land W\oplus s^{-1}S,D)$ be its Sullivan minimal model. We can suppose $(\land V,d)$ equipped with a new gradation preserved by $\varphi$. In particular, $(\land V^{(1)},d)$ is a sub cdga and its homotopy Lie algebra is isomorphic to $E$. The injection $(\land V^{(1)},d)\to (\land V,d)$ induces the surjection $L_X\to E$. 

Denote by $(\land V^{(1)}\otimes \land U,d)$ be the acyclic closure of $(\land V^{(1)},d)$ equipped with the new gradation given by $U = U^{(0)}$. Then let $\land Z$ be the Sullivan minimal model of $(\land V\otimes \land U,d)$. Then $Z = \oplus_{p\geq 2} V^{(p)}$, and since $H(\land V\otimes \land U,d) = \mathbb Q \oplus H^{(2)}(\land V\otimes \land U,d)$, the cdga $(\land Z,d)$ is quasi-isomorphic to $(\mathbb Q \oplus (Z\cap\mbox{ker}\, d),0)$. Thus $L_Z$ is profree and we get the required extension
$$0\to L_Z\to L_Y\to E\to 0.$$
\hfill$\square$

 \section{Wedges of circles}
 
 \subsection{Inertness}
 
 For simplicity we say that $[f]\in \pi_1(X)$ is rationally inert if $f: S^1\to X$ is rationally inert.

\begin{Theorem} \label{t13} If $X$ is a wedge of at least two circles then any non-zero $[f]\in \pi_1(X)$     is rationally inert  or, equivalently, $(X\cup_fD^2)_{\mathbb Q}$ is aspherical. \end{Theorem}
  
  \vspace{3mm}\noindent {\sl proof:} 
 It follows from Proposition 10.6 and its proof that the homotopy Lie algebra, $L_X$ is profree and concentrated in degree zero. Thus by Theorem 10, every non-zero element in $L_X$ is inert. It is therefore sufficient to show that $\pi_1(X)\to \pi_1(X_{\mathbb Q})$ is injective. Moreover, any element of $\pi_1(X)$ is in the image of $\pi_1(Y)$ for some finite sub wedge of circles $Y\subset X$. Then $Y$ (and hence $Y_{\mathbb Q}$) are retracts of $X$ and $X_{\mathbb Q}$ and so it is sufficient to prove that each element $[g]\in\pi_1(Y)$ is rationally inert.
  
For this, denote $G = \pi_1(Y)$, so that $G_{\mathbb Q}= \pi_1(Y_{\mathbb Q})$. According to \cite[Theorem 7.5]{RHTII}, $G^n/G^{n+1} \otimes \mathbb Q \stackrel{\cong}{\to} G^n_{\mathbb Q}/G^{n+1}_{\mathbb Q}$. But by \cite{hall}, $G^n/G^{n+1}$ is a free abelian group, and hence $G^n/G^{n+1}\to G^n_{\mathbb Q}/G^{n+1}_{\mathbb Q}$ is injective. Since $G$ is a free group, $G\to \varprojlim_n (G/G^n)_{\mathbb Q}$ is injective and the image of $[g] $ in $G_{\mathbb Q}$ is non-zero.

Finally since any non-trivial $[f]\in\pi_*(X)$ is rationally inert it follows by definition that
$$\pi_1(X_{\mathbb Q}) \to \pi_*((X\cup_f D^2)_{\mathbb Q})$$
is surjective. Since $\pi_*(X_{\mathbb Q})$ is concentrated in degree 1, $(X\cup_fD^2)_{\mathbb Q}$ is aspherical.

\hfill$\square$

\subsection{A connected finite CW complex, $Y$, whose homotopy Lie algebra, $L_Y$ is not nilpotent and satisfies dim$\,(L_Y/L_Y^2)_1=\infty$.}
 
 \vspace{3mm}\noindent {\sl Step One. Construction of $Y$ and a first dgl model}
 
 \vspace{2mm} In \cite{JML} Lemaire constructs a simply connected finite CW complex, $X$, whose homotopy Lie algebra, $L_X$, is not finitely generated. Since $L_X$ is a graded vector space of finite type, the same is true for $L_X/L_X^2$, and so the generators occur in infinitely many degrees.
 
 Lemaire's example is constructed by attaching seven cells to a wedge of five spheres. Here we use a similar process to construct a CW complex, $Y$ by attaching seven $2$-cells $D_{a_i}^2$ to a wedge of five circles $S^1_{b_j}$ via a map
 $$g : \vee_{i=1}^7 S^1_{a_i} \longrightarrow \vee_{j=1}^5 S^1_{b_j}.$$
 The respective fundamental groups are free groups $G_a$ and $G_b$ freely generated by generators $a_1, \dots , a_7$ and $b_1,\dots , b_5$.
 
 We denote the commutators in $G_b$ by
 $$(c,d)= cdc^{-1}d^{-1}.$$
 Then $g$ is defined by the homomorphism $\pi_1(g)$ given explicitly by
 
  $$a_1\mapsto  (b_1,b_3), a_2 \mapsto (b_1, b_4), a_3\mapsto (b_2, b_3), a_4 \mapsto (b_2, b_4),\hspace{3mm}\mbox{and }$$
  $$ a_5 \mapsto (b_5, b_1b_3^{-1}), a_6\mapsto(b_5, b_1b_4^{-1}), a_7\mapsto (b_5, b_2b_3^{-1}).$$
Now set $X = \vee_{j=1}^5 S^1_{b_j}$ and define
$$Y = X\cup_g \vee_{i=1}^7 D^2_{a_i}.$$
As usual the attachment is denoted by $\iota_X : X\to X\cup_g \vee_{i=1}^7 D^2_{a_i} = Y.$

We now combine Proposition \ref{p9.4} and the Example in \S 17.1 to construct a first dgl model of $Y$. As in \S 9.1, denote by $\log_{G_b}$ the map
$$\log_{G_b} : G_b\longrightarrow G_{L_X} \stackrel{\log}{\longrightarrow} L_X,$$
recalling that $\log$ is a bijection of sets. Then (Proposition \ref{p9.4}(i)) $L_X$ is a profree Lie algebra $\overline{\mathbb L}_{T(b)}$ in which the elements $\log_{G_b}(b_j)$ are a basis of $T(b)$. Then (\S 17.1) a dgl representative of $\iota_X$ is given by 
$$(L_X, 0) \to (L_X\, \widehat{\amalg}\, \overline{\mathbb L}_{sT(a)}, \partial),$$
in which the elements $\log_{G_a}(a_i)$ are a basis of $T(a)$ and $$\partial \, s\log_{G_a}(a_i) = \log_{G_b}(\pi_1(g)(a_i)).$$
\emph{In particular, $(L_X \, \widehat{\amalg}\, \overline{\mathbb L}_{sT(a)}, \partial)$ is a dgl model for $Y$. }

To simplify notation we set
$$\log_{G_b}b_j = x_j\hspace{5mm}\mbox{ and } \log_{G_a}(a_i) = y_i.$$
Then the dgl morphism above can be written as
$$(L_X, 0) \to L_X\, \widehat{\amalg}\, \overline{\mathbb L}_{sT(a)} = \overline{\mathbb L}(x_1, \dots x_5)\, \widehat{\amalg}\, \overline{\mathbb L}(sy_1, \dots , sy_7)$$
in which $\partial (sy_i) = \log_{G_b}(\pi_1(g)(a_i))$. Specifically,
$$\partial sy_1 = \log_{G_b}(b_1, b_3), \, \partial sy_2= \log_{G_b}(b_1,b_4),\, \partial sy_3= \log_{G_b}(b_2, b_3),\, \partial sy_4= \log_{G_b}(b_2, b_4)$$
and
$$\partial sy_5 = \log_{G_b}(b_5, b_1b_3^{-1}),\, \partial sy_6 = \log_{G_b}(b_5, b_1b_4^{-1}), \hspace{3mm}\mbox{and }\partial sy_7 = \log_{G_b}(b_5, b_2b_3^{-1}).$$

\vspace{2mm}\noindent {\sl Step Two. A second dgl model of $Y$.}

\vspace{2mm} We define a second dgl model, $(L_X\, \widehat{\amalg}\, \overline{\mathbb L}(sy_1, \dots , sy_7), \partial')$ by setting
$$\partial'sy_1= [x_1, x_3], \partial'sy_2= [x_1, x_4], \partial'sy_3= [x_2, x_3], \partial'sy_4= [x_2, x_4], \hspace{3mm}\mbox{and}$$
$$\partial'sy_5 = [x_5, x_1-x_3], \partial'sy_6 = [x_5, x_1-x_4] \hspace{3mm}\mbox{and } \partial'sy_7= [x_5, x_2-x_3].$$
We then identify this as a dgl model for $Y$ by constructing a commutative diagram of dgl morphisms
$$
\xymatrix{ 
&& (L_X\, \widehat{\amalg}\, \overline{\mathbb L}(sy_1, \dots , sy_7), \partial)\ar[dd]_\cong ^\psi\\
(L_X,0)\ar[rru]\ar[rrd] \\
&& (L_X\, \widehat{\amalg}\, \overline{\mathbb L}(sy_1, \dots , sy_7), \partial').}$$
For this, as in \S 9.1, denote by $\overline{J}$ the augmentation ideal of $\overline{UL_X}$. Then the left and right adjoint representations of $L_X$ extend to representations of $\overline{UL_X}$. In particular, if $u\in L_X$ and $\alpha \in \overline{J}$ then
$$u\cdot \alpha = \alpha^t\cdot u,$$
where if $\alpha = \alpha_1\dots \alpha_k$ with $\alpha_i\in L_X$ then
$$\alpha^t = (-1)^k \alpha_k\dots \alpha_1.$$
More generally, if $\alpha = \sum \alpha_k$, with $\alpha_k\in \overline{J}^k$ then $\alpha^t = \sum \alpha_k^t$. 

Similarly, if $\alpha = \alpha_1\dots \alpha_k\in \overline{J}$, with $\alpha_i\in L_X$, we define $\alpha \cdot sy_j \in L_X\, \widehat{\amalg}\, \overline{\mathbb L}(sy_j)$ by setting
$$\alpha \cdot sy_j = [\alpha_1, [\alpha_2, \dots [\alpha_k, sy_j]\dots ]].$$
It is then immediate that for any $u\in L_X$, $\alpha \in \overline{J}$ and $1\leq i\leq 7,$
$$[u\cdot \alpha, sy_i] = [u, \alpha\cdot sy_i].$$

We now define $\psi$. It is immediate from Proposition \ref{p9.4}(ii) that for some $\alpha_1, \dots , \alpha_4\in \overline{J}$,
$$\partial sy_1= [x_1, x_3] + [x_1, x_3]\cdot \alpha_1, \hspace{5mm}\partial sy_2= [x_1, x_4] + [x_1, x_4]\cdot \alpha_2, $$
and
$$\partial sy_3 = [x_2, x_3] + [x_2, x_3]\cdot \alpha_3, \hspace{5mm}\partial sy_4= [x_2, x_4] + [x_2, x_4]\cdot \alpha_4.$$
We define 
$$\psi (sy_i)= sy_i + \alpha_i\cdot sy_i, \hspace{5mm} 1\leq i\leq 4.$$

Next, observe that for some $\beta$ and $\alpha_5\in \overline{J}$,
$$\log_{G_b}(b_1b_3^{-1}) = x_1-x_3- \frac{1}{2}[x_1,x_3]+ [x_1, x_3]\cdot \beta,$$
and
$$ (x_5, \log_{G_b}(b_1b_3^{-1})) = [x_5, \log_{G_b}(b_1b_3^{-1})] + [x_5, \log_{G_b}(b_1, b_3^{-1})]\cdot \alpha_5.$$
It follows that 
$$[x_5, \log_{G_b}(b_1b_3^{-1})] = [x_5, x_1-x_3] + \frac{1}{2}[x_1, x_3],x_5] - [x_1, x_3]\cdot \beta x_5.$$
This gives $$\partial sy_5 = [x_5, x_1-x_3] + [x_1,x_3]\cdot \beta_5'= [x_5, x_1-x_3] + \beta_5\cdot [x_1, x_3],$$
for some $\beta'_5$ and $\beta_5$ in $\overline{J}$. 
Similarly,
$$\partial sy_6 = [x_5,x_1-x_4] + \beta_6\cdot [x_1, x_4] \hspace{3mm}\mbox{and } \partial sy_7= [x_5, x_2-x_3] + \beta_7 \cdot [x_2, x_3].$$
Now we define
$$\psi (sy_5)= sy_5 + \beta_5\cdot sy_1, \psi (sy_6) = sy_6 + \beta_6sy_2, \psi (sy_7) = sy_7 + \beta_7\cdot sy_3.$$
Then $\partial'\psi = \psi \partial$, and it is immediate that $\psi$ is an isomorphism. 

\vspace{3mm}\noindent {\sl Step Three. dim$\, (L_Y/L_Y^2)_1= \infty$.}

\vspace{2mm}
We first consider the uncompleted dgl $(L,d')= (\mathbb L(x_i, sy_j), d')$. Following the computations of \cite[page 118, 119]{JML},
 the dimension of $(H(L,d')/ H(L,d')^2)_1$ is infinite.

Denote by $L_s^r$ the vector space generated by the Lie brackets containing exactly $r$ elements $sy_j$ and $s$ elements $x_i$. 
Since $d(sy_j)\subset \mathbb L^2(x_i)$, the dgl $L$ is the sum of the finite dimensional complexes $C_n = \oplus_{2r+s= n} L_s^r $.  Moreover $\sum_{r\geq r_0}C_r$ is an ideal and we denote by $L(r_0)$ the quotient dgl $L/(\sum_{r\geq r_0}C_r)$.   Thus for any integer $k$ there is a $N$ such that   dim$\, (L(N)/L(N)^2)_1 >k$.

 It follows that its completion $\overline{\mathbb L}(x_i, sy_j)$ decomposes as the product $\prod_n C_n$, so that its homology is $\prod_n H_*(C_n)$ is the completion of $H(L)$.

Therefore 
$$L_Y= H(\overline{\mathbb L}(x_i, sy_j),d) \cong  \prod_n H(C_n).$$
Remark now that $\prod_{r\geq r_0} H(C_r)$ is an ideal in $L_Y$ whose quotient is $L(r_0)$. 
Therefore if dim $(L_Y/L_Y^2)_1 \leq k$ for some integer $k$, then for any $r_0$, dim $(L(r_0)/L(r_0)^2)\leq k$, which is not true.
It follows  that dim$\,(L_Y/L_Y^2)_1=\infty$.

\vspace{3mm}\noindent {\sl Step Four.  $L_Y$ is not a pronilpotent Lie algebra.}

 This follows from   the next Proposition.

\vspace{3mm}
Let $(L, \partial)$ be an enriched dgl satisfying dim$\,L/L^2<\infty$. Denote by $\widehat{H(L)}$ the pronilpotent completion of $H(L)$,
 $$\widehat{H(L)} = \varprojlim_n H(L)/H(L)^n.$$

\begin{Prop} There is a natural injection $j : \widehat{H(L)} \to H(L)$ that is an isomorphism if and only if dim$\, H(L)/H(L)^2 <\infty$.
\end{Prop}
 
 \vspace{3mm}\noindent {\sl proof.} Recall that the complete enriched structure of $H(L)$ is given by the ideals $I^n = \mbox{ker}\, H(L)\to H(L/L^n)$. The pronilpotent completion $\widehat{H(L)}$ can be viewed as the set of the series $\sum \alpha_i$ with $\alpha_i \in H(L)^i$. Since $\alpha_i = [a_i]$ with $a_i\in L^i$, $\alpha_i \in I^i$ and the series $\sum \alpha_i$ is a well defined element in $H(L)$.  The induced morphism $j : \widehat{H(L)}\to H(L)$ is injective because it is an injection when followed by the projections $H(L)\to H(L)/I^n$.
 
 When dim$\, H(L)/H(L)^2 <\infty$, then for each integer $n$ there is $p(n)$ such that $I^{p(n)} \subset H(L)^n$. It follows that each series $\sum \alpha_i$, with $\alpha_i \in I^i$ belongs to $\widehat{H(L)}$, which proves the surjectivity in that case.
 
 Now suppose dim$\, H(L)/H(L)^2 =\infty$. There exists then a sequence of indecomposable elements $\alpha_n \in H(L)$ with $\alpha_n \in I^n$. It follows that $\sum \alpha_n \in H(L)$, but not in $\widehat{H(L)}$.
 \hfill$\square$
 
 \vspace{3mm}\noindent {\bf Example.} Recall that the fundamental group of $S^1_a\vee S^1_b$ is a free group on two generators, $a,b$. Denote by $X$ the space obtained by adding 2 cells of dimension 2 to $S^1_a\vee S^1_b$ to oblige $a$ and $b$ to commute with the commutator $aba^{-1}b^{-1}$. Then $\pi_1(X)$ is a metabelian group and the Lie algebra associated to its Malcev completion is the quotient of the completion of the free Lie algebra on 2 elements $a,b$ by the relations $[a,[a,b]], [b,[a,b]]$. 

Denote by $(\land W,d)$ the minimal model of $S^1_a\vee S^1_b$ and let $(\land V,d)\stackrel{\simeq}{\to} (\land W\oplus s^{-1}P,D)$ be the minimal model of $X$ as in Proposition \ref{p17.5}(ii).  It follows that $V^1$ has a basis with 3 elements $x,y,t$ with $dx=dy=0$, $dt= xy$, and $(\land V,d)$ is quasi-isomorphic to $\land (x,y,t)/xyt$. Decompose as above $\land V = \land V^1\otimes \land Z$, and denote by $\land V^1\otimes \land U$ the acyclic closure of $\land V^1$. The short exact sequence
$$0\to \mathbb Q xyt \to \land (x,y,t)\to \land (x,y,t)/(xyt)\to 0$$ shows that 
$$(\land Z, \overline{d} ) \stackrel{\simeq}{\longrightarrow} (\mathbb Q \oplus  s^{-1}(xyt\otimes \land U), 0).$$
It follows that $L_Z$ is a profree Lie algebra.

\section{Sequential cell attachments}

A sequence $g_k : S^{n_k}\to X$ defines a map $g : \vee S^{n_k}\to X$. This in turn yields a morphism
 $$(L, \partial) \longrightarrow (L\, \widehat{\amalg}\, \overline{\mathbb L}_Q, \partial + \partial_\gamma)$$
 of enriched dgl's, the dgl model for $X\cup_f \vee_k D^{n_k+1}$. Here $sQ = H^{\geq 1}(\vee_k S^{n_k})^\vee.$ Thus if each $S^{n_k}$ is a circle the $sQ = H_1(\vee_k S^{n_k};\mathbb Q)^{\vee\vee}$. If there are infinitely many circles then $Q$ is the dual of a vector space of uncountable dimension, and this may not be a useful model for the construction of computable topological spaces. Therefore we introduce the
 
 \vspace{3mm}\noindent {\bf Definition.} \emph{The sequential cell attachment} associated with the sequence $(g_k)$ is the sequence of inclusions
 $$\varphi_k : X\to X_k:= X\cup_{g(k)} \, \vee_{\ell\leq k} D^{n_\ell+1},$$
 in which $g(k) = \{g_\ell \,\vert\, \ell\leq k\}$. 
 
 Here an enriched dgl representative for $\varphi_k$ is provided by the morphism
 $$(L, \partial) \to (L\, \widehat{\amalg}\, \overline{\mathbb L}_{T_k}, \partial + \partial_\gamma),$$
 where $sT_k = H^{\geq 1}(\vee_{\ell = 1}^k\, S^{n_k})^\vee= H_{\geq 1}(\vee_{\ell = 1}^k S^{n_\ell};\mathbb Q)$ is finite dimensional. The spatial realization of this morphism is then (Proposition \ref{p16.3}) the map
 $$X_{\mathbb Q} \to (X_k)_{\mathbb Q}.$$
 In particular corresponding to the sequence $\to X_k\to X_{k+1}\to $ we have the sequence
 $$\to L\, \widehat{\amalg}\, \overline{\mathbb L}_{T_k} \to L\, \widehat{\amalg}\, \overline{\mathbb L}_{T_{k+1}}\to .$$
 This in turn yields the sequence
 $$ X_{\mathbb Q} \to \dots \to (X_k)_{\mathbb Q} \to (X_{k+1})_{\mathbb Q}\to \dots \to (\cup_k X_k)_{\mathbb Q},$$
 and we set
 $$ {Y} = \varinjlim (X_k)_{\mathbb Q}.$$
 
 Note that $ {Y}$ may not be the Sullivan completion of a topological space ! However we do have
 
 \begin{Prop}
 \label{p29.1} With the hypotheses and notation above
 \begin{enumerate}
 \item[(i)] $\pi_*( {Y}) = \varinjlim \pi_*(X_k)_{\mathbb Q}$ and
 \item[(ii)] $H_*( {Y};\mathbb Q) \cong T= \varinjlim_k T_k = H_*(\vee_k S^{n_k};\mathbb Q)$.
 \end{enumerate}
 \end{Prop}
 
 \vspace{3mm}\noindent {\sl proof.} (i) follows because $ {Y}$ is the limit of an increasing sequence of CW complexes while (ii) follows from Proposition \ref{p13.2} because $\partial_\gamma : T_k\to L.$\hfill$\square$
 
  \subsection{A connected CW complex, $Y$, for which $\pi_*(Y_{\mathbb Q})=\pi_1(Y_{\mathbb Q})$, dim$\, H_1(Y;\mathbb Q) =3$, dim $H_2(Y;\mathbb Q)= \infty$ and $H_p(Y;\mathbb Q)= 0$ for $p>2$. }

\vspace{3mm} Denote by $X= S^1_a\vee S^1_b\vee S^1_c$ a wedge of three circles. Then $\pi_1(X)$ is the free group, $G$, on generators $a,b,c$ defined by the inclusions of the circles. Thus (Proposition \ref{p9.4}(i)) the profree dgl $(L_X,0) = (\overline{\mathbb L}(x,y,z), 0)$ is a dgl model for $X$ where $x,y,z= \log_Ga, \log_Gb,\log_Gc$.

Now use $\overline{ad}(-)$ to denote commutation with  elements   in $G$,
$$\overline{ad}(u)(v)= (u,v)= uvu^{-1}v^{-1}.$$   Then a sequential cell attachment
$$X\longrightarrow Y= \varinjlim_r (X\cup_g \vee_{n=1}^r D_n^2)$$
is defined by $\pi_1(g)[S_n^1]= \overline{ad}^n(a)\overline{ad}^n(b)(c)$. A dgl representative for this attachment is then provided (Proposition \ref{p9.4}(i) and \S 17.1) by the sequence
$$(L_X,0) \to \varinjlim_r (L_X\, \widehat{\amalg}\, \overline{\mathbb L}_{sT(r)}, \partial)$$
where $T(r) = \oplus_{n=1}^r \mathbb Q z_n$ and $\partial sz_n = \log_G(\overline{ad}^n(a)\overline{ad}^n(b)(c))$. Since $\partial : sT(r) \to L_X^2$ it follows from Proposition \ref{p13.2} that
$$\mbox{dim}\, H_p(X\cup_g\vee_{n=1}^r D_n^2;\mathbb Q ) = \left\{\begin{array}{ll}
3 & \mbox{if }p=1\\r&\mbox{if }p=2\\
0 & \mbox{if } p>2\end{array}\right.$$
It follows from the obvious direct limit argument that 
 $$\mbox{dim}\, H_p(Y;\mathbb Q ) = \left\{\begin{array}{ll}
3 & \mbox{if }p=1\\r&\mbox{if }p=2\\
0 & \mbox{if } p>2\end{array}\right.$$

Denote   by $ad(-)$ the conjugation in $L_X: ad_x(y) = [x,y]$. Then   for some element $\gamma_n \in L_X^{>2n+1}$    we have
$$\log_G(\overline{ad}^n(a)\overline{ad}^n(b)(c)) =  {ad}^n(x) {ad}^n(y)(z)  + \gamma_n.$$
Now consider the dgl morphism
$$(L_X,0)\to (L_X\, \widehat{\amalg}\, \overline{\mathbb L}_{sT(r)}, \partial'),$$
where 
$$\partial'sz_n = ad^n(x)ad^n(y)(z).$$

Write $L= L_X\, \widehat{\amalg}\, \overline{\mathbb L}_{sT(r)}$ and suppose  we have   proved that   $H_q(L, \partial')= 0$ for $q>1$. We introduce an extra gradation in $L$ by putting $x,y$ and $z$ in degree $1$ and $z_n$ in degree $2n+1$. We denote $L_{(n)}$ the component of $L$ in the new degree $n$. Then $\partial': L_{(n)}\to L_{(n)}$   and $(\partial-\partial')(L_{(n)}) \subset \oplus_{q>n} L_{(q)}$.

Suppose $\omega$ is a $\partial$-cycle in $L^{\geq 2}$, and write $\omega= \sum_{q\geq q_0} \omega_q$ with $\omega_q\in L_{(q)}$. Then $\omega_{q_0}$ is a $\partial'$-cycle and so a $\partial'$-boundary. We write $\omega_{q_0}= \partial'(\beta_0)$ with $\beta_0\in L_{(q_0)}$. Then denote $\omega'= \omega -\partial (\beta_0)$. We have $\omega'= \sum_{q> q_0}\omega'_q$, and $\omega'_{q_0+1}$ is a $\partial'$-cycle, so a boundary, and there is $\beta_1\in L_{(q_0+1)} $ with $\partial'\beta_1= \omega'_{q_0+1}$. We continue and construct in this way a sequence of elements $\beta_q\in L_{(q_0+q)}$. The series $\sum \beta_q$ is a well defined element in $L$ and $$\partial (\sum \beta_q)= \omega.$$
This shows that $H_{\geq 2}(L, \partial)=0$.

Now we prove that $H_{\geq 2}(L,\partial')=0$. Let $(\land W,d)$ be the quadratic Sullivan model of $L_X= \overline{\mathbb L}(x,y,z)$. Proposition \ref{p6.4} equips $(\land W,d)$ with a second gradation $W=\oplus_{n\geq 0}W(n)$ satisfying
$$W(0) = W\cap \mbox{ker}\, d \hspace{5mm}\mbox{and } d: W(n)\to (\land^2W)(n-1).$$
Then denote by $w_x,w_y$ and $w_z$ the basis of $W(0)$ which is the dual basis for $sx,sy$ and $sz$.

Next set $\alpha_n = ad^n(x)ad^n(y)(z)$. Then $\alpha_n$ determines a cdga morphism $\rho_n : \land W\to \mathbb Q \oplus \mathbb Q a_n$ and satisfying
$$\rho_n(W)\neq 0\hspace{5mm}\mbox{and } \rho_n(W(k))= 0\hspace{2mm}\mbox{ if }k\neq 2n.$$
Set $S(r) = \oplus_{n=1}^r \mathbb Q a_n$. Then
$$\rho(r) =\oplus_{n=1}^r \rho_n : \land W\to \mathbb Q \oplus S(r)$$
is well defined and surjective.

Now we introduce a third gradation in $\land W$ by assigning $W(k)$ the new degree $2k+3$. Then $d$ increases this new degree by $1$, so that $\land W$ is identified with a quadratic Sullivan algebra $\land Z$ in which $W(k) = Z^{2k+3}$ and $Z=Z^{\geq 3}$. In particular, $a_n$ has new degree $2(2n)+3$. Thus $\rho(r)$ with this new grading becomes a surjective morphism 
$$
\rho'(r) : \land Z\to \mathbb Q\oplus S'(r).$$

Now $\land Z\stackrel{\simeq}{\to} \mathbb Q\oplus (Z\cap \mbox{ker}\, d)$, and $Z\cap\mbox{ker}\, d= W(0)$ but with new degree 3. It follows that the homotopy Lie algebra of $\land Z$ is the free graded Lie algebra on the elements $x,y,z$, but regarded as having degree 3. We now show that the sequence $\alpha_n$, regarded as elements  in $L_{\land Z}$, is inert.

 Recall from \cite[Theorem 1.4]{Anick} that the sequence $\alpha_n$ is inert in $L_{\land T}$ if and only if it is inert in its universal enveloping algebra. Moreover following \cite[Theorem 2.1]{Anick} if for some ordering of $a,b,c$ the family $w_i$ of the   terms of maximum height of the $\alpha_i$ satisfies the two following conditions, then the family $\alpha_i$ is inert. The conditions are
 
 (1) No $w_i$ is a submonomial of any other $w_j$.
 
 (2) No $w_i$ "overlaps with" any other $w_j$, i.e. is $w_i = uv$ and $w_j= bw$ then $v= 1$. 
 
 \vspace{2mm} In our situation let $x>y>z$. Then the maximum height in    $\alpha_n $ is $x^ny^nz$ and this family satisfies  the two conditions above.
 
 It follows that the sequence $\alpha_n$ is inert in $\land Z$. This implies that we can suppose   the induced morphism $\lambda': \land V' \to \land Z$ satisfies $\lambda'(V')\subset Z$.
 
 Now return to the original gradings. Since $\land V\cong \land Z$ and $\lambda= \lambda'$, it follows that    $\lambda : V\to W$ is injective, and by Proposition \ref{p19.1} this shows that the sequence $\alpha_n$ is inert in $L_X$. In particular, for each $r$,
 $$L_X\to H(L_X\, \widehat{\amalg}\, \overline{\mathbb L}_{sT(r)}, \partial')\hspace{3mm}\mbox{and } 
 L_X\to H(L_X\, \widehat{\amalg}\, \overline{\mathbb L}_{sT(r)}, \partial)$$
 are surjective. 
 
 Since dim$\, T(r)<\infty$ it follows from Proposition \ref{p16.3} that 
 $$X_{\mathbb Q }\to (X\cup_g \vee_{n=1}^rD^2_n)_{\mathbb Q}$$
 is the spatial realization of $L_X\to H(L_X\, \widehat{\amalg}\, \overline{\mathbb L}_{sT(r)})$. This identifies the surjection $L_X\to H(L_X\, \widehat{\amalg}\, \overline{\mathbb L}_{sT(r)})$ with 
 $$\pi_*(X_{\mathbb Q}) \to \pi_*((X\cup_g \vee_{n=1}^r D_n^2)_{\mathbb Q}).$$
Since $L_X = (L_X)_0$, $\pi_*(X_{\mathbb Q}) = \pi_1(X_{\mathbb Q})$ and so $\pi_{\geq 2} (X\cup_g \vee_{n=1}^r D_n^2)_{\mathbb Q}= 0$.
Passage to direct limits shows that $\pi_*(Y_\mathbb Q)= \pi_1(Y_{\mathbb Q})$. \hfill$\square$

 \vspace{5mm}\noindent Institut de Math\'ematique et de Physique, Universit\'e Catholique de Louvain, 2, Chemin du cyclotron, 1348 Louvain-La-Neuve, Belgium, yves.felix@uclouvain.be

 \vspace{1mm}\noindent Department of Mathematics, Mathematics  Building, University of Maryland, College Park, MD 20742, United States, shalper@umd.edu

\newpage
\begin{thebibliography}{xx}
\bibitem{A1} D. Anick, \emph{Non-commutative algebras and their Hilbert series}, J. Algebra 78 (1982), 120-140.
 \bibitem{An} D. Anick, \emph{A rational homotopy analogue of Whitehead's problem}, in:Algebra, Algebraic Topology and its Applications, Lecture Notes in Math.   1183, Springer (1986), 28-31.
 \bibitem{Anick} D. Anick, \emph{Inert sets and the Lie algebra associated to a group}, J. Algebra 111 (1987), 154-165.
 \bibitem{Milnor} M.G. Barratt and J. Milnor, \emph{An example of anomalous singular homology}, Proceedings AMS 13 (1962), 293-297.
\bibitem{Bou} N. Bourbaki, \emph{El\'ements de Math\'ematique, Th\'eorie des Ensembles, tome 3}, 1970, Hermann.
\bibitem{BG} A.K. Bousfield and V.K. Gugenheim, \emph{On PL de Rham theory and rational homotopy type}, Mem. Amer. Math. Soc. {\bf 179} (1976).
\bibitem{BK} A.K. Bousfield and D.M. Kan, \emph{Homotopy Limits, completions and localizations}, Lecture Notes in Math.  304, Springer, 1972.
\bibitem{brown} K. Brown, \emph{Cohomology of groups}, Graduate Texts in Math. 87, Springer Verlag
\bibitem{Bub} P. Bubenik, \emph{Free and Semi-Inert Cell Attachments}, Trans. Amer. Math. Soc. 357 (2005), 4533-4553.
\bibitem{4A} U.Buijs, Y. F\'elix, A. Murillo and D. Tanr\'e, \emph{Lie models in Topology}, Birkhauser 2020.
 \bibitem{DV} E. Dyer and A. Vasquez, \emph{Some small aspherical spaces}, J. Austr. Math. Soc. 16 (1973), 332-352.
\bibitem{Sdepth} Y. F\'elix and S. Halperin, \emph{The Sdepth of a homotopy Lie algebra}, Geometry, Topology and Mathematical Physics Journal, 1 (2018), 4-26.
\bibitem{Malcev} Y. F\'elix and S. Halperin, \emph{Malcev completions, LS category and depth}, Bol. Soc. Math. Mex. 23 (2017), 267-288.
\bibitem{MScand}  Y. F\'elix and S. Halperin, \emph{The depth and LS category of a topological space}, Math. Scand., 123 (2018), 220-238.
\bibitem{Artin} Y. F\'elix and S. Halperin, \emph{The depth of a Riemann surface and of a right-angled
Artin group},   Journal of Homotopy and Related Structures 15 (2020, 223-248.
\bibitem{Sur}  Y. F\'elix and S. Halperin, \emph{Rational homotopy theory via Sullivan spaces: a survey}, ICCM Notices of the International Congress of Chinese Mathematicians, 5 (2017), 14-36.
\bibitem{FHTI} Y. F\'elix, S. Halperin and J.-C. Thomas, \emph{Rational Homotopy Theory}, Graduate Texts in Math. 201, Springer verlag, 2001.
\bibitem{ranks} Y. F\'elix, S. Halperin and J.C. Thomas, \emph{The ranks of the homotopy groups of a space of finite  LS category}, Expos. Math. 25 (2007), 67-76.
\bibitem{RHTII} Y. F\'elix, S. Halperin and J.-C. Thomas, \emph{Rational Homotopy Theory II}, World Scientific, 2015.
\bibitem{FL} Y. F\'elix and J.-M. Lemaire, \emph{On the Pontryagin algebra of the loops on a space with a cell attached}, Internat. J. Math. 2 (1991), 429-438.
\bibitem{FT} Y. F\'elix and J.-C. Thomas, \emph{Module d'holonomie d'une fibration}, Bull. Math. Soc. France 113 (1985), 255-258.
\bibitem{hall}  M. Hall, \emph{A basis for free Lie rings and higher commutators in free groups}, Proc. Amer. Math. Soc. 1(1950) 575-581. 
\bibitem{HGT} S. Halperin, A. Gomez-Tato and D. Tanr\'e, \emph{Rational homotopy theory for non-simply connected spaces}, Trans. Amer. Math. Soc. 352 (2000), 1493-1525.
\bibitem{HL} S. Halperin and J.M. Lemaire, \emph{Suites inertes dans les alg\`ebres de Lie gradu\'ees}, Math. Scand 61 (1987), 39-67.
\bibitem{HL2} S. Halperin and J.-M. Lemaire, \emph{The fibre of a cell attachment}, Proc. Roy. Soc. Edinburgh, 38 (1995), 295-311.
\bibitem{HessL} K. Hess and J.-M. Lemaire, \emph{Nice and lazy cell attachments}, J. of Pure and Applied Algebra 112 (1996), 29-39.
\bibitem{HS} S. Halperin and J.D. Stasheff, \emph{Obstructions to homotopy equivalences}, adv. in Math. 32 (1979), 233-279.
 \bibitem{IM} S. Ivanov and R. Mikhailov, \emph{A finite $\mathbb Q$-bad space}, Geometry and Topology 23 (2019), 1237-1249.
 \bibitem{Jech} Thomas Jech, \emph{Set Theory, the third millennium edition, revised and expanded}, 2003, Springer Monographs in Mathematics.
\bibitem{La} M. Lazard, \emph{Sur les groupes nilpotents et les anneaux de Lie}, Ann. Sc. Ec. Norm. Sup. 71 (1959), 101-190.
\bibitem{LM} A. Lazarev and M. Markl, \emph{Disconnected rational homotopy theory}, Adv. Math. 283 (2015), 303-361. 
\bibitem{JML} J.-M. Lemaire, \emph{A finite complex whose rational homology is not finitely generated}, in Lecture Notes in Mathematics 196 (1971), 114-120. 
\bibitem{Lemaire} J.-M. Lemaire, \emph{"Autopsie d'un meurtre" dans l'homologie d'une alg\`ebre de cha\^ines}, Ann. Scient. Ec. Norm. Sup. 11 (1978), 93-100.
 \bibitem{Ly} R.C. Lyndon, \emph{Cohomology theory of groups with a single defining relation}, Ann. of Math. 52 (1950), 650-656. 
  \bibitem{May} J.P. May, \emph{Simplicial objects in algebraic topology}, Chicago Univ. Press (1967).
  \bibitem{MM} J. Milnor and J.C. Moore, \emph{On the structure of Hopf algebras}, Annals of Math. 81 (1965), 211-264.
\bibitem{PS} S. Papadima and A. Suciu, \emph{Algebraic invariants for right-angled Artin groups}, Math. Annalen 334 (2006), 533-555.
\bibitem{Pu} A. Putman, \emph{One-relator groups}, preprint 2018.
\bibitem{Pridham} J.P. Pridham, \emph{Pro-algebraic homotopy types}, Arxiv 2008
\bibitem{Q} D. Quillen, \emph{Rational homotopy theory}, Ann. of Math. 90 (1969), 205-295.
\bibitem{Serre} J.P. Serre, \emph{Lie algebras and Lie groups}, Benjamin Inc. (1965).
\bibitem{Sh} A.I. Shirshov, \emph{Subalgebras of free Lie algebras}, Mat. Sbornik 33 (1975), 441-452.
\bibitem{Solo} R. Solovay, \emph{$2^{\aleph_0}$ can be anything it ought to be}, Symposium on the Theory of Models, North-Holland Publish Co (1965), pg 435. 
\bibitem{S} D. Sullivan, \emph{Infinitesimal computations in topology}, Publ. IHES 47 (1977), 269-331.
\bibitem{Wat} C. Watkiss, unpublished pager (1976).
 \bibitem{W} J.H.C. Whitehead, \emph{On adding relations to homotopy groups}, Ann. Math. 42 (1941), 409-428.

\end{thebibliography}
\end{document}